\documentclass[11pt]{article}

\usepackage{amsmath}
\usepackage{amsfonts}
\usepackage{amssymb}
\usepackage{epsfig}
\usepackage{dsfont}
\usepackage{mathrsfs}
\usepackage{bbm}
\usepackage{bm}
\usepackage{multirow}
\usepackage{float}
\usepackage{makecell} 
\usepackage{boldline}
\usepackage{diagbox}
\usepackage{threeparttable}
\usepackage[dvipsnames]{xcolor}
\usepackage{adjustbox}
\usepackage{booktabs}
\usepackage{algorithm}
\usepackage{algorithmic}
\allowdisplaybreaks[4]
\usepackage{caption}
\usepackage{subfigure}
\usepackage{mathabx} 
\usepackage{authblk}
\usepackage{float}
\usepackage[title]{appendix}
\usepackage{boldline}
\usepackage{diagbox}
\usepackage{threeparttable}
\usepackage[dvipsnames]{xcolor}
\usepackage{adjustbox}
\usepackage{booktabs}

\RequirePackage{xcolor}[2.11]
\colorlet{inlinkcolor}{green!60!black}
\colorlet{exlinkcolor}{red!50!black}
\colorlet{reviewcolor}{black!50}

\usepackage[colorlinks=true]{hyperref}
\hypersetup{urlcolor=exlinkcolor, citecolor=inlinkcolor,linkcolor=inlinkcolor}
\usepackage[capitalise,nameinlink]{cleveref}
\crefname{equation}{}{}
\crefname{figure}{Figure}{Figures}
\crefname{assumption}{Assumption}{Assumptions}
\crefname{subappendix}{Appendix}{Appendices}

\makeatletter
\newenvironment{breakablealgorithm}
{
\refstepcounter{algorithm}
\hrule height1pt depth0pt \vspace{-0.2cm} \kern2pt
\renewcommand{\caption}[2][\relax]
{
{\raggedright\textbf{\ALG@name~\thealgorithm} ##1\par}%
\ifx\relax##1\relax 
\addcontentsline{loa}
{algorithm}
{\protect\numberline{\thealgorithm}##2}%
\else 
\addcontentsline{loa}{algorithm}{\protect\numberline{\thealgorithm}##1}%
\fi
\kern2pt\hrule\kern0pt
}
\small
}
{

\kern2pt\hrule \relax
}
\makeatother

\newtheorem{theorem}{Theorem}[section]
\newtheorem{lemma}{Lemma}[section]
\newtheorem{proposition}{Proposition}[section]
\newtheorem{example}{Example}[section]
\newtheorem{remark}{Remark}[section]
\newtheorem{definition}{Definition}[section]
\newtheorem{assumption}{Assumption}[section]
\newtheorem{corollary}{Corollary}[section]

\newcommand{\dd}{\mathsf {d\kern -0.07em l}} 

\newcommand{\bgeqn}{\begin{eqnarray}}
\newcommand{\edeqn}{\end{eqnarray}}
\newcommand{\bgeq}{\begin{eqnarray*}}
\newcommand{\edeq}{\end{eqnarray*}}
\newcommand{\bec}{\begin{center}}
\newcommand{\enc}{\end{center}}
\newcommand{\R}{{\rm I\!R}}

\newcommand{\inmat}[1]{\mbox{\rm {#1}}}

\renewcommand{\Box}{\framebox{\rule{0.3em}{0.0em}}}

\newcommand{\ra}{\rightarrow}

\newcommand{\be}{\begin{equation}}
\newcommand{\ee}{\end{equation}}

\newcommand{\vt}{{\vartheta}}

\newcommand{\bdx}{\bm{x}}

\newcommand{\bdxi}{\bm{\xi}}

\newcommand{\bdz}{\bm{z}}
\newcommand{\bdt}{\bm{t}}

\newcommand{\bdf}{\bm{f}}

\numberwithin{equation}{section}
\numberwithin{definition}{section}

\newcommand{\calu}{\mathcal{U}}
\newcommand{\tldu}{\tilde{u}}

\newcommand{\calx}{\mathcal{X}}
\newcommand{\caly}{\mathcal{Y}}

\newcommand{\scru}{\mathscr{U}}
\newcommand{\scrg}{\mathscr{G}}
\newcommand{\lt}{\left}
\newcommand{\rt}{\right}
\newcommand{\disp}{\displaystyle}
\newcommand{\st}{{\rm s.t.}}

\newcommand{\1}{\mathds{1}}

\newcommand{\bdps}{\bm{\psi}}
\newcommand{\bdc}{\bm{C}}
\newcommand{\bdca}{\bm{A}}
\newcommand{\bdcb}{\bm{B}}

\newcommand{\bdg}{\bm{g}}
\newcommand{\bdcz}{\bm{Z}}

\def\la{{\langle}}
\def\ra{{\rangle}}

\def\bbe{{\mathbb{E}}}

\title{Multi-Attribute Utility Preference Robust Optimization: 
A Continuous Piecewise Linear Approximation Approach\thanks{This project is supported by RGC grant 14500620.
}}
\author[1]{Qiong Wu}
\author[1]{Sainan Zhang}
\author[2]{Wei Wang}
\author[1]{Huifu Xu\thanks{Corresponding Author.}}
\affil[1]{Department of Systems Engineering and Engineering Management,The Chinese University of Hong Kong, Shatin, N.T., Hong Kong.}
\affil[2]{School of Business, University of Southampton, Southampton SO17 1BJ, UK.}
\affil[ ]{\textsf{\{qiwu,snzhang,hfxu\}@se.cuhk.edu.hk,ww1e17@soton.ac.uk}.}
\vspace{-1cm}
\date{\today}
\begin{document}
\maketitle

\begin{abstract}
In this paper, we consider a multi-attribute decision making problem 
where the decision maker's (DM's) objective is to maximize the expected utility of
outcomes but the true utility function which captures the DM's risk preference is ambiguous. We propose a maximin multi-attribute
utility preference robust optimization (UPRO) model where the optimal decision is based on the worst-case utility function in an ambiguity set of plausible utility functions constructed 
using partially available information 
such as the DM's specific preferences between some lotteries. Specifically, we consider a
UPRO model with two attributes, where
the DM's risk attitude is multivariate 
risk-averse and the ambiguity set is defined by a linear system of inequalities 
represented by the Lebesgue–Stieltjes (LS) integrals of the DM's  utility functions.
To solve the maximin problem, 
we propose an explicit piecewise linear approximation (EPLA) scheme 
to approximate the DM's true unknown utility 
so that the inner minimization problem reduces to a linear program, 
and we solve the approximate maximin problem 
by a derivative-free (Dfree) method. 
Moreover, by introducing binary variables to locate the position of the reward function 
in a family of simplices,
we propose an implicit piecewise linear approximation (IPLA) representation 
of the approximate  UPRO and solve it using the Dfree method. Such IPLA technique prompts us to reformulate the approximate UPRO as a single mixed-integer program (MIP) and extend the tractability of the approximate UPRO to the multi-attribute case. Under some moderate conditions, we derive error bounds between the UPRO  and the approximate UPRO  in terms of the  ambiguity set, the optimal value and the optimal solution.
Furthermore, we extend the model to the expected utility maximization problem with expected utility constraints where the worst-case utility functions 
in the objective and constraints are considered simultaneously. Finally, we report the numerical results about 
performances of the proposed models and the computational schemes, and show that the schemes work efficiently 
and the UPRO model is stable against 
data perturbation.
\end{abstract}


\section{Introduction}

{\color{black}
Utility preference robust optimization (UPRO) model concerns the optimal decision making 
where the decision maker (DM) aims to maximize the expected utility but 
the true utility function which captures the DM's preference is ambiguous. 
Instead of finding an approximate utility function using partially available information as in the literature of behavioural economics (see, e.g., \cite{clemen2013making} and \cite[Chapter 10]{gonzalez2018utility}), 
the UPRO models construct a set of plausible utility functions and
base the optimal decision on the worst-case utility function from the set to mitigate such ambiguity.
This 
type
of approach can be traced back to Maccheroni~\cite{maccheroni2002maxmin}
who considers the worst-case utility evaluation among 
a number of available 
utilities when a conservative DM faces uncertain outcomes of lotteries.
He derives necessary and sufficient conditions 
for the existence of 
a set of utility functions such that the {\color{black}worst-case} 
in 
the set can be used to characterize the conservative decision making   
framework.
Armbruster and Delage \cite{AmD15}
give a comprehensive treatment 
of the problem from 
optimization perspective by formally 
proposing a maximin  UPRO paradigm.
Specifically, they consider  a class of  utility functions which 
are  concave or S-shaped and discuss how a DM's preference 
may be elicited through pairwise comparisons.  
Moreover, 
they demonstrate that solving the UPRO 
model
reduces to
solving a linear program (LP) under some mild conditions.
Over the past few years, the research of 
UPRO related models 
has received increasing attentions, see for instance \cite{hu2015robust,haskell2016ambiguity,hu2018robust,GXZ21,DGX22,WuX22}.

The above UPRO models are all about single-attribute decision making problem. 
In practice, there is a multitude of {\color{black}literature} focusing on the multi-attribute case.
For instance,
in healthcare it is typical to use several
metrics rather than just one to measure the quality of life (\cite{feeny2002multiattribute,torrance1982application}).
Similar problems can be found in network management \cite{azaron2008multi,chen2010stochastic},
scheduling \cite{Zakariazadeh2014}, 
multiobjective design optimization problem 
\cite{Dino2017,tseng1990minimax}, and portfolio
optimization \cite{fliege2014robust}. Indeed, over the past few decades,
there has been significant research on multi-attribute expected utility
\cite{fishburn1992multiattribute,miyamoto1996multiattribute,tsetlin2006equivalent,tsetlin2007decision,tsetlin2009multiattribute,von1988decomposition}.

Zhang et al.~\cite{zhang2020preference} seem to be the first to 
propose 
a preference robust optimization (PRO) model for multi-attribute decision making. Specifically, they 
consider 
a multivariate shortfall risk minimization problem where
there is an ambiguity in an investor's true disutility function of losses and they consider
the worst-case disutility function 
in an ambiguity set
to calculate the risk measure.
Wu et al.~\cite{wu2020preference}
propose a general PRO model 
for 
 multi-attribute decision making. Instead of 
considering expected utility,
they consider a 
quasi-concave choice function 
to measure the DM's multi-attribute rewards
which subsumes the expected utility model as a special case,
and propose
a support function-based approach 
to solve
the resulting 
preference robust choice problem. Since 
the model is very general,
the computational scheme does not benefit 
from the specific structure that it would do in expected utility maximization problems.
For example, it is unclear whether we can use piecewise linear utility functions to approximate the true unknown utility function
in the multi-attribute UPRO models 
as in the single-attribute case (\cite{GXZ21}).
We are interested in the piecewise linear approximation (PLA) approach for several reasons. First, a DM's utility preference is usually elicited at some discrete points. Connecting the 
utility function values at these points to form a piecewise linear utility function is the easiest way to obtain an approximate utility function.
Second, the  PLA approach works for a broader class of
UPRO models without specific requirements on convexity, S-shapedness or quasiconvexity of the true utility function.
Third, despite the PLA approach does not solve UPRO models precisely as the support function-based approach, it allows 
us to derive an  error bound under some moderate conditions.
 
In this paper, we endeavour to carry out a comprehensive 
study on the multi-attribute UPRO from modelling to computational schemes and underlying theory. Unlike 
single-attribute case, a conservative DM's utility function is not necessarily concave which means that
Armbruster-Delage's support function-based approach 
is not applicable in this case.  
This prompts us to 
adopt the PLA approach. The extension of the 
PLA approach from single-attribute UPRO to 
multi-attribute
UPRO would be trivial if the utility function is additive or concave. 
However, 
when we consider 
a general multi-attribute utility function without
specific independence condition,
the construction,
representation of PLA and subsequent computation of the approximate 
UPRO require 
a lot of new work.
One of the 
challenges that we have to tackle is to find an appropriate
representation
of a piecewise linear utility function which is 
easy to construct, and 
to embed
in the objective function and in the ambiguity set.
The main contributions of this paper can be summarized as follows. 

First, we 
propose a maximin robust optimization model for
bi-attribute decision making 
where the DM is
multivariate
risk-averse,
there is incomplete information to identify 
the DM's true utility function, and the optimal decision is based on the worst-case utility function in an ambiguity set. We discuss in detail how the ambiguity set 
of bivariate utility functions may be constructed by 
standard preference elicitation 
methods such as pairwise comparisons.
To solve 
the maximin problem, we propose a two-dimensional continuous 
PLA scheme
to approximate 
the true unknown utility function
so that the inner minimization problem 
can be reduced to a finite-dimensional program.
Differing from the one-dimensional case, 
we divide the domain of the utility functions
into a set of mutually exclusive triangles
and define an approximate 
utility function which is linear 
over
each of the triangles. 
Moreover, 
we reformulate the ambiguity
set defined by 
the expected utility values of 
the DM's preferences 
between
lotteries
into
the one where the expected utility values are 
represented by the Lebesgue–Stieltjes (LS) integrals with respect to (w.r.t.) the utility function. 
The PLA approach allows us to derive the approximate utility function explicitly using indicator functions 
{\color{black}and} characterize the Lipschitz continuity of the utility function by individual variables,   
and 
{\color{black}enables} us to calculate the  LS integrals conveniently.




Second, by exploiting the piecewise linearity of the approximate utility function,  
we use the well-known polyhedral method
(\cite{DLM10,KDN04,LeW01,VAN10,vielma2015mixed}) 
to represent the 
multi-attribute reward functions
using  a convex combination of 
the vertices of the 
simplex
containing the vector 
in the domain of 
the  multivariate utility function,
and subsequently
reformulate the inner 
approximate 
utility minimization problem as a mixed-integer program (MIP). Differing from the PLA approach described above, the approximate utility function cannot be represented explicitly, rather it is determined by solving an MIP. We call this implicit PLA (IPLA) whereas the former is explicit PLA (EPLA).  
A clear benefit of IPLA is that 
it works for multidimensional cases and 
also allows us to reformulate the whole maximin problem as a single MIP. 

Third, we extend the preference robust approach to the expected utility maximization problem with expected utility constraints. Instead of considering the worst-case utility in the objective and 
the worst-case utility in the constraints separately, 
we propose a UPRO model where the optimal decision is based on the same worst-case utility function in both the objective and the constraints. We derive conditions under which the two robust formulations are equivalent and
carry out comparative analysis through numerical studies to identify 
the differences that the two models may render.

Fourth, to justify the PLA scheme, we derive {\color{black}error bounds} for the optimal value and the optimal solutions, which is built on a newly derived Hoffman's lemma for the linear system in the infinite-dimensional space under the pseudo-metric. 
We also quantify the difference between {\color{black}the ambiguity sets} before and after the PLA
and indicate the special cases when these two ambiguity sets coincide. Moreover, to facilitate the application of the UPRO model in a data-driven environment, we carry out stability analysis on the optimal value and the optimal solutions of the UPRO model against 
data 
perturbation/contamination. 

Finally, we undertake extensive numerical tests on the proposed UPRO models and 
computational schemes and obtain the following main findings. The EPLA scheme (see (\ref{eq:PRO-N-reformulate})) and 
the IPLA scheme (see (\ref{eq:PRO_MILP_eqi}))
generate the same results in terms of the convergence of the worst-case utility functions and the optimal values, 
but the former works much faster because 
the IPLA requires solving a MILP as opposed to an LP in the EPLA and the number of integer variables 
increases rapidly with the increase of scenarios of
the underlying exogenous uncertainty.
The 
IPLA works also well in tri-attribute case 
although
the CPU time is long. 
For the constrained expected utility maximization problem, the two robust models 
may coincide in some cases but differ in other cases depending on 
the constraints. The approximate maximin model is stable in the presence of small perturbations arising 
during
the preference elicitation process and
{\color{black}resulting from} exogenous uncertainty data.

The rest of the paper is organized as follows. 
Section~\ref{sec:2Bi-Attribute} introduces the multi-attribute UPRO model 
and the definition 
of the ambiguity set. 
Section~\ref{sec:numer-methods}
gives the details of
the EPLA 
approach
and 
tractable formulation of 
approximate UPRO in bi-attribute case.
%
Section~\ref{sec:multi-atrribute} discusses the IPLA approach for the UPRO 
in multi-attribute case.  
Section~\ref{sec-errorbound} investigates the error bound of the approximate ambiguity set as well as the impact on the optimal value and the optimal solutions to the UPRO model. 
Section~\ref{sec:constrained} extends the UPRO model to the constrained optimization problem.
Section~\ref{sec:numerical results} reports the numerical tests of the UPRO model. Concluding remarks are given in Section~\ref{sec:Concluding remarks}. 




\section{The bi-attribute UPRO model}
\label{sec:2Bi-Attribute}

We consider the following one-stage expected bi-attribute utility maximization problem 
\bgeqn 
\label{eq:UMP-bi}
\max_{\bdz\in Z} \; \bbe_P[u(\bdf(\bdz,\bdxi))],
\edeqn
where $\bdf:\R^{n}\times \R^{m} \to \R^2$ is a continuous vector-valued function representing the rewards from two attributes,
$\bdz$ is a decision vector which is restricted to taking values over a specified feasible set $Z\subset \R^n$, 
$\bdxi$ is a random vector
representing exogenous uncertainties in the decision making problem mapping from probability space $(\Omega,\mathcal{F},\mathbb{P})$ to $\R^m$,
the expectation is taken w.r.t.~the probability of $\bdxi$, i.e., $P:=\mathbb P\circ\bdxi^{-1}$, 
and $u:\R^2\to\R$ is a real-valued 
componentwise non-decreasing continuous utility function, which maps each value of $\bdf$ to a utility value of the DM's interest. 
To facilitate our discussions, we make the following assumption throughout the paper.

\begin{assumption}
\label{assu-original}
$\bdf$ is a continuous function 
with its range covered by
$T:=X\times Y$ with $X:=[\underline{x},\bar{x}]$ and $Y:=[\underline{y},\bar{y}]$,
$Z$ is a compact convex subset of $\R^n$ and the support set $\Xi$ of $\bdxi$ is compact. 
\end{assumption}

Assumption~\ref{assu-original} allows us to restrict the domain of the unknown true utility function to a rectangle $T$.
We follow \cite{hu2015robust} and the literature of behavioural economics to 
normalize the utility function with $u(\underline{x},\underline{y})=0$ and $u(\bar{x},\bar{y})=1$. 
In most of the existing research on multi-attribute decision making,  
utility functions are assumed to be known (\cite{greco2016multiple,liesio2021nonadditive}) 
or can be elicited and estimated through a tolerable amount of questions
(\cite{andre2007non}). 
In practice, however, a DM's utility function is often unknown 
either from the DM's perspective or from the modeller's perspective (\cite{AmD15}).

In this paper, our focus 
is on the situation where the DM does not have complete information to identify the true utility function $u^*$, i.e., 
risk preference, 
but it is possible to elicit partial information to construct an ambiguity set of utility functions, denoted by $\mathcal{U}$, 
such that the true utility function 
which represents the DM's preference lies within $\mathcal{U}$ with high likelihood.
Under this circumstance, it might be sensible to consider 
the following bi-attribute utility preference robust optimization model 
to mitigate the model risk arising from the ambiguity in the true utility function 
\begin{equation}
\label{eq:MAUT-robust}
    \inmat{(BUPRO)} \quad
    \vt:=\max_{\bdz\in Z} \; \min_{u\in {\cal U}} \; \bbe_P[u(\bdf(\bdz,\bdxi))].
\end{equation}
The structure of the BUPRO model is largely determined by the structure of the ambiguity set $\calu$ as well as the nature of the utility functions in this set. 
Various approaches have been proposed
to construct 
an ambiguity set of utility functions in the literature of PRO depending on the availability of information (see \cite{AmD15, liu2021multistage,guo2022robust}). 
They are usually based on two types of information about a DM's preference: 
generic information such as risk aversion or risk taking 
and specific information such as preferring one prospect to another (see \cite{WaX23}).  

In single-attribute decision making, 
a DM is risk-averse if and only if the DM's utility function is concave (see \cite{tsanakas2003risk}).
Unfortunately, the equivalent relation 
does not hold in the multi-attribute case. 
Let $x_0,x_1\in X$, $y_0,y_1\in Y$ with  $x_0< x_1$ and $y_0<y_1$. 
Consider the following two lotteries:
Lottery one ($L_1$) gives the DM a 0.5 chance of receiving $(x_0,y_0)$ and a 0.5 chance of receiving $(x_1,y_1)$. 
Lottery two ($L_2$) gives the DM a 0.5 chance of receiving $(x_0,y_1)$ and a 0.5 chance of receiving $(x_1,y_0)$.
The DM is said to be \emph{multivariate risk-averse} (MRA) 
if the DM prefers $L_2$ to $L_1$ 
for all $x_0,x_1,y_0$ and $y_1$ described above
(see e.g.,\cite{richard1975multivariate}).
This type of behaviour means that the DM prefers  taking 
a mix of the best and worst in the two respective attribute
to getting 
either the ``best'' or the ``worst'' with equal probability.
Using the expected utility theory, we can write down the DM's preference mathematically as
$0.5u(x_0,y_0)+0.5u(x_1,y_1) \leq 0.5u(x_0,y_1)+0.5u(x_1,y_0)$,
which is equivalent to
\begin{equation}
\label{eq:conservative}
u(x_0,y_1)+u(x_1,y_0)\geq u(x_0,y_0)+u(x_1,y_1)
\end{equation}
for all $x_0,x_1,y_0$ and $y_1$. 
(\ref{eq:conservative}) is known as 
\emph{conservative property}.
In the case when the utility function is twice continuously differentiable,
the property is equivalent to 
$u_{x y}:=\frac{\partial^2 u}{\partial x \partial y}\leq 0$ for all $(x,y)\in T$,
see \cite[Theorem 1]{richard1975multivariate}.
This kind of definition is given in \cite{richard1975multivariate},  
{\color{black}and} there are some other definitions of MRA, see e.g.~\cite{ duncan1977matrix,karni1979multivariate,levy1991arrow} and references therein. 

From the definition, we can see immediately that a
risk-averse DM's utility function is not necessarily concave (e.g.~
%
$u(x,y)=x+y-(x y)^{1/4}$ for $x>0$ and $y>0$). 
This is a fundamental difference between the multi-attribute and single-attribute utility functions.
In the forthcoming discussions, we will consider 
utility functions satisfying
(\ref{eq:conservative})
since 
risk-averse is widely considered in the literature (\cite{abbas2005attribute,abbas2009multiattribute}),
 e.g., $u(x,y)=\frac{1-e^{-\gamma(x+\beta y)}}{1-e^{-\gamma(1+\beta)}}$ with $\gamma>0$, $\beta>0$,
and $u(x,y)=e^x-e^{-y}-e^{-x-2y}$.
Specific information about a DM's preference is often obtained by a modeller during 
a preference elicitation process.
The most widely used elicitation method is pairwise comparison
(\cite{AmD15}).
For instance, a DM is given a pair of lotteries ${\bm A}$ and ${\bm B}$ defined over $(\Omega,{\cal F},\mathbb{P})$ with different outcomes and asked for preference.
If the DM prefers ${\bm A}$, then we can use the expected utility theory to characterize the preference,
i.e., 
\begin{equation*}
\bbe_{\mathbb{P}}[u(\bdcb(\omega))] = \int_{T} u(x,y) d F_{\bdcb}(x,y)\leq \int_{T} u(x,y) d F_{\bdca}(x,y) = \bbe_{\mathbb{P}}[u(\bdca(\omega))]
\end{equation*}
or equivalently
\begin{equation}
\label{eq:ambi-U-ex}
\int_{T} u(x,y) d \psi(x,y):=\int_{T} u(x,y) d (F_{\bdcb}(x,y)-F_{\bdca}(x,y))\leq 0, 
\end{equation} 
where
$F_{{\bm A}}$ and $F_{{\bm B}}$ are the cumulative distribution functions of ${\bm A}$ and $\bdcb$, 
$u$ is the true utility function which represents the DM's preference but is unknown. 
The outcomes of the pairwise comparisons enable us to narrow down the scope of the utility function by inequalities. 
As more and more questions are asked,   
we 
can derive more inequalities
as such {\color{black}which 
lead to} a smaller ambiguity set.  
To facilitate discussions, we give a formal definition of the ambiguity set constructed as such.

\begin{definition}
\label{defi:ambguity-set}
Let $\scru$ be the set of continuous, componentwise non-decreasing, and 
normalized utility functions mapping from $T$
to $[0,1]$ 
satisfying conservative property (\ref{eq:conservative}). 
Define the ambiguity set of utility functions  
as
\begin{equation}
\label{eq:ambiguity_set}
    \mathcal{U}:=\left\{ u\in \scru 
    \,:\,\int_{T} u(x,y)d\psi_l(x,y)\leq c_l, l=1,\ldots,M  \right\},
\end{equation}
where $\psi_l:T\rightarrow \R$ is a real-valued function and
$c_l$ is a given constant  for $l=1,\ldots, M$, and the integrals are in the sense of Lebesgue-Stieltjes 
integration.
\end{definition}
In this definition, we make a blanket assumption that 
the LS
integrals are well-defined,
we refer readers to \cite{clarkson1933definitions}, 
\cite[page 129]{hildebrandt1963introduction}
and \cite{Mcs47} for the concept and properties of the integration. 
${\cal U}$ in (\ref{eq:ambiguity_set}) is defined by a system of inequalities 
which are linear in both $u$ and $\psi_l$. 
Thus the ambiguity set ${\cal U}$ defined as such is a convex set. 
Moreover, we assume that the DM's preferences shown during the elicitation process are consistent, which means that ${\cal U}$ is non-empty.
In practice, preferences observed over an elicitation process may be inconsistent due to observation/measurement errors, noise in data 
or the DM's wrong answers. We refer readers to 
\cite{AmD15} and
\cite{BeO13} for approaches to handle the inconsistency.

\section{Explicit piecewise linear approximation of BUPRO}
\label{sec:numer-methods}

We now move on to discuss how to solve the maximin problem (\ref{eq:MAUT-robust}).
Since the true utility function
is not necessarily concave, we cannot
adopt the 
{\color{black}support function-based approach} 
used in single-attribute UPRO models (see \cite{AmD15})  
and in multi-attribute UPRO models (see \cite{zhang2020preference}).
Instead, we use the PLA
approach considered in \cite{GXZ21}.
The main challenge is that 
constructing a PLA of a bivariate utility function is much more complex than that of a univariate utility function. In this section, we discuss the details.

Let ${\cal X}:=\{x_i,i=1,\ldots,N_1\}\subset X$ and 
${\cal Y}:=\{y_j,j=1,\ldots,N_2\}\subset Y$ with 
$\underline{x}=x_1<\ldots<x_{N_1}=\bar{x}$ and $\underline{y}=y_1<\ldots<y_{N_2}=\bar{y}$.
We define $\calx\times \caly:=\{(x_i,y_j), x_i\in {\cal X},y_j\in {\cal Y}\}$ as a set of $N_1N_2$ gridpoints. 
Let $X_1:=[x_1,x_2]$, $X_i:=(x_i,x_{i+1}]$ for $i=2,\ldots,N_1-1$ and $Y_1:=[y_1,y_2]$, $Y_j:=(y_j,y_{j+1}]$ for $j=2,\cdots,N_2-1$. We divide $T$ into $(N_1-1)(N_2-1)$ mutually exclusive 
cells $T_{i,j}:=X_i\times Y_j$, $i=1,\cdots,N_1-1$, $j=1,\cdots,N_2-1$
and
$T=\bigcup_{i=1}^{N_1-1}\bigcup_{j=1}^{N_2-1} T_{i,j}$.

There are two ways to define a continuous piecewise linear function over a cell
$T_{i,j}$.
One is to define two linear pieces over the two triangle areas separated 
using the {\em main diagonal} (Type-1 PLA) connecting $(x_i,y_j)$ and $(x_{i+1},y_{j+1})$ and the other is using the {\em counter diagonal} (Type-2 PLA) connecting $(x_i,y_{j+1})$ and $(x_{i+1},y_{j})$, see Figure~\ref{fig-division-all} for an illustration.
Consider Type-1.
For any $(x,y)\in T_{i,j}$, if $\frac{y_{j+1}-y_j}{x_{i+1}-x_i}\leq\frac{y-y_j}{x-x_i}$, then $(x,y)$ lies in the upper triangle 
and the upper linear piece of the utility function is defined as
\begin{equation}
\label{eq-up}
    u^{1u}_{i,j}(x,y) := \frac{y_{j+1}-y}{y_{j+1}-y_j} u_{i,j}
    +\lt(\frac{y-y_j}{y_{j+1}-y_j}-\frac{x-x_i}{x_{i+1}-x_i}\rt) u_{i,j+1}
    +\frac{x-x_i}{x_{i+1}-x_i} u_{i+1,j+1}.
\end{equation}
If $0\leq\frac{y-y_j}{x-x_i}\leq \frac{y_{j+1}-y_j}{x_{i+1}-x_i}$,
then $(x,y)$ lies in the lower triangle 
and
\begin{equation}
\label{eq-lo}
    u^{1l}_{i,j}(x,y) := \frac{x_{i+1}-x}{x_{i+1}-x_i} u_{i,j}
    +\lt( \frac{x-x_i}{x_{i+1}-x_i}-\frac{y-y_j}{y_{j+1}-y_j} \rt) u_{i+1,j}
    +\frac{y-y_j}{y_{j+1}-y_j} u_{i+1,j+1},
\end{equation}
where $u_{i,j} := u(x_i,y_j)$, $i=1,\cdots,N_1,j=1,\cdots,N_2$.

Note that this kind of definition is based on 
the interpolation method using the utility values 
at the three vertices of the triangles. 
It differs
significantly from Guo and Xu \cite{GXZ21} and Hu~et~al.~\cite{hu2022distributionally}, 
where each linear piece is defined by
a slope intercept form.
We do not
adopt their approaches because
in multi-attribute case  they 
require the utility values of two neighbouring active linear pieces to 
coincide on the boundary of each cell, 
which would significantly complicate the representation of PLA.

We now turn to discuss the construction of Type-2 PLA.
The upper 
and lower linear 
pieces 
can be defined respectively as
\begin{equation}
\label{eq-up-2}
    u_{i,j}^{2u}(x,y) := \frac{x_{i+1}-x}{x_{i+1}-x_i} u_{i,j+1}
    +\lt( \frac{x-x_i}{x_{i+1}-x_i}-\frac{y_{j+1}-y}{y_{j+1}-y_j} \rt) u_{i+1,j+1}
    +\frac{y_{j+1}-y}{y_{j+1}-y_j} u_{i+1,j}
\end{equation}
and 
\begin{equation}
\label{eq-lo-2}
    u_{i,j}^{2l}(x,y) := \frac{y-y_j}{y_{j+1}-y_j} u_{i,j+1}
    +\lt( \frac{x_{i+1}-x}{x_{i+1}-x_i}-\frac{y-y_j}{y_{j+1}-y_j} \rt) u_{i,j}
    + \frac{x-x_i}{x_{i+1}-x_i} u_{i+1,j}.
\end{equation}
Notice that the conservative property for the utility function plays an important role, that is, 
\begin{equation}
\label{eq:consevative-condition}
u_{i,j+1}+u_{i+1,j}\geq u_{i,j}+u_{i+1,j+1} \quad \forall i=1,\cdots,N_1-1,j=1,\cdots,N_2-1.
\end{equation}
If the conservative property holds at each cell, 
then the graph of the Type-2 piecewise linear function majorizes 
that of the Type-1, see Figure~\ref{fig-diagall-b}.
In this case, the diagonal line 
connecting points $1$ and $4$ looks like a ``valley", 
while the segment connecting points $2$ and $3$ looks like a ``ridge".
\begin{figure}[!ht]
  \vspace{-0.4cm}
  \centering
   \hspace{-1em}
  \subfigure[\scriptsize{Main $\&$ counter diagonals}]{ 
    \label{fig-division-all} 
    \includegraphics[width=0.3\linewidth]{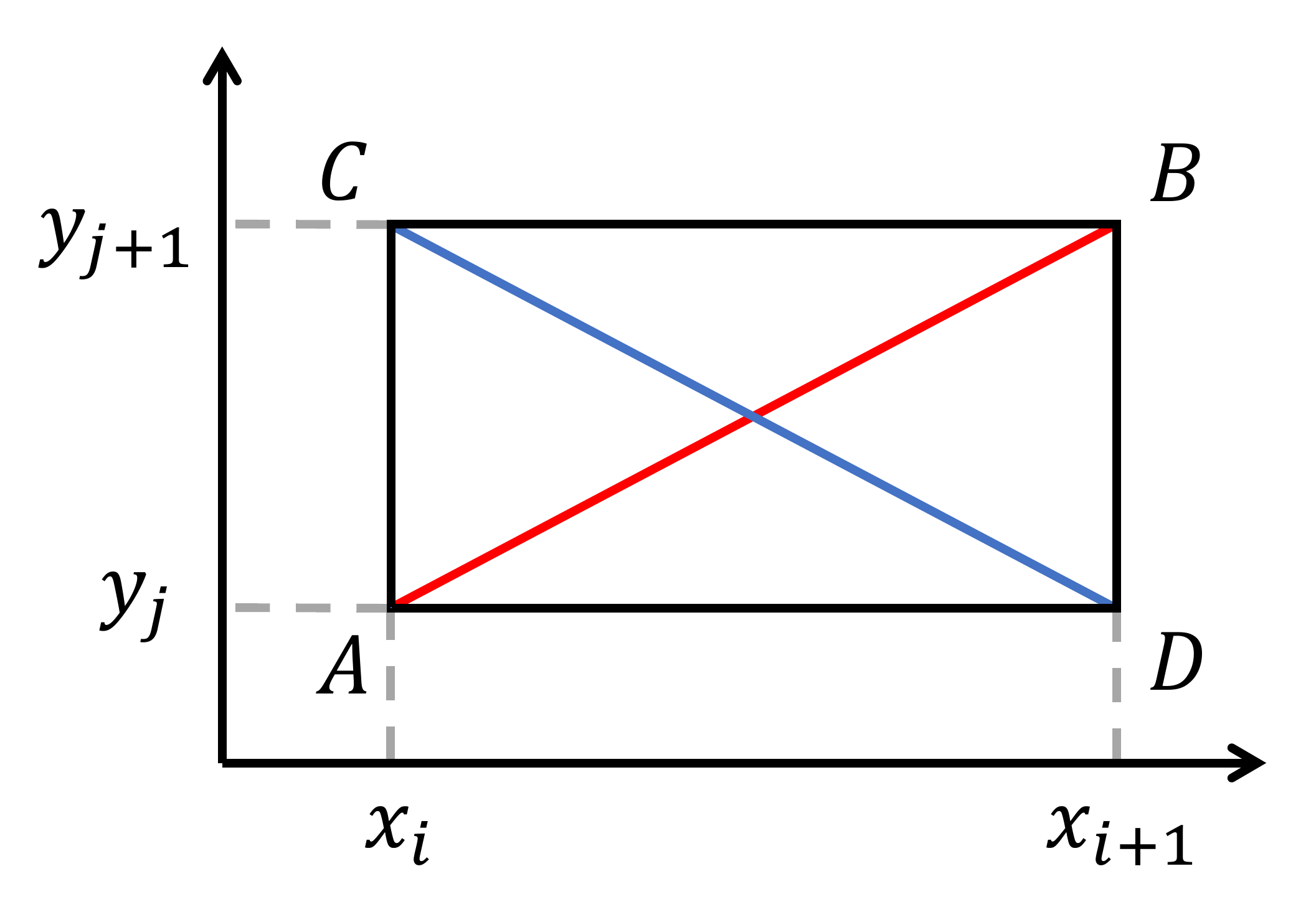}
  }
  \subfigure[\scriptsize{Conservative conditions hold}]{
    \label{fig-diagall-b} 
    \includegraphics[width=0.3\linewidth]{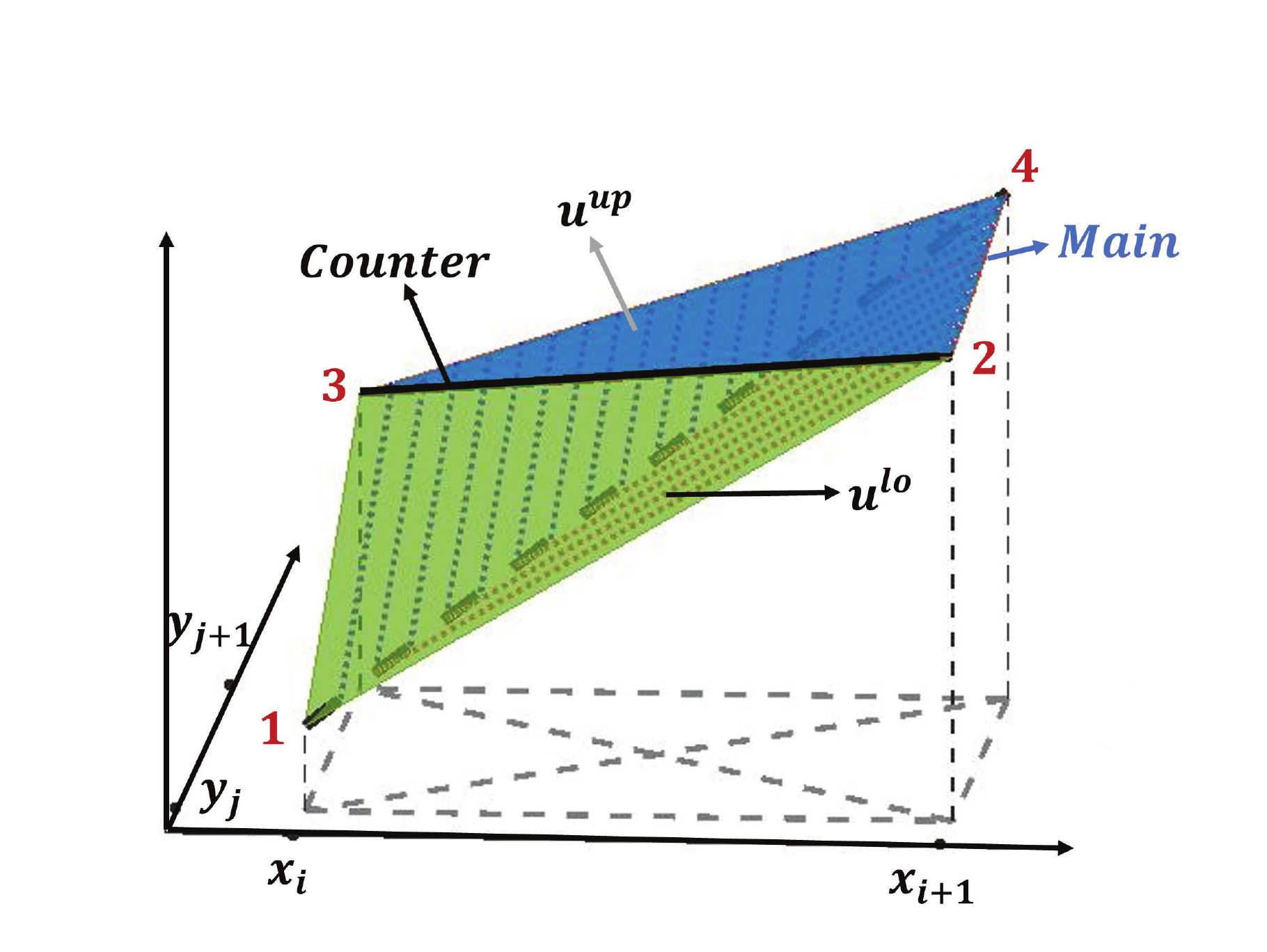}
    }
  \hspace{-1em}
  \subfigure[\scriptsize{Conservative conditions fail}]{
    \label{fig-diagall-c} 
       \includegraphics[width=0.3\linewidth]{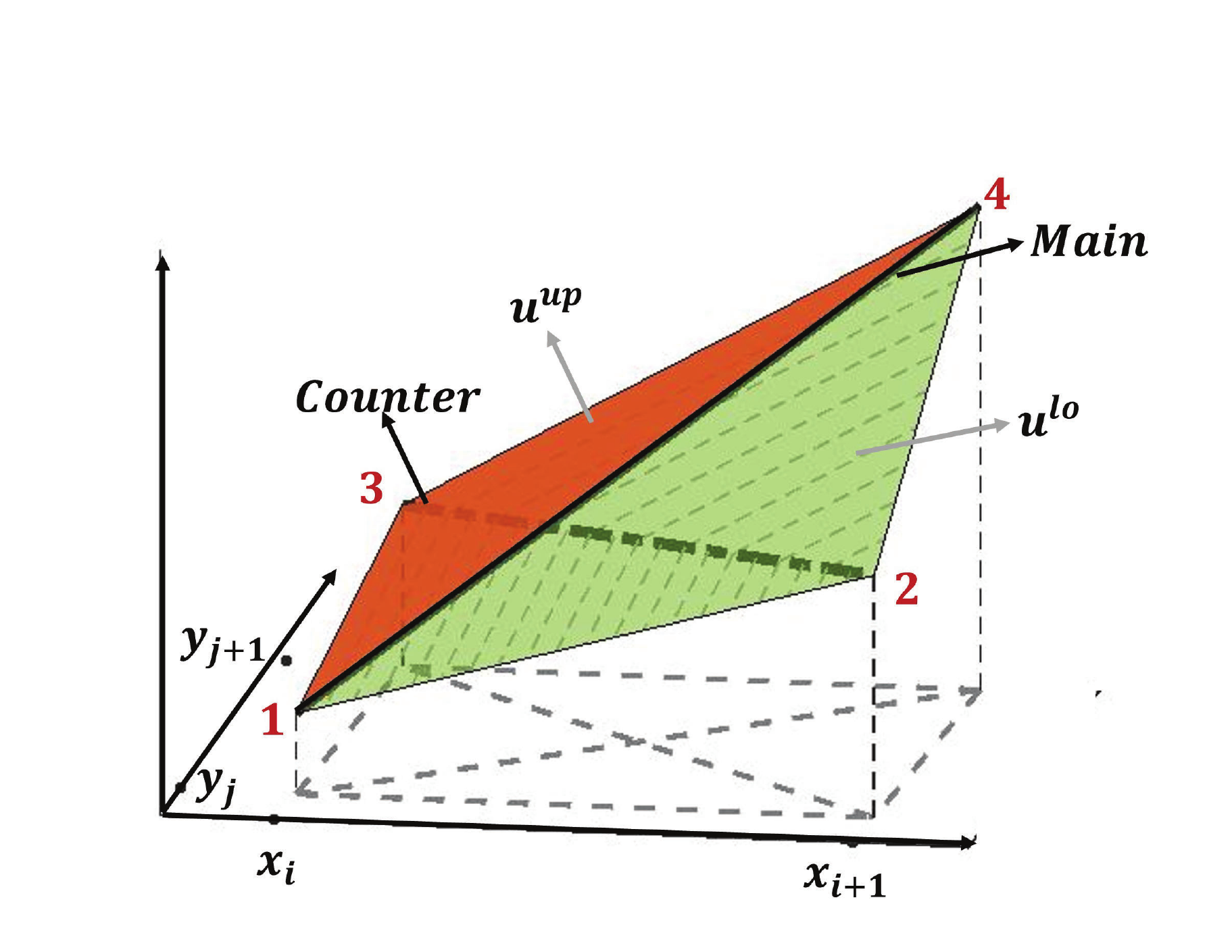}
    }
    \vspace{-1em}
  \captionsetup{font=footnotesize}
  \caption{ 
  (a) The red line is the main diagonal and
  the blue line is the counter diagonal.
  (b) When the conservative condition holds, 
  the graph of the 
  Type-2 piecewise linear function (PLF) (blue and green planes) lies above 
  that of the Type-1 PLF (represented by
dotted lines).
  (c) When the conservative condition fails,
   the graph of the 
  Type-1 
  PLF (orange and green planes) lies above 
  that of the Type-2 PLF (represented by dotted lines).
  }
  \vspace{-0.5cm}
  \label{fig-division} 
\end{figure}

\begin{definition}[Ambiguity set of piecewise linear utility functions]
\label{def-ambi-N}
Let $\scru_N\subset\scru$ be the set of all Type-1 (or Type-2) piecewise linear utility functions 
over $T_{i,j}$
for $i=1,\cdots,N_1-1, j=1,\cdots,N_2-1$.
Define the ambiguity set of piecewise linear utility functions  as
\begin{equation}
\label{eq:U_N-PLA}
    \mathcal{U}_N:=\lt\{ u_N\in \mathscr{U}_N \,:\,
    \int_{T} u_N(x,y) d\psi_l(x,y)
    \leq c_l, \; l=1,\ldots,M  \rt\}.
\end{equation}
\end{definition}

We propose to use $\calu_{N}$ to approximate $\calu$.
Since $\scru_{N}\subset\scru$, then $\calu_{N}\subset\calu$.
Conversely, for any $u\in\calu$, 
we can construct a piecewise linear utility function $u_{N}\in\scru_{N}$ by connecting the utility values at 
gridpoints $(x_i,y_j)$, {\color{black}$i=1,\cdots,N_1$, $j=1,\cdots,N_2$}. 
In general, $u_{N}\notin\calu_{N}$ but the inclusion may hold 
{\color{black}in} some special cases.

\begin{proposition}
\label{prop-uti-N}
Let $\psi_l(x,y)$ be a simple function over $T$ for $l=1,\cdots,M$,
which takes constant values over cells $T_{i,j}$
for $i=1,\cdots,N_1-1, j=1,\cdots,N_2-1$.
Then for any $u\in\calu$, there exists a function $u_{N}\in\scru_{N}$ with $u_N(x_i,y_j)=u(x_i,y_j)$ for $i=1,\ldots,N_1$, $j=1,\ldots,N_2$ such that $u_{N}\in\calu_{N}$. 
Specifically, for $(x,y)\in T$
such $u_{N}$ can be constructed as Type-1 or Type-2 piecewise linear functions defined as
\begin{equation}
\label{eq-utility-N-1}
\begin{split}
    & \inmat{(Type-1)\quad} u_{N}(x,y) = \sum_{i=1}^{N_1-1} \sum_{j=1}^{N_2-1} \1_{T_{i,j}} (x,y) \times \\
    & \lt[
    u^{1u}_{i,j}(x,y) \1_{\lt(\frac{y_{j+1}-y_j}{x_{i+1}-x_i},+\infty\rt)} \lt(\frac{y-y_j}{x-x_i}\rt) 
    + u^{1l}_{i,j}(x,y) \1_{\lt[0,\frac{y_{j+1}-y_j}{x_{i+1}-x_i}\rt]} \lt(\frac{y-y_j}{x-x_i}\rt) \rt]
\end{split}
\end{equation}
or 
\begin{equation}
\label{eq-utility-N-2}
\begin{split}
    & \inmat{(Type-2)\quad}  u_{N}(x,y) = \sum_{i=1}^{N_1-1} \sum_{j=1}^{N_2-1} \1_{T_{i,j}} (x,y) \times \\
    & \lt[u^{2u}_{i,j}(x,y) \1_{\lt[0,\frac{y_{j+1}-y_j}{x_{i+1}-x_i}\rt]} \lt(\frac{y_{j+1}-y}{x-x_i}\rt) +
    u^{2l}_{i,j}(x,y) 
    \1_{\lt(\frac{y_{j+1}-y_j}{x_{i+1}-x_i},+\infty\rt)} \lt(\frac{y_{j+1}-y}{x-x_i}\rt) \rt],
\end{split}
\end{equation}
where $u^{1u}_{i,j}(x,y)$, $u^{1l}_{i,j}(x,y)$, $u^{2u}_{i,j}(x,y)$, and $u^{2l}_{i,j}(x,y)$ are defined as in (\ref{eq-up})-(\ref{eq-lo-2}),
and $\1_A(\cdot)$ denotes the indicator function of set $A$.

\end{proposition}

The proof is deferred to Appendix~\ref{app:proof-uN}.
Using $\calu_{N}$, we propose to solve the BUPRO problem (\ref{eq:MAUT-robust}) 
by solving the following approximate problem: 
\begin{equation}
\label{eq:MAUT-robust-N}
\inmat{(BUPRO-N)} \quad
\vt_N:=\max_{\bdz\in Z}\min_{u_{N}\in\calu_{N}} \bbe_P[u_{N}(\bdf(\bdz,\bdxi))].
\end{equation}

In the rest of the section, we discuss numerical schemes for solving the BUPRO-N problem. 
To this end, we need to restrict 
our discussion to the case that $\bdxi$ is discretely distributed.
\begin{assumption}
\label{assu-discrete}
    $P$ 
    is a discrete distribution with $P(\bdxi=\bdxi^k)=p_k$ for $k=1,\ldots,K$.
\end{assumption}

Under Assumption~\ref{assu-discrete},
we can write the BUPRO-N model
as 
\begin{equation}
\label{eq:MAUT-robust-N-dis}
    \max_{\bdz\in Z}\min_{u_{N}\in\calu_{N}} \sum_{k=1}^K p_k u_{N}(\bdf(\bdz,\bdxi^k)).
\end{equation}
The maximin problem can be decomposed into an inner minimization problem
\begin{equation}
\label{eq:MAUT-robust-N-dis-min}
    v_{N}(\bdz):=\min_{u_{N}\in\calu_{N}} \sum_{k=1}^K p_k u_{N}(\bdf(\bdz,\bdxi^k))
\end{equation}
and an outer maximization problem $ \vt_N=\max_{\bdz\in Z} v_{N}(\bdz)$.
Our strategy is to formulate (\ref{eq:MAUT-robust-N-dis-min}) as an LP 
and solve 
the outer maximization problem
by 
derivative-free (Dfree) methods.
We will discuss 
the performance of PLA 
in Section~\ref{sec-errorbound}.
Note that if $\bdxi$ is continuously distributed, then we may regard (\ref{eq:MAUT-robust-N-dis}) as
a
discrete
approximation to the BUPRO-N model.

We now move on to 
derive the tractable formulation of 
(\ref{eq:MAUT-robust-N-dis-min}) 
when $\scru$ is a class of componentwise non-decreasing and Lipschitz continuous utility functions which is concave in 
each single variate. 

\begin{assumption}
\label{A:concave-in-x-and-y}
For any $u\in\scru$,
the single-variate
utility functions $u(x,\hat{y})$ and $u(\hat{x},y)$ are concave at any instantiations $\hat{y}\in Y$ and $\hat{x}\in X$.
\end{assumption}

The 
single-variate
utility functions 
in Assumption~\ref{A:concave-in-x-and-y} 
can be regarded as non-normalized
single-attribute utility functions.
The concavity condition is used widely  in the literature of expected utility theory, 
which implies the 
DM is risk-averse for each 
individual  attribute (\cite{tsanakas2003risk}).

\begin{assumption}
\label{assu-lip}
Each function $u\in\scru$ is Lipschitz continuous over $T$
with the modulus being bounded by $L$ in the sense that 
\begin{equation}
\label{eq-Lip-condition}
    |u(x,y)-u(x',y')|\leq L \|(x-x',y-y')\|_1 \quad \forall \, (x,y),(x',y')\in T.
\end{equation}
\end{assumption}

The normalization condition and the Lipschitz condition imply that  $L\geq 1/(\bar{x}-\underline{x}+\bar{y}-\underline{y})$. This Lipschitz condition means that the DM's utility change is not drastic at any level of the attributes. It is satisfied
when $u$ is locally Lipschitz continuous over an open set containing $T$.



Notice that in the case when $\psi_l$ is not a simple function for $l=1,\cdots,M$, the LS integrals in (\ref{eq:U_N-PLA}) cannot be calculated directly. Fortunately, we can tackle the issue by swapping the positions
between $u_N$ and $\psi_l$.
Specifically, 
using multivariate integration by parts for the 
LS
integrals (see, e.g., \cite{young1917multiple} and \cite{Ans22}), 
we can rewrite 
ambiguity set 
(\ref{eq:ambiguity_set}) as
\begin{equation}
\label{eq:u_N-int}
\begin{split}
    \calu_N=\left\{u_N \in \scru_N: \int_T \right.
    \hat{\psi}_l(x,y) & d u_N(x,y)  + \int_{X}\psi_{1,l}(x)d u_N(x,\underline{y})\\ 
    &\left. +\int_{Y} \psi_{2,l}(y)d u_N(\underline{x},y) \leq c_l,\;l=1,\ldots,M  \right\},
\end{split}
\end{equation}
 where $\hat{\psi}_l(x,y):=\psi_l(\bar{x},\bar{y})-\psi_l(x,\bar{y})-\psi_l(\bar{x},y)+\psi_l(x,y)$, 
 $\psi_{1,l}(x):=\psi_l(\bar{x},\bar{y})-\psi_l(x,\bar{y})-\psi_l(\bar{x},\underline{y})+\psi_l(x,\underline{y})$, and 
 $\psi_{2,l}(y):=\psi_l(\bar{x},\bar{y})-\psi_l(\underline{x},\bar{y})-\psi_l(\bar{x},y)+\psi_l(\underline{x},y)$ for $l=1,\ldots,M$. 
 Likewise, we can reformulate the ambiguity set ${\cal U}$  defined in (\ref{eq:ambiguity_set})
 as 
 \begin{equation}
\label{eq-equipres}
\begin{split}
    \calu=\left\{u \in \scru: \int_T \right.
    \hat{\psi}_l(x,y) &d u(x,y)  + \int_{X}\psi_{1,l}(x)d u(x,\underline{y})\\
    &\left. +\int_{Y} \psi_{2,l}(y)d u(\underline{x},y) \leq c_l,\;l=1,\ldots,M  \right\}.
\end{split}
\end{equation}
In the case that 
the decision making problem has only one variable (e.g. $y$ disappears),  the first term and the third term at the left hand side of the inequalities 
will disappear and consequently 
the two-dimensional conditions defined as in (\ref{eq:u_N-int}) reduce to the one-dimensional moment-type conditions in \cite{GXZ21}.
The next proposition states how
the two-dimensional LS integrals 
in (\ref{eq:u_N-int})
may be converted into one-dimensional 
Riemann integrals. 
The proof is deferred to Appendix~\ref{app:proof-LS}.

\begin{proposition}
\label{prop-int-pl}
Let 
$F: 
[\underline{a}, \overline{a}]\times [\underline{b}, \overline{b}]
\rightarrow \R$ be a 
continuous 
function.
Assume: 
(a) $F$ is a piecewise linear function
with two linear pieces divided   
by line segment connecting points $A(\underline{a}, \underline{b})$ and
$B(\bar{a}, \bar{b})$;
(b) $\psi$ is
a real-valued measurable function 
w.r.t.~a measure induced by $F$,
and is 
 Riemann integrable 
over the line segment connecting points $A$ and $B$. 
Then 
\begin{equation}
\label{eq-int-pl}
    \int_{\underline{a}, \bar{a}}^{\underline{b}, \bar{b}}\psi(x,y)d F(x,y)=
    \frac{F(\bar{a},\bar{b})-F(\underline{a},\bar{b})-F(\bar{a},\underline{b})+F(\underline{a},\underline{b})}{\bar{a}-\underline{a}}\int_{\underline{a}}^{\bar{a}} \psi(x,y(x))d x, 
\end{equation}
where $y(x)$ is the linear function 
representing the 
segment $AB$.
\end{proposition}

With this, we are ready to
reformulate the inner minimization problem (\ref{eq:MAUT-robust-N-dis-min})
as an LP.

\begin{proposition}
Under Assumptions \ref{assu-original}-\ref{assu-lip}, the inner minimization problem (\ref{eq:MAUT-robust-N-dis-min}) with Type-1 PLA can be reformulated as the following LP:
\begin{subequations}
\label{eq:PRO-N-reformulate}
\begin{align}
    \disp{ \min_{{\bm u}}} \quad
    & \sum_{k=1}^K p_k \sum_{i=1}^{N_1-1} \sum_{j=1}^{N_2-1} \1_{T_{i,j}} (\bdf^k) 
    \lt[ u_{i,j}^{1l}(\bdf^k) \1_{\lt[0,\frac{y_{j+1}-y_j}{x_{i+1}-x_i}\rt]} 
    \lt( \frac{f_2^k-y_j}{f_1^k-x_i} \rt) \rt. \nonumber \\
    & \lt. + u_{i,j}^{1u}(\bdf^k) \1_{\lt( \frac{y_{j+1}-y_j}{x_{i+1}-x_i},+\infty \rt)} \lt(\frac{f_2^k-y_j}{f_1^k-x_i}\rt) \rt] \label{eq-traform-obj} \\
    \st \quad
    & 
    \sum_{i=1}^{N_1-1}\sum_{j=1}^{N_2-1} \frac{u_{i,j+1}-u_{i+1,j+1}-u_{i,j}+u_{i+1,j}}{x_{i+1}-x_i} \int_{x_i}^{x_{i+1}} \hat{\psi}_l (x,y(x)) d x  \nonumber \\
    & + 
    \sum_{i=1}^{N_1-1} \frac{u_{i+1,1}-u_{i,1}}{x_{i+1}-x_i} \int_{x_i}^{x_{i+1}} \psi_{1,l}(x) d x \nonumber \\
    & +
    \sum_{j=1}^{N_2-1} \frac{u_{1,j+1}-u_{1,j}}{y_{j+1}-y_j} \int_{y_j}^{y_{j+1}} \psi_{2,l}(y) d y \leq c_l, l=1,\ldots,M, \label{eq-traform-paircom}\\
    & \frac{u_{i+1,j}-u_{i,j}}{x_{i+1}-x_i} \leq \frac{u_{i,j}-u_{i-1,j}}{x_i-x_{i-1}}, i=2,\ldots,N_1-1, j=1,\ldots,N_2, \label{eq-traform-concave1} \\
    & \frac{u_{i,j+1}-u_{i,j}}{y_{j+1}-y_j} \leq \frac{u_{i,j}-u_{i,j-1}}{y_j-y_{j-1}}, i=1,\ldots,N_1, j=2,\ldots,N_2-1, \label{eq-traform-concave2} \\
    & u_{i+1,j}-u_{i,j}\leq L(x_{i+1}-x_i), i=1,\ldots,N_1-1, j=1,\ldots,N_2, \label{eq-traform-lip1} \\
    & u_{i,j+1}-u_{i,j}\leq L(y_{j+1}-y_j), i=1,\ldots,N_1,j=1,\ldots,N_2-1, \label{eq-traform-lip2} \\ 
    & u_{i+1,j}\geq u_{i,j}, i=1,\ldots,N_1-1,j=1,\ldots,N_2, \label{eq-traform-mon1} \\
    & u_{i,j+1}\geq u_{i,j}, i=1,\ldots,N_1,j=1,\ldots,N_2-1, \label{eq-traform-mon2} \\
    & u_{i,j}+u_{i+1,j+1} \leq
    u_{i,j+1}+u_{i+1,j}, i=1,\ldots,N_1-1,j=1,\ldots,N_2-1, \label{eq-traform-conservative} \\
    & u_{1,1}=0, u_{N_1,N_2}=1, \label{eq-traform-norm1}
\end{align}
\end{subequations}
where {\color{black}${\bm u}:={\rm vec}\left((u_{i,j})_{1\leq i\leq N_1}^{1\leq j\leq N_2}\right)=(u_{1,1},\ldots,u_{N_1,1},\ldots,u_{1,N_2},\ldots,u_{N_1,N_2})^T\in \R^{N_1N_2}$,}
$\bdf^k:=\bdf(\bdz,\\ \bdxi^k)=(f_1^k,f_2^k)$  with $f_1^k:=f_1(\bdz,\bdxi^k)$, $f_2^k:=f_2(\bdz,\bdxi^k)$,
$\hat{\psi}_l$,
$\psi_{1,l}$
$\psi_{2,l}$ are defined as in (\ref{eq:u_N-int}),
$u_{i,j}^{1l}$ and $u_{i,j}^{1u}$ are defined as in (\ref{eq-up}) and (\ref{eq-lo}).
\end{proposition}

\noindent 
\textbf{Proof.}
Using the Type-1 PLA as defined in (\ref{eq-utility-N-1}), we may reformulate the objective as (\ref{eq-traform-obj}).
Moreover,
constraint (\ref{eq-traform-paircom})
represents
the integral inequalities conditions defined as in (\ref{eq:u_N-int}) from
Proposition~\ref{prop-int-pl}.
Constraints (\ref{eq-traform-concave1}) and (\ref{eq-traform-concave2})
characterize
concavity of single-variate utility functions
assumed in
Assumption~\ref{A:concave-in-x-and-y}.
Constraints (\ref{eq-traform-lip1}) and (\ref{eq-traform-lip2}) 
capture
the Lipschitz continuity for the utility function.
Constraints (\ref{eq-traform-mon1}) and (\ref{eq-traform-mon2}) 
reflect
componentwise monotonicity of utility functions.
Constraint (\ref{eq-traform-conservative}) 
states
the conservative property.
Constraint (\ref{eq-traform-norm1})
is
the normalization condition for the utility function.
\hfill \Box

\begin{remark}
\label{rem:BUPRO-DF}
(i)
 Note that (\ref{eq:PRO-N-reformulate}) is reformulated based on the Type-1 PLA.
A similar formulation can be obtained for Type-2 PLA. 
By solving (\ref{eq:PRO-N-reformulate}),
we can 
obtain the worst-case utility function $u_N^{\rm worst}$.
The information on $u_N^{\rm worst}$ 
gives us 
a guidance 
as to how the inner minimization problem 
approximates the true expected utility.
The problem size depends on the number of gridpoints and is independent of 
the scenarios of $\bdxi$.

Note also that the single-attribute utility functions are assumed to be concave in Assumption~\ref{A:concave-in-x-and-y}.
This is in accordance with single-attribute decision making in the risk-averse case.
Likewise, we can also assume that one (or both) single-attribute utility at any instantiation is (are) convex.
In that case, it suffices to input constraints (\ref{eq-traform-concave1}) or (and) (\ref{eq-traform-concave2}) in the reverse direction.

The Lipschitz continuity is also reflected by the Type-1 PLA  (\ref{eq-utility-N-1}) which can be formulated as
\begin{equation*}
    u^{1u}_{i,j}(x,y) = \lt( \frac{u_{i+1,j+1}-u_{i,j+1}}{x_{i+1}-x_i}, \frac{u_{i,j+1}-u_{i,j}}{y_{j+1}-y_j} \rt) (x,y)^T +b^{1u}_{i,j}
\end{equation*}
and
\begin{equation*}
    u^{1l}_{i,j}(x,y) = \lt( \frac{u_{i+1,j}-u_{i,j}}{x_{i+1}-x_i}, \frac{u_{i+1,j+1}-u_{i+1,j}}{y_{j+1}-y_j} \rt) (x,y)^T+b^{1l}_{i,j},
\end{equation*}
where $b^{1u}_{i,j}$
and $b^{1l}_{i,j}$ are constants representing the 
intercepts respectively.
The Lipschitz continuity defined as in Assumption~\ref{assu-lip} corresponds to 
$$
\max\lt\{ \lt\|\lt( \frac{u_{i+1,j+1}-u_{i,j+1}}{x_{i+1}-x_i}, \frac{u_{i,j+1}-u_{i,j}}{y_{j+1}-y_j} \rt)\rt\|_{\infty},\lt\|\lt( \frac{u_{i+1,j}-u_{i,j}}{x_{i+1}-x_i}, \frac{u_{i+1,j+1}-u_{i+1,j}}{y_{j+1}-y_j} \rt)\rt\|_{\infty}  \rt\}\leq L,
$$ 
{\color{black} over each cell $T_{i,j}$,} which implies constraints (\ref{eq-traform-lip1}) and (\ref{eq-traform-lip2}).

(ii)
The aforementioned PLA utility function $u_N$ is constructed either in Type-1 or in Type-2.
It is possible to allow both.
Specifically, we can define 
\begin{equation*}
\begin{split}
    & u_N (x,y)  =\sum_{i=1}^{N_1-1}\sum_{j=1}^{N_2-1} \1_{T_{i,j}}(x,y) \times \\
    & \lt[
    h_{i,j}
    \lt(
    u^{1u}_{i,j}(x,y) \1_{\lt(\frac{y_{j+1}-y_j}{x_{i+1}-x_i},+\infty\rt)} \lt(\frac{y-y_j}{x-x_i}\rt) 
    + u^{1l}_{i,j}(x,y) \1_{\lt[0,\frac{y_{j+1}-y_j}{x_{i+1}-x_i}\rt]} \lt(\frac{y-y_j}{x-x_i}\rt) \rt) \rt.\\
    & \lt. + (1-h_{i,j})
    \lt(u^{2u}_{i,j}(x,y) \1_{\lt[0,\frac{y_{j+1}-y_j}{x_{i+1}-x_i}\rt]} \lt(\frac{y_{j+1}-y}{x-x_i}\rt) +
    u^{2l}_{i,j}(x,y) 
    \1_{\lt(\frac{y_{j+1}-y_j}{x_{i+1}-x_i},+\infty\rt)} \lt(\frac{y_{j+1}-y}{x-x_i}\rt) \rt)
    \rt],
\end{split}
\end{equation*}
where 
$\{h_{i,j}, i=1,\ldots,N_1-1,j=1,\ldots,N_2-1\}$ is a set of binary variables 
taking values $0$ or $1$.
In the case that $h_{i,j}=1$, Type-1 PLA is invoked over $T_{i,j}$.
Otherwise, Type-2 PLA is active.
Obviously, this approach significantly extends 
the class of piecewise linear utility functions 
and consequently the optimal value of the inner minimization problem is smaller than that of Type-1 and Type-2.
With regard to the tractable formulation, we will have $(N_1-1)(N_2-1)$ additional binary variables 
and the inner minimization becomes an MILP.

(iii) 
Type-1 PLA
$u_N(\bdf(\bdz,\bdxi^k))$
can be alternatively represented in the following form:
\begin{equation*}
    u_{N}(\bdf(\bdz,\bdxi^k)) = \sum_{i=1}^{N_1-1}\sum_{j=1}^{N_2-1}  \left[ \alpha_{i,j}^k(\bdz) u_{i,j}
    +\alpha_{i,j+1}^k(\bdz) u_{i,j+1}
    +\alpha_{i+1,j}^k(\bdz) u_{i+1,j}
    +\alpha_{i+1,j+1}^k(\bdz) u_{i+1,j+1}\right],
\end{equation*}
where
\begin{align*}
    &\alpha_{i,j}^k(\bdz):=\mathds{1}_{T_{i,j}}(\bdf^k) \lt[ 
     \frac{y_{j+1}-f_2^k}{y_{j+1}-y_j} \mathds{1}_{\left[0,\frac{y_{j+1}-y_j}{x_{i+1}-x_i}\right]}\left(\frac{f_2^k-y_j}{f_1^k-x_i} \right) 
    +  \frac{x_{i+1}-f_1^k}{x_{i+1}-x_i} \mathds{1}_{\left(\frac{y_{j+1}-y_j}{x_{i+1}-x_i},\infty \right)}  \left(\frac{f_2^k-y_j}{f_1^k-x_i}\right) \rt], 
    \nonumber \\
    &\alpha_{i,j+1}^k(\bdz):= \left(\frac{f_2^k-y_j}{y_{j+1}-y_j}-\frac{f_1^k-x_i}{x_{i+1}-x_i}\right) \mathds{1}_{T_{i,j}}(\bdf^k) \mathds{1}_{\left[0,\frac{y_{j+1}-y_j}{x_{i+1}-x_i}\right]} \left(\frac{f_2^k-y_j}{f_1^k-x_i}\right) , \\
    & \alpha_{i+1,j}^k(\bdz):= \left( \frac{f_1^k-x_i}{x_{i+1}-x_i}-\frac{f_2^k-y_j}{y_{j+1}-y_j}\right) \mathds{1}_{T_{i,j}}(\bdf^k) \mathds{1}_{\left(\frac{y_{j+1}-y_j}{x_{i+1}-x_i},\infty \right)}\left(\frac{f_2^k-y_j}{f_1^k-x_i}\right), \nonumber \\
    & \alpha_{i+1,j+1}^k(\bdz):=\mathds{1}_{T_{i,j}} (\bdf^k) \lt[  \frac{f_1^k-x_i}{x_{i+1}-x_i} \mathds{1}_{\left[0,\frac{y_{j+1}-y_j}{x_{i+1}-x_i}\right]}\left(\frac{f_2^k-y_j}{f_1^k-x_i}\right) 
    +  \frac{f_2^k-y_j}{y_{j+1}-y_j} \mathds{1}_{\left(\frac{y_{j+1}-y_j}{x_{i+1}-x_i},\infty \right)}\left(\frac{f_2^k-y_j}{f_1^k-x_i}\right) \rt]. \nonumber
\end{align*}
For fixed $\bdz$,
the inner minimization problem (\ref{eq:MAUT-robust-N-dis-min}) 
is also an LP with this $u_N$.
Let $L({\bm u},{\bm \lambda};\bdz)$ be the Lagrange function of the inner problem
and ${\bm \lambda}$ be the vector of Lagrange multipliers.
Then the inner problem can be reformulated as 
$\min_{{\bm u}}\max_{{\bm \lambda}} L({\bm u},{\bm \lambda};\bdz)$.
In this way, we can reformulate the maximin problem (\ref{eq:PRO-N-reformulate}) as a single maximization problem
$
\max_{\bdz\in Z,{\bm \lambda}} \{\min_{{\bm u}}L({\bm u},{\bm \lambda};\bdz)\}.
$
Unfortunately, 
this is not helpful 
since the coefficients of
$u_{i,j}$,
$u_{i,j+1}$,
$u_{i+1,j}$,
$u_{i+1,j+1}$
are composed of  indicator functions of 
$\bdf^k$.
In the next section,
we will propose a new approach to 
handle the issue 
and extend the discussions to 
the multi-attribute case.

\end{remark}

\section{Implicit piecewise linear approximation of UPRO -- from bi-attribute to multi-attribute case}
\label{sec:multi-atrribute}

{\color{black}
In this section,
we look into the PLA approach from 
a slightly different perspective: instead of deriving an explicit form of 
piecewise linear function as we discussed in
the previous section, we propose to use
the well-known polyhedral method
(see e.g. \cite{LeW01,KDN04,DLM10,VAN10,vielma2015mixed}),
where the PLA
function 
at each cell is implicitly determined by solving a minimization or a maximization program.
There are two 
advantages 
for 
doing this.
One is that the implicit approach allows us to
extend the PLA for UPRO problem from bi-attribute decision making problems to the multi-attribute case and
this would be extremely 
complex under the   
explicit PLA framework. 
The other 
is that 
the implicit approach enables us to
reformulate 
the approximate UPRO problem
into a single MILP
when ${\bm f}(\bdz,\bdxi)$ is linear in $\bdz$.
}

\subsection{Bi-attribute case}
\label{sec:two-dim-u}

Inspired by the polyhedral method, we 
can obtain the coefficients $\alpha_{i,j}$ of $u_N(\bdf(\bdz,\bdxi^k))$ in terms of $u_{i,j}$ 
under Type-1 PLA
in Remark~\ref{rem:BUPRO-DF}~(iii)
by solving a system of linear equalities and inequalities:
\begin{subequations}
\begin{align}
    & \sum_{i=1}^{N_1} \sum_{j=1}^{N_2} \alpha_{i,j}^{k}=1,\; k=1,\cdots,K,
    \label{eq:mixed-integer-R2-b}\\
    & \sum_{i=1}^{N_1} \sum_{j=1}^{N_2} \alpha_{i,j}^{k} x_{i}=f_1^k,\;\; \sum_{i=1}^{N_1} \sum_{j=1}^{N_2} \alpha_{i,j}^{k} y_{j}= f_2^k,\;\; k=1,\cdots,K,\label{eq:mixed-integer-R2-c}\\
    & \sum_{i=1}^{N_1-1} \sum_{j=1}^{N_2-1} \lt(h_{i,j,k}^{u}+h_{i,j,k}^{l}\rt)=1,\;\;
    k=1,\cdots,K,\label{eq:mixed-integer-R2-d}\\
    & {\bm h}^u_k,{\bm h}^l_k\in \{0,1\}^{(N_1-1)(N_2-1)},\;k=1,\cdots,K,
    \label{eq:mixed-integer-R2-e}\\
    &0\leq \alpha_{i,j}^{k}\leq 
    h_{i,j,k}^{u}
    +h_{i,j,k}^{l} 
    +h_{i,j-1,k}^{u}
    +h_{i-1,j-1,k}^{l}
    +h_{i-1,j-1,k}^{u}
    +h_{i-1,j,k}^{l},\nonumber \\
    & \qquad \qquad \qquad  \qquad \qquad \quad
i=1,\cdots,N_1,\; j=1,\cdots,N_2, \; k=1,\cdots,K,
\label{eq:mixed-integer-R2-f}
\end{align}
\end{subequations}
where {\color{black} ${\bm h}^u_k:= {\rm vec}\left((h^u_{i,j,k})_{1\leq i\leq N_1}^{1\leq j\leq N_2}\right)$, ${\bm h}^l_k:= {\rm vec}\left((h^l_{i,j,k})_{1\leq i\leq N_1}^{1\leq j\leq N_2}\right)$ for $k=1,\ldots,K$,}
${\bm f}(\bdz,\bdxi^k)=(f_1^k,f_2^k)^T$ with $f_i^k:=f_i(\bdz,\bdxi^k)$ for $i=1,2$,
$h_{0,*,*}^*=h_{*,0,*}^*=h_{N_1,*,*}^*=h_{*,N_2,*}^*=0$.
Here $*$ represents all indexes possibly taken at the subscripts and superscripts.
Constraint (\ref{eq:mixed-integer-R2-b}) and $\alpha_{i,j}^k\geq 0$
result from the coefficients of the convex combinations of $u_{i,j}$
for $u_N({\bm f}(\bdz,\bdxi^k))$.
Constraint (\ref{eq:mixed-integer-R2-c}) 
arises because 
the linearity of $u_N$ over $T_{i,j}$ guarantees that
the convex combination coefficients of ${\bm f}(\bdz,\bdxi^k)$ and $u_N({\bm f}(\bdz,\bdxi^k))$ are 
identical.
Since $h_{i,j,k}^u, h_{i,j,k}^l\in \{0,1\}$,
constraint (\ref{eq:mixed-integer-R2-d}) imposes 
a restriction 
that only one is used 
for the convex combination among all triangles.
The constraint (\ref{eq:mixed-integer-R2-f}) imposes that the
only nonzero $\alpha_{i,j}$ 
can be those associated with the
three vertices of a such triangle.
For example, if $h_{i,j,k}^l=1$, then 
$\bdf(\bdz,\bdxi^k)$ lies in the
lower triangle
of the cell $T_{i,j}$. This is 
indicated by the fact that
$\alpha_{i,j}^k\leq h_{i,j,k}^l=1$,
$\alpha_{i+1,j+1}^k\leq h_{i,j,k}^l=1$,
$\alpha_{i+1,j}^l\leq h_{i,j,k}^l=1$,
and $\alpha_{i',j'}^k=0$ for $(i',j')\notin \{(i,j),(i+1,j+1),(i+1,j)\}$,
see Figure~\ref{fig:Type1IPLA} where
the six triangles are related to point $(x_i,y_j)$ 
and we indicate the corresponding binary variables $h_{i,j,k}^u$ and $h_{i,j,k}^l$ in each triangle to facilitate readers understanding.
Consequently,
under Assumption~\ref{assu-discrete},
we can formulate
the bi-attribute utility maximization problem  
$\max_{\bdz\in Z}\sum_{k=1}^Kp_k[u_N({\bm f(\bdz,\bdxi^k)})]$ 
 as
an MIP:
\begin{subequations}
\label{eq:mixed-integer-R2-2}
\begin{align}
 \max\limits_{\bdz\in Z,{\bm \alpha},{\bm h^l},{\bm h}^u} \; &  \sum_{k=1}^K p_k \sum_{i=1}^{N_1} \sum_{j=1}^{N_2} \alpha_{i,j}^{k} u_{i,j}\\
 {\rm s.t. }\quad\;\;\; & \inmat{constraints } (\ref{eq:mixed-integer-R2-b})-
 (\ref{eq:mixed-integer-R2-f}),
\end{align}
\end{subequations}
where ${\bm \alpha}:=({\bm \alpha}^1,\cdots,{\bm \alpha}^K)\in \R^{(N_1N_2)\times K}$,
${\bm \alpha}^{k}:={\rm vec}\left((\alpha_{i,j}^k)_{1\leq i\leq N_1}^{1\leq j\leq N_2}\right)$
for $k=1,\cdots,K$,
${\bm h}^l:=({\bm h}_1^l,\cdots,{\bm h}_K^l)\in \R^{(N_1-1)(N_2-1)\times K}$,
${\bm h}^u:=({\bm h}_1^u,\cdots,{\bm h}_K^u)\in \R^{(N_1-1)(N_2-1)\times K}$.
If $f(\bdz,\bdxi)$ is linear in $\bdz$, then the problem (\ref{eq:mixed-integer-R2-2}) is an MILP.
This idea can be applied to the BUPRO-N model.
To ease the exposition, we consider 
the case that the ambiguity set is constructed 
by pairwise comparisons, 
that is,
$\psi_l=F_{{\bm B}_l}-F_{{\bm A}_l}$,
and
\begin{equation*}
\begin{split}
    {\cal U}_N
    &=\left\{u_N\in \mathscr{U}_N:
    \int_{T} u_N(x,y) d (F_{\bdcb_l}(x,y)-F_{\bdca_l}(x,y))\leq 0,\;l=1,\cdots,M\right\}. 
\end{split}
\end{equation*}
Under Assumption~\ref{assu-lip},
suppose  the set of gridpoints $\{(x_i,y_j):i=1,\cdots,N_1,j=1,\cdots,N_2\}$ contains all the outcomes of lotteries $ {\bm A}_l$ and ${\bm B}_l$ for $l=1,\cdots,M$,
then we can reformulate 
BUPRO-N problem
(\ref{eq:MAUT-robust-N-dis})
as: 
\begin{subequations}
\label{eq:PRO_MILP_eqi}
\begin{align}
\max\limits_{\bdz\in Z,{\bm \alpha}, {\bm h}^u, {\bm h}^l} 
\min\limits_{{\bm u}} \;\;&
\sum_{k=1}^K p_k \sum_{i=1}^{N_1} \sum_{j=1}^{N_2} \alpha_{i,j}^k u_{i,j}\\
{\rm s.t.} \;\; &
 \sum_{i=1}^{N_1} \sum_{j=1}^{N_2} (\mathbb{P}({\bm B}_l=(x_i,y_j))- \mathbb{P}({\bm A}_l=(x_i,y_j)))u_{i,j}
   \leq 0, \nonumber \\
   & \hspace{16em} l=1,\cdots,M,
   \label{eq:PRO_MILP_mina2-c}\\
  & \inmat{constraints } (\ref{eq-traform-lip1})-(\ref{eq-traform-norm1}),  \\
 & \inmat{constraints } (\ref{eq:mixed-integer-R2-b})
  -(\ref{eq:mixed-integer-R2-f}),
\end{align}
\end{subequations}
where constraints (\ref{eq-traform-lip1})-(\ref{eq-traform-norm1}) characterize
the restrictions on ${\bm u}$
and constraints (\ref{eq:mixed-integer-R2-b})
-(\ref{eq:mixed-integer-R2-f}) 
stipulate the coefficients of ${\bm u}$ implicitly as discussed earlier.
{\color{black}Problem (\ref{eq:PRO_MILP_eqi}) is equivalent to problem (\ref{eq:PRO-N-reformulate}) without constraints (\ref{eq-traform-concave1})-(\ref{eq-traform-concave2}) and with $\psi_l=F_{{\bm B}_l}-F_{{\bm A}_l}$, $l=1,\cdots,M$.}
  It is possible to change the maximization
  w.r.t. ${\bm \alpha}$, ${\bm h}^l$  and ${\bm h}^u$ into 
  minimization. The next proposition explains this.


\begin{proposition}
The BUPRO-N problem (\ref{eq:PRO_MILP_eqi}) is equivalent to 
\begin{subequations}
\label{eq:PRO_MILP_mina2}
\begin{align}
 \displaystyle \max_{\bdz\in Z}  \displaystyle \min_{{\bm \alpha},{\bm h}^u,{\bm h}^l,
 {\bm u}} \;\;&
\sum_{k=1}^K p_k \sum_{i=1}^{N_1} \sum_{j=1}^{N_2} \alpha_{i,j}^{k} u_{i,j}\\
 \inmat{s.t.}\quad\;\; & \inmat{constraints } (\ref{eq-traform-lip1})-(\ref{eq-traform-norm1}),\; 
   (\ref{eq:PRO_MILP_mina2-c}), \label{eq:PRO_MILP_mina2-b}\\
    & \inmat{constraints } (\ref{eq:mixed-integer-R2-b})
  -(\ref{eq:mixed-integer-R2-f}).
\end{align}
\end{subequations}
\end{proposition}
\noindent{\bf Proof.}
We begin by 
writing part of the outer maximization (w.r.t. ${\bm \alpha},{\bm h}^u,{\bm h}^l$)
and  the inner minimization
problem of
(\ref{eq:PRO_MILP_eqi}) as
\begin{subequations}
\label{eq:PRO_MILP_mina}
\begin{align}
 \displaystyle \max_{{\bm \alpha},{\bm h}^u,{\bm h}^l} \; &  \min_{{\bm u}}\left\{
\sum_{k=1}^K p_k \sum_{i=1}^{N_1} \sum_{j=1}^{N_2} \alpha_{i,j}^{k} u_{i,j}:  (\ref{eq-traform-lip1})-(\ref{eq-traform-norm1}) \right \}\\
{\rm s.t.} \quad &
\inmat{\,constraints } (\ref{eq:mixed-integer-R2-b})
 - (\ref{eq:mixed-integer-R2-f}),
 (\ref{eq:PRO_MILP_mina2-c}). \label{eq:PRO_MILP_mina-b}
\end{align}
\end{subequations}
Since 
the representation of point $\bdf(\bdz,\bdxi^k)$ by the convex combination of the vertices of a simplex is unique, 
the feasible set of the outer 
maximization problem (\ref{eq:PRO_MILP_mina})
(specified by
(\ref{eq:PRO_MILP_mina-b}))
is a singleton 
for each fixed $\bdz$.
Thus we can replace operation ``$\max_{{\bm \alpha},{\bm h}^u,{\bm h}^l}$'' with ``$\min_{{\bm \alpha},{\bm h}^u,{\bm h}^l}$'' without affecting the optimal value and the optimal solutions of (\ref{eq:PRO_MILP_mina}).
The replacement effectively reduces 
(\ref{eq:PRO_MILP_eqi}) to (\ref{eq:PRO_MILP_mina2}). \hfill $\Box$

Note that the outer maximization problem (\ref{eq:PRO_MILP_mina2}) can be solved by a Dfree method,
where the inner problem can be seen as an MILP when ${\bm f}(\bdz,\bdxi^k)$ is linear in $\bdz$.

\begin{remark}
\begin{itemize}
\item[(i)] Inequality (\ref{eq:mixed-integer-R2-f}) corresponds to Type-1 PLA.
For Type-2 case
(see Figure~\ref{fig:Type2IPLA}),
we can replace (\ref{eq:mixed-integer-R2-f}) by 
\begin{equation}
\label{eq:constraint-alpha}
\begin{split}
0\leq \alpha_{i,j}^k\leq  h^u_{i-1,j,k}  + h^l_{i-1,j,k} & +  h^u_{i-1,j-1,k}+h^l_{i,j-1,k}+h^u_{i,j-1,k}+h^l_{i,j,k}, \\
    & i=1,\cdots,N_1,j=1,\cdots,N_2,k=1,\cdots,K.
\end{split}
\end{equation}

\item[(ii)] 
For the mixed-type PLA,
we 
can also obtain the coefficients $\alpha_{i,j}$ of $u_N(\bdf(\bdz,\bdxi^k))$ in terms of $u_{i,j}$ 
by solving a system of linear equalities and inequalities:
\begin{subequations}
\begin{align}
    & \sum_{i=1}^{N_1} \sum_{j=1}^{N_2} \alpha_{i,j}^{k\tau}=1,\;\tau=1,2,\;k=1,\cdots,K,
    \label{eq:mixed-integer-Type2-b}\\
    & \sum_{i=1}^{N_1} \sum_{j=1}^{N_2} \alpha_{i,j}^{k\tau} x_{i}=f_1^k,\;\; \sum_{i=1}^{N_1} \sum_{j=1}^{N_2} \alpha_{i,j}^{k\tau} y_{j}= f_2^k,\;\; k=1,\cdots,K,\;\tau =1,2, \label{eq:mixed-integer-Type2-c}\\
    &  \sum_{i=1}^{N_1-1}\sum_{j=1}^{N_2-1} \lt(h_{i,j,k}^{1u}+h_{i,j,k}^{1l}+h_{i,j,k}^{2u}+h_{i,j,k}^{2l}\rt)=1,\;k=1,\cdots,K,
    \label{eq:mixed-integer-Type2-d}\\
    & {\bm h}^{\tau u}_{k},{\bm h}^{\tau l}_{k}\in \{0,1\}^{(N_1-1)(N_2-1)},\;k=1,\cdots,K,\;\tau=1,2,
    \label{eq:mixed-integer-Type2-e}\\
    &0\leq \alpha_{i,j}^{k1}\leq 
    h_{i,j,k}^{1u}
    +h_{i,j,k}^{1l} 
    +h_{i,j-1,k}^{1u}
    +h_{i-1,j-1,k}^{1l }
    +h_{i-1,j-1,k}^{1u }
    +h_{i-1,j,k}^{1l},\nonumber \\
    & \qquad \qquad \qquad  \qquad
i=1,\cdots,N_1,\; j=1,\cdots,N_2, \; k=1,\cdots,K,\\
&0\leq \alpha_{i,j}^{k2}\leq  h^{2u}_{i-1,j,k}  + h^{2l}_{i-1,j,k} +  h^{2u}_{i-1,j-1,k}+h^{2l}_{i,j-1,k}+h^{2u}_{i,j-1,k}+h^{2l}_{i,j,k}, \nonumber \\
& \qquad \qquad \qquad  \qquad i=1,\cdots,N_1,j=1,\cdots,N_2,k=1,\cdots,K,
\label{eq:mixed-integer-Type2}
\end{align}
\end{subequations}
where variables $\alpha_{i,j}^{k1}$, $h_{i,j,k}^{1u}$, $h_{i,j,k}^{1l}$ represent the Type-1 PLA case,
and $\alpha_{i,j}^{k2}$, $h_{i,j,k}^{2u}$, $h_{i,j,k}^{2l}$ 
represent the Type-2 case.
The constraint (\ref{eq:mixed-integer-Type2-d})
indicates that only one type partition is used for each cell.
\end{itemize}
\end{remark}

\vspace{-2em}
\begin{figure}[!hbpt]
  \centering
   \subfigure[]
   {
  \label{fig:Type1IPLA}
    \includegraphics[width=0.4\linewidth]{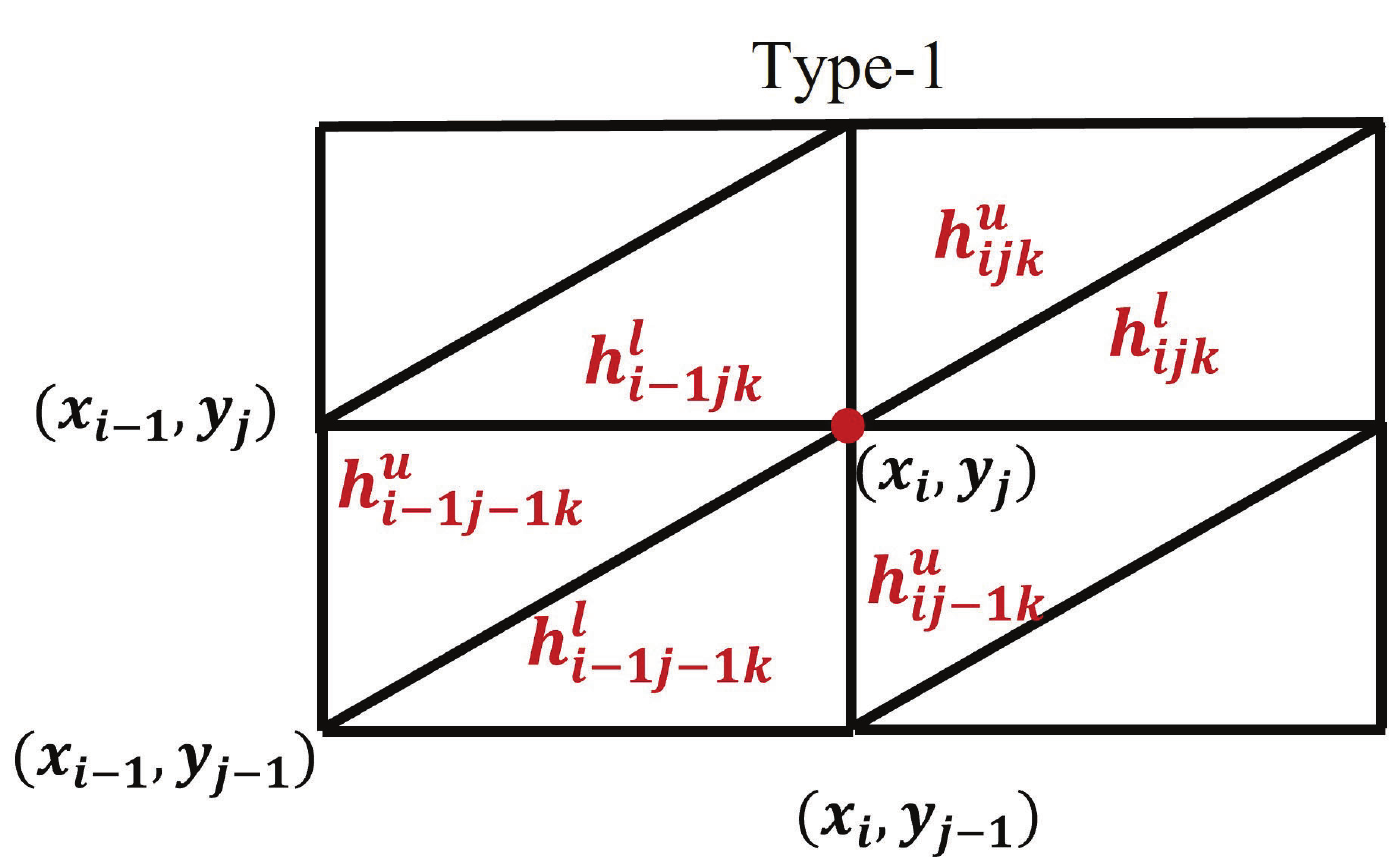}
    \hspace{0.2cm}
  }
   \subfigure[]
   {
    \label{fig:Type2IPLA} 
    \includegraphics[width=0.4\linewidth]{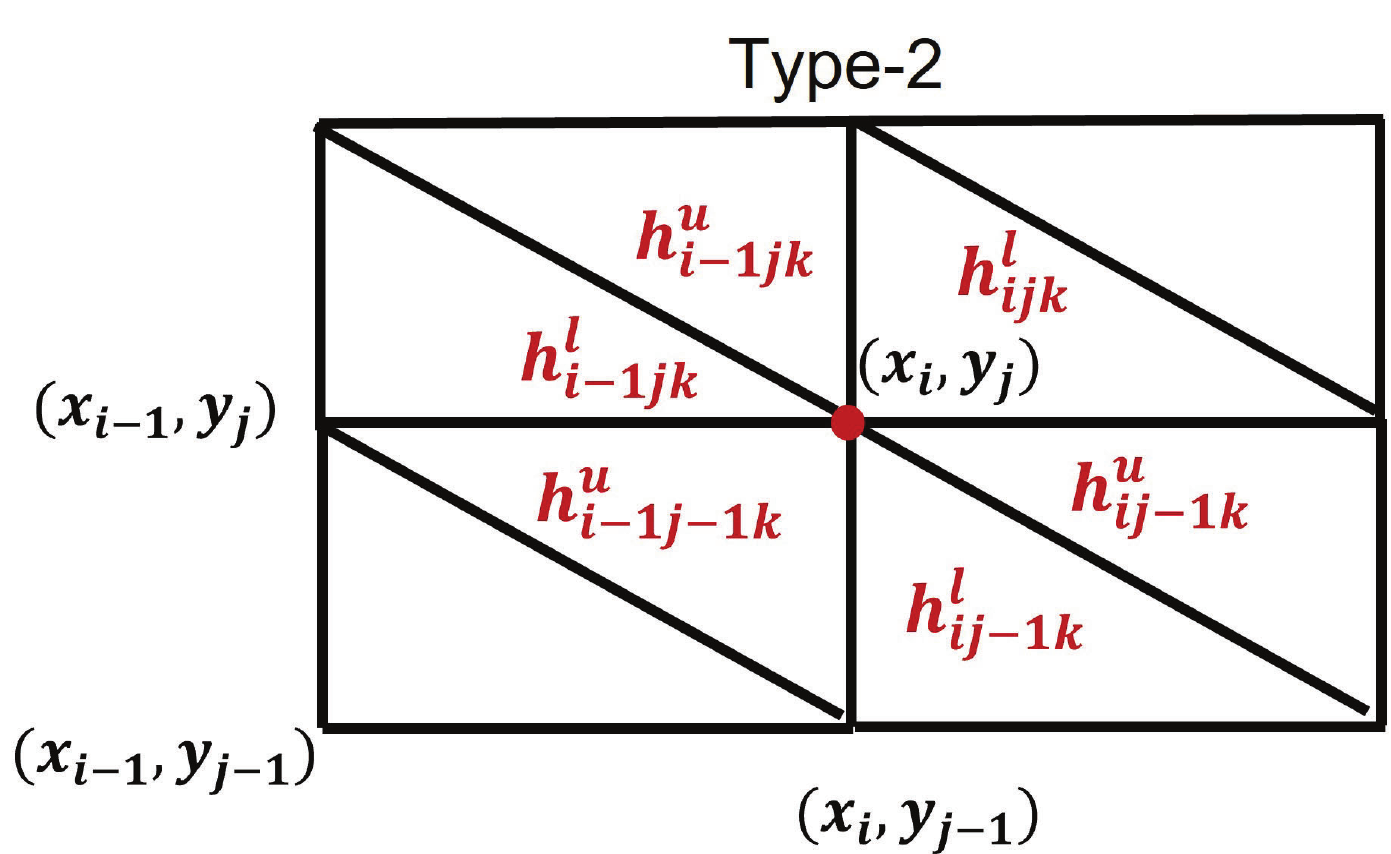}
  }
  \vspace{-1em}
  \captionsetup{font=footnotesize}
  \caption{\footnotesize 
  (a) \& (b) represent the bi-attribute case 
  over $T_{i,j}$.
  They show the six triangles related to point $(x_i,y_j)$.
  }
  \vspace{-0.5cm}
  \label{fig:alpha_ij} 
\end{figure}


\subsubsection{Single mixed-integer reformulation of \texorpdfstring{(\ref{eq:PRO_MILP_eqi})}{}}

By deriving the Lagrange dual of the inner minimization problem of (\ref{eq:PRO_MILP_eqi})
which is established under Type-1 PLA, we can recast the maximin problem as a 
single MILP when ${\bm f}(\cdot,\bdxi)$ is linear.

\begin{proposition}
[Reformulation of (\ref{eq:PRO_MILP_eqi})]
\label{prop:single-MILP}
Problem (\ref{eq:PRO_MILP_eqi})
can be reformulated as a single MILP when ${\bm f}(\bdz,\bdxi)$ is linear in $\bdz$,
\begin{subequations}
\label{eq:PRO_MILP_single}
\hspace{-0.5cm}
\begin{align}
\displaystyle \max_{\substack{{\bm z}\in Z, {\bm \alpha},
{\bm h}^u, {\bm h}^l\\
{\bm \lambda}^1,
{\bm \lambda}^2,
{\bm \eta}^1\\
{\bm \eta}^2,
{\bm \tau}, {\bm \zeta}}} \; &  -\sum_{i=1}^{N_1-1}\sum_{j=1}^{N_2}\eta_{i,j}^1L(x_{i+1}-x_i)
-\sum_{i=1}^{N_1}\sum_{j=1}^{N_2-1} \eta_{i,j}^2 L(y_{j+1}-y_j) 
+ \sum_{k=1}^Kp_k \alpha_{N_1,N_2}^k
\nonumber  \\
& 
-\lambda_{N_1-1,N_2}^1
-\lambda_{N_1,N_2-1}^2
+\eta_{N_1-1,N_2}^1
+\eta_{N_1,N_2-1}^2  +\tau_{N_1-1,N_2-1}
 +{\bm \zeta}^T {\bm Q}_{N_1,N_2} \\
{\rm s.t.} \quad\;\;\,  
&
\sum_{k=1}^Kp_k\alpha_{i,j}^k
 +\lambda_{i,j}^1-\lambda_{i-1,j}^1+\lambda_{i,j}^2-\lambda_{i,j-1}^{2}+\eta_{i-1,j}^1 -\eta_{i,j}^1+\eta_{i,j-1}^2-\eta_{i,j}^2 
 \nonumber\\
 & 
 \quad  + \tau_{i,j}+\tau_{i-1j-1}
 -\tau_{i,j-1} -\tau_{i-1,j}
+{\bm \zeta}^T{\bm Q}_{i,j}
\geq  0,
i\in {\cal I}, j\in {\cal J},\\
& \sum_{k=1}^Kp_k\alpha_{N_1,j}^k
-\lambda_{N_1-1,j}^1
+\lambda_{N_1,j}^2
-\lambda_{N_1,j-1}^2
+\eta_{N_1-1,j}^1
+\eta_{N_1,j-1}^2
-\eta_{N_1,j}^2
\nonumber \\
& 
\quad  +\tau_{N_1-1,j-1}
-\tau_{N_1-1,j} +{\bm \zeta}^T {\bm Q}_{N_1,j}\geq  0, j\in {\cal J},\\
& \sum_{k=1}^Kp_k\alpha_{1,j}^k
+\lambda_{1,j}^1
+\lambda_{1,j}^2
-\lambda_{1,j-1}^2
-\eta_{1,j}^1
+\eta_{1,j-1}^2
-\eta_{1,j}^2
+\tau_{1,j}
-\tau_{1,j-1} \nonumber \\
&
\quad  +{\bm \zeta}^T {\bm Q}_{1,j} \leq  0, j\in {\cal J}, \\
& \sum_{k=1}^Kp_k \alpha_{i,N_2}^k
+ \lambda_{i,N_2}^1
-\lambda_{i-1,N_2}^1
-\lambda_{i,N_2-1}^2
+\eta_{i-1,N_2}^1
-\eta_{i,N_2}^1
+\eta_{i,N_2-1}^2 
 \nonumber \\
& 
\quad  +\tau_{i-1,N_2-1}
-\tau_{i,N_2-1}
 +{\bm \zeta}^T {\bm Q}_{i,{N_2}}
\geq  0, i\in {\cal I},\\
& \sum_{k=1}^Kp_k \alpha_{i,1}^k
+ \lambda_{i,1}^1
-\lambda_{i-1,1}^1
+\lambda_{i,1}^2
+\eta_{i-1,1}^1
-\eta_{i,1}^1
-\eta_{i,1}^2 
+\tau_{i,1}
-\tau_{i-1,1}\nonumber \\
&
\quad  +{\bm \zeta}^T {\bm Q}_{i,{N_2}}
\geq  0, i\in {\cal I}\\
& \sum_{k=1}^Kp_k \alpha_{1,1}^k
+ \lambda_{1,1}^1
+\lambda_{1,1}^2
-\eta_{1,1}^1
-\eta_{1,1}^2 
+\tau_{1,1}
 +{\bm \zeta}^T {\bm Q}_{i,{N_2}}
\geq  0,\\
& \sum_{k=1}^Kp^k\alpha_{N_1,1}^k-\lambda_{N_1-1,1}^1+\lambda_{N_1,1}^2+\eta_{N_1-1,1}^1-\eta_{N_1,1}^2-\tau_{N_1-1,1}\geq 0,\\
& \sum_{k=1}^Kp^k\alpha_{1,N_2}^k+\lambda_{1,N_2}^1-\lambda_{1,N_2-1}^2-\eta_{1,N_2}^1+\eta_{1,N_2-1}^2-\tau_{1,N_2-1}\geq 0,\\
&\inmat{constraints } (\ref{eq:mixed-integer-R2-b}),
  (\ref{eq:mixed-integer-R2-d})-
  (\ref{eq:mixed-integer-R2-f}),\label{eq:PRO_MILP_single-j}\\
&  \sum_{i=1}^{N_1} \sum_{j=1}^{N_2} \alpha_{i,j}^{k} x_{i}=f_1(\bdz,\bdxi^k),\;\; \sum_{i=1}^{N_1} \sum_{j=1}^{N_2} \alpha_{i,j}^{k} y_{j}= f_2(\bdz,\bdxi^k), \nonumber \\
& \hspace{16em} k=1,\cdots,K,\\
& {\bm \lambda}^1 \geq 0, 
{\bm \lambda}^2 \geq 0,
{\bm \eta}^1\geq 0,
{\bm \eta}^2 \geq 0,
{\bm \tau} \geq 0,
\end{align}
\end{subequations}
where ${\cal I}:=\{2,\cdots,N_1-1\}$,
${\cal J}:=\{2,\cdots,N_2-1\}$,
${\bm Q}_{ij}:=(Q_{ij}^1,\cdots,Q_{ij}^M)^T\in \R^M$,
$Q_{ij}^l:= \mathbb{P}({\bm B}_l=(x_i,y_j))-\mathbb{P}({\bm A}_l=(x_i,y_j))$,
${\bm \lambda}^1\in \R^{(N_1-1)\times N_2}_+$,
${\bm \lambda}^2\in \R^{N_1\times (N_2-1)}_+$,
${\bm \eta}^1\in \R^{(N_1-1)\times N_2}_+$,
${\bm \eta}^2 \in \R^{N_1\times (N_2-1)}_+$,
${\bm \tau}\in \R^{(N_1-1)\times (N_2-1)}_+$,
${\bm \zeta}\in \R^M$.
\end{proposition}
We can also reformulate BUPRO-N with the Type-2 PLA as a single MIP.
We only need to
replace (\ref{eq:PRO_MILP_single-j}) with 
(\ref{eq:mixed-integer-R2-b}),
  (\ref{eq:mixed-integer-R2-d})-(\ref{eq:mixed-integer-R2-e}) and (\ref{eq:constraint-alpha}).
  Note that 
Hu et al.~\cite{hu2022distributionally}
consider a distributionally robust model for 
the random utility maximization problem 
in multi-attribute decision making 
and reformulate a maximin PRO 
as a single MILP. 
The main difference is that they considered 
the true utility function to be 
in additive form (sum of the single-attribute utility functions).
Here we consider a general   
multivariate true utility function.
Thus, we believe this is a step forward
from computational perspective
in handling BUPRO-N. 
Note that in this formulation, we have not incorporated  
Assumption~\ref{A:concave-in-x-and-y}
because the dual formulation of the problem with 
the convexity/concavity constraints 
would be very complex.

\subsection{Tri-attribute case}
\label{sec:three-dim-u}

We now extend our discussions 
on the implicit PLA of UPRO to the tri-attribute case.

\subsubsection{Triangulation of a cube and interpolation}


We follow the well-known triangulation method (see e.g.
\cite{chien1977solving,meyer2005convex,misener2010piecewise}) to 
divide each
cube into six non-overlapping simplices. 
There are six ways to divide, 
{\color{black}and}
here we 
use the second way (called Type B in the references).
Specifically, we consider $u:[\underline{x},\bar{x}]\times[\underline{y},\bar{y}]\times [\underline{z},\bar{z}]\rightarrow \R$
with 
$\underline{x}=x_1< x_2 < \cdots< x_{N_1}=\bar{x}$,
$\underline{y}=y_1< y_2 < \cdots< y_{N_2}=\bar{y}$
and $\underline{z}=z_1 < z_2< \cdots< z_{N_3}=\bar{z}$.
Let $X_i:=(x_i,x_{i+1}]$,
$Y_j:=(y_j,y_{j+1}]$
and $Z_l:=(z_l,z_{l+1}]$.
For any given point $(x,y,z)\in X_i\times Y_j\times Z_l$, 
consider the cube with vertices 
$1$: $(x_i,y_j,z_l)$, $2$: $(x_{i+1},y_j,z_l)$, $3$: $(x_{i},y_{j+1},z_l)$ $4$: $(x_{i+1},y_{j+1},z_l)$,
$5$: $(x_i,y_j,z_{l+1})$, $6$: $(x_{i+1},y_j,z_{l+1})$, $7$: $(x_{i},y_{j+1},z_{l+1})$, $8$: $(x_{i+1},y_{j+1},z_{l+1})$.
We first divide a cube $[\underline{x},\bar{x}]\times[\underline{y},\bar{y}]\times [\underline{z},\bar{z}]$ in $\R^3$ into two parts,
denoted by 
Part $1$-$2$-$4$-$5$-$6$-$8$
and Part $1$-$3$-$4$-$5$-$7$-$8$.
Then we can produce six simplices by three planes,
see 
Figures~\ref{fig-divisioin-6-1} $\&$ \ref{fig-divisioin-6-2}.
%
\begin{figure}[!ht]
\vspace{-0.1cm}
  \centering
   \subfigure[{\color{red}Red} simplex $1$-$2$-$4$-$8$]{
     \label{fig:3a} 
    \includegraphics[width=0.3\linewidth]{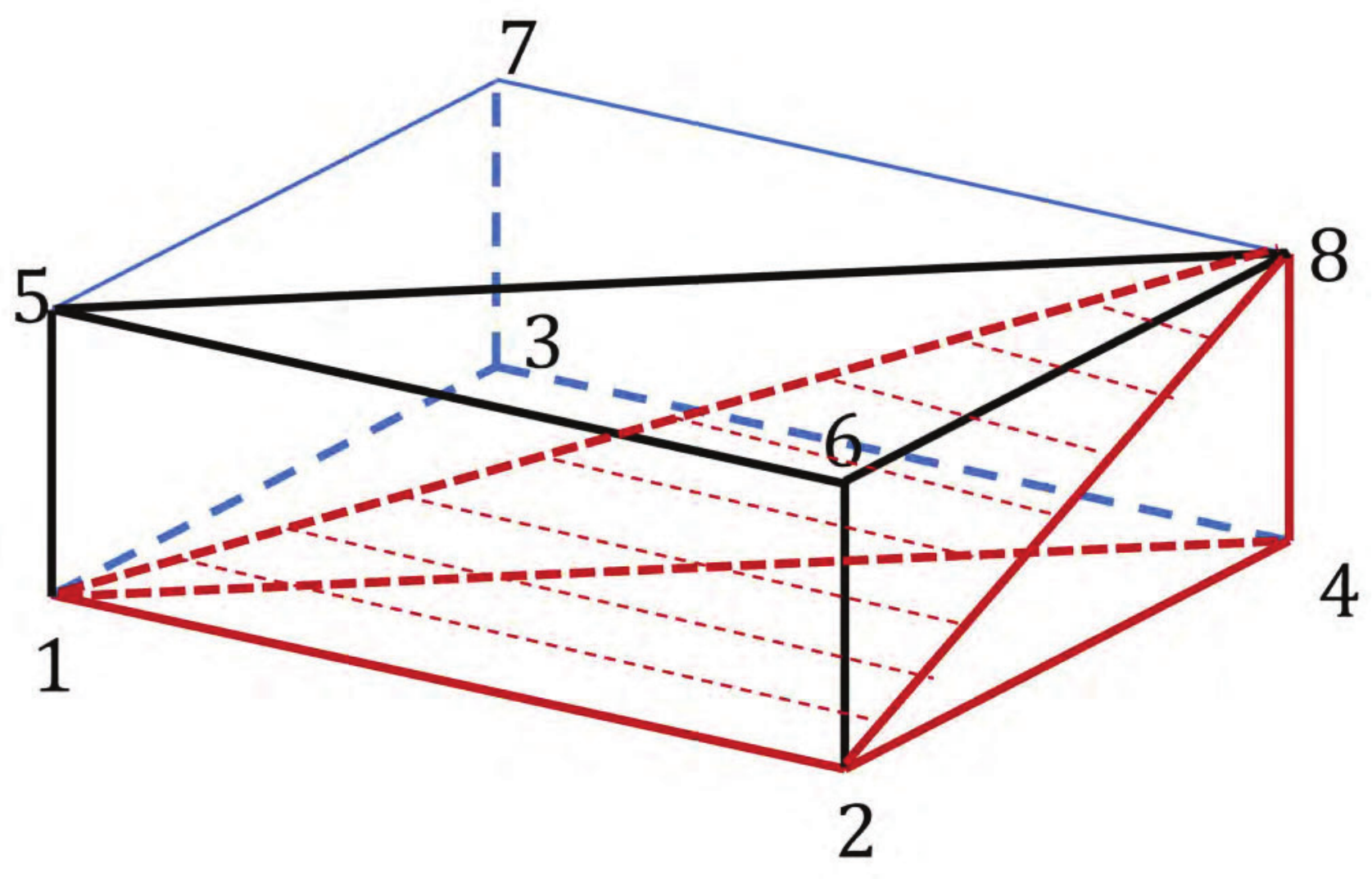}
  }
   \subfigure[{\color{PineGreen}Green} simplex $1$-$2$-$6$-$8$]{
     \label{fig:3b} 
    \includegraphics[width=0.3\linewidth]{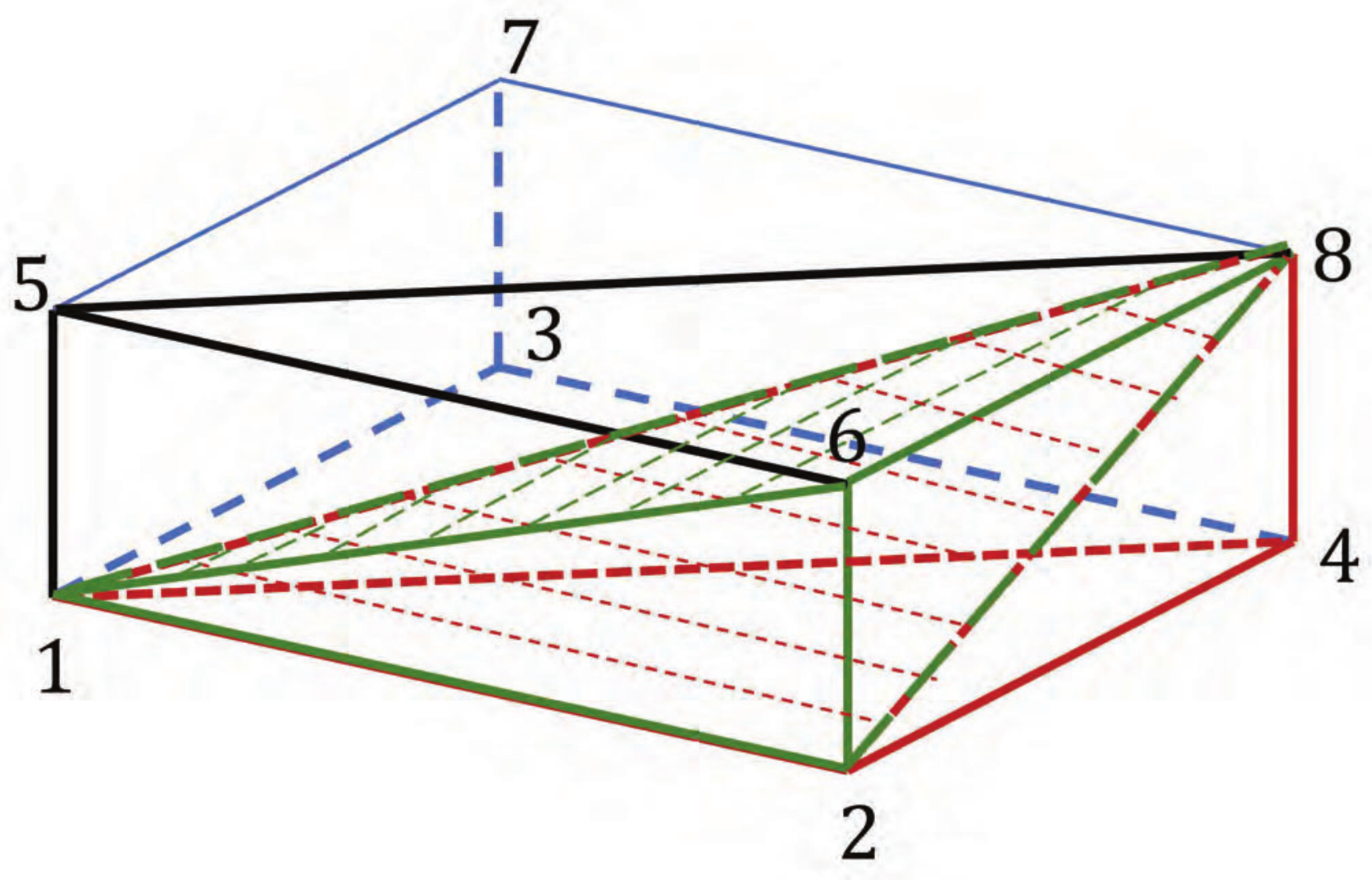}
  }
   \subfigure[{\color{purple}Purple} simplex $1$-$5$-$6$-$8$]{
     \label{fig:3c}
    \includegraphics[width=0.3\linewidth]{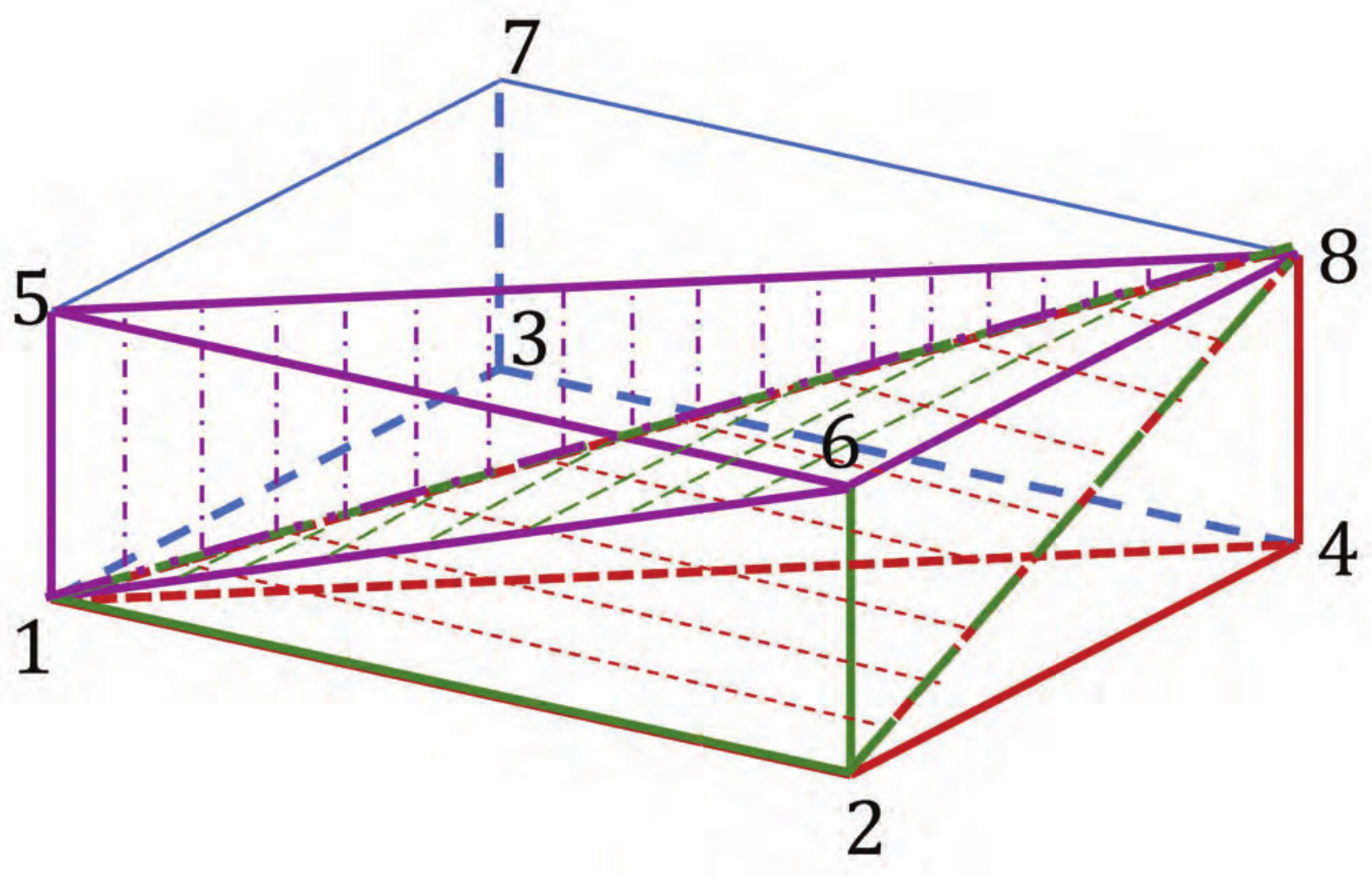}
  }
   \vspace{-0.3cm}
  \captionsetup{font=footnotesize}
  \caption{\footnotesize
  Divide the part $1$-$2$-$4$-$5$-$6$-$8$ into three simplices in $\R^3$.
(a) cuts the part $1$-$2$-$4$-$5$-$6$-$8$ by plane with vertices $1$-$2$-$8$,
and get the first simplex $1$-$2$-$4$-$8$ (red color).
(b) \& (c) go on to cut the rest part by plain $1$-$2$-$8$,
and obtain the second simplex $1$-$2$-$6$-$8$ (green color)
and the third simplex $1$-$5$-$6$-$8$ (purple color).}
\label{fig-divisioin-6-1} 
\end{figure}

\begin{figure}[!ht]
\vspace{-0.5cm}
  \centering
   \subfigure[{\color{red}Red} simplex $1$-$5$-$7$-$8$]{
     \label{fig:4a} 
    \includegraphics[width=0.3\linewidth]{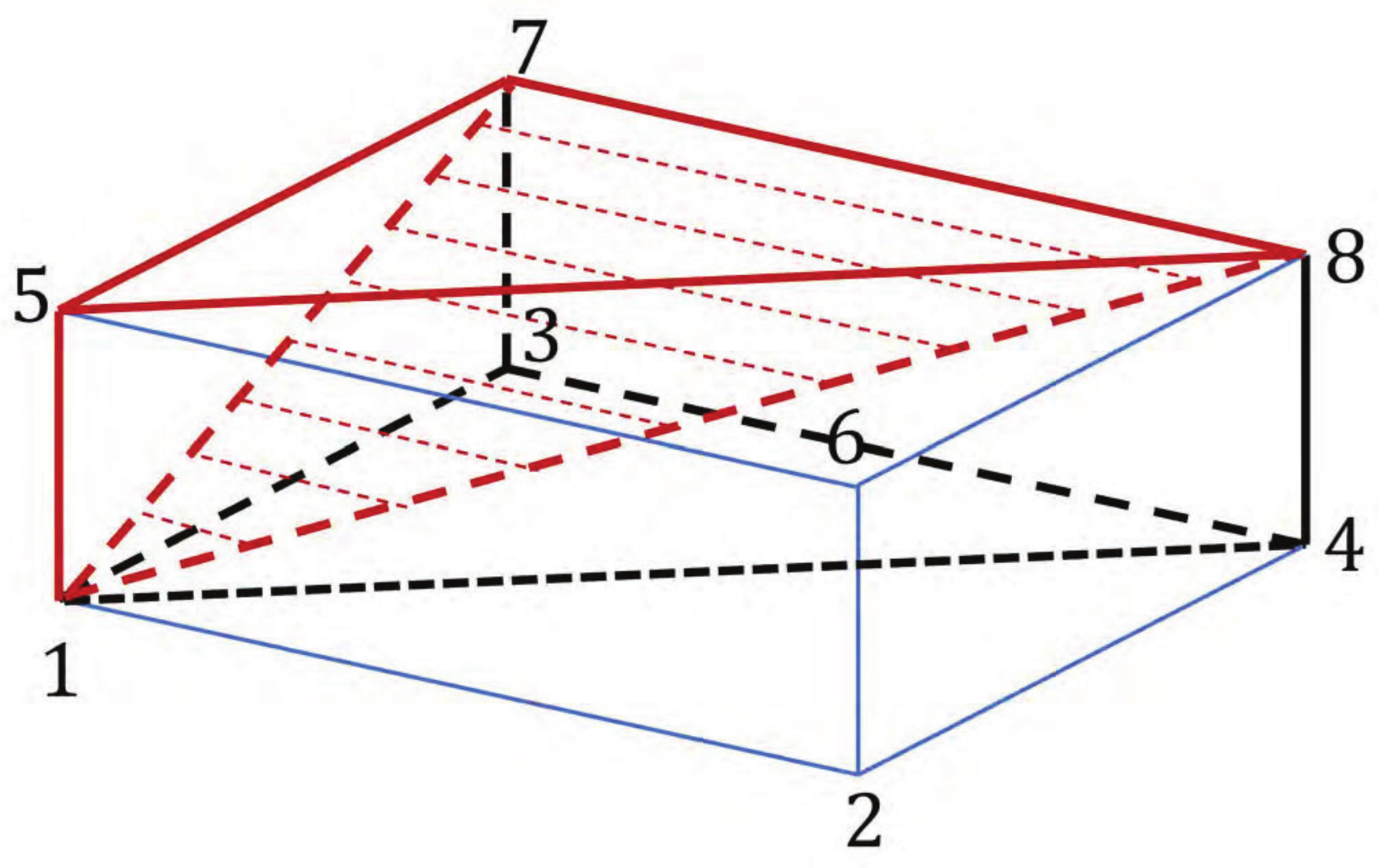}
  }
   \subfigure[{\color{PineGreen}Green} simplex $1$-$3$-$7$-$8$]{
     \label{fig:4b} 
    \includegraphics[width=0.3\linewidth]{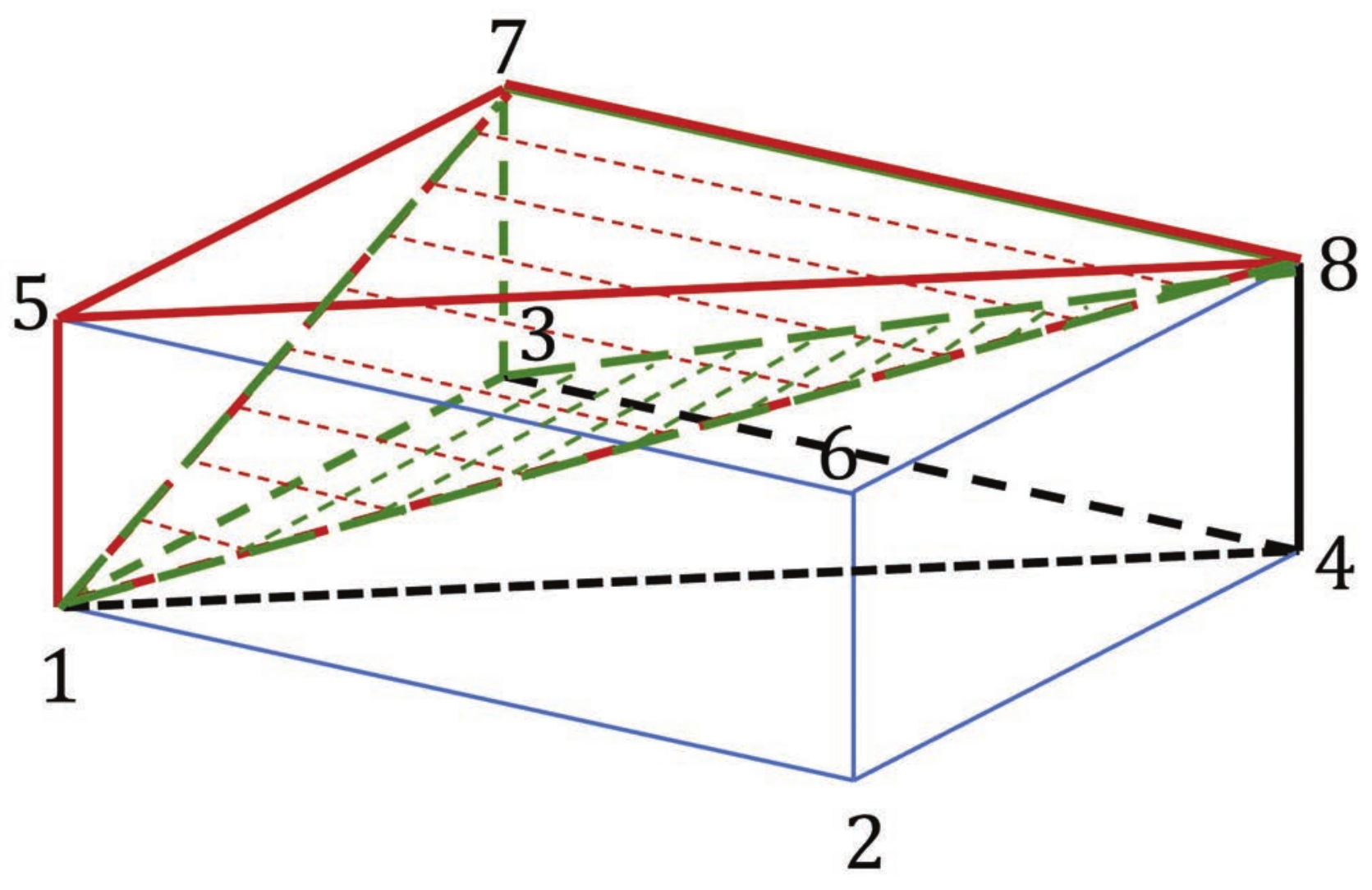}
  }
   \subfigure[{\color{purple}Purple} simplex $1$-$3$-$4$-$8$]{
     \label{fig:4c}
    \includegraphics[width=0.3\linewidth]{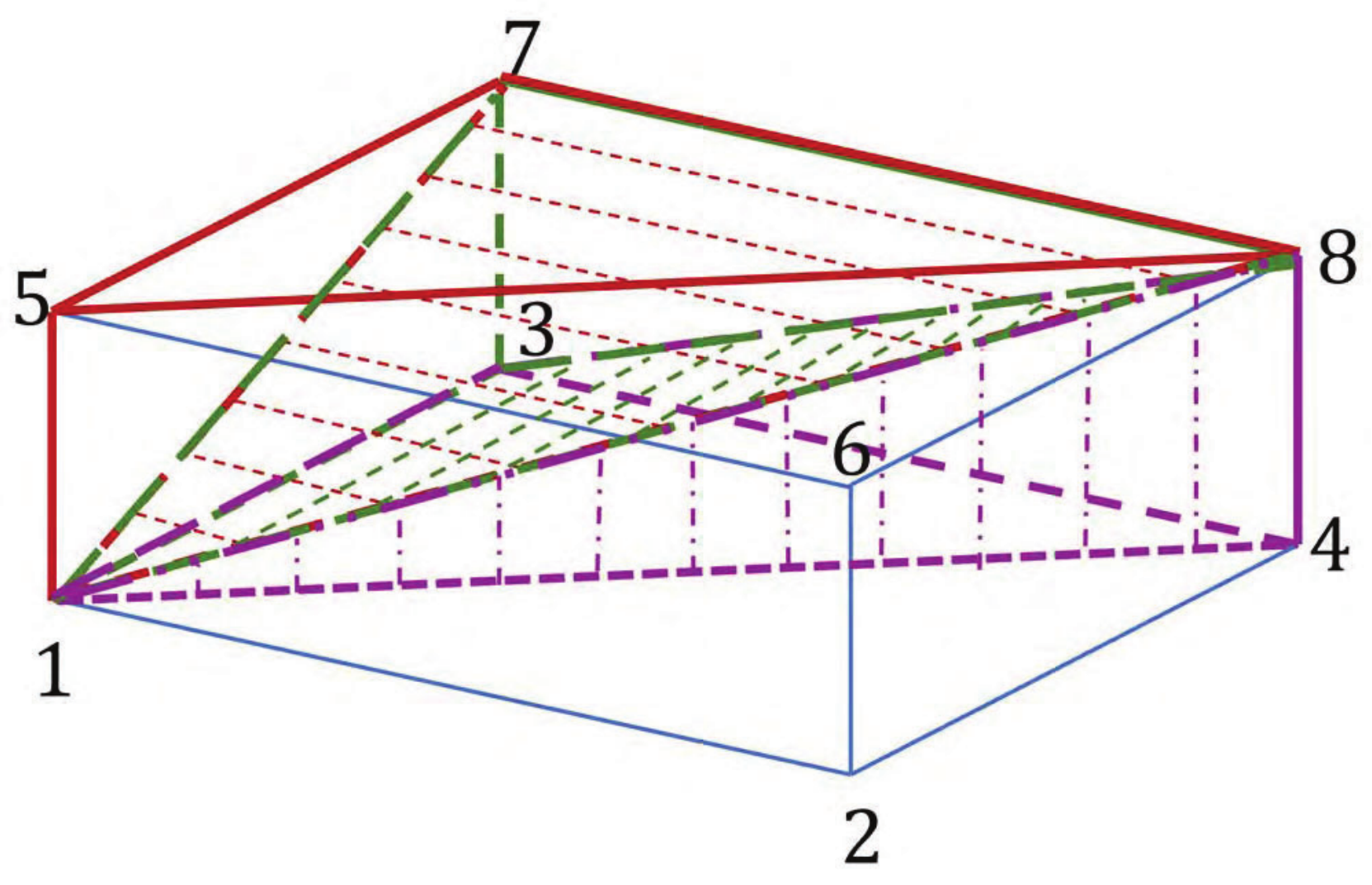}
  }
  \vspace{-0.2cm}
  \captionsetup{font=footnotesize}
  \caption{\footnotesize Divide the part $1$-$3$-$4$-$5$-$7$-$8$ into three simplices in $\R^3$.
(a) cuts the part $1$-$3$-$4$-$5$-$7$-$8$ by plane of vertices $1$-$7$-$8$,
and get the first simplex $1$-$5$-$7$-$8$ (red color).
(b) \& (c) go on to cut the rest part by plain $1$-$3$-$8$,
and obtain the second simplex $1$-$3$-$7$-$8$ (green color)
and the third simplex $1$-$3$-$4$-$8$ (purple color).}
  \vspace{-0.5cm}
  \label{fig-divisioin-6-2} 
\end{figure}

Let $1$-$4$-$5$-$8$$\searrow$
denote the front half subspace of  the plane constructed by points $1$-$4$-$5$-$8$,
i.e.,
\bgeq
\inmat{$1$-$4$-$5$-$8$$\searrow$}:=\{(x,y,z)^T\in \R^3:
(y_{j+1}-y_j)x-(x_{i+1}-x_i)y+x_{i+1}y_j-x_iy_{j+1}\geq 0\},
\edeq
and let
$1$-$2$-$7$-$8$ $\uparrow$ denote the upper subspace of the plane constructed by points $1$-$2$-$7$-$8$,
that is,
\bgeq
&& \inmat{ $1$-$2$-$7$-$8$} \uparrow:=\{(x,y,z)^T\in \R^3: (z_{l+1}-z_l) y-(y_{j+1}-y_j)z +y_{j+1}z_l-y_jz_{l+1} \geq 0\},\\
&& \inmat{ $1$-$3$-$6$-$8$} \uparrow:=\{(x,y,z)^T\in \R^3: 
(z_{l+1}-z_l) x -(x_{i+1}-x_i) z +x_{i+1}z_l-x_iz_{l+1}\leq 0\}.
\edeq
The function value $u({x},{y},{z})$ is approximated by a convex combination of the function values evaluated at the vertices of the simplex containing $({x},{y},{z})$,
that is,
$$
u({x},{y},{z})=\lambda u_{i,j,l}
+\mu u_{{i+1},{j+1},{l+1}}+
\bar{u},
$$
where  $\lambda, \mu\in [0,1]$ and 
{\small \bgeq
&& \bar{u}=\left\{\begin{array}{ll}
\eta u_{{i+1},j,l}
+(1-\lambda-\mu-\eta)u_{{i+1},{j+1},l} & \inmat{if }\;
(x,y)\in \inmat{$1$-$4$-$5$-$8$}
\searrow 
\bigcap 
\inmat{$1$-$2$-$7$-$8$ $\downarrow$ (Fig.~\ref{fig:3a})},\\
\eta u_{{i+1},j,{l+1}} +(1-\lambda-\mu-\eta) u_{{i+1},j,l} & \inmat{if }\;
(x,y)\in \inmat{$1$-$4$-$5$-$8$}
\searrow
\bigcap
\inmat{$1$-$3$-$6$-$8$} \downarrow 
\bigcap 
\inmat{$1$-$2$-$7$-$8$ $\uparrow$ (Fig.~\ref{fig:3b})},\\
\eta u_{{i+1},j,{l+1}} +(1-\lambda-\mu-\eta) u_{i,j,{l+1}} & \inmat{if }\;
(x,y)\in \inmat{$1$-$4$-$5$-$8$}
\searrow
\bigcap
\inmat{$1$-$3$-$6$-$8$} \uparrow 
\bigcap \inmat{$1$-$2$-$7$-$8$ $\uparrow$ (Fig.~\ref{fig:3c})},
\end{array}
\right.\\
&& \bar{u}=\left\{\begin{array}{ll}
\eta u_{i,j,{l+1}}
+(1-\lambda-\mu-\eta)u_{i,{j+1},{l+1}} & 
\inmat{if }\;
(x,y)\in \inmat{$1$-$4$-$5$-$8$}
\nearrow 
\bigcap
\inmat{$1$-$2$-$7$-$8$ $\uparrow$ (Fig.~\ref{fig:4a})},\\
\eta u_{{i},{j+1},l}
+(1-\lambda-\mu-\eta) u_{{i},{j+1},{l+1}} & 
\inmat{if }\;
(x,y)\in \inmat{$1$-$4$-$5$-$8$}
\nearrow 
\bigcap
\inmat{$1$-$3$-$6$-$8$} \uparrow
\bigcap
\inmat{$1$-$2$-$7$-$8$ $\downarrow$(Fig.~\ref{fig:4b})},\\
\eta u_{{i},{j+1},l}
+(1-\lambda-\mu-\eta) u_{{i+1},{j+1},l}& 
\inmat{if }\; (x,y)\in \inmat{$1$-$4$-$5$-$8$}
\nearrow 
\bigcap
\inmat{$1$-$3$-$6$-$8$} \downarrow 
\bigcap
\inmat{$1$-$2$-$7$-$8$ $\downarrow$(Fig.~\ref{fig:4c})}.
\end{array}
\right.
\edeq
}

\subsubsection{Implicit PLA}

As in the two-dimensional case,
since $u_N$ is linear over each simplex,
a target point ${\bm f}(\bdz,\bdxi^k)\in \R^3$ and its approximate utility value 
$u_N({\bm f}(\bdz,\bdxi^k))$ have the same 
convex combination. 
We use binary variables $h_{i,j,l,k}^{1u}$,
$h_{i,j,l,k}^{1m}$,
$h_{i,j,l,k}^{1l}$
to characterize 
whether  point 
${\bm f}(\bdz,\bdxi^k)$ lies 
in the upper or middle,
or lower simplex in Part $1\inmat{-}2\inmat{-}4\inmat{-}5\inmat{-}6\inmat{-}8$ or beyond,
see Figure~\ref{fig-divisioin-6-1}.
Likewise, we use binary variables $h_{i,j,l,k}^{2u}$,
$h_{i,j,l,k}^{2m}$,
$h_{i,j,l,k}^{2l}$
to characterize 
whether  ${\bm f}(\bdz,\bdxi^k)$ 
lies in the upper, or middle,
or lower simplex in Part $1\inmat{-}3\inmat{-}4\inmat{-}5\inmat{-}7\inmat{-}8$  or beyond,
see Figure~\ref{fig-divisioin-6-2}.
As in the bi-attribute case, we can identify the coefficients of the convex combinations 
by solving a system of linear equalities and inequalities:
\begin{subequations}
\label{eq:mixed-integer-R3}
\begin{align}
& \sum_{i=1}^{N_1} \sum_{j=1}^{N_2}
\sum_{l=1}^{N_3}
\alpha_{i,j,l}^k=1,\;\;k=1,\cdots,K, \label{eq:mixed-integer-R3-b}\\
& \sum_{i=1}^{N_1} \sum_{j=1}^{N_2} \sum_{l=1}^{N_3} \alpha_{i,j,l}^k x_{i}= f_1^k,\;
\sum_{i=1}^{N_1} \sum_{j=1}^{N_2} \sum_{l=1}^{N_3}\alpha_{i,j,l}^k y_{j}=f_2^k,\;
\sum_{i=1}^{N_1} \sum_{j=1}^{N_2} \sum_{l=1}^{N_3}\alpha_{i,j,l}^k z_{l}=f_3^k, \nonumber \\
& \hspace{20em} k=1,\cdots,K,
\label{eq:mixed-integer-R3-d} \\
& \sum_{i=1}^{N_1-1} \sum_{j=1}^{N_2-1}
\sum_{l=1}^{N_3-1}
h_{i,j,l,k}^{1u}+h_{i,j,l,k}^{1m}+h_{i,j,l,k}^{1l}
+h_{i,j,l,k}^{2u}+h_{i,j,l,k}^{2m}+h_{i,j,l,k}^{2l}=1, \nonumber \\
& \hspace{20em} k=1,\cdots,K,
\label{eq:mixed-integer-R3-e} \\
 & {\bm h}_k^{\tau u},
 {\bm h}_k^{\tau m},
 {\bm h}_k^{\tau l}
 \in \{0,1\}^{(N_1-1) (N_2-1) (N_3-1)},\;
 \tau=1,2,\;k=1,\cdots,K,
 \label{eq:mixed-integer-R3-f}\\
& 0\leq \alpha_{i,j,l}^k \leq \sum_{\nu = {\rm I}}^{\rm VIII} H_{i,j,l,k}^{\nu},\;
i=1,\cdots,N_1,\; j=1,\cdots,N_2, \;
l=1,\cdots,N_3, \nonumber \\
& \hspace{20em} k=1,\cdots,K, 
\label{eq:mixed-integer-R3-g}
\end{align}
\end{subequations}
where 
${\bm h}_{k}^{\tau u}:= (h_{1,1,1,k}^{\tau u},\cdots, h_{N_1-1,N_2-1,N_3-1,k}^{\tau u})^T$,
${\bm h}_{k}^{\tau m}:=(h_{1,1,1,k}^{\tau m},\cdots, h_{N_1-1,N_2-1,N_3-1,k}^{\tau m})^T$,
${\bm h}_{k}^{\tau l}:=(h_{1,1,1,k}^{\tau l}$,
$\cdots, h_{N_1-1,N_2-1,N_3-1,k}^{\tau l})^T$
for $\tau =1,2$,
 $k=1,\cdots,K$,
and
\bgeq
&& H_{i,j,l,k}^{\rm I}:= h_{i,j,l,k}^{1u}+h_{i,j,l,k}^{1m}+h_{i,j,l,k}^{1l}
+h_{i,j,l,k}^{2u}+h_{i,j,l,k}^{2m}+h_{i,j,l,k}^{2l},  
 \\
&& H_{i,j,l,k}^{\rm II}
:=h_{i-1,j,l,k}^{1m}+h_{i-1,j,l,k}^{1l}, \qquad~~~
H_{i,j,l,k}^{\rm III}
:=h_{i-1,j-1,l,k}^{1l}+h_{i-1,j-1,l,k}^{2l},\\
&& 
H_{i,j,l,k}^{\rm IV}
:=h_{i,j-1,l,k}^{2m}
+h_{i,j-1,l,k}^{2l},
\qquad~~~ 
H_{i,j,l,k}^{\rm V}
:=h_{i,j,l-1,k}^{1u}+h_{i,j,l-1,k}^{2u},\\
&&   
H_{i,j,l,k}^{\rm VI}:
=h_{i-1,j,l-1,k}^{1u}
+h_{i-1,j,l-1,k}^{1m}, 
\quad
H_{i,j,l,k}^{\rm VIII}:=
h_{i,j-1,l-1,k}^{2u}+h_{i,j-1,l-1,k}^{2m},\\
&& 
H_{i,j,l,k}^{\rm VII}:=
h_{i-1,j-1,l-1,k}^{1u}
+h_{i-1,j-1,l-1,k}^{1m}
+h_{i-1,j-1,l-1,k}^{1l}
+h_{i-1,j-1,l-1,k}^{2u}\\
&& \qquad \qquad 
+h_{i-1,j-1,l-1,k}^{2m}
+h_{i-1,j-1,l-1,k}^{2l},
\edeq
 ${\bm f({\bm w},\bdxi^k)}=(f_1^k,f_2^k,f_3^k)^T$,
$f_i^k:=f_i({\bm w},\bdxi^k)$ for
$i=1,2,3$,
$h_{0,*,*,*}^*=h_{*,0,*,*}^*=h_{*,*,0,*}^*=h_{N_1,*,*,*}^*=h_{*,N_2,*,*}^*=h_{*,*,N_3,*}^*=0$.
Constraint (\ref{eq:mixed-integer-R3-e})
imposes the restriction that only one is active
for the convex combination among all 
$6$ simplices.
Constraint (\ref{eq:mixed-integer-R3-g}) imposes that the
only nonzero $\alpha_{i,j,l}^k$ can be those associated with the
four vertices of such simplex,
see 
Figure~\ref{fig:I-VIII}.
$H_{i,j,l,k}^\nu$
represents the sum of $h_{i,j,l,k}^*$
with $*\in \{1u,1m,1l,2u,2m,2l\}$
in Octant $\nu$ that are related to point $(x_i,y_j,z_l)$ for $\nu={\rm I},\ldots,{\rm VIII}$.
Specifically,
there are $6$ simplices in Octant I that are related to point $(x_i,y_j,z_l)$,
and the corresponding binary variables are
$h_{i,j,l,k}^{1u}$,
$h_{i,j,l,k}^{1m}$,
$h_{i,j,l,k}^{1l}$,
$h_{i,j,l,k}^{2u}$,
$h_{i,j,l,k}^{2m}$,
$h_{i,j,l,k}^{2l}$.
There are two triangles in Octant II that are related to $(x_i,y_j,z_l)$,
and
the corresponding binary variables are $h_{i-1,j,l,k}^{1m}$ and
$h_{i-1,j,l,k}^{1l}$.
The related binary variables in Octant III-VIII
can also be observed.
Such $h_{i,j,l,k}^*$ can be used to identify which vertices are used to represent 
${\bm f}(\bdz,\bdxi^k)$.
For example, if $h_{i,j,l,k}^{1l}=1$, then 
$\bdf(\bdz,\bdxi^k)$ lies in the
lower simplex in the former part of
cube $X_i\times Y_j\times Z_l$.
This is 
indicated
by the fact 
that
$\alpha_{i,j,l}^k\leq h_{i,j,l,k}^{1l}=1$,
$\alpha_{i+1,j+1,l+1}^k\leq h_{i,j,l,k}^{1l}=1$,
$\alpha_{i+1,j,l}^k\leq h_{i,j,l,k}^{1l}=1$,
$\alpha_{i+1,j+1,l}^k\leq h_{i,j,l,k}^{1l}=1$,
and $\alpha_{i',j',l'}^k=0$ for $(i',j',l')\notin \{(i,j,l),(i+1,j+1,l+1),(i+1,j,l),(i+1,j+1,l)\}$,
see Figure~\ref{fig:I-VIII} for
the 24 simplices that are related to point $(x_i,y_j,z_l)$.

\begin{figure}[!ht]
  \centering
   \subfigure[A cube]{
    \includegraphics[width=0.25\linewidth]{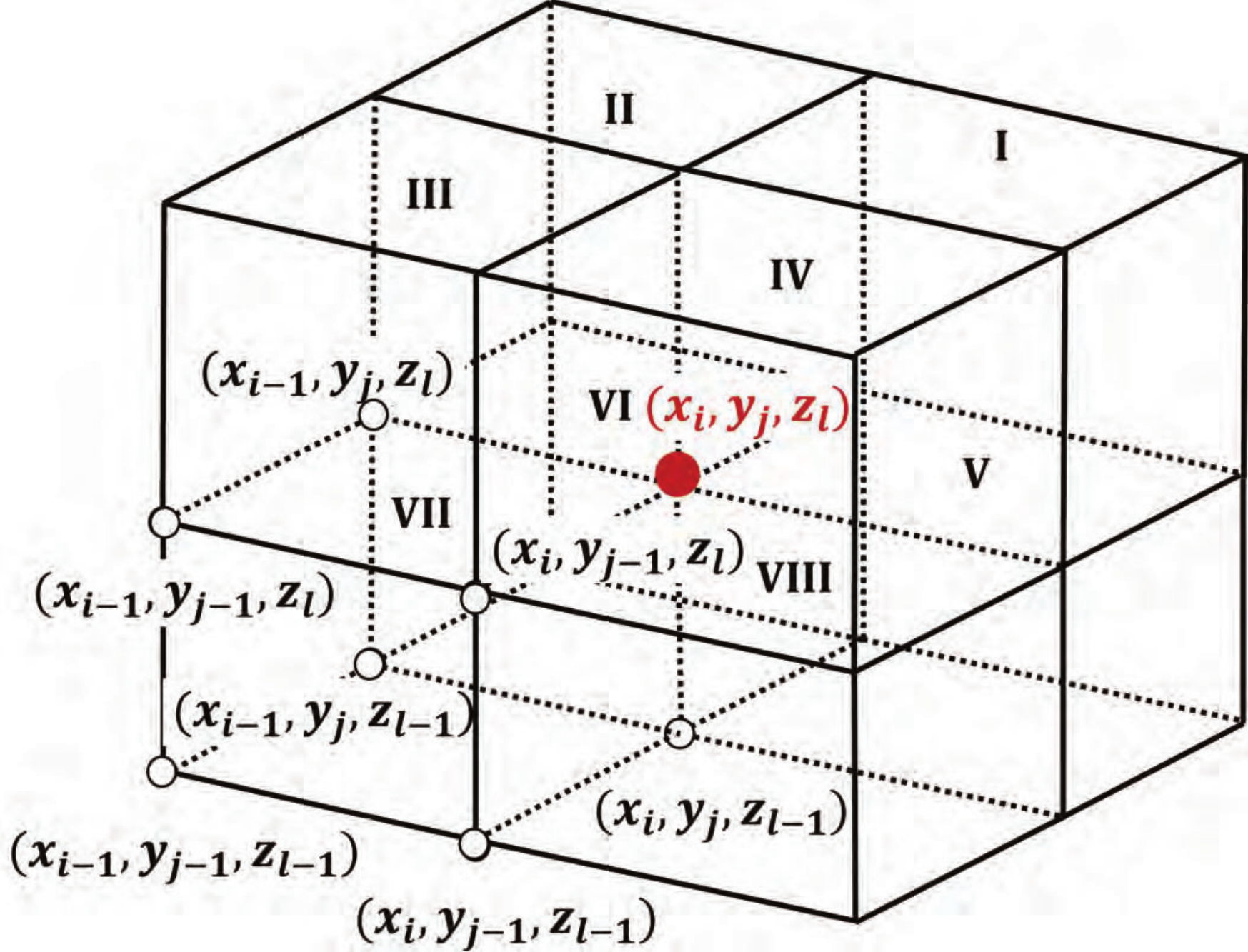}
    }
     \hspace{-0.3em}
   \subfigure[Octants I-IV]{
    \includegraphics[width=0.32\linewidth]{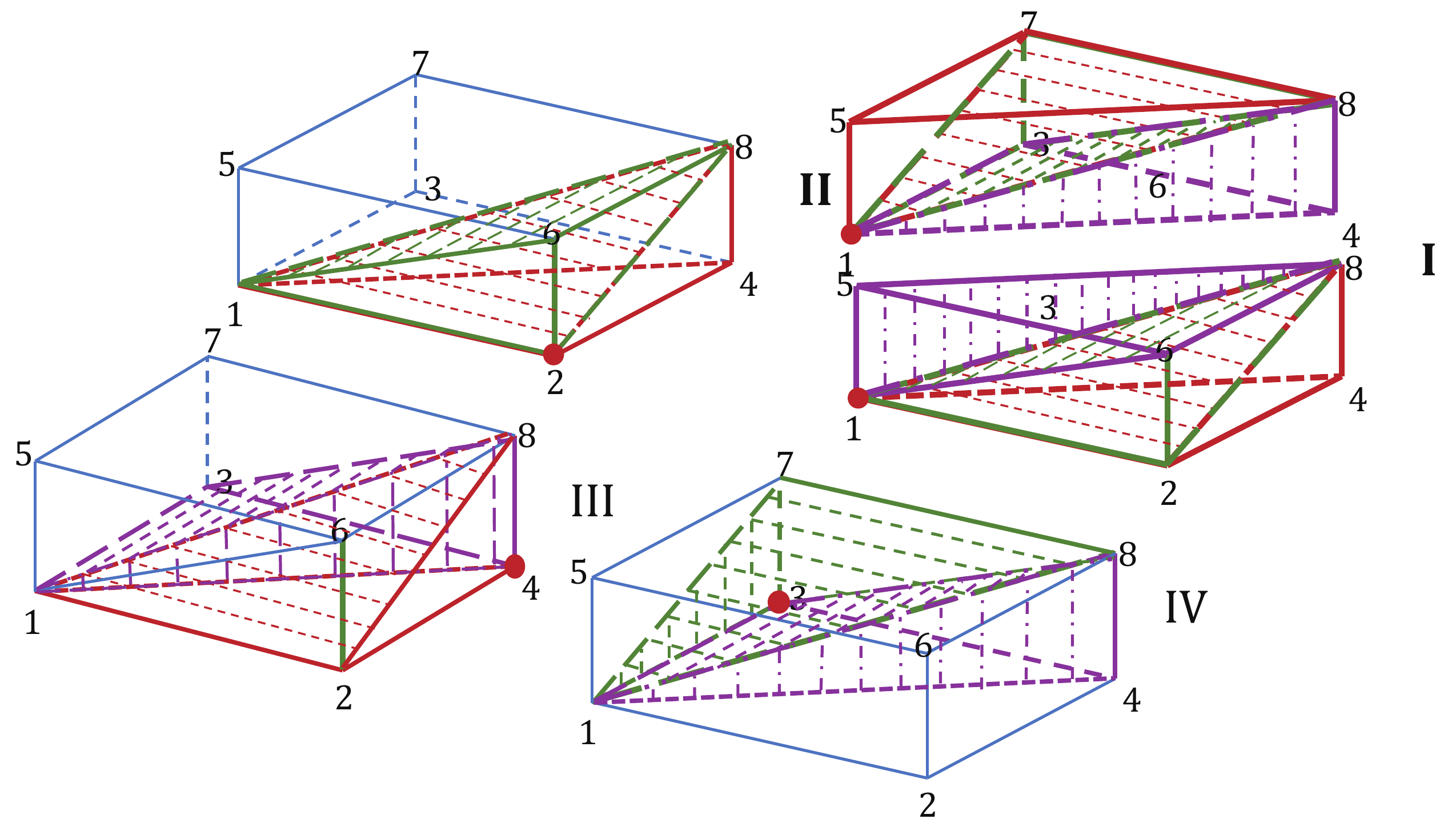}
  }
  \hspace{-0.3em}
   \subfigure[Octants V-VIII]{
    \includegraphics[width=0.32\linewidth]{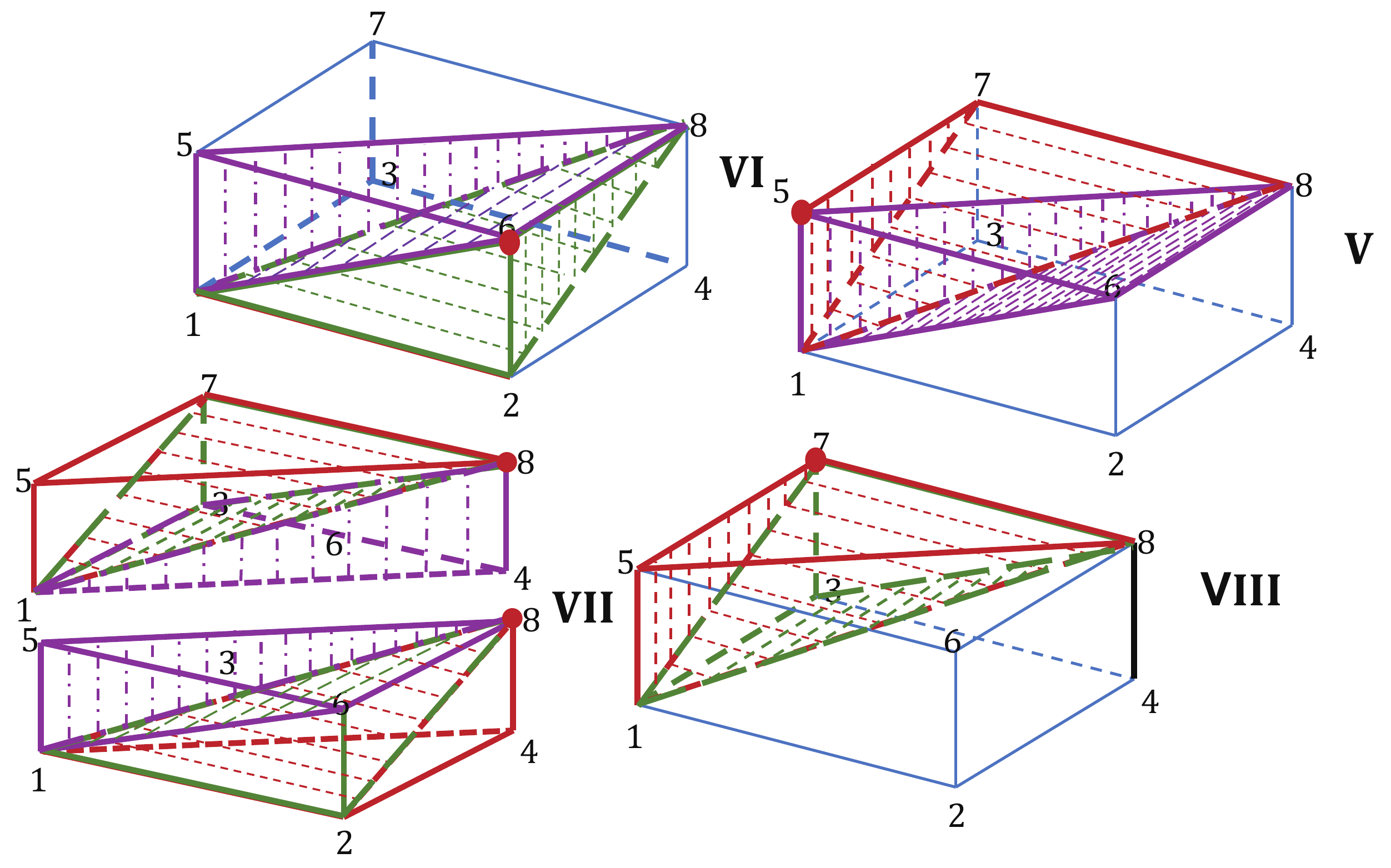}
  }
  \vspace{-0.35cm}
  \caption{\footnotesize 
  (a) divides a cube into $8$ sub-cubes denoted by octants I, II, III, IV, V, VI, VII and VIII.
  The red point in (a)-(c) is $(x_i,y_j,z_l)$.
  (b) $\&$ (c) illustrate all the $24$ simplices related to the point $(x_i,y_j,z_l)$,
which 
prompts
the last constraint in problem (\ref{eq:mixed-integer-R3}).
(b) represents the cases in octants I-IV.
In octant I,
the vertex $1$ is $(x_i,y_j,z_l)$,
and there are six simplices containing the red point $(x_i,y_j,z_l)$.
In octant II,
the vertex $1$ is $(x_{i-1},y_j,z_l)$,
and there are two simplices containing the red point.
In octant III,
the vertex $1$ is $(x_{i-1},y_{j-1},z_l)$,
and there are two simplices containing the red point.
In octant IV,
the vertex $1$ is $(x_{i},y_{j-1},z_l)$,
and there are two simplices containing the red point. 
(c) represents the cases 
in octants V-VIII.
In octant V,
the vertex $1$ is $(x_i,y_j,z_{l-1})$,
and there are two simplices containing the red point.
In octant VI,
the vertex $1$ is $(x_{i-1},y_j,z_{l-1})$,
and there are two simplices containing the red point.
In octant VII,
the vertex $1$ is $(x_{i-1},y_{j-1},z_{l-1})$,
and there are six simplices containing the red point. 
In octant VIII,
the vertex $1$ is $(x_{i},y_{j-1},z_{l-1})$,
and there are two simplices containing the red point.}
  \label{fig:I-VIII} 
\end{figure}
Consequently, we can reformulate the tri-attribute utility maximization problem \linebreak
$\max_{{\bm w}\in Z}\sum_{k=1}^K p_k[u_N({\bm f({\bm w},\bdxi^k)})]$ as:
\begin{subequations}
\label{eq:mixed-integer-R3-2}
\begin{align}
 \max\limits_{{\bm w}\in Z,{\bm \alpha},{\bm h^l},{\bm h}^m,{\bm h}^u} \; &  \sum_{k=1}^K p_k \sum_{i=1}^{N_1} \sum_{j=1}^{N_2}\sum_{l=1}^{N_3}\alpha_{i,j,l}^k u_{i,j,l}\\
{\rm s.t.} \qquad\;\; & \inmat{constraints } (\ref{eq:mixed-integer-R3-b})-
(\ref{eq:mixed-integer-R3-g}),
\end{align}
\end{subequations}
where 
${\bm \alpha}:=({\bm \alpha}^1,\cdots,{\bm \alpha}^K)\in \R^{(N_1N_2N_3)\times K}$,
${\bm \alpha}^k:=(\alpha_{1,1,1}^k,\cdots,\alpha_{N_1,N_2,N_3}^k)^T$ for $k=1,\cdots,K$,
${\bm h}^u:=({\bm h}^{1u}_1,\cdots,{\bm h}^{1u}_K,{\bm h}^{2u}_1,\cdots,{\bm h}^{2u}_K)\in \R^{(N_1-1)(N_2-1)(N_3-1)\times 2K}$,\;
${\bm h}^m:=({\bm h}^{1m}_1,\cdots,{\bm h}^{1m}_K,{\bm h}^{2m}_1,\cdots, \\ {\bm h}^{2m}_K)$,\;
${\bm h}^l:=({\bm h}^{1l}_1,\cdots,$
${\bm h}^{1l}_K,{\bm h}^{2l}_1,\cdots,{\bm h}^{2l}_K)$,\;
${\bm u}:=(u_{1,1,1},\cdots,$
$u_{N_1,N_2,N_3})^T\in \R^{N_1 N_2N_3}$.
If $f({\bm w},\bdxi)$ is linear in ${\bm w}$,
then (\ref{eq:mixed-integer-R3-2}) is an MILP.}
Extending this to the UPRO model, we consider the ambiguity set ${\cal U}_N$ constructed by pairwise comparison of questions $({\bm A}_m,{\bm B}_m)$.
Under Assumption~\ref{assu-lip},
suppose that the set of gridpoints $\{(x_i,y_j,z_l):i=1,\cdots,N_1,j=1,\cdots,N_2,l=1,\cdots,N_3\}$ contains all the outcomes of lotteries $ {\bm A}_m$ and ${\bm B}_m$ for $m=1,\cdots,M$,
then we can solve the tri-attribute utility preference robust optimization (TUPRO) problem by solving the approximate TUPRO-N problem 
$\max_{\bdz\in Z} \min_{u_N\in {\cal U}_N} \sum_{k=1}^Kp_k[u_N({\bm f(\bdz,\bdxi^k)})]$ 
as:
\begin{subequations}
\label{eq:PRO_MILP_3m}
\begin{align}
\max\limits_{\bdz\in Z}
\min\limits_{
\substack{{\bm \alpha}, {\bm h}^u,{\bm h}^m\\
{\bm h}^l,\bm u}}\;&
\sum_{k=1}^K p_k \sum_{i=1}^{N_1} \sum_{j=1}^{N_2}
\sum_{l=1}^{N_3}
\alpha_{i,j,l}^k u_{i,j,l}\\
{\rm s.t.} \;\;\;\; & 
 u_{i+1,j,l}\geq u_{i,j,l}, i=1,\cdots,N_1-1, j=1,\cdots,N_2,l=1,\cdots,N_3,  \label{eq:PRO_MILP_3m-f}\\
& u_{i,j+1,l}\geq u_{i,j,l}, i=1,\cdots,N_1, j=1,\cdots,N_2-1,l=1,\cdots,N_3, \label{eq:PRO_MILP_3m-h}\\
& u_{i,j,l+1}\geq u_{i,j,l}, i=1,\cdots,N_1, j=1,\cdots,N_2,l=1,\cdots,N_3-1, \label{eq:PRO_MILP_3m-i}\\
& u_{i+1,j,l}-u_{i,j,l}\leq L(x_{i+1}-x_i), \nonumber \\
& \qquad \qquad i=1,\cdots,N_1-1, j=1,\cdots,N_2,l=1,\cdots,N_3,\label{eq:PRO_MILP_3m-j}\\
& u_{i,j+1,l}-u_{i,j,l}\leq L(y_{j+1}-y_j), \nonumber \\
& \qquad \qquad i=1,\cdots,N_1, j=1,\cdots, N_2-1, l=1,\cdots,N_3, \label{eq:PRO_MILP_3m-k}\\
& u_{i,j,l+1}-u_{i,j,l}\leq L(z_{l+1}-z_{l}), \nonumber \\
& \qquad \qquad  i=1,\cdots,N_1,j=1,\cdots,N_2, l=1,\cdots,N_3-1, \label{eq:PRO_MILP_3m-l} \\
& u_{1,1,1}=0, \; u_{N_1,N_2,N_3}=1,  \label{eq:PRO_MILP_3m-m}\\
&\sum_{i=1}^{N_1} \sum_{j=1}^{N_2}\sum_{l=1}^{N_3} ( \mathbb{P}({\bm B}_m=(x_i,y_j,z_l))-\mathbb{P}({\bm A}_m=(x_i,y_j,z_l)))u_{i,j,l} \leq 0,
\nonumber \\
& \hspace{16em} m=1,\cdots,M, \label{eq:PRO_MILP_3m-n}\\
& \inmat{constraints }(\ref{eq:mixed-integer-R3-b})-(\ref{eq:mixed-integer-R3-g}),\label{eq:PRO_MILP_3m-b}
\end{align}
\end{subequations}
where
 ${\bm u}:=(u_{1,1,1},\cdots,u_{N_1,N_2,N_3})^T\in \R^{N_1 N_2N_3}$.
We can solve problem (\ref{eq:PRO_MILP_3m})
by a Dfree method,
where the  inner problem is an MILP when ${\bm f}(\bdz,\bdxi)$ is linear in $\bdz$.
It is also possible to reformulate 
the problem further as a single MILP, we leave this for interested readers.

\subsection{Multi-attribute case}
\label{sec:m-dim-u}


Since a large number of simplices are needed to partition hypercubes of dimension greater than three
(see \cite{hughes1996simplexity}),
we give a general framework 
for the $m$ attributes case.
For ${\bm x}\in \R^m$,
we  can divide the domain of utility 
$\bigtimes_{i=1}^{m}[\underline{x}_i,\overline{x}_i]$ into $(N_1-1)\times (N_2-1)\times \cdots \times (N_m-1)$ subsets $\{\bigtimes_{i=1}^{m}[x_{i_j},x_{i_{j+1}}]:j=1,\cdots,N_{i}-1\}$.
We denote the values of $u_N$ at $(x_{1_{j_1}},\cdots,x_{m_{j_m}})$ by $u_{1_{j_1},\cdots,m_{j_m}}$ for
$j_i=1,\cdots,N_i$,
$i=1,\cdots,m$.
We {\color{black}
reshape $(u_{1_{j_1},\cdots,m_{j_m}})_{N_1\times \cdots \times N_m}\in \R^{N_1\times \cdots \times N_m}$
as a vector ${\bm u}=(u_1,\cdots,u_V)^T\in \R^{V}$}
with $V:=N_1\times\cdots \times N_m$,
and label the 
corresponding vertices 
by $1,\cdots,V$.
{\color{black}
We divide the domain $\bigtimes_{i=1}^{m}[\underline{x}_i,\overline{x}_i]$ into mutually exclusive simplices and  label them
by $1,\cdots,S$.}
The $v$-th vertice is ${\bm x}_v:=(x_{1_v},\cdots,x_{i_v},\cdots,x_{m_v})^T \in \R^m$ for $v=1,\ldots,V$.
Let ${\cal V}_s$ denote the set of vertices of the $s$-th simplex.
As in the bi-attribute and tri-attribute cases,
for given
${\bm f}(\bdz,\bdxi^k)=(f_1(\bdz,\bdxi^k),f_2(\bdz,\bdxi^k),\cdots,f_m(\bdz,\bdxi^k))^T$,
we can
identify 
 the simplex containing 
${\bm f}(\bdz,\bdxi^k)$
and obtain
the coefficients 
of the representation of ${\bm u}$ 
at ${\bm f}(\bdz,\bdxi^k)$
in terms of the utility values at the 
vertices of the simplex
by solving a system of linear equalities and
inequalities:
\begin{subequations}
\label{eq:mixed-integer-Rm}
\begin{align}
& \sum_{v=1}^{V}
\alpha_{v}^k=1,\;\;k=1,\cdots,K,
\label{eq:mixed-integer-Rm-b}\\
& \sum_{v=1}^V \alpha_{v}^k x_{i_v}=f_i(\bdz,\bdxi^k),\;\;i=1,\cdots,m,\;
k=1,\cdots, K,
\label{eq:mixed-integer-Rm-c}\\
& \sum_{s=1}^S h_s^k=1,\;\;h_s^k\in \{0,1\}, \;s=1,\cdots,S, \;k=1,\cdots,K,
\label{eq:mixed-integer-Rm-d}\\
&  0\leq \alpha_v^k\leq \sum_{s:{\bm x}_v\in {\cal V}_s} h_s^k, \; v=1,\cdots,V,\;k=1,\cdots,K,
\label{eq:mixed-integer-Rm-e}
\end{align}
\end{subequations}
where $s:{\bm x}_v\in {\cal V}_s$ means all $s\in \{1,\cdots,S\}$ satisfying that 
the vertice ${\bm x}_v$ belongs to the set ${\cal V}_s$. 
Constraint
(\ref{eq:mixed-integer-Rm-e})
implies that
only $\alpha_{v}$ values different from $0$ 
are those associated with the
vertices of 
the simplex.
Then we can reformulate the multi-attribute utility maximization problem 
$\max_{\bdz\in Z}\sum_{k=1}^K p_k[u_N({\bm f(\bdz,\bdxi^k)})]$ as 
an MIP
(see e.g. \cite{VAN10}),
\begin{subequations}
\label{eq:mixed-integer-Rm-2}
\begin{align}
\max\limits_{\bdz \in Z, {\bm \alpha}, {\bm h}} \; &  \sum_{k=1}^K p_k \sum_{v=1}^{V}\alpha_{v}^k u_v\\
{\rm s.t.}\;\;\; & 
\inmat{constraints }(\ref{eq:mixed-integer-Rm-b})-(\ref{eq:mixed-integer-Rm-e}),
\end{align}
\end{subequations}
where 
${\bm \alpha}:=({\bm \alpha}^1,\cdots,{\bm \alpha}^K)\in \R^{V\times K}$,
${\bm h}:=({\bm h}^1,\cdots,{\bm h}^K)\in \R^{S \times K}$.
If $f({\bdz,\bdxi})$ is linear in $\bdz$,
then problem (\ref{eq:mixed-integer-Rm-2}) is an MILP.
We continue to assume that the ambiguity set ${\cal U}_N$ 
is constructed by pairwise comparisons of lotteriess ${\bm A}_l$ and ${\bm B}_l$ with $l=1,\cdots,M$.
Let 
$
\tilde{\cal U}_N={\cal U}_N \bigcap \{u_N: u_N \inmat{ is Lipschtz with its modulus }L\}.
$
Consequently,
we can solve 
the multi-attribute utility preference robust  problem
by solving 
MUPRO-N problem
$\max_{\bdz\in Z} \min_{u_N\in \tilde{\cal U}_N}\sum_{k=1}^K p_k[u_N({\bm f(\bdz,\bdxi^k)})]$ 
as 
an MIP:
\begin{subequations}
\begin{align}
\label{eq:PRO_m-dim}
\disp{ \max_{\bdz\in Z} \min_{{\bm u}, {\bm \alpha},{\bm h} } }\; &
\sum_{k=1}^K p_k \sum_{v=1}^{V}\alpha_{v}^k u_v \\
{\rm s.t.} \;\; &
\inmat{constraints }(\ref{eq:mixed-integer-Rm-b})-(\ref{eq:mixed-integer-Rm-e}), \;\\
& u_{1_{j_1},\cdots,i_{j_{i}+1},\cdots,m_{j_m}}
\geq u_{1_{j_1},\cdots,i_{j_{i}},\cdots,m_{j_m}}, \nonumber \\
& \hspace{12em}
j_i=1,\cdots,N_i,
i=1,\cdots,m,
\label{eq:PRO_m-dim-c}\\
& u_{1_{j_1},\cdots,i_{j_{i}+1},\cdots,m_{j_m}}
-
u_{1_{j_1},\cdots,i_{j_{i}},\cdots,m_{j_m}} \leq L(x_{i_{j_i+1}}-x_{i_{j_i}}), \nonumber \\
& \hspace{12em}
j_i=1,\cdots,N_i,i=1,\cdots,m, \qquad~~\;
\label{eq:PRO_m-dim-e}\\
& u_1=0,\; u_V=1,\label{eq:PRO_m-dim-f}\\
&\sum_{v=1}^{V}  \mathbb{P}({\bm B}_l=\bdx_v)u_v\leq \sum_{v=1}^{V} \mathbb{P}({\bm A}_l=
\bdx_v)u_v, 
\;l=1,\cdots,M,
\end{align}
\end{subequations}
where ${\bm u}\in \R^{V}$,
${\bm \alpha}=({\bm \alpha}^1,\cdots,{\bm \alpha}^K)\in \R^{V\times K}$,
${\bm h}=({\bm h}^1,\cdots,{\bm h}^K)\in \R^{S \times K}$.
Constraint (\ref{eq:PRO_m-dim-c}) represents the 
non-decreasing
property of the utility function.
Constraint (\ref{eq:PRO_m-dim-e}) represents the Lipschitz continuity of $u_N$
and (\ref{eq:PRO_m-dim-f}) characterizes the normalization of $u_N$.

\section{Error bounds for  the PLA}
\label{sec-errorbound}

In the previous section, we outline
computational schemes 
to solve BUPRO-N problem. In this section, we investigate 
the error bounds of the optimal value and the optimal solutions obtained from solving BUPRO-N problem when we use them to approximate the optimal value and optimal solutions of BUPRO problem.
Notice that 
the only difference between 
the two maximin optimization problems 
is the feasible set of the inner minimization problem,
thus we proceed with our investigation 
by quantifying the difference between
$\calu_N$ and $\calu$ and then apply classical 
stability results in parametric programming to
derive the error bounds of the optimal value and optimal solutions.
Proofs of all technical results are deferred to the appendix. 

To ease the exposition, we 
write $\langle u, \psi_l \rangle$ for 
$\int_T u(x,y) d \psi_l(x,y)$,
and subsequently 
(\ref{eq:ambiguity_set}) as
\begin{equation}
\label{eq-pseume}
    \calu=\{u\in\scru: \la u,\bdps \ra \leq \bdc\},
\end{equation}
where $\bdps := (\psi_1(x,y), \ldots, \psi_M(x,y))^T \in \R^M$, $\bdc:=(c_1,\ldots,c_M)^T\in \R^M$.
Note that $\langle u, \psi_l \rangle$ should 
not be read as a kind of inner product as 
we cannot swap the positions between $u$ and $\psi_l$.
We adopt the 
notation since
(\ref{eq-pseume})
clearly      indicates 
${\cal U}$ as the set of the solutions of the 
inequality system $\la u,\bdps \ra \leq \bdc$
relative to $\scru$.
To quantify the difference between two utility functions,
we  define, 
for any $u,v\in\scru$,
the pseudo-metric between $u$ and $v$
under the function set
$\scrg$ by
\begin{equation*}
    \dd_{\scrg}(u,v):=\sup_{g\in\scrg} |\la g,u \ra - \la g,v \ra|.
\end{equation*}
It is easy to observe that $\dd_{\scrg}(u,v)=0$ if and only if $\la g,u \ra=\la g,v \ra$ for all $g\in\scrg$. 
In practice,
we may regard $\scrg$ as a set of “test functions” associated with some prospects and interpret $u$ as a measure induced by utility.
The pseudo-metric means that if $u$ and $v$ give the same average value
for each $g\in\scrg$, then they are regarded as “equal” under $\dd_{\scrg}$ 
although
they may not be identical.
Thus $\dd_{\scrg}$ is a kind of pseudo-metric defined over the space of utility-induced measures $\scru$.
This definition is in parallel to a similar definition in probability theory, where $u$ and $v$ are 
in a position of probability measures and the corresponding pseudo-metric is known as $\zeta$-metric, see \cite{Rom03}.
Here we continue to adopt the terminology although
the background is different.



\begin{example}
\label{exm-g}
Recall that $T= [\underline{x},\bar{x}]\times[\underline{y},\bar{y}]$.

(a) Let 
$$
\scrg=\scrg_M:=\lt\{ g: T \rightarrow \R \,\lt|\; 
\inmat{g is measurable,} \sup_{\bdt\in T}|g(\bdt)| \leq 1 \rt.\rt\}.
$$ 
Then $\dd_{\scrg_M}(u,v)$ corresponds to the total variation metric and $\dd_{\scrg}(u,v)\leq 1$.

(b) Let 
\begin{equation}
\label{eq-kantorovich}
    \scrg=\scrg_K:=\{g:T\to\R \,|\; \inmat{g is Lipschitz continuous with the modulus bounded by 1}\}.
\end{equation}
Then $\dd_{\scrg_K}(u,v)$ corresponds to the Kantorovich metric in which case we have
$\dd_{\scrg_K}(u,v)=\int_T \|\bdt-\bdt'\|  d\pi(\bdt,\bdt') \leq \sqrt{(\bar{x}-\underline{x})^2+(\bar{y}-\underline{y})^2)}
$,
where $\int_T \pi(\bdt,\bdt') d \bdt'=u(\bdt)$, $\int_T \pi(\bdt,\bdt') d \bdt=v(\bdt')$
and $\|\bdt\|$ denotes the Euclidean norm.

(c) Let $\scrg:=\scrg_L\cap\scrg_M$. Then $\dd_{\scrg}(u,v)$ corresponds to the bounded Lipschitz metric and $\dd_{\scrg}(u,v)\leq\min\left\{1,\sqrt{(\bar{x}-\underline{x})^2+(\bar{y}-\underline{y})^2)}\right\}$.

(d) Let 
\begin{equation*}
    \scrg=\scrg_I:=\{g:T\to\R \,|\; g=\1_{[\underline{x},x]\times[\underline{y},y]}(\cdot), (x,y)\in T \}.
\end{equation*}
Then $\dd_{\scrg_I}(u,v)$ corresponds to the Kolmogorov metric in which case we have $\dd_{\scrg_I}(u,v)\leq~1$.
\end{example}

For any two sets $U, V\subset \mathscr{U}$, 
let
$
\mathbb{D}_{\mathscr{G}}(U,V)  := \sup_{u\in U}\inf_{v\in V} \dd_{\mathscr{G}}(u,v),
$
which quantifies the deviation of $U$ from $V$
and
$
\mathbb{H}_{\mathscr{G}}(U,V) := \max\left\{\mathbb{D}_{\mathscr{G}}(U,V), \mathbb{D}_{\mathscr{G}}(V,U)\right\},
$
which denotes the Hausdorff distance between the two sets under the pseudo-metric.
By convention, when $U=\{u\}$ is a singleton, we write
the distance $\dd_{\mathscr{G}}(u,V)$ from $u$ to set $V$ 
rather than
$\mathbb{D}_{\mathscr{G}}(U,V)$.

Using the pseudo-metric, we 
can
derive an error bound of any utility function $u\in\scru$ deviating from $\calu$ in terms of the residual of the linear system defining ${\cal U}$. 
This 
type
of result is known as Hoffman’s lemma.
We state this in the next lemma.

\begin{lemma}[Hoffman's lemma]
\label{lem-hof}
Consider (\ref{eq-pseume}).
Assume: (a)
 $\mathscr{G}$ is chosen so that 
the resulting pseudo-distance 
between any two utility functions is finite-valued,
and 
(b) there exist a constant $\alpha$ and a function $u^0\in\calu$ such that 
\begin{equation}
\label{eq-sla}
    \la u^0, \bdps \ra -\bdc+\alpha \mathbb{B}^M \subset \R_-^M.
\end{equation}
Then 
\begin{equation}
\label{eq-hof}
    \dd_{\scrg} (u,\calu) \leq \frac{\dd_{\scrg}(u,u^0)}{\alpha} 
    \|(\la u, \bdps \ra -\bdc)_+\| \quad \forall u\in\scru,
\end{equation}
where $({\bm a})_+:=\max\{0,{\bm a}\}$ which is taken componentwise.
\end{lemma}

Condition (\ref{eq-sla}) is known as 
Slater’s condition. It implies that there is at least one utility function $u^0$ such that 
 $\langle u,\bdps\rangle $ lies in the interior of $\R^M_-$. 
 This kind of condition 
 is widely used in the literature of 
 Hoffman’s lemma for linear and convex systems, see \cite{Rob75} and references therein. 
Since the proof of Hoffman’s lemma in the case that utility function in $\R^2$ is similar to the case with utility function in $\R$ (\cite{GXZ21}), we omit the details.



\subsection{Error bound on the ambiguity set}

We move on to quantify the difference between $\calu$ and $\calu_N$.
First, we give the following technical result.

\begin{proposition}
\label{prop-d}
Let $u\in\scru$ and $u_N$ be the PLA of $u$ defined as in Proposition~\ref{prop-uti-N}, then the following assertions hold:

(i) If $\scrg=\scrg_K$, 
then
\begin{equation}
    \label{}
    \dd_{\scrg_K} (u,u_N)\leq 2(\beta_{N_1}^2+\beta_{N_2}^2)^{1/2},
\end{equation}
where $\scrg_K$ is defined as in Example~\ref{exm-g}~(b)
 and
\begin{equation}
\label{eq-be}
    \beta_{N_1}:=\max_{i=2,\ldots,N_1}(x_i-x_{i-1}), \; \beta_{N_2}:=\max_{j=2,\ldots,N_2}(y_j-y_{j-1}).
\end{equation}

(ii) If $\scrg=\scrg_I$ and 
$u$ is Lipschitz continuous over
$T$
with the modulus $L$,
then 
\begin{equation}
\label{eq-d}
    \dd_{\scrg_I} (u,u_{N}) \leq 
    2L\lt( \beta_{N_1}+\beta_{N_2} \rt),
\end{equation}
where $\scrg_I$ is defined as in Example~\ref{exm-g}~(d).

\end{proposition}


With Lemma~\ref{lem-hof} and Proposition~\ref{prop-d}, we are ready to quantify the difference between $\calu_N$ and $\calu$.

\begin{theorem}[Error bound on $\mathbb{H}_{\scrg}(\calu_N,\calu)$]
\label{thm-erramb}
Assume: 
(a) 
Slater's condition
in Lemma~\ref{lem-hof} is satisfied;
(b) 
$\int_T d\psi_l(\bdt)$ is well-defined;
{\color{black}
(c)
$u$ is Lipschitz continuous over
$T$
with the modulus $L$.}
Then there exist a positive constant  $\hat{\alpha}<\alpha$, $N^0_1$ and $N^0_2$ such that the following assertions hold for specific $\scrg$ defined as in Example~\ref{exm-g}.

(i) If $\scrg=\scrg_K$,
then
\begin{equation}
\label{eq-erram-L}
\begin{split}
    \mathbb{H}_{\scrg_K} (\calu,\calu_{N}) \leq  & 2(\beta_{N_1}^2+\beta_{N_2}^2)^{1/2} \\
    &  + L(\beta_{N_1}+ \beta_{N_2}) 
    \frac{\lt((\bar{x}-\underline{x})^2+(\bar{y}-\underline{y})^2\rt)^{1/2}}{\hat{\alpha}} \lt(\sum_{l=1}^M \lt| \int_T d \psi_l (\bdt)\rt|^2\rt)^{1/2}
\end{split}
\end{equation}
for all $N_1\geq N_1^0$ and $N_2\geq N_2^0$.

(ii) If $\scrg=\scrg_I$, then
\begin{equation}
\label{eq-erram-I}
    \mathbb{H}_{\scrg_I}(\calu,\calu_{N})\leq
    L\lt( \beta_{N_1}+\beta_{N_2} \rt) \lt( 2+\frac{1}{\hat{\alpha}}  \lt(\sum_{l=1}^M \lt| \int_T d \psi_l (t)\rt|^2\rt)^{1/2} \rt)
\end{equation}
for all $N_1\geq N_1^0$ and $N_2\geq N_2^0$,
where $\beta_{N_1}$, $\beta_{N_2}$ are defined as in (\ref{eq-be}) and $\beta_{N_i}\to 0$ as $N_i\to\infty$ for $i=1,2$. 
\end{theorem}

The constant $\hat{\alpha}$ is related to 
Slater’s condition for the linear system when the utility function is restricted to space $\scru_N$, 
{\color{black} see \cite[pages 16]{GXZ21}.}
It is well-known that Kolmogorov metric $\dd_{\scrg_I}$ is tighter than Kantorovich metric $\dd_{\scrg_K}$ defined as in Example~\ref{exm-g}~(b) and (d) because the former is about the largest difference between two  utility functions whereas the latter is about the area between the graphs of the two utility functions, see \cite{GiS02}.
Consequently,
$\mathbb{H}_{\scrg_I}(\calu_N,\calu)$ is tighter than $\mathbb{H}_{\scrg_K}(\calu_N,\calu)$.
The following corollary shows that the second term disappears in both cases 
when $\psi_l$, $l=1,\ldots,M$
are simple functions. 


\begin{corollary}
\label{cor-err-optval-discrete}
Let $u\in\scru$. Assume that $u$ is Lipschitz continuous over
$T$
with the modulus $L$.
 If $\psi_l$ is a simple function 
taking constant values 
over each cell of  $T$
for 
$l=1,\ldots,M$, then $\mathbb{H}_{\scrg_K}(\calu_{N},\calu)\leq 2(\beta_{N_1}^2+\beta_{N_2}^2)^{1/2}$ and $\mathbb{H}_{\scrg_I}(\calu_{N},\calu)\leq 2L\lt( \beta_{N_1}+\beta_{N_2} \rt)$.
\end{corollary}

The corollary provides us with some useful insights: 
if $\psi_l$ is a simple function for $l=1,\cdots,M$
(which 
corresponds to the case when the DM's preference is elicited via pairwise comparison lotteries), 
then we can construct the grid of $T$ in such a way that $\psi_l$ 
is constant over $T_{i,j}$ (the vertices of the cells comprise
all outcomes of the lotteries).
In this way, we may effectively reduce 
the modelling error arising from  PLA of the utility function.
Note also that in this case, 
Slater’s condition is not 
required,
which means 
that the error bound holds
for all $N_1$ and $N_2$ rather than for them to be sufficiently large.

\subsection{Error bound on the optimal value and the optimal solution}

We are now ready to quantify the difference between the BUPRO-N and BUPRO models. 
Let $\vt_N$ and $\vt$ denote the respective optimal values, and $Z_N^*$
and $Z^*$ denote the corresponding sets of optimal solutions.

\begin{theorem}[Error bound on the optimal value and the optimal solution]
\label{thm-optval}
Assume the settings and conditions of  
Theorem~\ref{thm-erramb}.
Then the following assertions hold.

(i)
\begin{equation}
\label{eq-err-vt}
    |\vt_{N}-\vt| \leq L\lt( \beta_{N_1}+\beta_{N_2} \rt) \lt( 3+\frac{1}{\hat{\alpha}}  \lt(\sum_{l=1}^M \lt| \int_T d \psi_l (t)\rt|^2\rt)^{1/2} \rt)
\end{equation}
for all $N_1\geq N_1^0$ and $N_2\geq N^0_1$, where $L$, $\hat{\alpha}$, $\beta_{N_1}$, $\beta_{N_2}$, $N_1^0$ and $N^0_2$ are defined as in Theorem~\ref{thm-erramb}.

(ii) 
Let $v(\bdz):=\min_{u\in\calu} \bbe_P[u(\bdf(\bdz,\bdxi))]$. 
Define the growth function 
$\Lambda(\tau):=\min\{v(\bdz)-\vt^*:d(\bdz,Z^*)\geq \tau,\forall\, \bdz\in Z\}$
and $\Lambda^{-1}(\eta):=\sup\{\tau:\Lambda(\tau)\leq\eta\}$ where $d(\bdz,Z^*)=\inf_{\bdz'\in Z^*} \|\bdz-\bdz'\|$.
Then 
\begin{equation}
\label{eq-err-so}
    \mathbb{D}(Z_{N}^*,Z^*)\leq \Lambda^{-1} \lt( 2L\lt( \beta_{N_1}+\beta_{N_2} \rt) \lt( 3+\frac{1}{\hat{\alpha}}  \lt(\sum_{l=1}^M \lt| \int_T d \psi_l (t)\rt|^2\rt)^{1/2} \rt) \rt),
\end{equation}
where $\mathbb{D}(Z_{N}^*,Z^*):=\sup_{\bdz\in Z_{N}^*} \inf_{\bdz' \in Z^*} \|\bdz-\bdz'\|$.
\end{theorem}

\begin{remark}
\label{rem:distance}
(i) Note that $\vt$ is not computable whereas $\vt_N$ is.
The error bound established in (\ref{eq-err-vt}) gives
the DM an interval centred at $\vt_N$ which contains $\vt$.
We can say that for a specified precision $\epsilon$, we can use the inequality to estimate $\beta_N$ such that $|\vt_N-\vt|\leq\epsilon$.
In the case when $x_1,\ldots,x_{N_1}$ and $y_1,\ldots,y_{N_2}$ 
are evenly spread over $[\underline{x},\bar{x}]$ and $[\underline{y},\bar{y}]$, we know the specified precision is reached when 
$L\lt( \frac{\bar{x}-\underline{x}}{N_1}+\frac{\bar{y}-\underline{y}}{N_2} \rt) \lt( 3+\frac{1}{\hat{\alpha}}  \lt(\sum_{l=1}^M \lt| \int_T d \psi_l (\bdt)\rt|^2\rt)^{1/2} \rt)\leq\epsilon$.

(ii) The error bound (\ref{eq-err-vt}) is established without restricting
the utility functions to being concave and it is derived under the PLA scheme. We envisage that similar results may be obtained using spline approximation and leave interested readers to investigate. 
Note that these are mesh-dependent approximation schemes which means that the quality of approximation depends on the number of
gridpoints $N=N_1N_2$.

(iii) 
Let $u^{\rm worst}_N\in \arg\min_{u_N\in {\cal U}_N} \sum_{k=1}^K p_k u_N({\bm f}(\bdz^N,\bdxi^k))$,
where $\bdz^N$ denotes the optimal solution of (\ref{eq:MAUT-robust-N-dis}).
Then
$\dd_{\mathscr{G}_I}(u^*,u_N^{\rm worst})=\sup_{\bdt\in T} |u^*(\bdt)-u_N^{\rm worst}(\bdt)|$.
Let $u_N^*$ denote the PLA of $u^*$ with identical values at the gridpoints.
Then 
$$
\dd_{\mathscr{G}_I}(u^*,u^*_N)=\sup_{\bdt\in T}|u^*(\bdt)-u^*_N(\bdt)|=\sup_{\substack{i=1,\cdots,N_1-1, \\ j=1,\cdots,N_2-1}} \sup_{\bdt \in T_{i,j}}|u^*(\bdt)-u_N^*(\bdt)|\leq L (\beta_{N_1}+\beta_{N_2}),
$$
and 
$\dd_{\mathscr{G}_I}(u^*_N,u_N^{\rm worst})=\sup_{\bdt\in T} |u_N^*(\bdt)-u_N^{\rm worst}(\bdt)|=\max_{i=1,\cdots,N_1,j=1,\cdots,N_2}|u_N^*(\bdt_{i,j})-u_N^{\rm worst}(\bdt_{i,j})|$,
where $\bdt_{i,j}:=(x_i,y_j)$.
In Section~\ref{sec:numerical results},
we will examine how $u_{N}^{\rm worst}$ converges to $u^*$ as the number of queries increases.

(iv) The error bounds established 
under $\dd_{\mathscr{G}_I}$ and $\dd_{\mathscr{G}_K}$ require
conservative property of the utility function.
Specifically, the bound of Hausdorff distance between $\calu$ and $\calu_N$ is related to two terms $\dd_{\scrg}(u,u_N)$ and $\dd_{\scrg}(u_N,u_N^0)$ (see (\ref{eq-uU})),
where $u\in\calu$, $u_N$ is the PLA of $u$, and $u_N^0$ is defined in (\ref{eq-sla-0}).
It can be observed that in the case that $\scrg=\scrg_K$, the bound of $\dd_{\scrg}(u,u_N)$ relies on the conservative property as shown in (\ref{eq-u-u-N}), whereas in the case $\scrg=\scrg_I$, the bound of $\dd_{\scrg}(u_N,u_N^0)$ relies on the conservative property in Example~\ref{exm-g} (d). This makes it difficult to extend 
the theoretical results to multivariate utility case. We leave this for future research.




\end{remark}

\begin{example}
Consider the ambiguity set 
defined as in (\ref{eq:ambi-U-ex}). 
Since ${\bm A}$ is preferred,
there exists some $u^0\in\scru$ and a small positive number $\epsilon$ such that $\int_T u^0(x,y) d(F_{\bm A}(x,y)-F_{\bm B}(x,y))<-\epsilon$.
Let $\alpha=-\epsilon-\int_T u^0(x,y) d(F_{\bm A}(x,y)-F_{\bm B}(x,y))>0$.
Then 
Slater’s condition (\ref{eq-sla}) is satisfied.
Let $\hat{\alpha}\in(0,\alpha)$ be such that $\int_T u^0_N(x,y) d(F_{\bm A}(x,y)-F_{\bm B}(x,y))+\hat{\alpha}
\in 
\R_-$.
Observe that
$\psi(x,y) :=F_{\bm A}(x,y)-F_{\bm B}(x,y)$ 
satisfies 
 $|\psi(x,y)|\leq2$ for all $(x,y)\in T$.
By Theorem~\ref{thm-optval},
$
|\vt-\vt_N| \leq L\lt( \beta_{N_1}+\beta_{N_2} \rt) \lt( 3+\frac{2}{\hat{\alpha}} \rt).
$
Moreover, if ${\bm A}$ and ${\bm B}$ follow discrete distributions, then $\psi$ is a step function. 
In that case, we may select the gridpoints in ${\cal X}\times {\cal Y}$ (in the PLA) from the gridpoints
of $\psi$ and subsequently it follows by Corollary~\ref{cor-err-optval-discrete} that $|\vt-\vt_N| \leq \mathbb{H}_{\scrg_I} (\calu,\calu_{N})+ L(\beta_{N_1}+\beta_{N_2})\leq 3L(\beta_{N_1}+\beta_{N_2})$.
\end{example}

\section{
BUPRO models for 
constrained optimization problem}
\label{sec:constrained}

In this section, we extend the UPRO model to the
expected utility maximization problem with expected utility constraints. Specifically, we consider the following problem:
\begin{equation}
\label{eq:SPR-x}
\begin{split}
    {\vt}^*:=\max_{\bdz\in Z} \;\; & \bbe_P[ u(\bdf(\bdz,\bdxi))] \\
    \st \;\;\, &
    \bbe_P[u(\bdg(\bdz,\bdxi))] \geq c,
\end{split}
\end{equation}
where $\bdf$ and $\bdg$ are continuous functions and $c$ is a constant. We may interpret 
$\bdf$ as the total return of a portfolio 
and $\bdg$ is an important part of it or vice versa. 
Suppose that the true utility function is unknown but it is possible to construct an ambiguity set ${\cal U}$
using partially available information as we discussed earlier. 
Then we may consider the following maximin preference robust optimization problem
\begin{equation}
\label{eq:PRO-x}
\begin{split}
   \hat{\vt}:= \max_{\bdz\in Z} \;\;
    \min_{u\in {\cal U}} \;\; & \bbe_P[ u(\bdf(\bdz,\bdxi))] \\
    \st \;\;  &
    \bbe_P[u(\bdg(\bdz,\bdxi))] \geq c.
\end{split}
\end{equation}
In this formulation, we consider the same worst-case utility function in the objective and constraint. There is an alternative
way to develop a robust formulation of (\ref{eq:SPR-x}):
\begin{equation}
\label{eq:PRO-x-1}
\begin{split}
   \tilde{\vt}:= \max_{\bdz\in Z} \;\;
    \min_{u\in {\cal U}} \;\; & \bbe_P[ u(\bdf(\bdz,\bdxi))] \\
    \st \;\;\,
    \min_{u\in {\cal U}} \;\; &
    \bbe_P[u(\bdg(\bdz,\bdxi))] \geq c.
\end{split}
\end{equation}
Formulation (\ref{eq:PRO-x-1}) means that the worst-case 
utility in the objective and in the constraint  might differ.
It is easy to observe that $\tilde{\vt} \leq \hat{\vt}$
which means (\ref{eq:PRO-x-1}) is more conservative than
(\ref{eq:PRO-x}).
Moreover,
if the true utility $u^*$ lies within ${\cal U}$,
then $\tilde{\vt}\leq \vt^*$.
However, under some conditions,
the two formulations are equivalent. The next proposition states this.

\begin{proposition}
\label{Prop-equivalence}
Let $\hat{\bdz}$ denote the optimal solution of problem (\ref{eq:PRO-x}) and
\begin{equation}
    \tilde{Z}:= 
    \left\{\bdz\in Z \,:\, \inf_{u\in {\cal U}} \; \bbe_P[u(\bdg(\bdz,\bdxi))] -c
    \geq 0  \right\}.
\label{eq:x*-PRO-U}
\end{equation}
If $\hat{\bdz}\in \tilde{Z}$, then  
$\hat{\vt}=\tilde{\vt}$.
\end{proposition}

\noindent
\textbf{Proof.} Let $\hat{v}(\bdz)$ denote the optimal value of
the inner minimization problem of (\ref{eq:PRO-x})
and 
$$
\tilde{v}(\bdz) :=  \inf_{u\in {\cal U}} \; \bbe_P[ u(\bdf(\bdz,\xi))].
$$
Let
$\hat{\vt}$ and $\tilde{\vt}$ be defined as in
(\ref{eq:PRO-x}) and (\ref{eq:PRO-x-1}).
Define
$$
{\cal U}(\bdz) := \{u\in {\cal U}: \bbe_P[u(\bdg(\bdz,\bdxi))] \geq 
c
\}.
$$
Since ${\cal U}(\bdz)\subset {\cal U}$, then
$\hat{v}(\bdz)\geq \tilde{v}(\bdz)$ for all $\bdz\in Z$.
Moreover, since
$\tilde{Z}\subset Z$, then
$$
\hat{\vt} = \max_{\bdz\in Z} \hat{v}(\bdz) \geq \max_{\bdz\in \tilde{Z}}
\tilde{v}(\bdz)  = \tilde{\vt}.
$$
Conversely, for any $\bdz\in \tilde{Z}$,
$
{\cal U}(\bdz)={\cal U}.
$
Thus, the assumption that $\hat{\bdz}\in \tilde{Z}$ implies that
$
{\cal U}(\hat{\bdz})={\cal U}
$
and subsequently $\hat{v}(\hat{\bdz}) = \tilde{v}(\hat{\bdz})$.
This shows
$
\hat{\vt} =\hat{v}(\hat{\bdz}) = \tilde{v}(\hat{\bdz}) \leq \tilde{\vt}
$
because $\hat{\bdz}\in \tilde{Z}$.
\hfill $\Box$

From a practical point of view, Proposition~\ref{Prop-equivalence}
is not useful 
in that we do not know the optimal solution $\hat{\bdz}$ and hence are unable to verify the condition $\hat{\bdz}\in \tilde{Z}$.
Consequently, 
it might be sensible to consider (\ref{eq:PRO-x})
as (\ref{eq:PRO-x-1}) might be too conservative.
Using the definition of ${\cal U}(\bdz)$, we 
can
write (\ref{eq:PRO-x}) succinctly as
\begin{equation}
\label{eq:PRO-x-DD}
    \inmat{(BUPRO-D)\quad} \max_{\bdz\in Z} \; \min_{u\in {\cal U}(\bdz)} \; \bbe_P[ u(\bdf(\bdz,\bdxi))].
\end{equation}
Problem (\ref{eq:PRO-x-DD}) looks as if the ambiguity set 
${\cal U}(\bdz)$ is decision-dependent. 
We propose to use the PLA
approach
to solve problem (\ref{eq:PRO-x-DD}).
In this case, 
\[
\mathcal{U}_N(\bdz):=\{u_N\in\mathcal{U}_N \,|\; \bbe_P[u_N(\boldsymbol{g}(\bdz,\bdxi))] \geq c\},
\]
where $\mathcal{U}_N$ is defined as in 
(\ref{eq:U_N-PLA}).
The approximate BUPRO 
can be subsequently written as
\begin{equation}
\label{eq:SPR-x-approx}
\inmat{(BUPRO-DN)} \quad \max_{\bdz\in Z}\;\min_{u\in {\cal U}_N(\bdz)}\; \bbe_P[u(\bdf(\bdz,\bdxi))]. 
\end{equation}
The inner minimization problem 
based on EPLA 
can be reformulated as an LP: 
\begin{align}
\label{eq:PRO-x-inner}
    \disp{ \min_{{\bm u} }} \; & \sum_{k=1}^K p_k \sum_{i=1}^{N_1-1} \sum_{j=1}^{N_2-1} \1_{T_{i,j}} (\bdf^k) 
    \lt[ u^{1l}_{i,j}(f_1^k,f_2^k) \1_{\lt[0,\frac{y_{j+1}-y_j}{x_{i+1}-x_i}\rt]} 
    \lt( \frac{f_2^k-y_j}{f_1^k-x_i} \rt) \rt. \nonumber \\
    & \lt. + u^{1u}_{i,j}(f_1^k,f_2^k) \1_{\lt( \frac{y_{j+1}-y_j}{x_{i+1}-x_i},+\infty \rt)} \lt(\frac{f_2^k-y_j}{f_1^k-x_i}\rt) \rt]  \nonumber \\
    \st \; & 
    \sum_{k=1}^K p_k \sum_{i=1}^{N_1-1} \sum_{j=1}^{N_2-1} \1_{T_{i,j}} (\bdg^k) 
    \lt[ u^{1l}_{i,j}(g_1^k,g_2^k) \1_{\lt[0,\frac{y_{j+1}-y_j}{x_{i+1}-x_i}\rt]} 
    \lt( \frac{g_2^k-y_j}{g_1^k-x_i} \rt) \rt. \nonumber  \\
    & \lt. + u^{1u}_{i,j}(g_1^k,g_2^k) \1_{\lt( \frac{y_{j+1}-y_j}{x_{i+1}-x_i},+\infty \rt)} \lt(\frac{g_2^k-y_j}{g_1^k-x_i}\rt) \rt] \geq c,  \\
    & \inmat{constraints} \; (\ref{eq-traform-paircom})-(\ref{eq-traform-norm1}), \nonumber
\end{align}
where ${\bm u}={\rm vec}\left((u_{i,j})_{1\leq i\leq N_1}^{1\leq j\leq N_2}\right)\in \R^{N_1N_2}$,
$\bdg^k:=\bdg(\bdz,\bdxi^k)=(g_1^k,g_2^k)^T$ with $g_1^k:=g_1(\bdz,\bdxi^k)$, $g_2^k:=g_2(\bdz,\bdxi^k)$.
Since (\ref{eq:PRO-x-inner})
is an LP
for fixed $\bdz$, we can use a Dfree method to solve (\ref{eq:SPR-x-approx}).
Similar formulations can be 
derived based on the IPLA approach.

\section{Numerical results}
\label{sec:numerical results}

We have carried out numerical 
tests on the performances of the proposed models
and computational schemes discussed in the 
previous sections 
by applying them to a portfolio optimization problem. 
In this section, we report the test results.

\subsection{Setup}
\label{eq:setup-1}
As an example of a real-life portfolio selection problem with uncertain project outcomes, we consider an application of the UPRO models in healthcare resource allocation problem studied by Airoldi~et~al.~\cite{airoldi2011healthcare}.
In this application, public health officials (PHO) decide on a portfolio of projects that seek to improve the quality of life.
Specifically, the health benefits of $n=8$ projects (access to dental, workforce development, primary prevention, Obesity training, CAMHS School, early detection and diagnostics, palliative \&
EOL, active treatment)
are evaluated 
through two attributes, 
commissioning areas of children and cancer.
Moreover, the outcomes of the projects are uncertain and 
represented by
discretely distributed
random 
vector $\bdxi^k=(\xi^k_1,\ldots,\xi^k_8)^T$ {\color{black}supported by $\Xi\subset \R^8$} with equal probabilities $p_k:=1/K$ for $k=1,\ldots, K$.
Let $\bdz=(z_1,\ldots,z_8)^T$ be the proportions of a fixed fund.
For the convenience of calculation,
we generate 
samples of 
$\bdxi^k$ 
by the uniform distribution over $[0,1]^8$.
We consider a situation where 
the PHO's utility of 
the 
bi-attribute outcomes 
is ambiguous and the optimal allocation
is based on the worst-case utility in 
ambiguity 
set $\calu$
\begin{equation*}
    \max_{\bdz\in Z} \min_{u\in\calu} \; \sum_{k=1}^K p_ku(\bdf(\bdz,\bdxi^k)),
\end{equation*}
where $f_1(\bdz,\bdxi^k):=\sum_{i=1}^{5} z_i \xi_i^k
\in [0,1]$, $f_2(\bdz,\bdxi^k):=\sum_{i=6}^{8} z_i \xi_i^k
\in [0,1]$ and $Z:=\{\bdz\in\R^8_+ : \sum_{i=1}^8 z_i=1\}$.
To examine the performance of BUPRO-N, we carry out the tests with a specified true utility function and investigate how the optimal value and the worst-case utility function converge as 
information 
about the PHO's utility preference increases.
We 
consider the true utility
$
u(x,y)=e^x-e^{-y}-e^{-x-2y}
$
defined over $[0,1]\times [0,1]$
and normalize it 
by setting $u^*(x,y) :=
(u(x,y)-u(0,0))/(u(1,1)-u(0,0))$. 
This function
satisfies the 
conservative 
property
(\ref{eq:conservative}), and
is convex 
w.r.t.~$x$ and concave 
w.r.t.~$y$.
Although the PHO is unaware that the preference can be characterized as this function, we assume that the decision of PHO never contradicts with results suggested by such a function
unless specified otherwise (we will remove this assumption in Section~\ref{subsec:preference-incon}), 
see similar assumption in \cite{AmD15}.
We may refine 
$\scru$ to 
a set of 
normalized 
non-decreasing
utility functions mapping from $[0,1]^2$ to $[0,1]$, and $\scru_N$ the corresponding set of PLA functions.
All of the tests are carried out in MATLAB R2022a installed on a PC (16GB, CPU 2.3 GHz) with an Intel Core i7 processor.
We use GUROBI and YALMIP \cite{lofberg2004yalmip} to solve the inner minimization problem (LP or MILP) and single MILP, 
and SURROGATEOPT to solve the outer maximization problem (unconstrained problem
(\ref{eq:MAUT-robust-N-dis}) and constrained problem (\ref{eq:SPR-x-approx})).

\subsection{Design of the pairwise comparison lotteries}
\label{sec:PC-design}

As we discussed earlier, the ambiguity set of utility functions $\calu_N$ is characterized by  available 
information about the DM's preferences.
We ask PHO
questions by 
showing preference between 
a risky lottery
with two
outcomes and a 
lottery 
with certain outcome
(we call it ``certain lottery'' 
following the terminology of
\cite{AmD15}),
denoted respectively by
\begin{equation}
\label{eq:lottery}
    \bdcz_1 = \lt\{
\begin{array}{ll}
    (\underline{x},\underline{y}) & \inmat{\;w.p.\;} 1-p, \\
    (\bar{x},\bar{y}) & \inmat{\;w.p.\;} p,
\end{array} 
\rt.
\inmat{ and }
\bdcz_2=(x,y) \inmat{\;w.p.\;} 1,
\end{equation}
where $\underline{x},\bar{x}, \underline{y}$ and $\bar{y}$
are fixed and 
 $(x,y)\in [\underline{x},\bar{x}]\times[\underline{y},\bar{y}]$
is randomly generated.
Since we assume that 
$u(\underline{x},\underline{y})=0$ and $u(\bar{x},\bar{y})=1$,
the only parameters to be identified are $x,y,p$, 
so that the question is properly posed.
Observe that 
$$
\bbe_{\mathbb{P}}[u(\bdcz_1(\omega))]=(1-p) u(\underline{x},\underline{y})+ p u(\bar{x},\bar{y})=p \inmat{\quad and\quad} \bbe_{\mathbb{P}}[u(\bdcz_2(\omega))]=u(x,y). 
$$
Thus the question is down to checking 
whether inequality $u(x,y) \geq p$ holds or not.
Next,
we turn to discuss how to generate $M$ lotteries, or more specifically how to
set values for $x$, $y$ and $p$.
We generate randomly $M_1$ points of the first attribute including $\underline{x}$, $\bar{x}$, and $M_2$ points of the second attribute including $\underline{y}$ and $\bar{y}$. Thus
the number of 
the 
certain 
lotteries 
is at most $M=M_1M_2-2$.
Let 
$$
S:=\{(x_{i_l},y_{j_l}), i_l\in\{1,\ldots,M_1\}, j_l\in\{1,\ldots,M_2\}, l=1,\ldots,M \}
$$
be the set of all 
certain
lotteries except points $(\underline{x},\underline{y})$, $(\bar{x},\bar{y})$ and $\calu_{N}^{l-1}$ be the set of all piecewise linear utility functions which are consistent to the
previously
generated $l-1$ questions.
Assume that 
the $l$th 
lottery with the certain outcome
is $\bdcz_2^l=(x_{i_l},y_{i_l})$. 
Define 
\begin{equation}
\label{eq:I1-I2}
    I_1^l:=\min_{u\in\scru_N\cap\calu_{N}^{l-1}} u(\bdcz_2^l) 
    \inmat{\quad and\quad} 
    I_2^l:=\max_{u\in\scru_N\cap\calu_{N}^{l-1}} u(\bdcz_2^l).
\end{equation}
Since $u(x_{i_l},y_{i_l})\in[0,1]$, then $I_1^l,I_2^l\in[0,1]$. 
We set $p^l:=\frac{I_1^l+I_2^l}{2}$,
and 
use the true utility function $u^*$ to 
check whether inequality
\begin{equation}
\label{eq:lottery-l}
    u^*(x_{i_l},y_{i_l}) =\bbe_{\mathbb{P}}[u^*(\bdcz_2^l(\omega))] \geq \bbe_{\mathbb{P}}[u^*(\bdcz_1^l(\omega))] = p^l
\end{equation}
holds or not.
If it holds,
then $\bdcz_2^l$ is preferred 
to $\bdcz_1^l$.
The following algorithm describes the procedures 
for constructing
$\calu_{N}=\calu_N^M$.

\vspace{0.6em}
\begin{breakablealgorithm}
\caption{}
{
\noindent
\textbf{Initialization.}
Set $m_1:=1, m_2:=1, l:=1$, 
$\calu_N^0:=\scru_N$
and $S:=\emptyset$. 
}
\begin{algorithmic}[1]
\begin{small}
\STATE
Choose two positive integers $M_1$ and $M_2$ 
as the numbers of the 
gridpoints of the two attributes.
Generate $M_1-2$ points 
within $[\underline{x},\bar{x}]$ and $M_2-2$ points
within $[\underline{y},\bar{y}]$ randomly
using the uniform distribution; 
sort them out 
in 
increasing order
of their values
and label them by $x_i, i=1,\ldots,M_1-2$ and $y_j, j=1,\ldots,M_2-2$.
Let ${\cal X}:= \{\underline{x},x_1,\ldots,x_{M_1-2},\bar{x}\}$ and ${\cal Y}:= \{\underline{y},y_1,\ldots,y_{M_1-2},\bar{y}\}$, and 
let ${\cal X}\times{\cal Y}:=\{(x_i,y_j),x_i\in{\cal X}, y_j\in{\cal Y}\}$ be the set of the gridpoints.

\STATE
Let the $l$th certain lottery be $\bdcz_2^l=(x_{i_l},y_{j_l})$,
solve the 
problem (\ref{eq:I1-I2}) to obtain $I_1^l$ and $I_2^l$.
Let $I^l=[I_1^l,I_2^l]$, $p^l=\frac{I_1^l+I_2^l}{2}$ and 
$S=S\cup\{\bdcz_1^l,\bdcz_2^l\}$.
\STATE
If 
$p^l\leq u^*(x_{i_l},y_{j_l})$,
then 
$$
{\cal U}_N^{l}:={\cal U}_N^{l-1} \bigcap \lt\{u_N\in \mathscr{U}_N: p^l\leq u_N(x_{i_l},y_{j_l}) \rt\}.
$$
Otherwise, 
$$
{\cal U}_N^{l} :={\cal U}_N^{l-1} \bigcap \lt\{u_N\in \mathscr{U}_N: p^l\geq u_N(x_{i_l},y_{j_l})\rt\}.
$$
Set $l:=l+1$, and go
to Step 1.
\end{small}
\end{algorithmic}
\end{breakablealgorithm}
\vspace{0.5em}

Steps 1-2 generate a lottery for pairwise comparison. 
Note that the minimization problem in (\ref{eq:I1-I2}) can be formulated as 
\begin{subequations}
\label{eq-lottery-I_1}
\begin{align}
    I_1^l = \min_{{\bm u}}\;
\; & u_{i_l,j_l}
    \nonumber \\
    \st \;\; & 
    h_{l'} (p_{l'}- u_{i_{l'},j_{l'}})\leq 0, l'=0,\ldots,l-1, \label{eq-lottery} \\
    & \frac{u_{i+1,j}-u_{i,j}}{x_{i+1}-x_i} \geq \frac{u_{i,j}-u_{i-1,j}}{x_i-x_{i-1}}, i=2,\ldots,M_1-1, j=1,\ldots,M_2, \label{eq-single-concave} \\
    & \frac{u_{i,j+1}-u_{i,j}}{y_{j+1}-y_j} \leq \frac{u_{i,j}-u_{i,j-1}}{y_{j}-y_{j-1}}, i=1,\ldots,M_1, j=2,\ldots,M_2-1, \label{eq-single-convex} \\
    & \inmat{constraints\;} (\ref{eq-traform-mon1})- (\ref{eq-traform-norm1}), \notag
\end{align}
\end{subequations}
where  ${\bm u}:=(u_{1,1},\cdots,u_{N_1,1},\cdots,u_{1,N_2},\cdots,u_{N_1N_2})^T$,
(\ref{eq-lottery}) requires the answer to the $l$th question to be consistent with the previous $l-1$ questions 
(if $\bdcz_1^{l'}$ is preferred, then (\ref{eq:lottery-l}) holds for $l=l'$ and we set $h_{l'}=1$, otherwise we set $h_{l'}=-1$),
(\ref{eq-single-concave}) and (\ref{eq-single-convex}) comply with the assumption that the single-attribute utility function $u(\cdot,\hat{y})$ is concave and $u(\hat{x},\cdot)$ is convex for any fixed $\hat{x}\in X$ and $\hat{y}\in Y$.
Step 3 asks the DM to choose 
between the risky lottery and the certain lottery. 
Here the true utility function $u^*$ 
(defined in Section~\ref{eq:setup-1})
is used to ``act as the DM''. 
After the DM makes a choice, an expected utility inequality is created and added to the ambiguity set $\calu_N$.
Since $p^l$ is chosen as the midpoint of $I^l$, we deduce that the true utility function value at $(x_{i_l},y_{j_l})$ lies within the right or left half of the interval $I^l$ and the pairwise comparison effectively reduces the ambiguity set by ``half'' in the sense that those $u_N$ whose values (at point $(x_{i_l},y_{j_l})$) lie within the other half of the interval $I^l$ are excluded from the ambiguity set.

\begin{example}
We use a simple example to explain the above steps where the true utility function $u^*$ (defined in Section~\ref{eq:setup-1})
is defined over $[0,1]^2$ and the piecewise utility functions have $N=M_1M_2=6$ gridpoints including $(0,0)$ and $(1,1)$.
We randomly generate one point in $[0,1]$ for the second attribute as the non-end gridpoints.
Then ${\cal X}=\{0,1\}$ and ${\cal Y}=\{0,0.3706,1\}$.
The number of questions is $M=M_1 M_2-2=4$.

\noindent
\textbf{Lottery 1} $(l=1)$. Set $i_1:=1, j_1:=2$ and $(x_{i_1},y_{j_1})=(0,0.3706)$.
By solving (\ref{eq-lottery-I_1}) and the corresponding maximization problem, we obtain that $[I_1^1,I_2^1]=[0,1]$ and set $p^1=0.5$.
By checking $u^*(0,0.3706)=0.252\leq p^1$, we set $h_1:=-1$.

\noindent
\textbf{Lottery 2} $(l=2)$. Set $i_2:=1, j_2:=3$ and $(x_{i_2},y_{j_2})=(0,1)$.
Solve (\ref{eq-lottery-I_1}) and the corresponding maximization problem to obtain $[I_1^2,I_2^2]=[0,0.5]$,
so $p^2=0.25$.
By checking whether $u^*(0,1)=0.454\geq p^2$ or not, we set $h_2:=1$.

\noindent
\textbf{Lottery 3} $(l=3)$. 
Set $i_3:=2, j_3:=1$ and $(x_{i_3},y_{i_3})=(1,0)$.
We obtain $[I_1^3,I_2^3]=[0.5,1]$, 
and $p^3=0.75$.
By checking $u^*(1,0)=0.712\leq p^3$, we set $h_3:=-1$.

\noindent
\textbf{Lottery 4} $(l=4)$. 
Set $i_4:=2, j_4:=2$ and $(x_{i_4},y_{j_4})=(1,0.3706)$ to obtain $[I_1^4,I_2^4]=[0.75,1]$,
and $p^4=0.875$.
By checking $u^*(1,0.3706)=0.864\leq p^4$, we set $h_4:=-1$.
\end{example}

\subsection{Convergence results}

In this subsection, we 
investigate 
the convergence of 
the 
worst-case approximate utility functions of the unconstrained problem (\ref{eq:MAUT-robust-N}) and the constrained optimization problem (\ref{eq:SPR-x-approx}) under EPLA and IPLA schemes as $N$ increases.

\textbf{(i) EPLA and IPLA for unconstrained problem (\ref{eq:MAUT-robust-N}).}

\begin{figure}[!ht]
  \centering
  \vspace{-0.5cm}
  \subfigure{
    \includegraphics[width=0.3\linewidth]{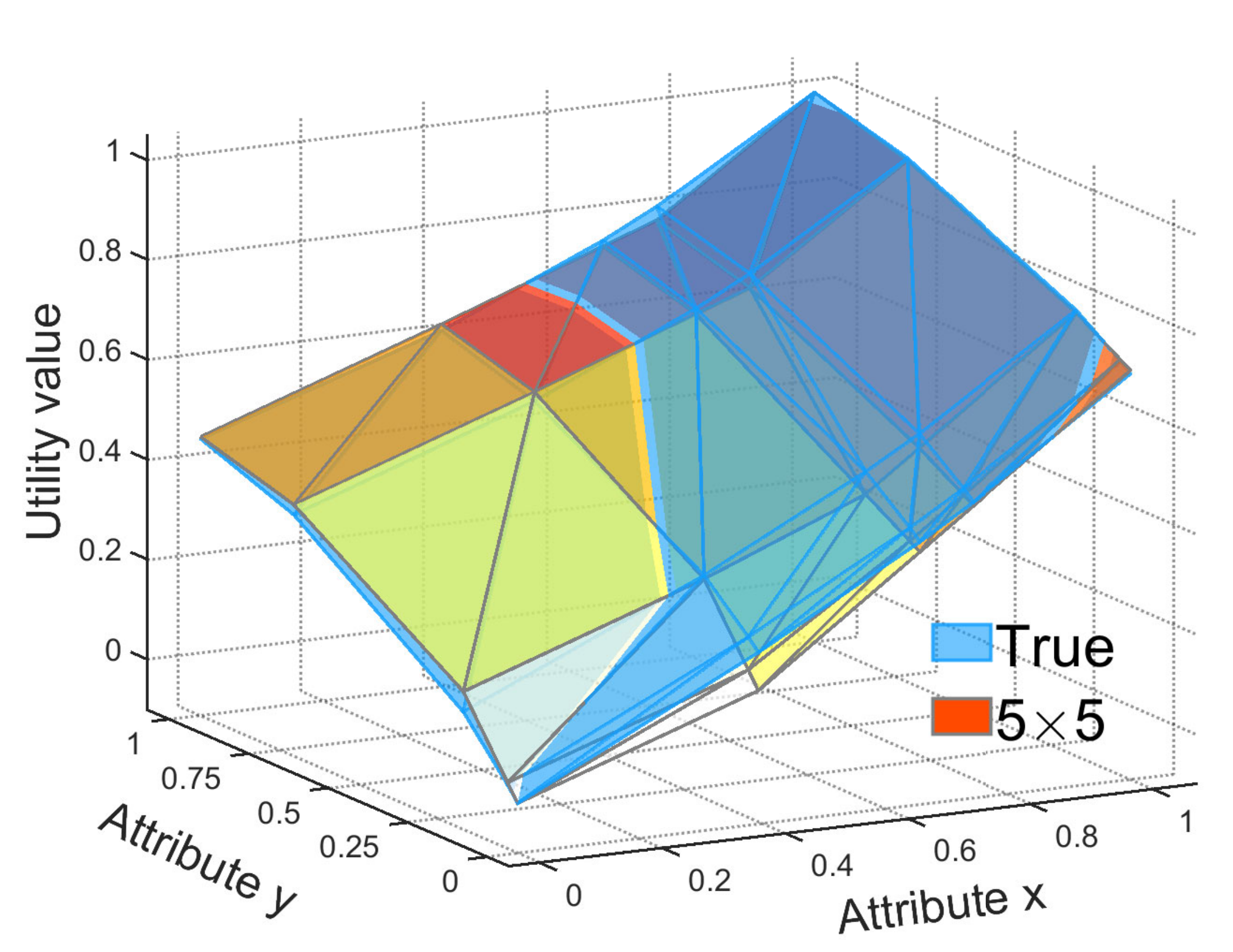}
  } 
  \hspace{-0.5em}
  \subfigure{
    \includegraphics[width=0.3\linewidth]{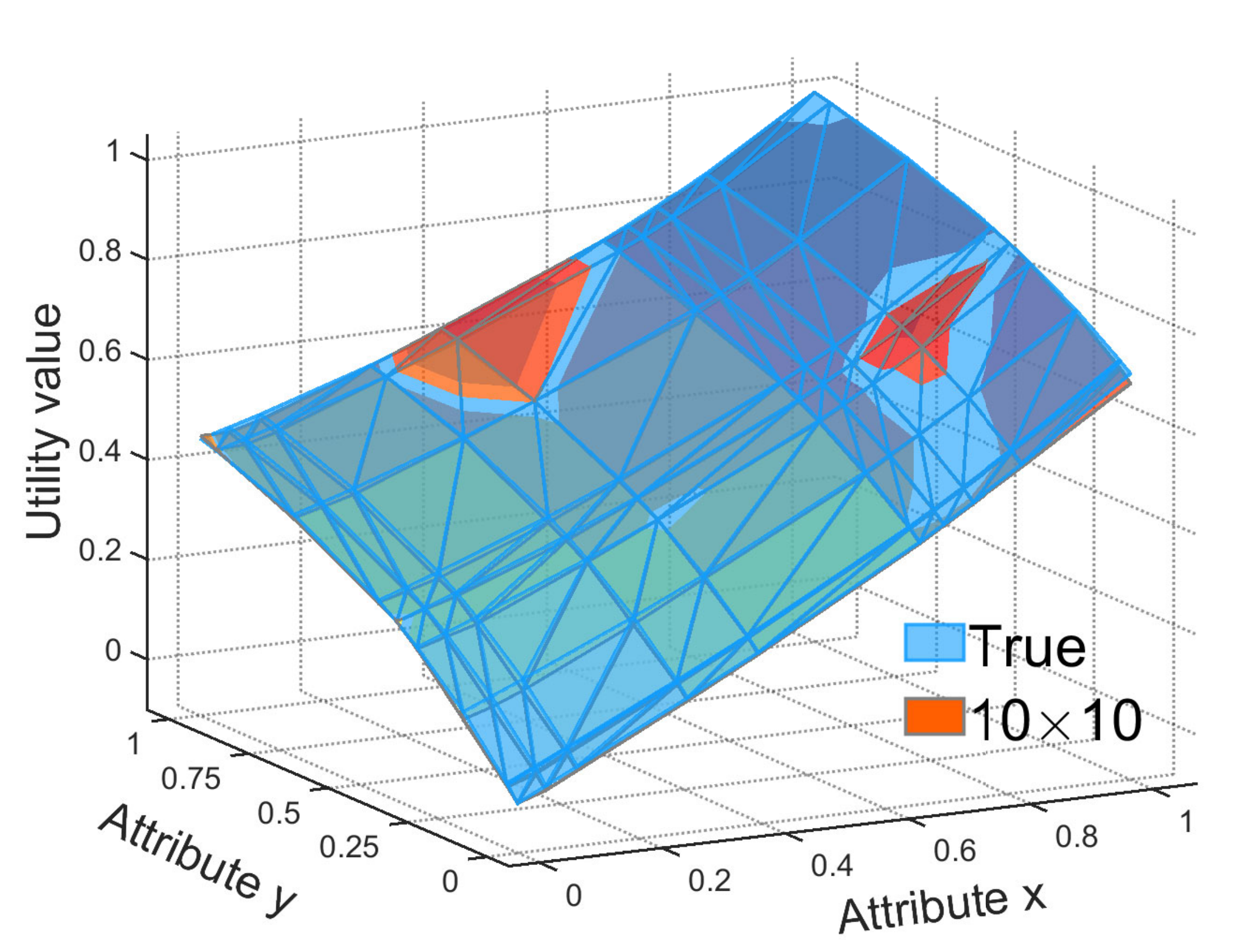}
  }
  \hspace{-0.5em}
  \subfigure{
    \includegraphics[width=0.3\linewidth]{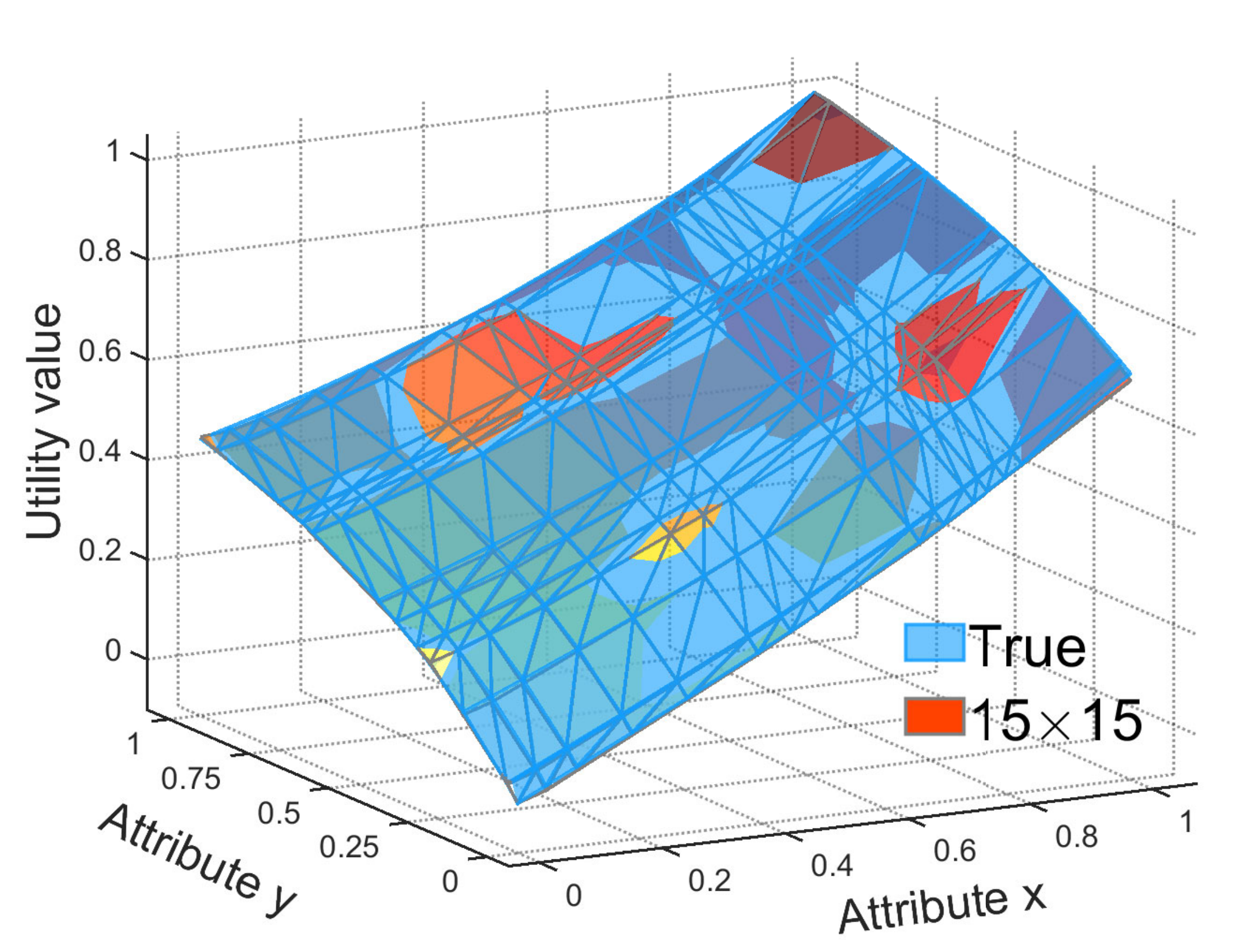}
  }
  \vspace{-0.2cm}
  \captionsetup{font=footnotesize}
  \caption{{\bf Type-1 EPLA}: the convergence of the worst-case utility function of EPLA model (\ref{eq:PRO-N-reformulate}) to the
  true utility function (in blue) as the number of questions 
  increases
  from $5\times 5$ to $15\times 15$.} 
  \label{fig-utility-main} 
  \vspace{-0.5cm}
\end{figure}

\begin{figure}[!ht]
  \centering
  \vspace{-0.5cm}
  \subfigure{
    \label{subfig-utility-b-counter}
    \includegraphics[width=0.3\linewidth]{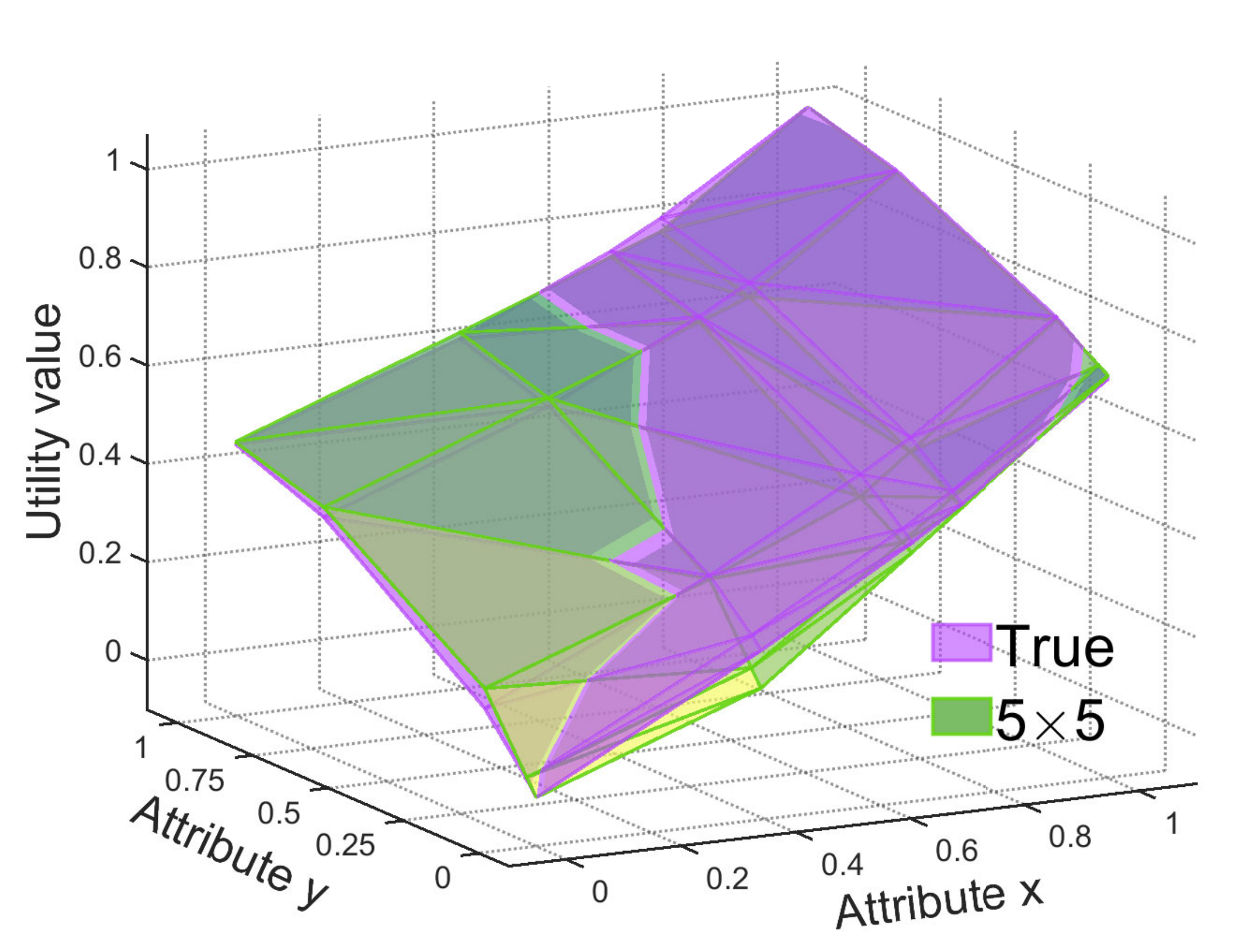}
  } \hspace{-0.7em}
  \subfigure{
    \label{subfig-utility-c-counter}
    \includegraphics[width=0.3\linewidth]{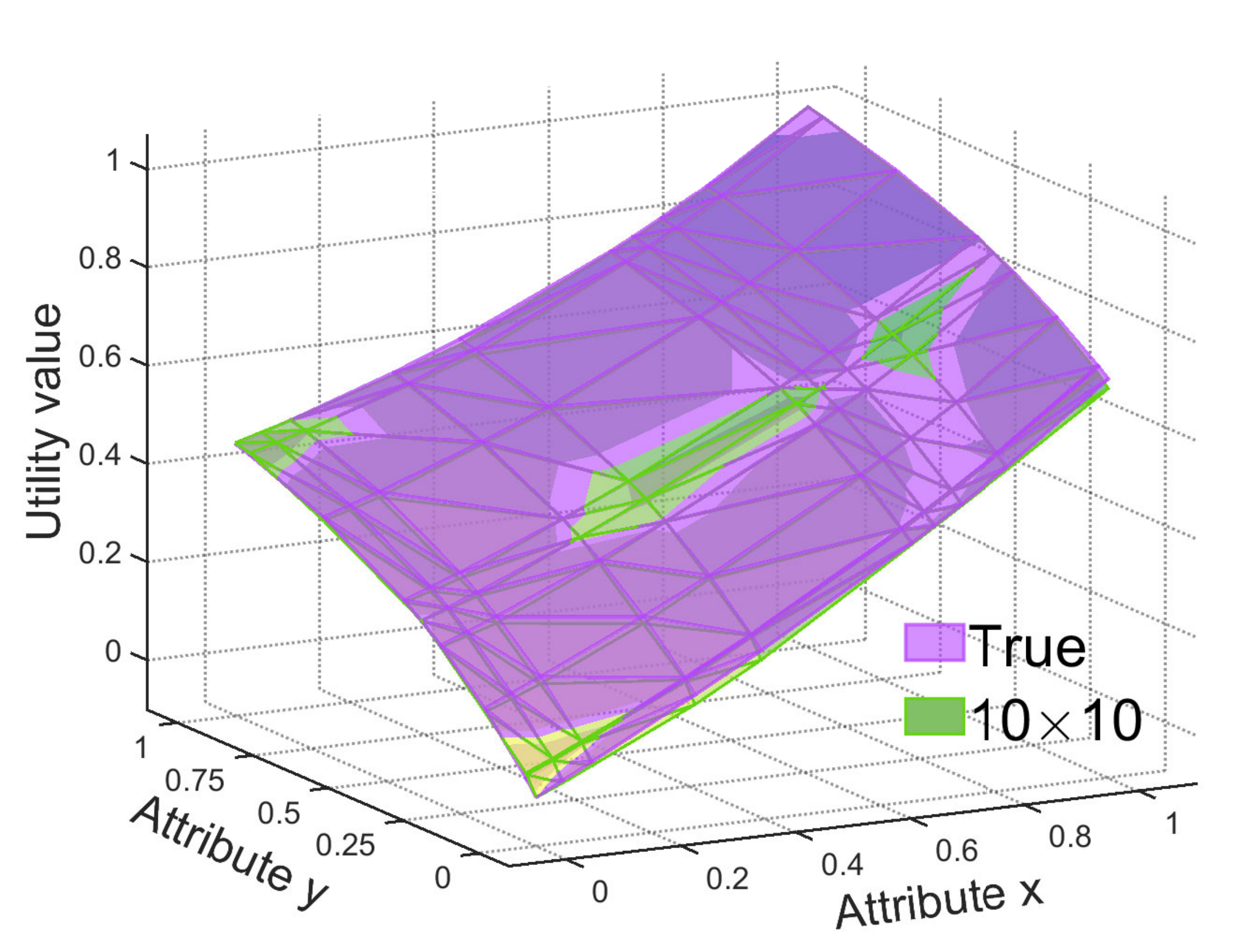}
  } \hspace{-0.7em}
  \subfigure{
    \label{subfig-utility-d-counter}
    \includegraphics[width=0.3\linewidth]{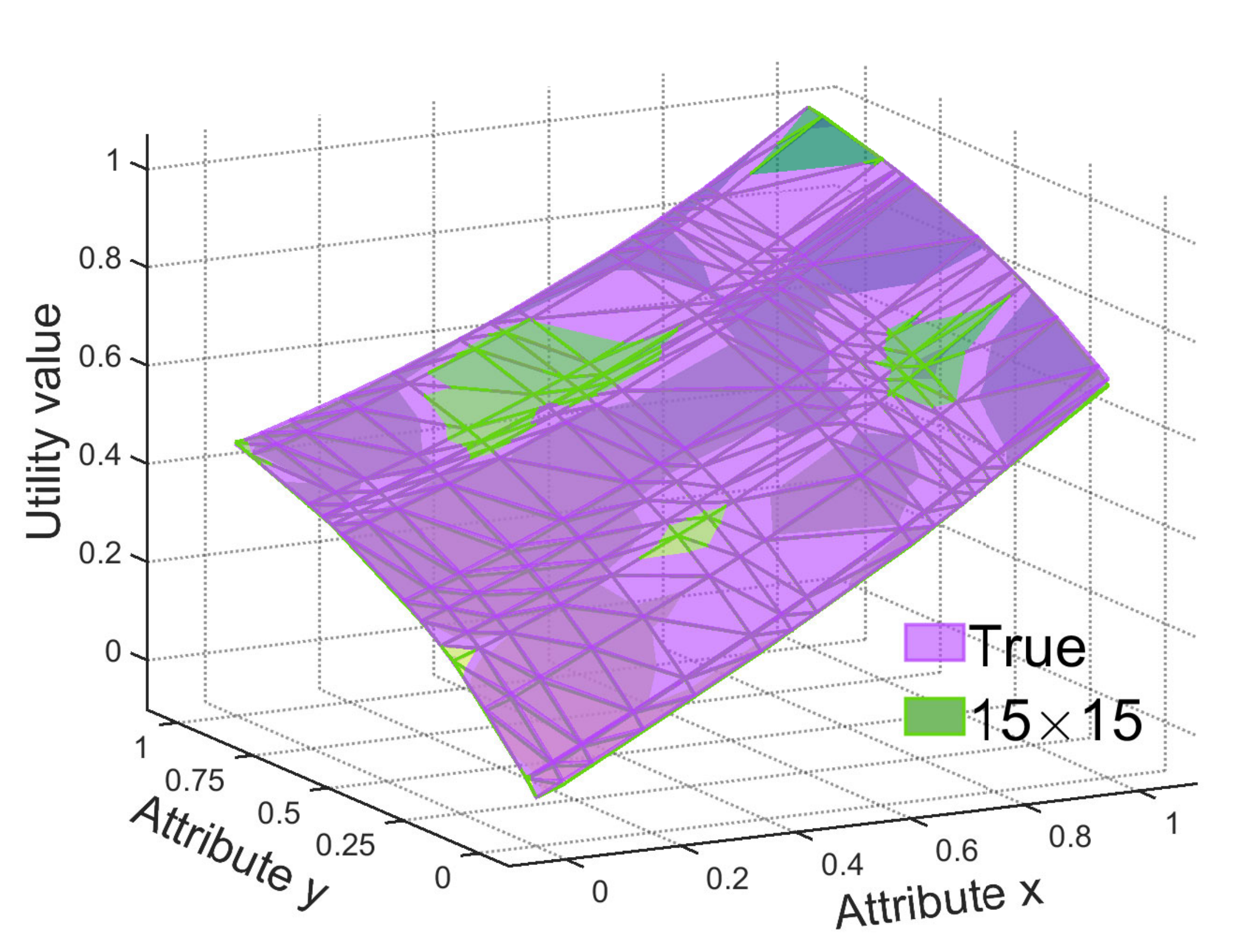}
  }
  \captionsetup{font=footnotesize}
  \caption{{\bf Type-2 EPLA}: the convergence of the worst-case utility function of Type-2 EPLA model
  to the true utility function (in purple) as the number of questions 
  increases from $5\times 5$ to $15\times 15$.}
  \label{fig-utility-counter} 
\end{figure}

\begin{figure}[!ht]
  \centering
  \vspace{-0.5cm}
  \subfigure{
    \label{subfig-utility-b-mixed}
    \includegraphics[width=0.3\linewidth]{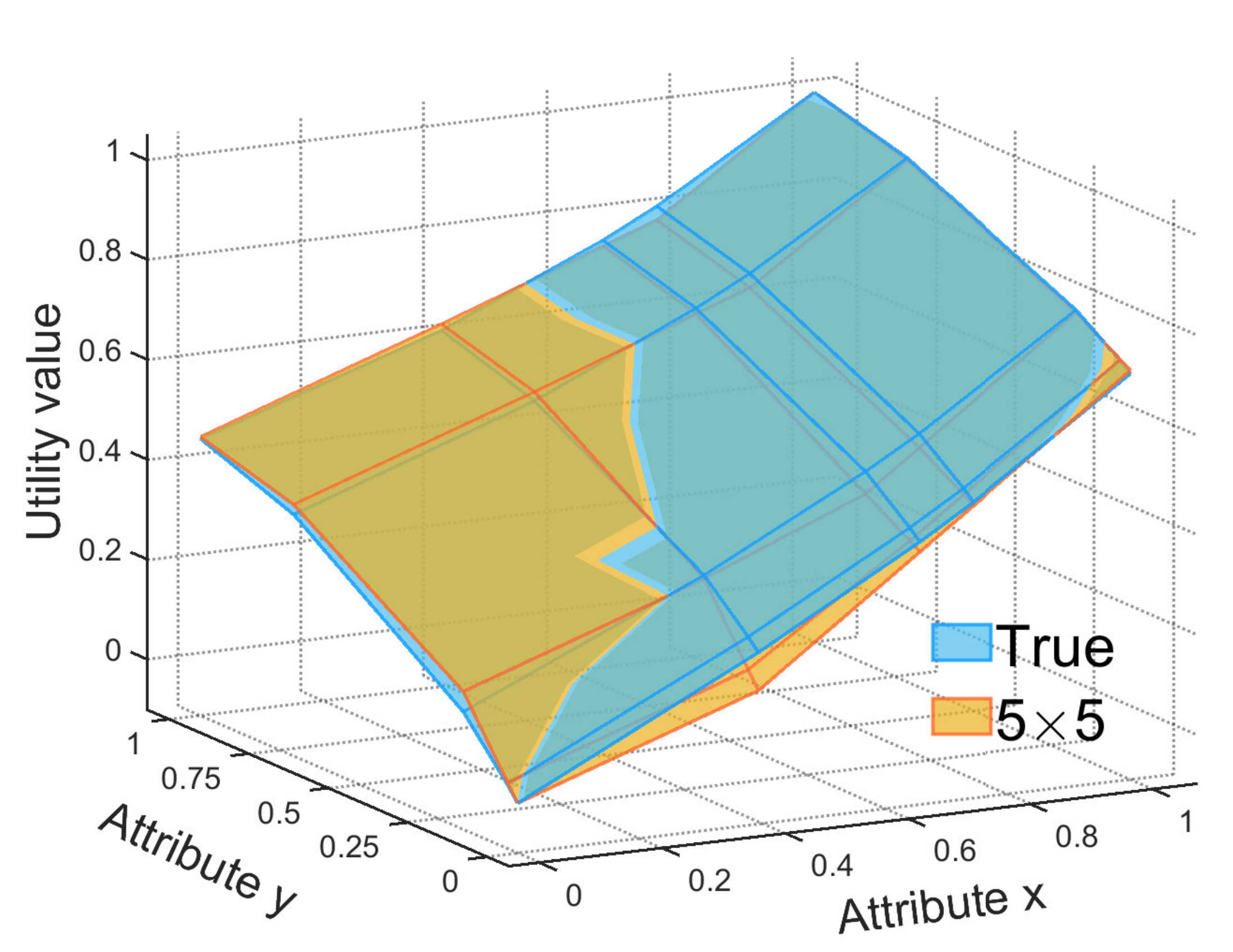}
  } \hspace{-0.7em}
  \subfigure{
    \label{subfig-utility-c-mixed}
    \includegraphics[width=0.3\linewidth]{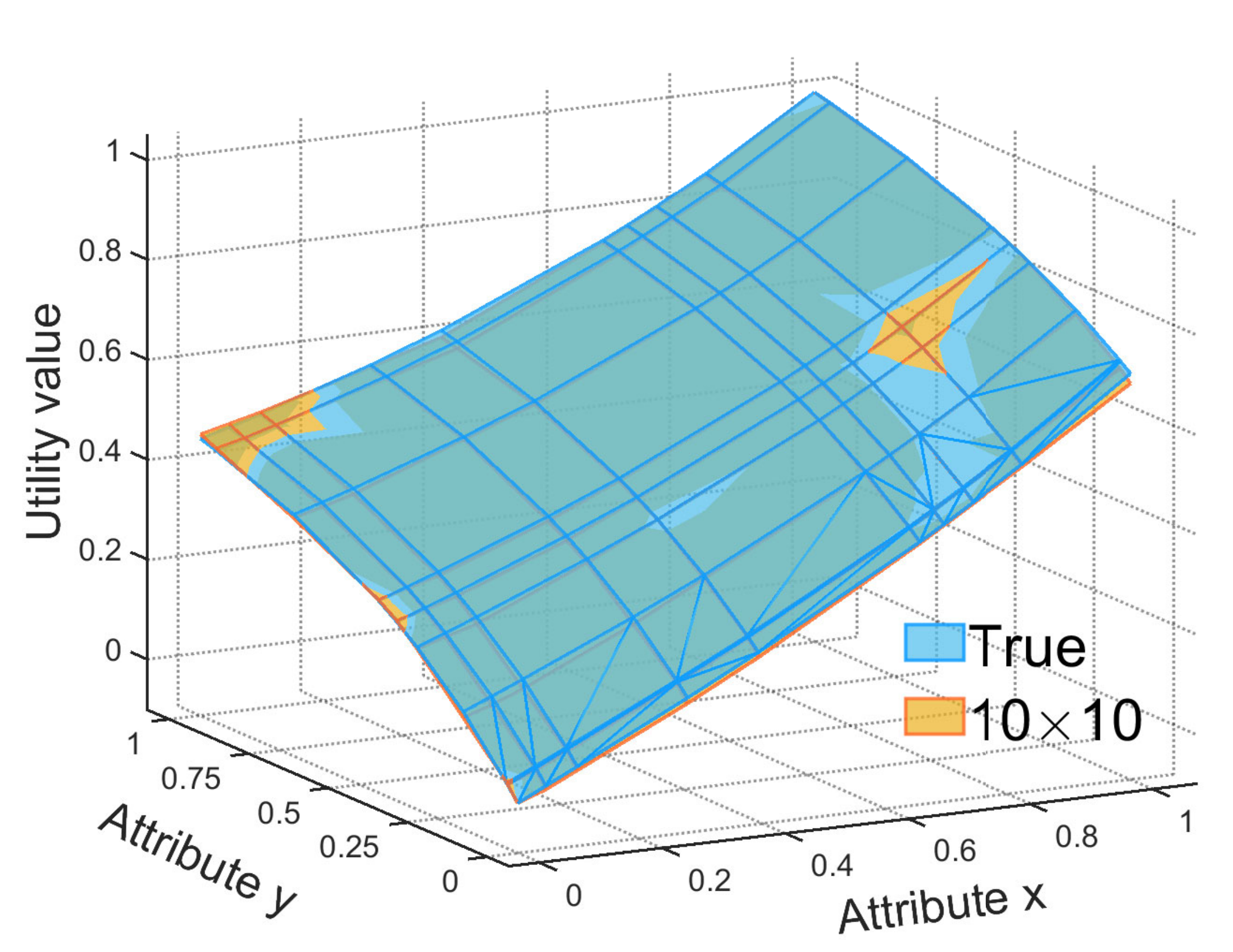}
  } \hspace{-0.7em}
  \subfigure{
    \label{subfig-utility-d-mixed}
    \includegraphics[width=0.3\linewidth]{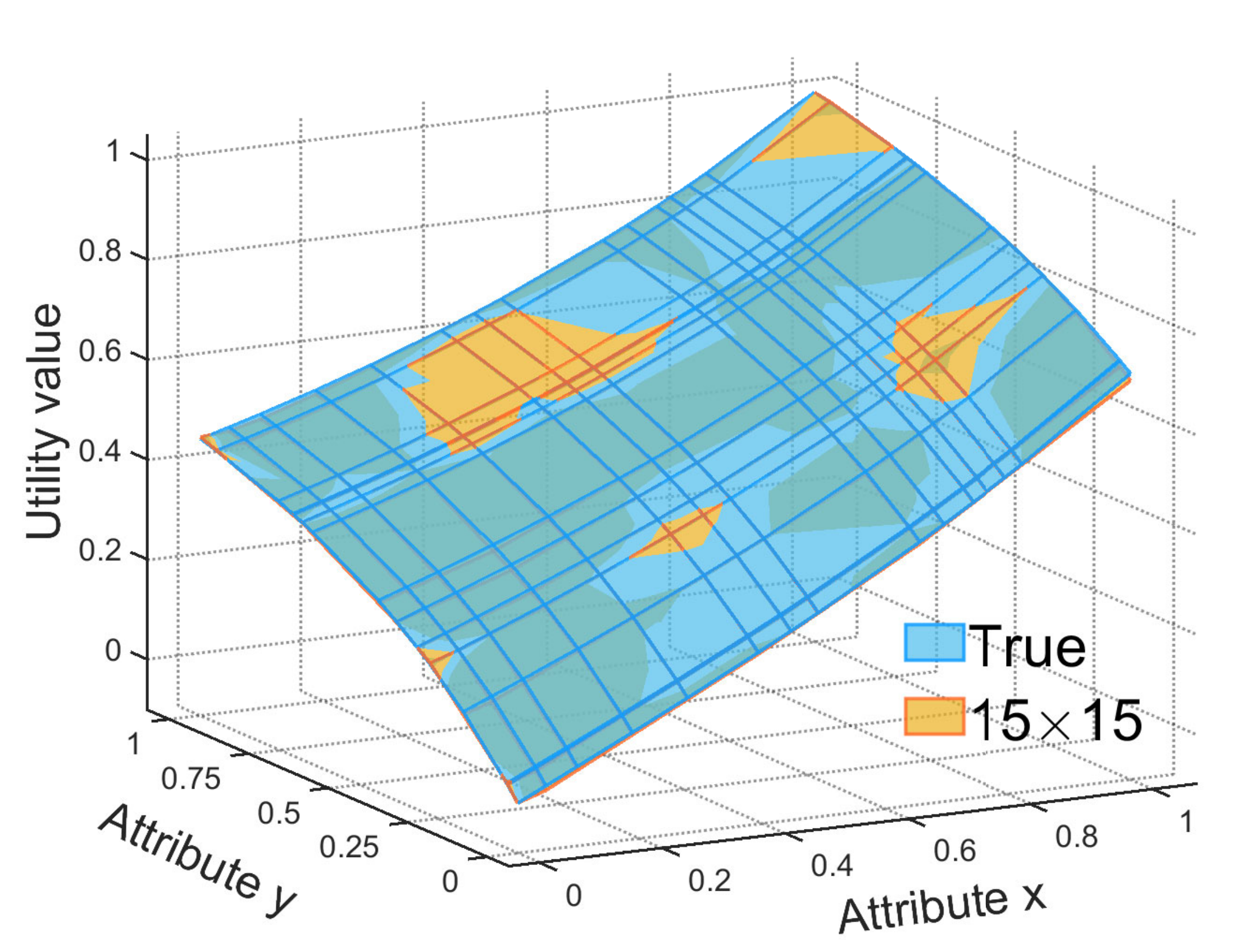}
  }
  \vspace{-0.2cm}
  \captionsetup{font=footnotesize}
  \caption{{\bf Mixed-type EPLA}:
 the convergence of the worst-case utility functions of mixed-type EPLA model to the true utility function (in green).
  The cell with no diagonal line means 
that Type-1 and Type-2 PLAs 
coincide because
in this case 
${\bm f}(\bdz,\bdxi^k)$ does not fall into the cell
for $k=1,\cdots,K$.
  }  
  \label{fig-utility-mixed} 
\end{figure}

\underline{EPLA approach.}
We begin by
examining
the performance of the EPLA approach 
with different 
types of partitions 
discussed in Section~\ref{sec:numer-methods}.
We assume
the number of the scenarios of  $\bdxi$ 
is  $K=1000$.
The convergence results are displayed in 
Figures~\ref{fig-utility-main}-\ref{fig-utility-mixed},
and Tables~\ref{tab-result-N}-\ref{tab:distance}.
Figures~\ref{fig-utility-main}-\ref{fig-utility-mixed} 
depict
the true utility function and the worst-case 
 utility functions for Type-1 PLA (Figure~\ref{fig-utility-main}), 
Type-2 PLA (Figure~\ref{fig-utility-counter}), 
and mixed-type PLA (see Remark~\ref{rem:BUPRO-DF}~(ii) for the definition) in Figure~\ref{fig-utility-mixed}.
We can see that the worst-case utility functions move closer and closer 
to the true utility function as more questions 
are asked in 
all the three cases,
which is in accordance
with our anticipation in Remark~\ref{rem:distance}~(iii) and Table~\ref{tab:distance}. 
Table~\ref{tab-result-N} displays the optimal solutions, the optimal values, 
the errors of the 
optimal values (which is defined as 
the difference between the true and the approximate optimal value), and computation time (CPU time).
We find that the 
optimal
values
increase as the number of queries increases.
This 
is
because
the ambiguity set $\calu_N$ becomes smaller as the number of queries increases.
Moreover, 
the errors
decrease as the number of 
questions  increases. 
The optimal values in the Type-1 PLA and mixed-type PLA are smaller than that of the Type-2 PLA
in that the conservative condition
makes the utility value of Type-2 
PLA
larger than the other cases, 
see Figure~\ref{fig-division}.

\begin{table}[!ht]
\tiny
    \centering
    \captionsetup{font=scriptsize}
    \caption{Computational results of BUPRO-N problem ($K=1000$, $\vt^*=0.3392$)}
    \vspace{-0.2em}
    \renewcommand\arraystretch{1.1}
    \begin{threeparttable}
    \resizebox{0.8\linewidth}{!}
    {
    \begin{tabular}{c|ccccc}
        \hline
         \textbf{EPLA} & 
        Lotteries & Optimal solutions & Optimal values & Error & CPU time (s) \\
        \hline
        \multirow{3}{*}{\makecell{Type-1 \\ }} & $5\times5$ & $[0,0,0,1,0,0,0,0]$ & 0.3122 & 0.0270 & 113.6 \\
        & $10\times10$ & $[0,0,0,0.955,0,0,0.016,0.029]$ & 0.3321 & 0.0071 & 151.6 \\
        & $15\times15$ & $[0,0,0,1,0,0,0,0]$ & 0.3349 & 0.0043 & 223.3 \\
        \hline
        \multirow{3}{*}{\makecell{Type-2 \\ }} &
        $5\times5$ & $[0,0,0,1,0,0,0,0]$ &  0.3122 & 0.0270 & 115.2 \\
        & $10\times10$ & $[0,0,0,0.961,0,0,0.008,0.031]$ & 0.3324 & 0.0068 & 164.7 \\
        & $15\times15$ & $[0,0,0,1,0,0,0,0]$ & 0.3349 & 0.0043 & 220.5 \\
        \hline
        \multirow{3}{*}{\makecell{Mixed-type \\ }} &
        $5\times5$ & $[0,0,0,1,0,0,0,0]$ &  0.3122 & 0.0270 & 886.2 \\
        & $10\times10$ & $[0,0,0,0.955,0,0.003,0.009,0.033]$ & 0.3321 & 0.0071 & - \\
        & $15\times15$ & $[0,0,0,1,0,0,0,0]$ & 0.3349 & 0.0043 & - \\
        \hline
    \end{tabular}
    }
    \begin{tablenotes}
    \raggedleft
        \item `-' implies runtime $>$ 3600s. 
    \end{tablenotes}
    \end{threeparttable}
    \vspace{-0.3cm}
    \label{tab-result-N}
\end{table}

\begin{table}[!htbp]
   \tiny
    \centering
    \captionsetup{font=scriptsize}
    \caption{{\bf EPLA:} upper bound for $\dd_{\mathscr{G}_I}(u^*,{u}_N^*)$ and distance $\dd_{\mathscr{G}_I}({u}_N^{*},u_{\rm worst}^N)$}
    \vspace{-0.2cm}
  \renewcommand\arraystretch{1.2}
    \begin{threeparttable}
    \resizebox{0.9\linewidth}{!}
    {
    \begin{tabular}{ccccccc}
        \hline
      Lotteries &
        $L(\beta_{N_1}+\beta_{N_2})$
        & $\dd_{\mathscr{G}_I}({u}_N^*,u_{\rm worst}^N)$ (Type-1) &  $\dd_{\mathscr{G}_I}({u}_N^*,u_{\rm worst}^N)$ (Type-2) & $\dd_{\mathscr{G}_I}({u}_N^*,u_{\rm worst}^N)$ (Mixed-type) \\
        \hline
       $5\times 5$ & $1.8541$ & $0.0763$ & $0.0763$ & $0.0763$\\
         $10\times 10$ & $1.0611$ & $ 0.0233$ &  $0.0306$ & $0.0226$\\
        $15\times 15$ & $0.7682$ & $0.0141$ &  $0.0230$ & $0.0141$ \\
        \hline
    \end{tabular}
    }
    \end{threeparttable}
    \vspace{-0.5cm}
    \label{tab:distance}
\end{table}

\begin{figure}[!ht]
  \centering
  \vspace{-1.5em}
      \subfigure
      {
    \includegraphics[width=0.23\linewidth]{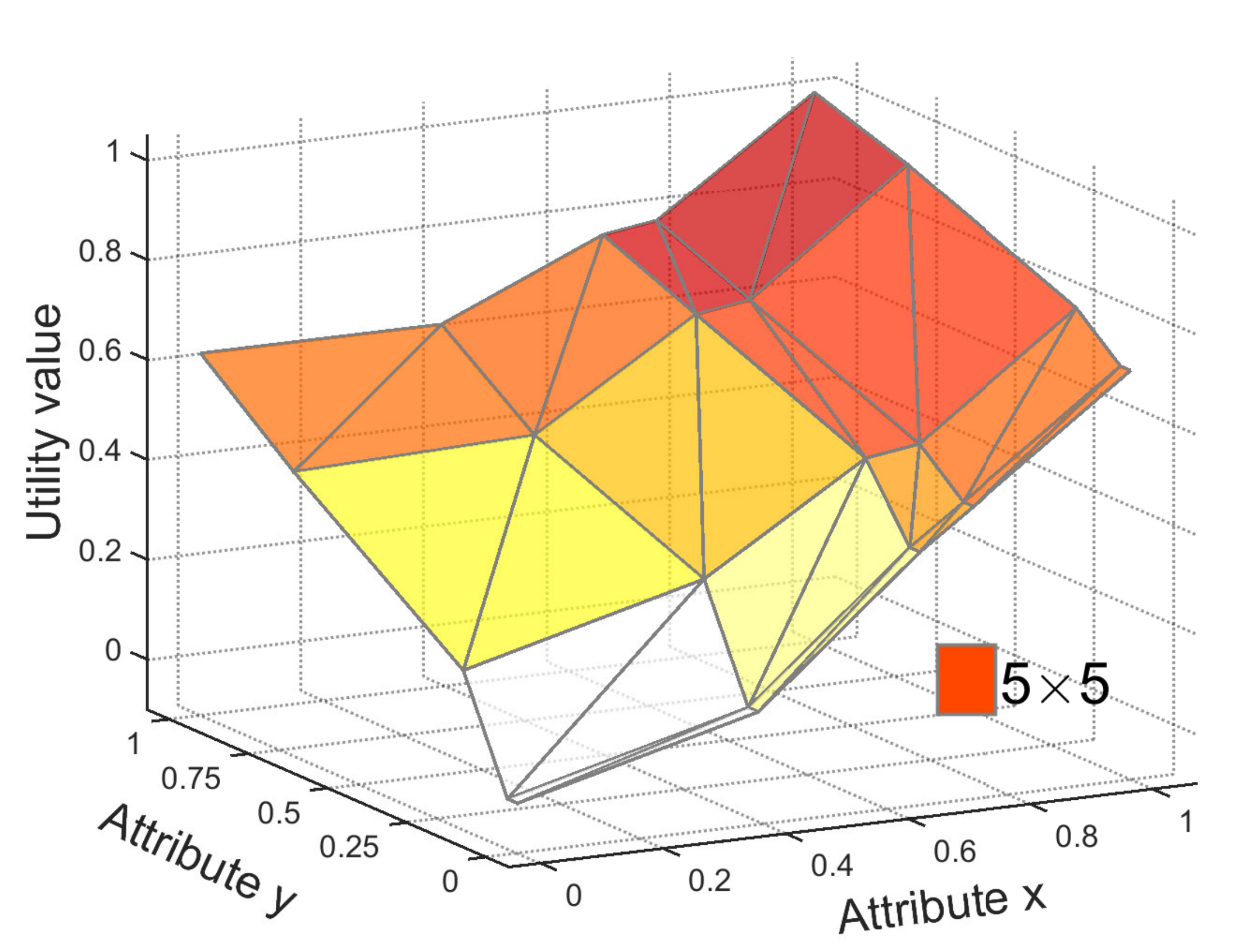}
  }
  \hspace{-0.5em}
    \subfigure
    {
    \includegraphics[width=0.23\linewidth]{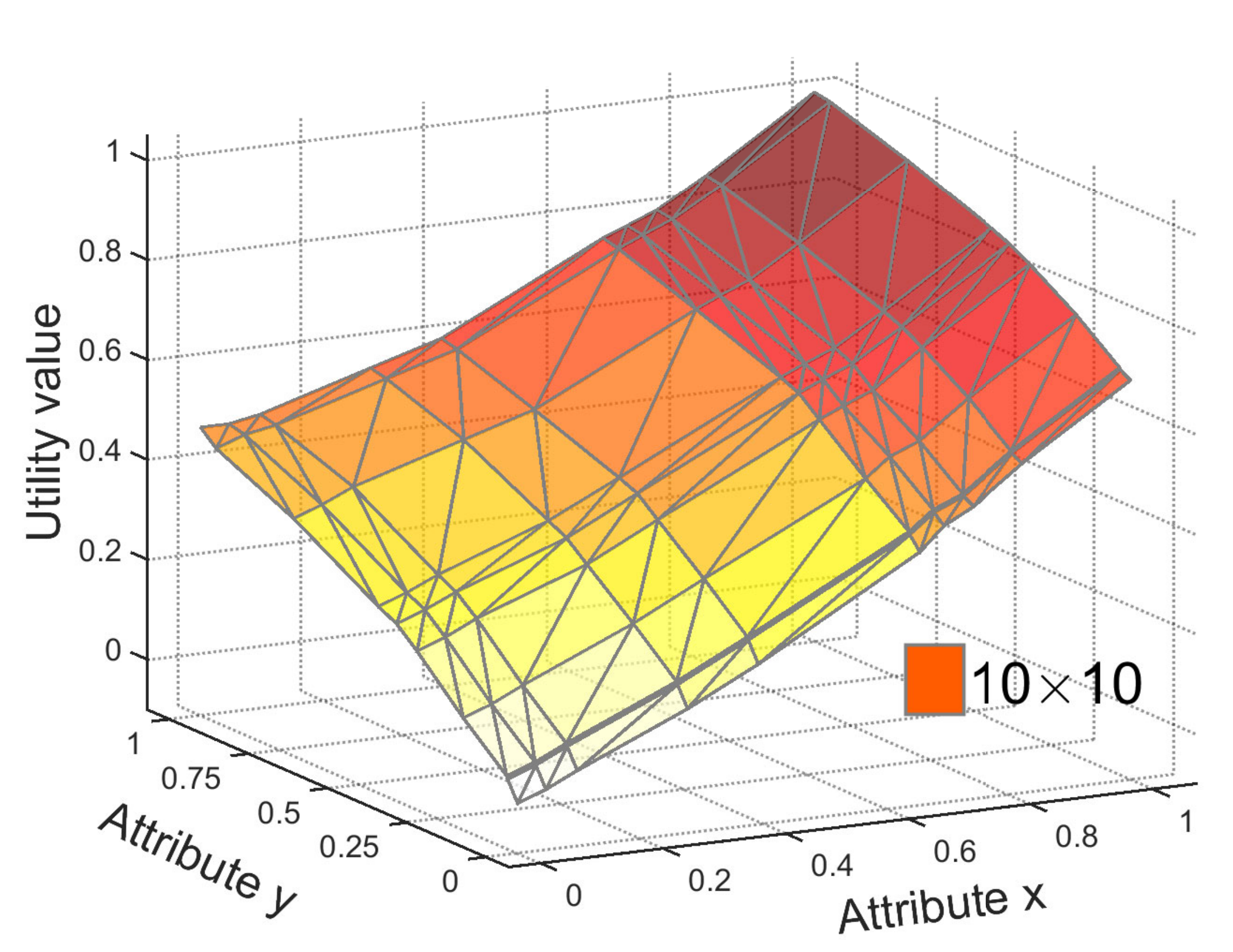}
  }
    \hspace{-0.5em}
    \subfigure
    {
    \includegraphics[width=0.23\linewidth]{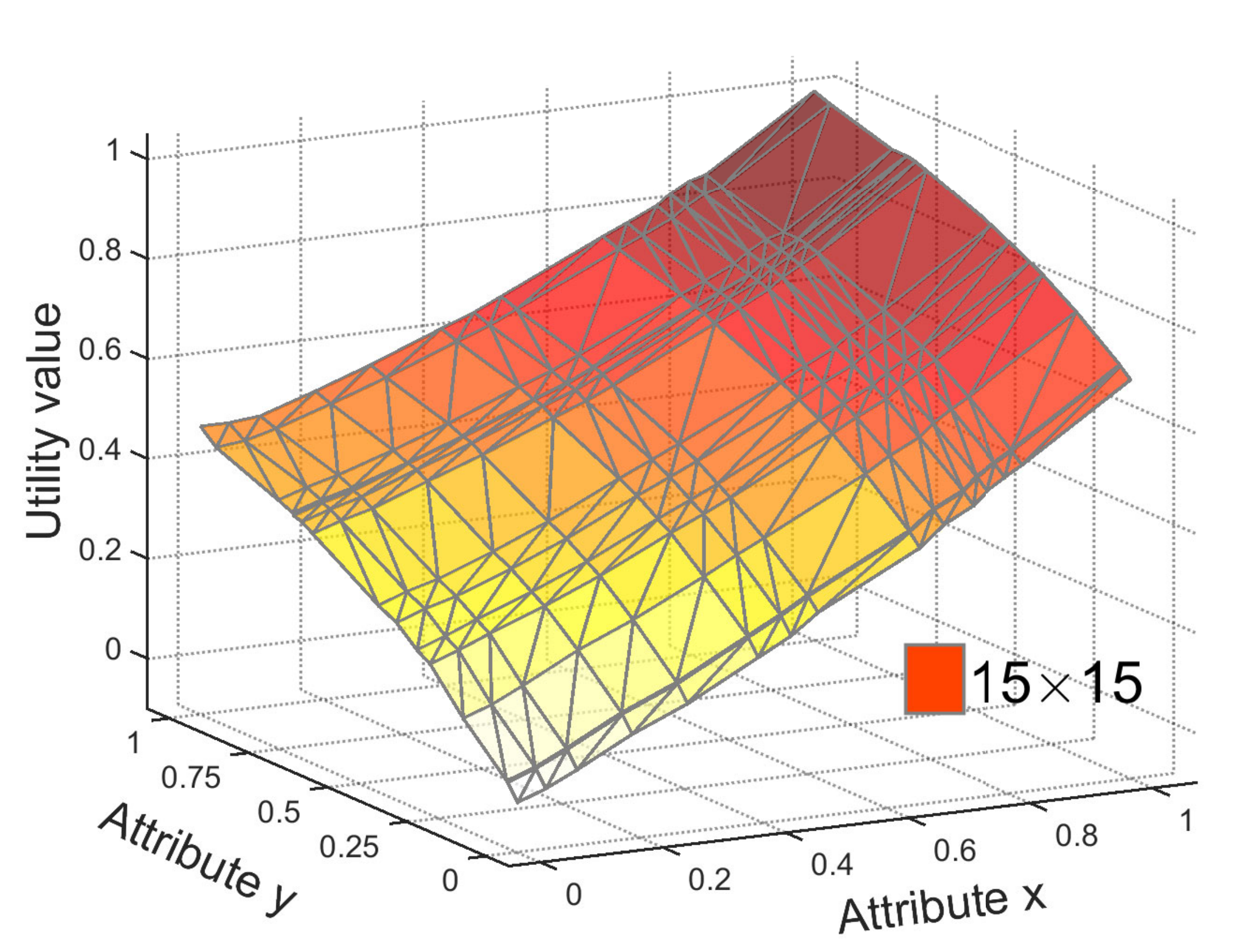}
  }
    \hspace{-0.5em}
      \subfigure
    {
    \includegraphics[width=0.23\linewidth]{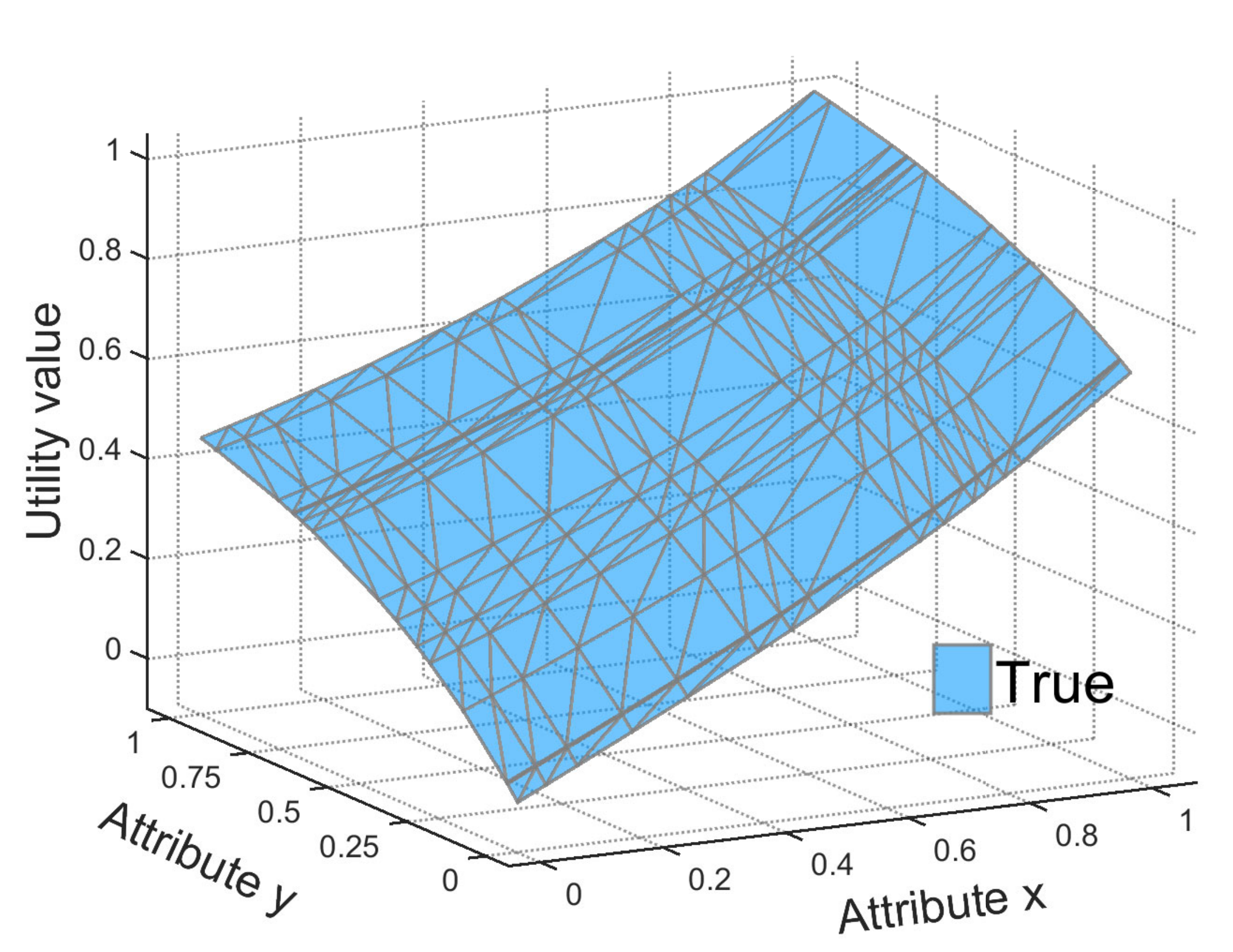}
  }
  \captionsetup{font=footnotesize}
  \vspace{-0.2cm}
  \caption{{\bf Type-1 IPLA}: the convergence of the worst-case utility function solved by the Dfree method for IPLA  model (\ref{eq:PRO_MILP_mina2}) without convex/concave constraints.
  }
  \vspace{-0.3cm}
  \label{fig:question-Utility_MILP} 
\end{figure}

\begin{figure}[!ht]
  \centering
  \vspace{-0.4cm}
      \subfigure{
    \includegraphics[width=0.23\linewidth]{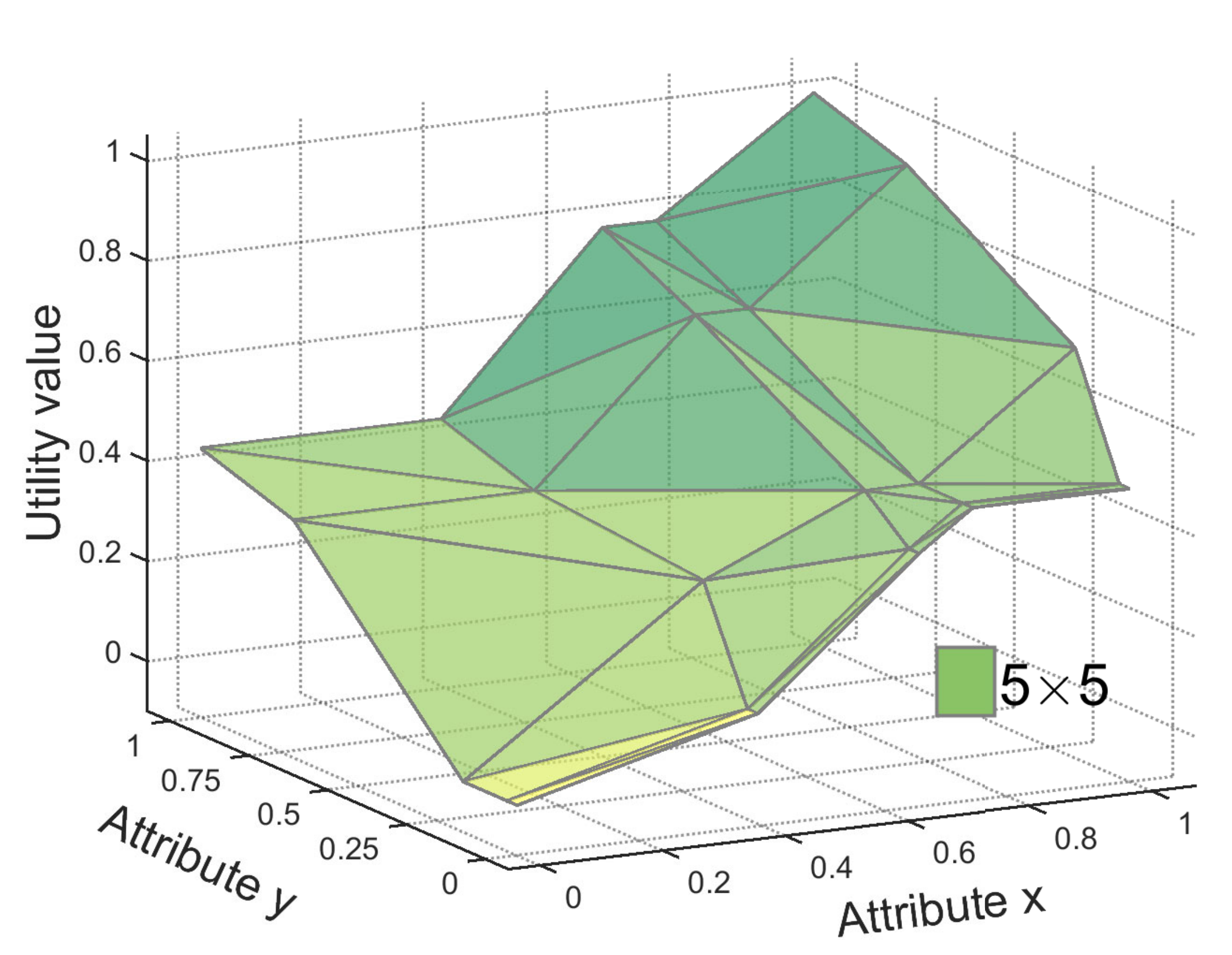}
  }
    \hspace{-0.5em}
    \subfigure{
    \includegraphics[width=0.23\linewidth]{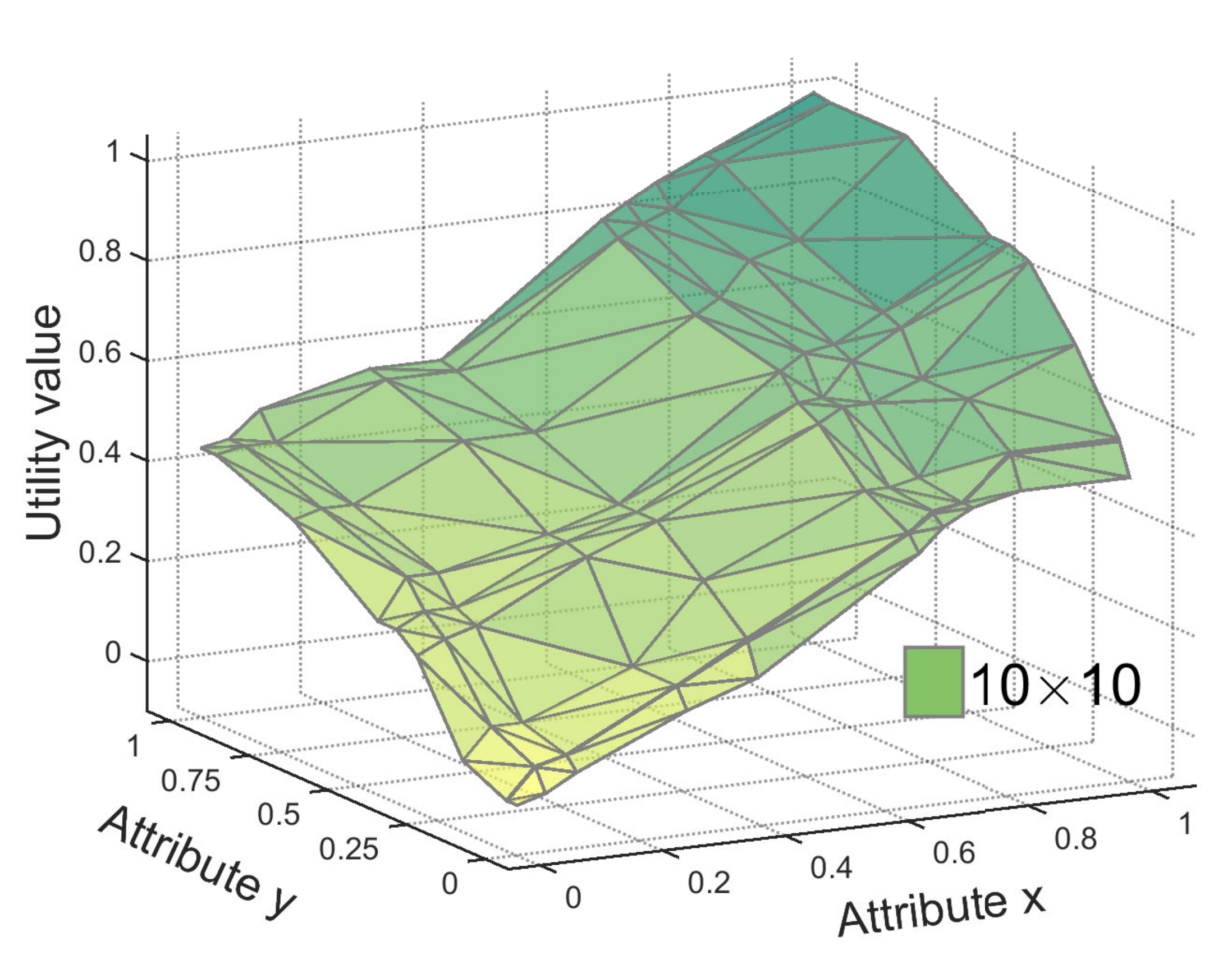}
   }
     \hspace{-0.5em}
    \subfigure{
    \includegraphics[width=0.23\linewidth]{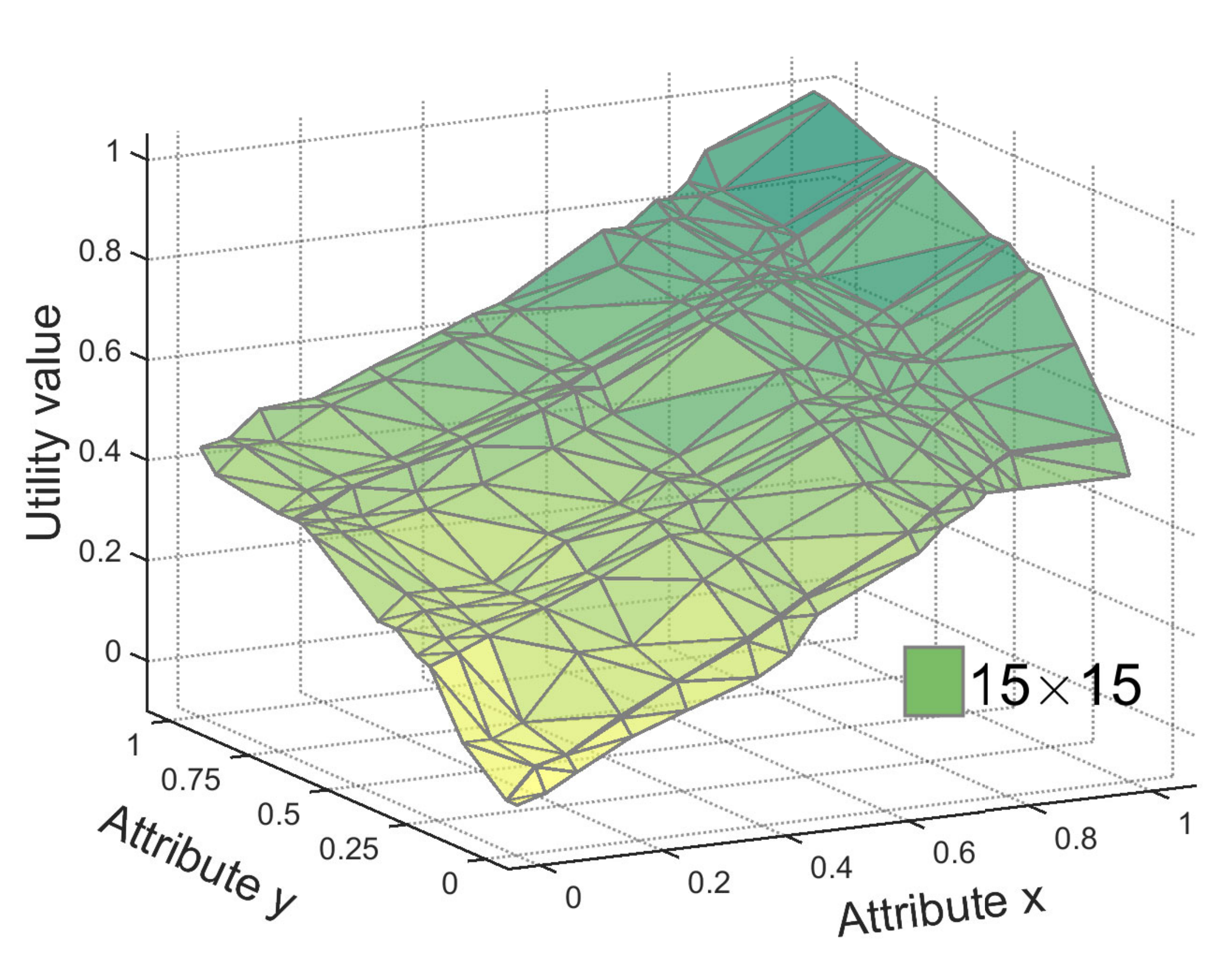}
  }
    \hspace{-0.5em}
      \subfigure{
    \includegraphics[width=0.23\linewidth]{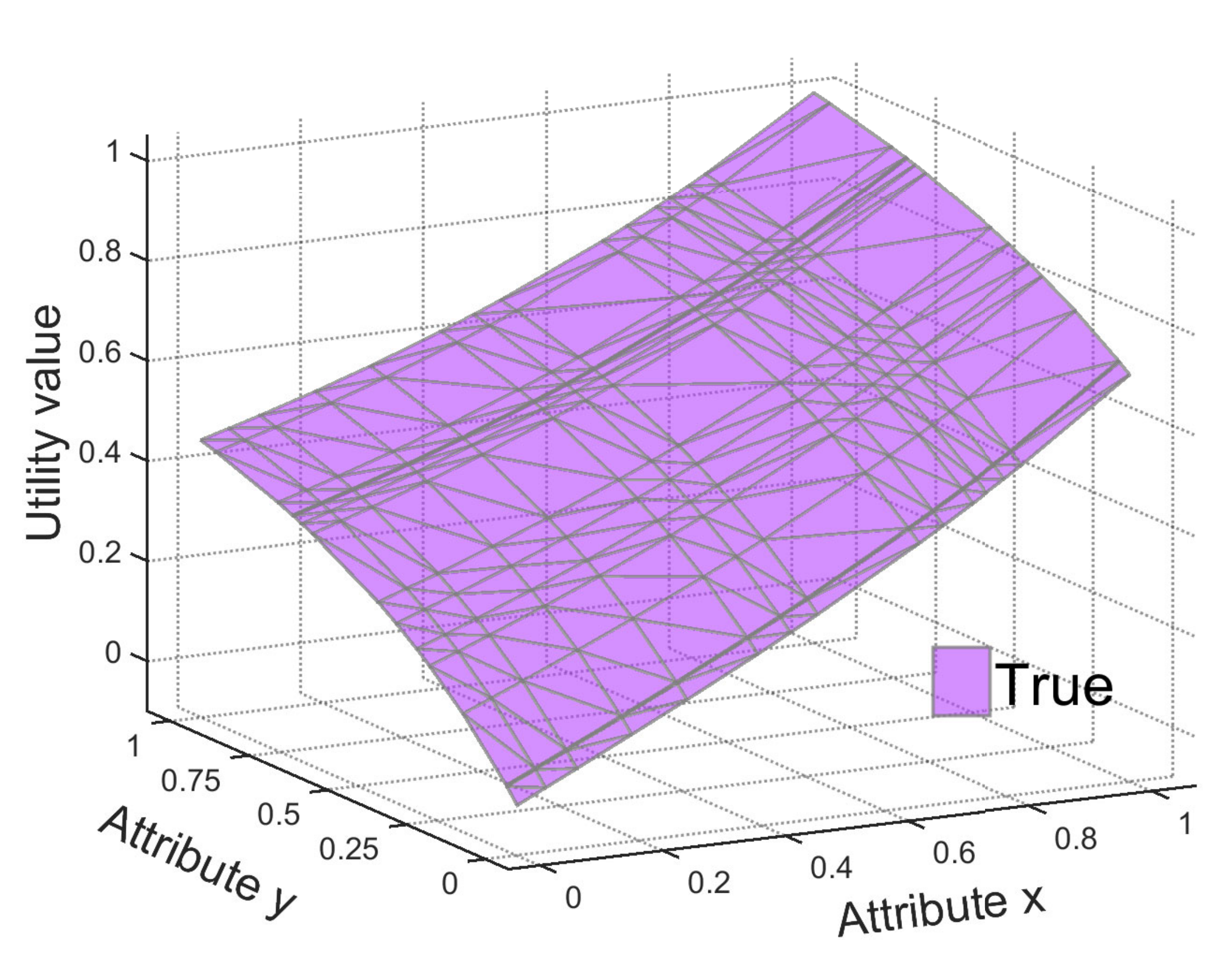}
  }
  \vspace{-0.2cm}
  \captionsetup{font=footnotesize}
  \caption{{\bf Type-2 IPLA}: the convergence of the worst-case utility function solved by Dfree method for IPLA  model (\ref{eq:PRO_MILP_mina2}) with (\ref{eq:mixed-integer-R2-f}) being replaced by 
(\ref{eq:constraint-alpha}).
}
  \label{fig:question-Utility-MILP2} 
\end{figure}

\underline{IPLA in bi-attribute case.}
Set $K=20$ (take the first $20$ from $1000$ samples, we do so because the problem size of (\ref{eq:PRO_MILP_mina2}) and (\ref{eq:PRO_MILP_single}) 
depends on the $K$ whereas 
problem size of (\ref{eq:PRO-N-reformulate}) under EPLA 
is independent of $K$),
the true utility $u^*$ is the same as in EPLA case.
In this set of tests,
the convexity/concavity of single-variate utility functions $u_N(\cdot,\hat{y})$ and $u_N(\hat{x},\cdot)$ for all $\hat{x}\in X$ and $\hat{y}\in Y$ is not considered to facilitate
comparison of the three models (maximin EPLA, maximin IPLA and single MILP using IPLA).
 because  in problem (\ref{eq:PRO_MILP_single}), we have not incorporated 
the constraints (see our comments there).
In Table~\ref{tab-result-N-MILP},
we compare the three models for the tractable reformulation of BUPRO-N:
EPLA 
(\ref{eq:PRO-N-reformulate}),
IPLA (\ref{eq:PRO_MILP_mina2})
and the single MILP (\ref{eq:PRO_MILP_single}) using IPLA,
for both Type-1 PLA and Type-2 PLA in terms of the optimal solution, the optimal value, error 
between BUPRO-N and utility maximization problem $\max_{\bdz \in Z} \sum_{k=1}^K p_k[u^*({\bm f}(\bdz,\bdxi^k))]$, and CPU time.
We 
find that the optimal values $\vt_N$ converge to the true optimal value $\vt^*$ in all cases.
We also find that for both types,
the EPLA
(\ref{eq:PRO-N-reformulate}) where the inner problem is an LP is most efficient,
 the single MILP (\ref{eq:PRO_MILP_single})
 obtains the best approximate optimal values but takes 
 longest CPU time.
 Note that although the three models are 
 equivalent theoretically,
 the actual computational results 
 differ  slightly  because of computational rounding errors.
Figures~\ref{fig:question-Utility_MILP}-\ref{fig:question-Utility-MILP2}  display the worst-case utility functions of IPLA maximin model (\ref{eq:PRO_MILP_mina2}) for Type-1 and Type-2 respectively.
We can see that the worst-case utility function 
displays some ``oscillations''
although it converges to the true.
The phenomenon disappears
when we  confine  
$u_N(\cdot,\hat{y})$ and $u_N(\hat{x},\cdot)$ 
to convex and concave functions respectively.

\begin{table}[!ht]
    \tiny
    \centering
    \captionsetup{font=scriptsize}
    \caption{{\bf The bi-attribute case:} comparison of the results of BUPRO-N problem (K=20, $\vt^*=0.3835$)}
    \vspace{-0.2cm}
    \renewcommand\arraystretch{1.1}
    \resizebox{0.9\linewidth}{!}
    {
    \begin{tabular}{c|ccccc}
        \hline
         & 
        Lotteries & Optimal solutions & Optimal values & Error & CPU time (s) \\
        \hline              
        \multirow{3}{*}{\makecell{Type-1 \\ Maximin \\({\bf EPLA})}} & $5\times5$ & $[ 0.112, 0.037,0,0.439,0.024, 0,0.054,0.335]$ &  0.2835  &  0.1000 & 47.4 \\
        & $10\times10$ & $[0,  0.599,0,0.316,0.008,0.037,0.012,0.027]$  & 0.3479 & 0.0356  & 82.8 \\
      & $15\times15$  & $[0,1,0,0,0,0,0,0]$  & 0.3754 & 0.0081 & 146.2 \\
        \hline
        \multirow{3}{*}{\makecell{Type-1 \\ Maximin
        \\ ({\bf IPLA}) } } &
      $5\times5$ & $[0.0996,0.0313,0.0297,0.4525,0, 0,0.0469,0.3400]$ &  0.2824 & 0.1011 & 240.5 \\
        & $10\times10$ & $[0, 0.9467,0,0, 0,0.0476,0,0.0057]$  & 0.3697 &  0.0138 & 916.0 \\
        & $15\times15$ &  $[0,0.9902,0,0.0038,0,0.0060,0,0]$  & 0.3748  & 0.0087 & 2442.2 \\
        \hline
        \multirow{3}{*}{\makecell{Type-1\\  Single MILP \\
        ({\bf IPLA})} } &
        $5\times5$ & $[0,1,0,0,0,0,0,0]$ &  0.3232 & 0.0603 & 1103.2 \\
        & $10\times10$ & $[0,0.946,0,0,0,0.043,0,0.011]$ & 0.3698 &  0.0137 & 4552.2 \\
        & $15\times15$ & $[0,1,0,0,0,0,0,0]$ & 0.3754 & 0.0081 &  3421.1 \\
        \hline
        \hline
         \multirow{3}{*}{\makecell{Type-2 \\ Maximin 
         \\({\bf EPLA})}
        } & $5\times5$ & $[0,0.1542,0.0117,0.1210,0.3554,0.1679,0,0.1898]$ &0.3113 &  0.0722 &  43.7 \\
        & $10\times10$ &  $[0, 0.5301,0,0.1857,0.2030,0.0410,0.0308,0.0094]$ & 0.3475 & 0.0360 &  56.4 \\
        & $15\times15$ &  $[0,1,0,0,0,0,0,0]$  & 0.3754 &  0.0081 & 98.2 \\
        \hline
        \multirow{3}{*}{\makecell{Type-2 \\ Maximin
        \\({\bf IPLA})}
        } &
        $5\times5$ & $[0, 0.8875,0.1125,0,0,0,0,0]$  & 0.3102  &  0.0733 &  213.6 \\
        & $10\times10$ & $[0,0.975,0.025,0,0,0,0,0]$  & 0.3410 & 0.0425 &  828.2 \\
        & $15\times15$ &  $[0,0.7129,0.0986,0,0.1884,0,0,0]$ & 0.3474 & 0.0392 &  2126.9 \\
        \hline
      \multirow{3}{*}{\makecell{Type-2\\ Single MILP \\
      ({\bf IPLA})} } &
        $5\times5$ & $[0,1,0,0,0,0,0,0]$ &  0.3232 & 0.0603 & 952.0 \\ 
        & $10\times10$ &$[0,0.9470,0,0,0, 0.0467,0,0.0062]$  &  0.3704 & 0.0131 & 1359.1\\
        & $15\times15$ & $[0,1,0,0,0,0,0,0]$  &  0.3754 & 0.0081 &  2680.9\\
        \hline
    \end{tabular}
    }
    \label{tab-result-N-MILP}
    \vspace{-0cm}
\end{table}




\underline{IPLA in tri-attribute case.}
The sample is the same as in the bi-attribute case with $K=20$.
The true utility function is
$u(x,y,z)=e^{x}-e^{-y}-e^{-z}-e^{-x-2y-z}:[0,1]^3 \to [0,1]$
and normalize it by setting $u^*(x,y,z)=(u(x,y,z)-u(0,0,0))/(u(1,1,1)-u(0,0,0))$.
We divide the eight projects into three groups in order of importance as the three attributes, that is,
$f_1^k:=\sum_{i=1}^{3} w_i \xi_i^k$, $f_2^k:=\sum_{i=4}^{6} w_i \xi_i^k$,
$f_3^k:=\sum_{i=7}^{8} w_i \xi_i^k$,
and ${\bm w} \in Z:=\{{\bm w}\in\R^8_+:\sum_{i=1}^8 w_i=1\}$.
Table~\ref{tab-result-N-MILP3} indicates that the IPLA model (\ref{eq:PRO_MILP_3m}) in  tri-attribute case is effective and the optimal values $\vt_N$ of the TUPRO-N problem converge to the true optimal value $\vt^*$ as the number of lotteries increases.

\begin{table}[!ht]
\vspace{-0.3cm}
    \tiny
    \centering
    \captionsetup{font=scriptsize}
    \caption{{\bf The tri-attribute case}: computational results of TUPRO-N problem in (K=20, $\vt^*=0.3193$)}
    \vspace{-1em}
    \renewcommand\arraystretch{1.08}
    \begin{threeparttable}
    \resizebox{0.8\linewidth}{!}
    {
    \begin{tabular}{c|ccccc}
        \hline
         & 
        Lotteries & Optimal solutions & Optimal values & Error & CPU time (s) \\
        \hline
        \multirow{3}{*}{ {\bf IPLA}
        } & $3\times3 \times 3$ & $[0,0,1,0,0,0,0,0]$ & 0.1994 & 0.1198 & 320.8 \\
        & $4 \times 4 \times 4$ & $[ 0.4992,0.5008,0,0, 0,0,0,0]$  & 0.2076 &  0.1117 & 944.1 \\
        & $5\times 5 \times 5$ & $[0,1,0, 0,0,0,0,0]$  & 0.2498 & 0.0694 & 2083.7 \\
        & $6\times 6 \times 6$  & $[0,1,0, 0,0,0,0,0]$ & 0.2774 & 0.0418 &  - \\
          \hline
    \end{tabular}
    }
    \begin{tablenotes}
    \raggedleft
        \item `-' implies runtime $>$ 3600s. 
    \end{tablenotes}
    \end{threeparttable}
    \vspace{-0.5cm}
    \label{tab-result-N-MILP3}
\end{table}

\textbf{(ii) EPLA for the constrained optimization problems (\ref{eq:PRO-x}) and 
(\ref{eq:PRO-x-1}).}
The second part of numerical tests is concerned with 
problems (\ref{eq:PRO-x}) and (\ref{eq:PRO-x-1}).
We set $g_1(\bdz,\bdxi^k):=\sum_{i=3}^5 z_i\xi_i^k$ and 
$g_2(\bdz,\bdxi^k):=\sum_{i=7}^8 z_i\xi_i^k$,
which represent the 
effects of part of the 
projects on 
mental health and cancer 
commissioning areas.
PHO expects 
this part of effects to
reach at least level $c$.
We consider two cases: (a)
$c=0.1$ and (b) $c=0.3$.

\underline{Case (a)}. 
The optimal values of problem (\ref{eq:PRO-x}) and problem (\ref{eq:PRO-x-1}) coincide (see
Table~\ref{tab-result-N}) because 
the optimal solution of the former
falls into set (\ref{eq:x*-PRO-U}),
which is consistent with our theoretical analysis in Proposition~\ref{Prop-equivalence}.
 \underline{Case (b)}. We repeat the tests but with different observations.
%
Recall that the optimal values of problems (\ref{eq:SPR-x}), (\ref{eq:PRO-x}) and (\ref{eq:PRO-x-1}) are denoted by $\vt^*$, $\hat{\vt}$ and $\tilde{\vt}$, respectively.

Observation 1. 
For problem (\ref{eq:PRO-x-1}),
we can see from  Table~\ref{tab-result-DN-3}
that $\tilde{\vt}<\vt^*$ and $\tilde{\vt}$  increases as $M$ increases. This is 
consistent with our theoretical analysis.
The increasing trend is underpinned 
by the fact that 
as $M$ increases,
$\calu_N$
becomes smaller
and consequently both the objective function $\min_{u\in {\cal U}} \bbe_{P}[u({\bm f}(\bdz,\bdxi))]$ and  the feasible set  
$\tilde{Z}$ (see  (\ref{eq:x*-PRO-U}))
become larger.

Observation 2. 
For
problem (\ref{eq:PRO-x}),
we can see from Table~\ref{tab-result-DN-2} that
$\vt^*<\hat{\vt}$ 
for the cases that $5\times 5$ and $10\times 10$ lotteries are used.
Note that by theory,
$\tilde{\vt}\leq \hat{\vt}$ and $\tilde{\vt}\leq \vt^*$. Moreover,
 when $\hat{\bdz}\in \tilde{Z}$,
we are guaranteed that $\tilde{\vt}=\hat{\vt}\leq \vt^*$.
The observed trend reflects the fact that $\hat{\vt}>\vt^*$ may occur when
$\hat{\bdz}\notin \tilde{Z}$.
Moreover,
$\vt^*>\hat{\vt}$ 
when $15\times 15$ lotteries are used since
$\hat{\bdz}\in \tilde{Z}$.

Observation 3. The optimal value $\hat{\vt}$ is decreasing from Table~\ref{tab-result-DN-2} as the number of questions increases. 
This phenomena is a bit difficult to explain. On one hand,
when the size of ${\cal U}_N$ decreases,
$\hat{v}(z)$ increases 
and on the other hand
the size of $\hat{Z}:=\{\bdz:\bbe_P[u(\bdg(\bdz,\bdxi))] \geq c\}$ decreases.
Note that $\hat{\vt} = \max_{z\in \hat{Z}}\hat{v}(z)$,
 it seems the reduction 
of the size of $\hat{Z}$ has more effect than
that of the increase of  $\hat{v}(z)$ in this test. 

We have not tested IPLA as our focus here is on the difference between 
model (\ref{eq:PRO-x}) and model (\ref{eq:PRO-x-1}) rather than 
different performances of EPLA and IPLA.

\begin{table}[!ht]
    \tiny
    \centering
    \captionsetup{font=scriptsize}
    \caption{Computational results of
    problem (\ref{eq:PRO-x-1}) ($K=1000$, $c=0.3$, $\vt^*=0.3387$)}
    \vspace{-0.5em}
    \renewcommand\arraystretch{1.1}
    \begin{threeparttable}
    \resizebox{0.8\linewidth}{!}
    {
    \begin{tabular}{c|ccccc}
        \hline
         \textbf{EPLA} & 
        Lotteries & Optimal solutions & $\tilde{\vt}$ & Error & CPU time (s) \\
        \hline
        \multirow{3}{*}{Type-1
        } & $5\times5$ & $[0,0,0,1,0,0,0,0]$ & 0.3122 & 0.0265 & 237.0 \\
        & $10\times10$ & $[0,0,0,0.955,0,0,0.011,0.034]$ & 0.3321 & 0.0066 & 379.9 \\
        & $15\times15$ & $[0,0,0,1,0,0,0,0]$ & 0.3349 & 0.0038 & 436.7 \\
        \hline
        \multirow{3}{*}{Type-2
        } &
        $5\times5$ & $[0,0,0,1,0,0,0,0]$ &  0.3122 & 0.0265 & 307.1 \\
        & $10\times10$ & $[0,0,0,0.959,0,0,0,0.041]$ & 0.3323 & 0.0064 & 321.2 \\
        & $15\times15$ & $[0,0,0,1,0,0,0,0]$ & 0.3349 & 0.0038 & 486.0 \\
        \hline
        \multirow{3}{*}{ Mixed-type
        } &
        $5\times5$ & $[0,0,0,1,0,0,0,0]$ &  0.3122 & 0.0265 & 2476.2 \\
        & $10\times10$ & - & - & - & - \\
        & $15\times15$ & - & - & - & - \\
        \hline
    \end{tabular}
    }
    \begin{tablenotes}
    \raggedleft
        \item `-' implies runtime $>$ 3600s.
    \end{tablenotes}
    \end{threeparttable}
    \vspace{-0.6cm}
    \label{tab-result-DN-3}
\end{table}

\begin{table}[!ht]
    \tiny
    \centering
    \captionsetup{font=scriptsize}
    \caption{Computational results of
    problem (\ref{eq:PRO-x}) ($K=1000$, $c=0.3$, $\vt^*=0.3387$)}
    \vspace{-0.2cm}
    \renewcommand\arraystretch{1.1}
    \begin{threeparttable}[b]
    \resizebox{0.9\linewidth}{!}
    {
    \begin{tabular}{c|ccccc}
        \hline
        \textbf{EPLA} & 
        Lotteries & Optimal solutions & $\hat{\vt}$ & Error & CPU time (s)\\
        \hline
        \multirow{3}{*}{Type-1
        } & $5\times5$ & $[0.118,0.115,0.178,0.179,0,0.130,0.112,0.169]$ & 0.3873 & -0.0486 & 216.0 \\
        & $10\times10$ & $[0,0.111,0,0.883,0,0.006,0,0]$ & 0.3413 & -0.0026 & 255.5 \\
        & $15\times15$ & $[0.007,0.098,0,0.875,0.020,0,0,0]$ & 0.3377 & 0.0010 & 363.4 \\
        \hline
        \multirow{3}{*}{Type-2
        } &
        $5\times5$ & $[0.176,0.129,0,0.077,0,0.178,0.084,0.357]$ &  0.4186 & -0.0799 & 242.6 \\
        & $10\times10$ & $[0.027,0.082,0,0.891,0,0,0,0]$ & 0.3384 & 0.0003 & 261.7 \\
        & $15\times15$ & $[0,0,0,1,0,0,0,0]$ & 0.3349 & 0.0038 & 335.2 \\
        \hline
        \multirow{3}{*}{Mixed-type
        } &
        $5\times5$ & $[0.073,0.131,0.118,0.022,0.151,0.194,0.107,0.205]$ & 0.4095 & -0.0708 & 2287.6 \\
        & $10\times10$ & - & - & - & - \\
        & $15\times15$ & - & - & - & - \\
        \hline
    \end{tabular}
    }
    \hspace{-0.1cm}
    \begin{tablenotes}
    \raggedleft
        \item `-' implies runtime $>$ 3600s.
    \end{tablenotes}
    \end{threeparttable}
       \vspace{-0.5cm}
    \label{tab-result-DN-2}
\end{table}

\subsection{Perturbation analysis}
This part of numerical tests is concerned with data perturbation including (i) elicitation data perturbation and (ii) sample average approximation (SAA) of the exogenous uncertainties.
SAA is needed when the true 
probability distribution $P$ in
(\ref{eq:MAUT-robust}) is continuously distributed. In this case, Assumption \ref{assu-discrete}
and the subsequent UPRO models
may be viewed as sample average approximations.
We skip the theoretical analysis about errors 
arising from SAA and refer interested readers to \cite{GXZ21}
in single-attribute case.

\begin{figure}[!htbp]
  \centering
  \vspace{-0.2cm}
  \subfigure{
    \label{subfig-main-ut-ptb5} 
    \includegraphics[width=0.3\linewidth]{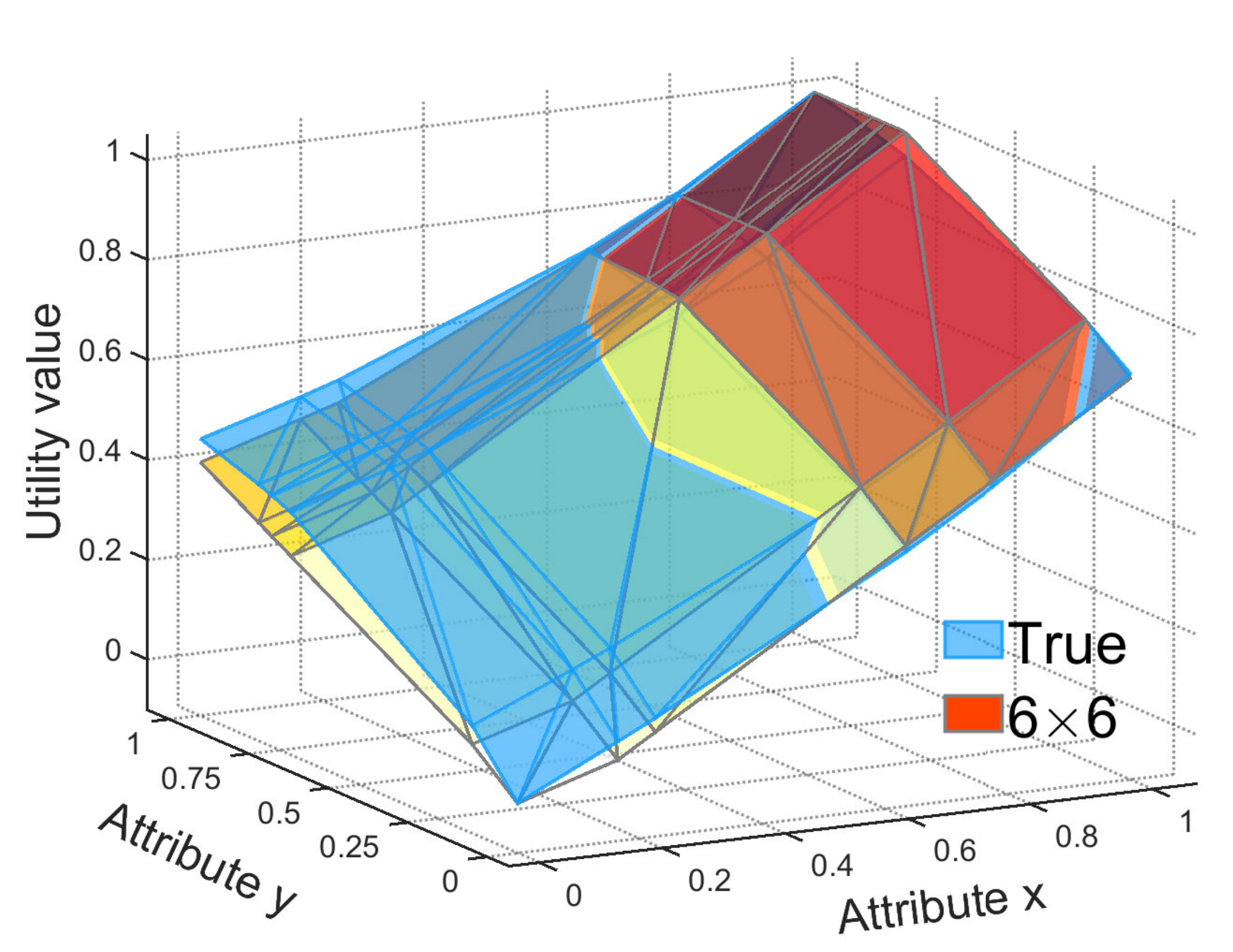}
  }
    \hspace{-0.5em}
  \subfigure{
    \label{subfig-main-ut-ptb10} 
    \includegraphics[width=0.3\linewidth]{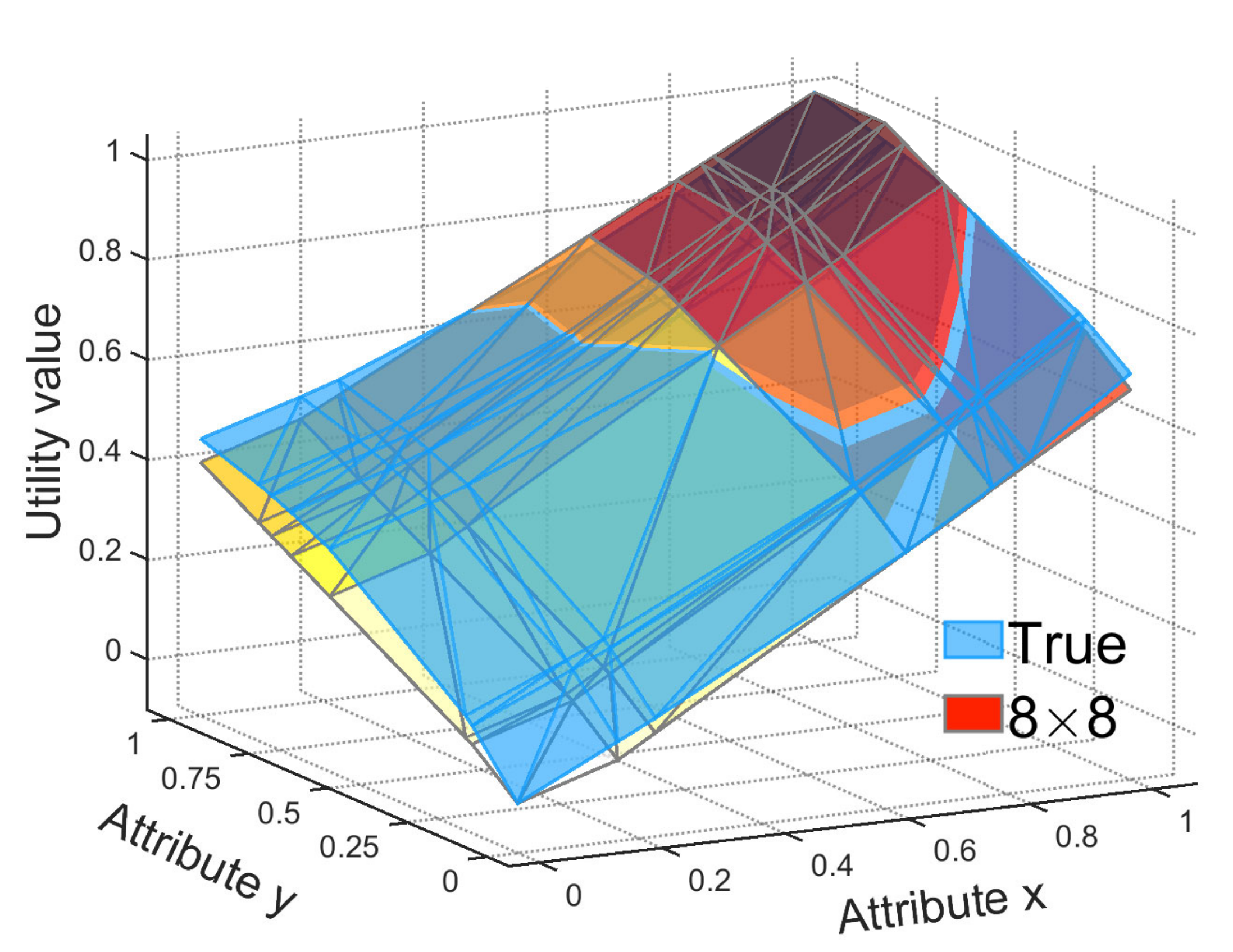}
  }
    \hspace{-0.5em}
  \subfigure{
    \label{subfig-main-ut-ptb15}
    \includegraphics[width=0.3\linewidth]{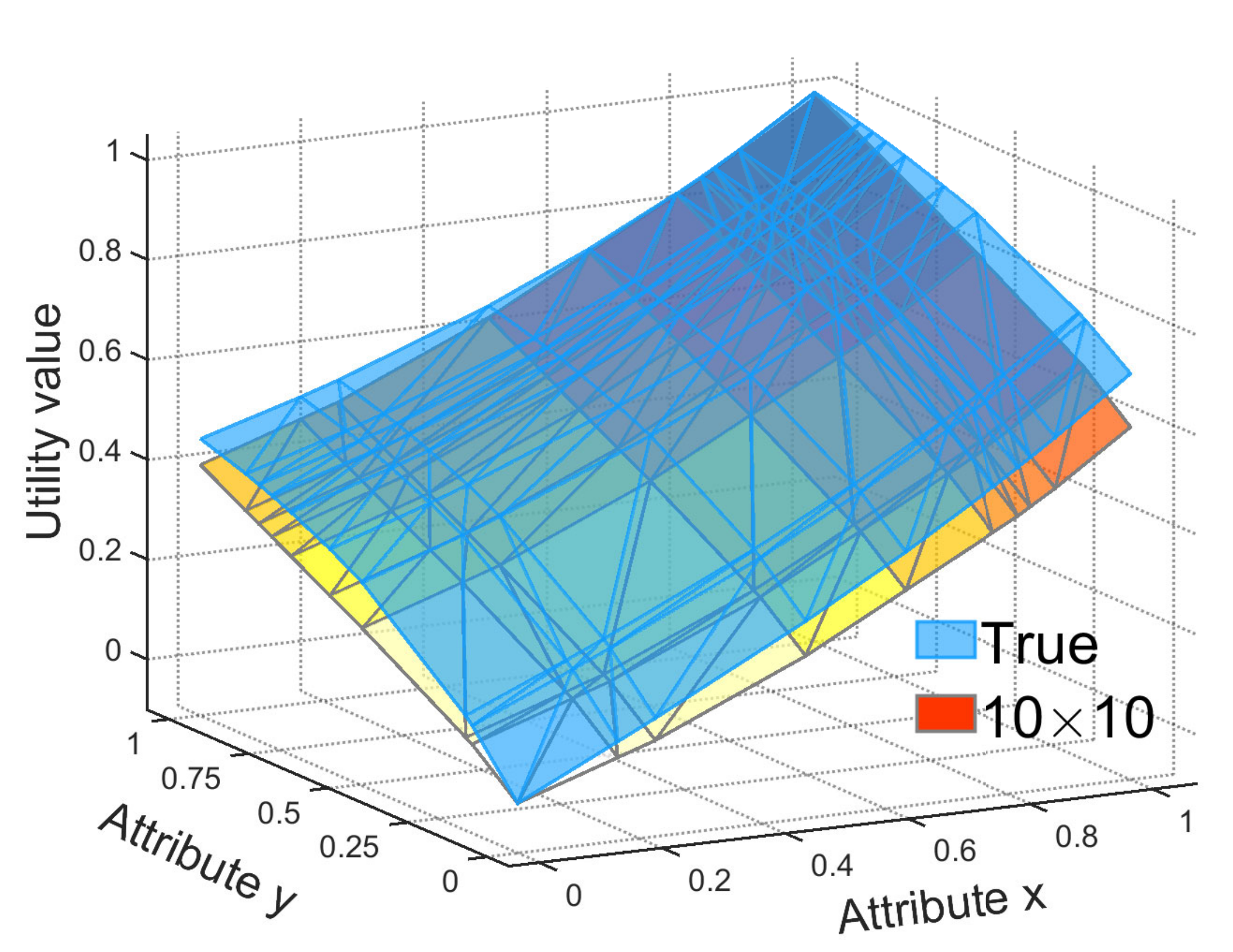}
  }
  \vspace{-1em}
  \captionsetup{font=footnotesize}
  \caption{{\bf Type-1 EPLA}: the worst-case utility function with $\delta_1=0.1$ }
  \vspace{-0.2cm}
  \label{fig-ptb-ut-main} 
\end{figure}

\begin{figure}[!htbp]
  \centering
  \vspace{-0.2cm}
  \subfigure{
    \label{subfig-counter-ut-ptb5} 
    \includegraphics[width=0.3\linewidth]{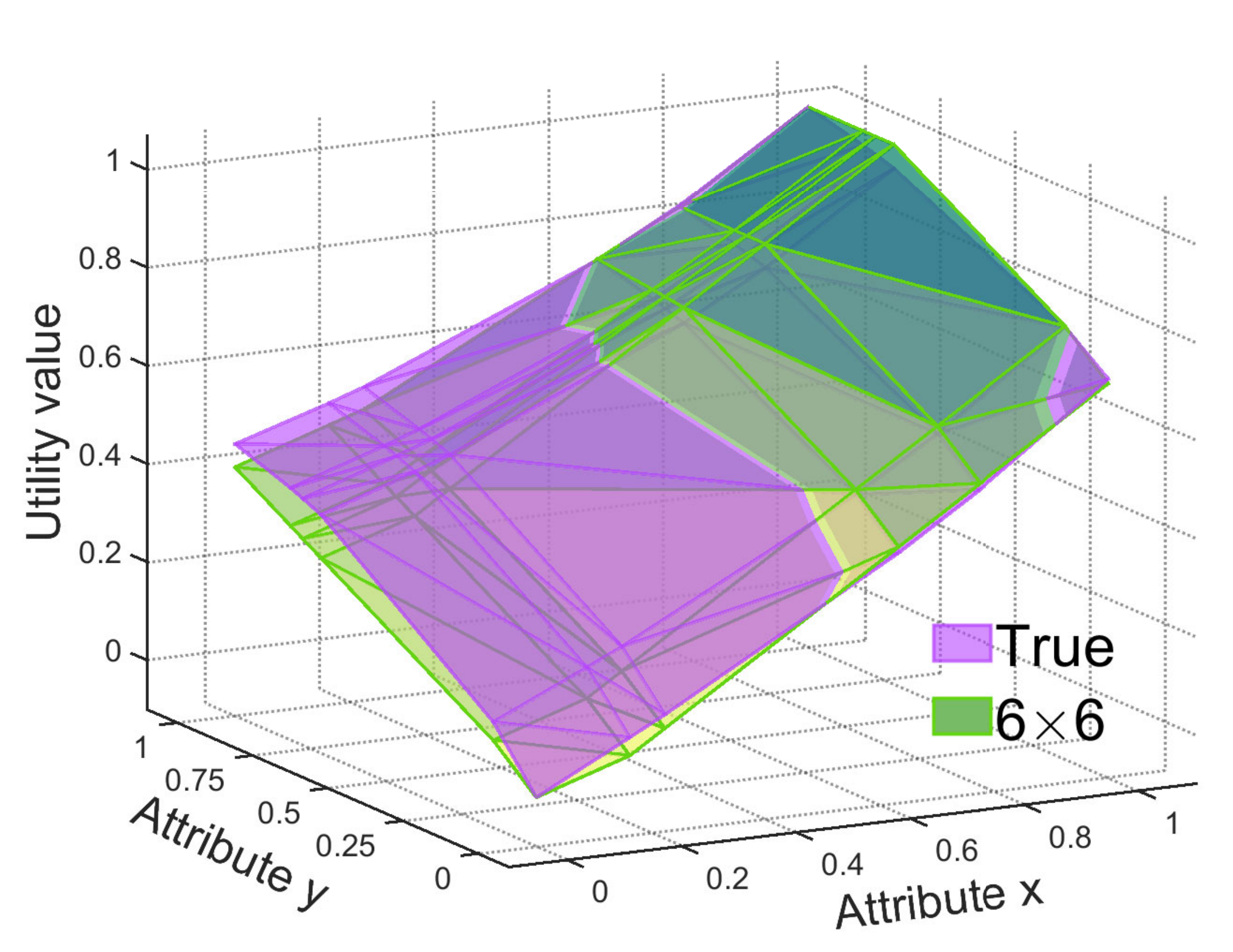}
  }
    \hspace{-0.5em}
  \subfigure{
    \label{subfig-counter-ut-ptb10} 
    \includegraphics[width=0.3\linewidth]{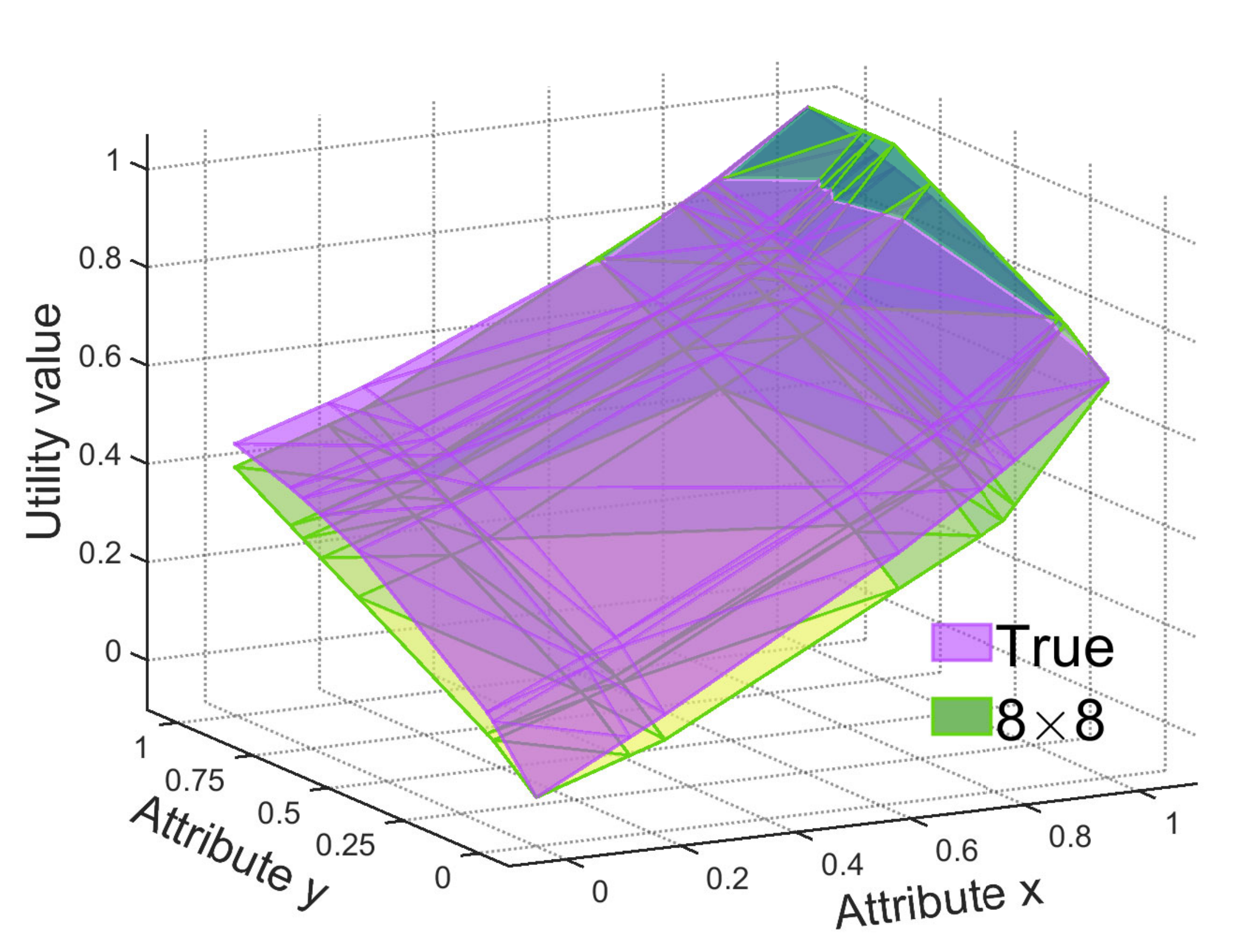}
  }
    \hspace{-0.5em}
  \subfigure{
    \label{subfig-counter-ut-ptb15}
    \includegraphics[width=0.3\linewidth]{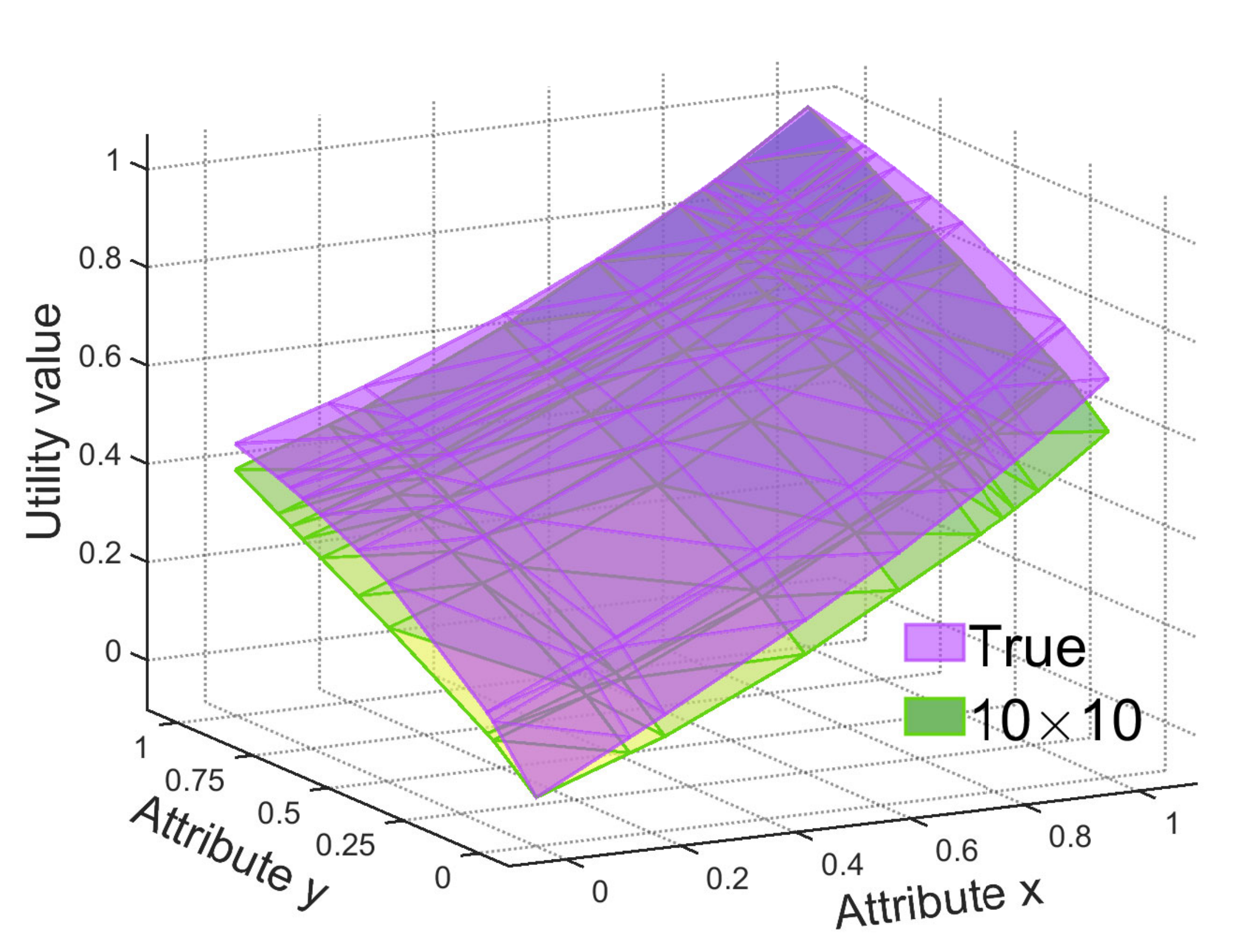}
  }
  \vspace{-1em}
  \captionsetup{font=footnotesize}
  \caption{{\bf Type-2 EPLA}: the worst-case utility function with $\delta_1=0.1$ }
  \vspace{-0.2cm}
  \label{fig-ptb-ut-counter} 
\end{figure}

\begin{figure}[!ht]
  \centering
  \subfigure[]{
    \label{subfig-main-ov} 
    \includegraphics[width=0.23\linewidth]{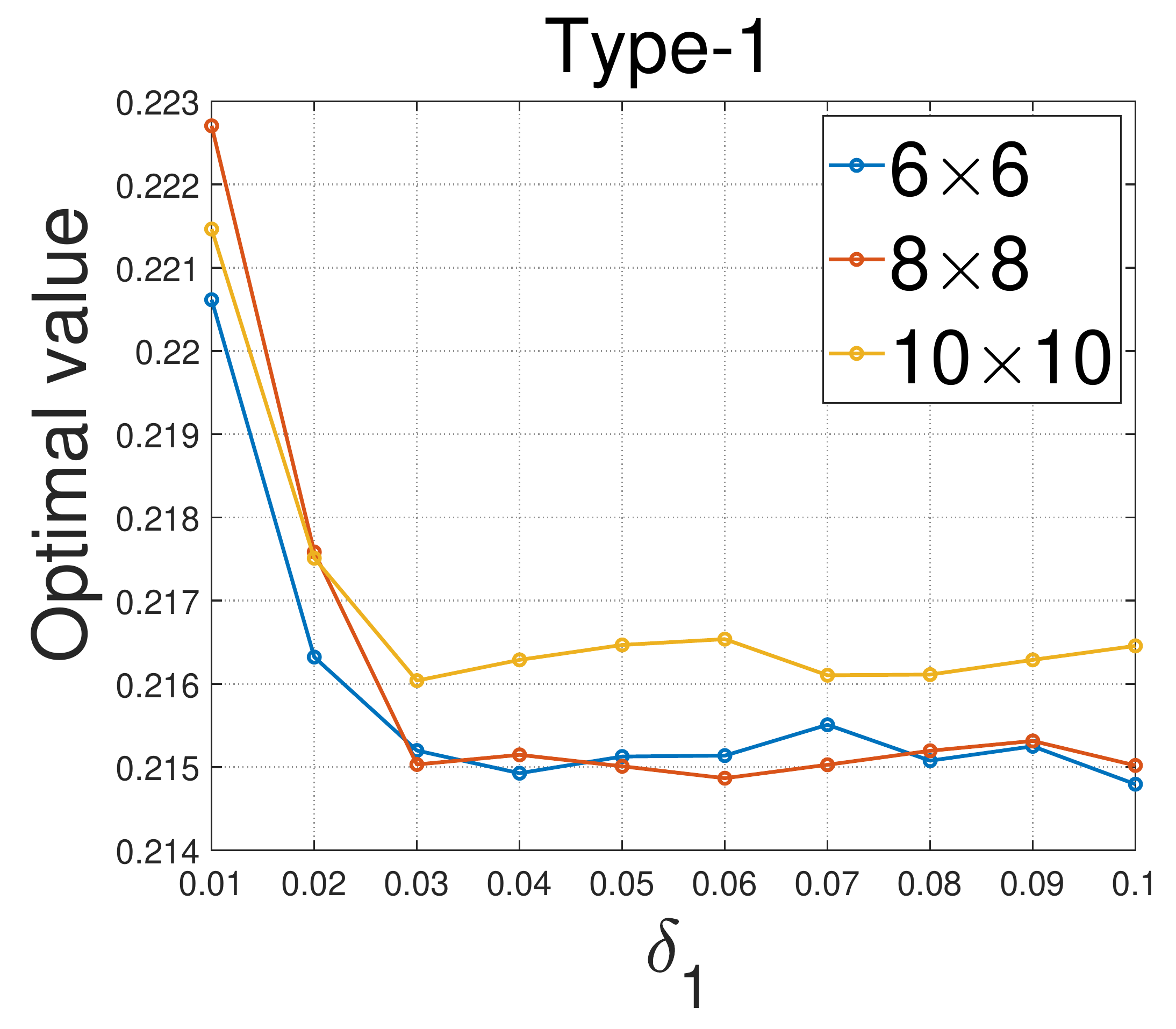}
  }
  \hspace{-1.5em}
  \subfigure[]{
    \label{subfig-counter-ov} 
    \includegraphics[width=0.23\linewidth]{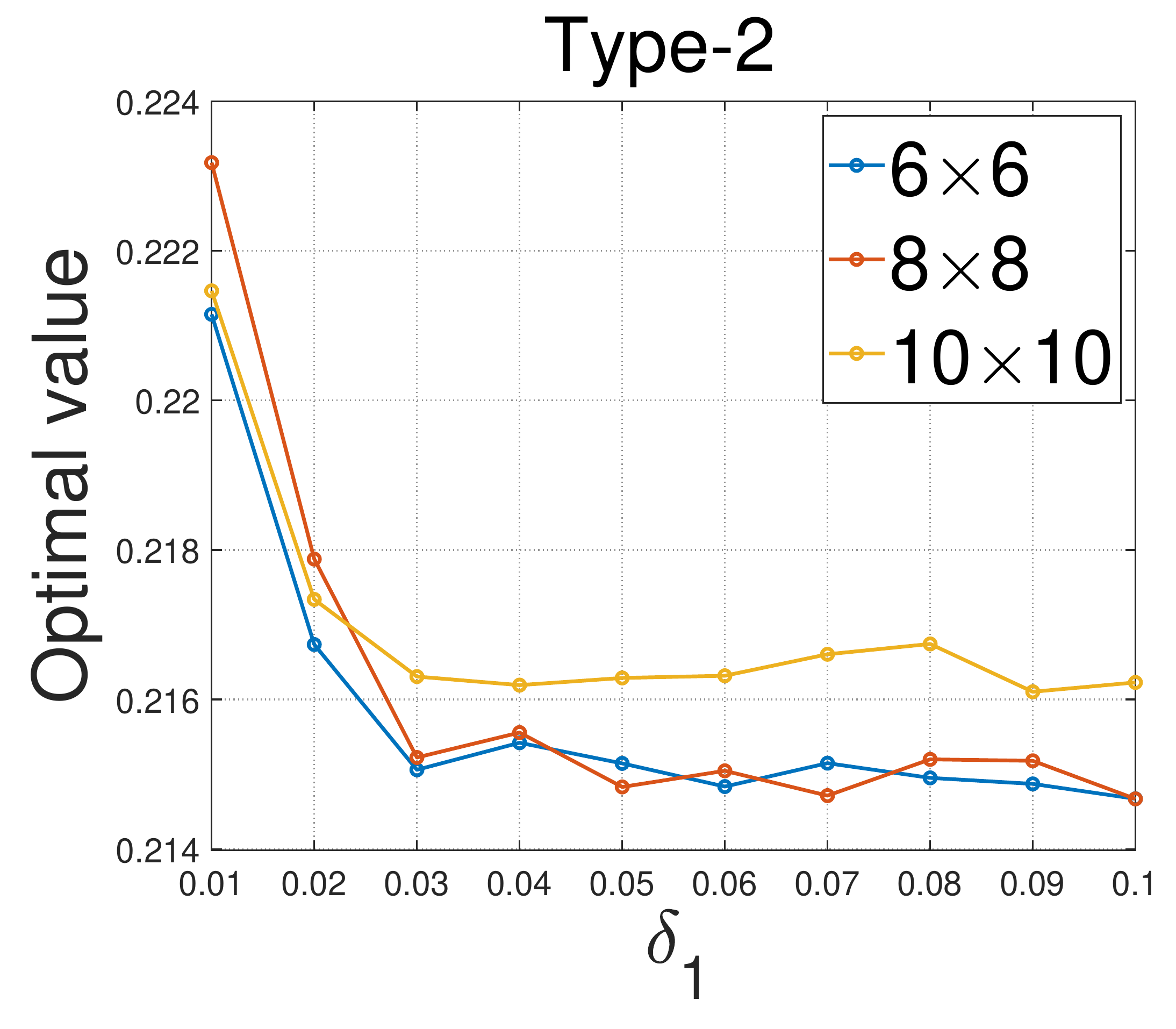}
  }
  \hspace{-1.5em}
  \subfigure[]{
    \label{subfig-SAA-main}
    \includegraphics[width=0.23\linewidth]{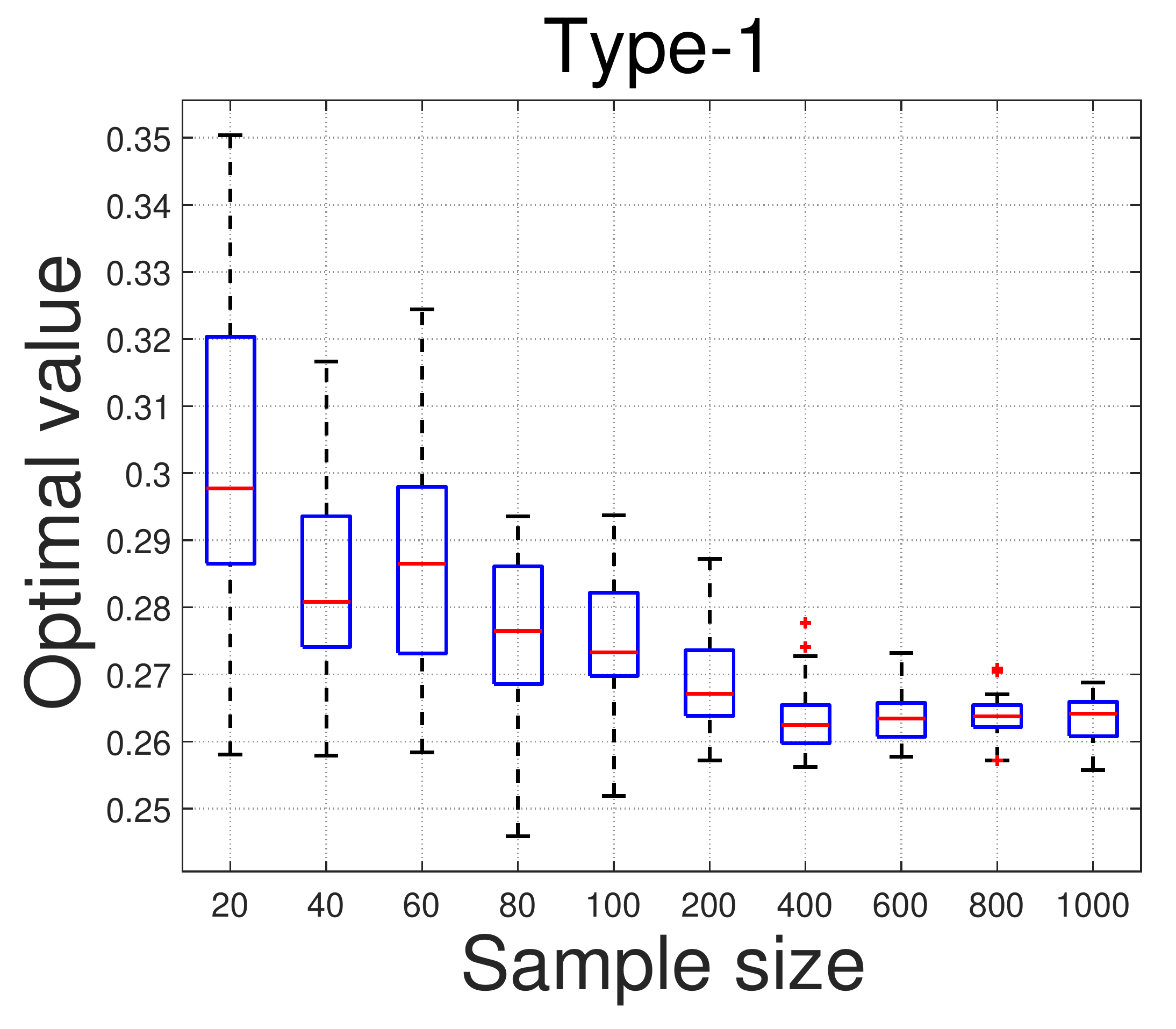}
  }
  \hspace{-1.5em}
  \subfigure[]{
    \label{subfig-SAA-counter}
    \includegraphics[width=0.23\linewidth]{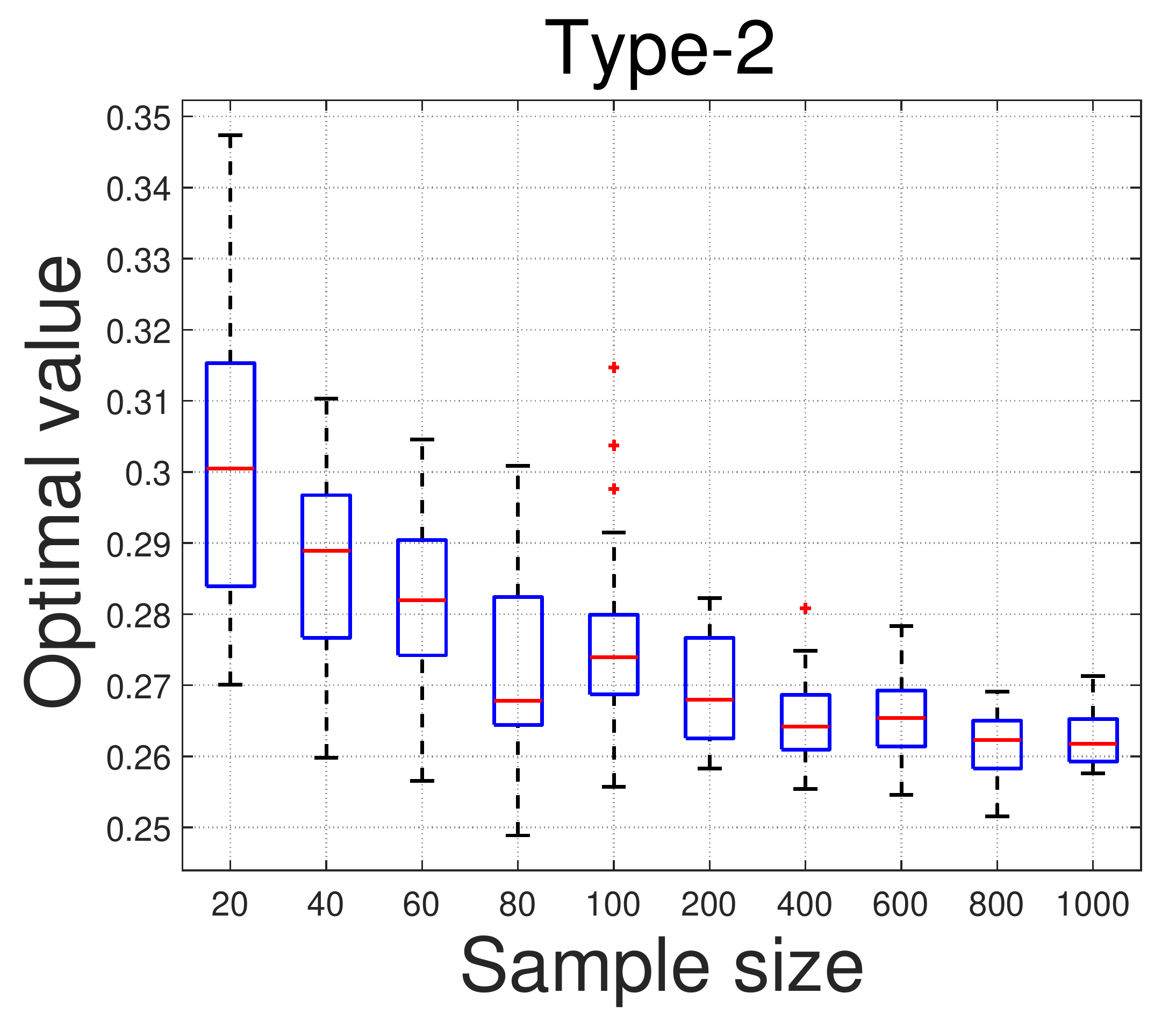}
  }
  \vspace{-1.2em}
  \captionsetup{font=footnotesize}
  \caption{{\bf EPLA:} the optimal values with $\delta_1=0.01$ and SAA problem as sample size increases}
  \vspace{-0.2cm}
  \label{fig-ptb-ov}
\end{figure}

\textbf{(i) Perturbation 
in the data in the ambiguity set.}
In this set of experiments, we will test the performance of the PLA scheme
when the ambiguity sets $\calu$ and $\tilde{\calu}$ are replaced by  $\calu_N$ and $\tilde{\calu}_N$ respectively.
We begin by considering a situation where the underlying functions $\psi_l, l=1,\ldots,M$ in the ambiguity set are perturbed by the observation error of the random data in pairwise comparison questions, i.e., 
\begin{equation*}
    \tilde{\psi}_l(x,y):= \1_{[\hat{x}^l+\delta_1,1]\times [\hat{y}^l+\delta_2,1]}(x,y)-(1-p^l) \1_{[0,1]\times[0,1]\setminus (1,1)}(x,y)- \1_{(1,1)}(x,y), l=1,\ldots,\hat{M},
\end{equation*}
where $\hat{M}$ is the number of perturbed functions $\psi_l$.
Notice that some lotteries are on the boundary of rectangle $T$ which can only be perturbed inwards.
Thus we assume that these lotteries are not perturbed for the convenience of discussion.
Let $\calu_N=\{u_N\in\scru_N: \la u_N,\psi_l \ra\leq c_l, l=1,\ldots,M\} $
and 
\begin{equation*}
    \tilde{\calu}_N=\{u_N\in\scru_N: \la u_N,\tilde{\psi}_l \ra\leq c_l, l=1,\ldots,M\}.
\end{equation*}
We can  solve problem (\ref{eq:PRO-N-reformulate}) with $\psi_l$ being replaced by $\tilde{\psi}_l$ to obtain the optimal value and the corresponding worst-case utility function.
Specifically, we assume $\delta_2=0$, which means we only consider the case that the first attribute is slightly perturbed but the second attribute is not.
Figures~\ref{fig-ptb-ut-main}-\ref{fig-ptb-ut-counter} depict the convergence of the worst-case utility functions as the number of questions increases for fixed $\delta_1=0.1$ with Type-1 PLA and Type-2 PLA.
Figures~\ref{subfig-main-ov}-\ref{subfig-counter-ov} depict the changes of the optimal values as $\delta_1$ varies from $0.01$ to $0.1$ with different $M$.

\textbf{(ii) SAA 
of exogenous uncertainty. }
In this set of experiments, we 
use
sample data to approximate the true probability distribution $P$ (of $\bdxi$), which is also known as  SAA.
We include this in the category of data perturbation in the sense that empirical distribution constructed with sample data may be regarded as a perturbation of $P$. 
We investigate how the variation of sample size affects the optimal values and the optimal solutions.
We solve problem (\ref{eq:PRO-N-reformulate}) with different sample size $K$ and run $20$ simulations for each fixed sample size $K$. 
We plot a boxplot diagram 
to examine the convergence of the optimal values as 
$K$ increases in Figures~\ref{subfig-SAA-main}-\ref{subfig-SAA-counter}.
We can see that as the sample size reaches 400, the optimal values of the SAA problem are close to the true optimal value in both Type-1 PLA and Type-2 PLA.

\subsection{Preference inconsistency}
\label{subsec:preference-incon}

In Section~\ref{sec:PC-design}, we consider pairwise comparisons to elicit the DM’s preference.
In practice, various errors may occur 
during
the elicitation process such as measurement errors
and DM’s wrong responses, 
all of which may lead to preference inconsistency. 
In this part, we examine the effects of the inconsistencies on 
the worst-case utility functions
and the optimal value  
in the following
two types of inconsistency 
during the preference elicitation process.

\begin{figure}[!ht]
  \centering
  \vspace{-0.5cm}
  \subfigure{
    \label{subfig-incon-main-ut1}
    \includegraphics[width=0.3\linewidth]{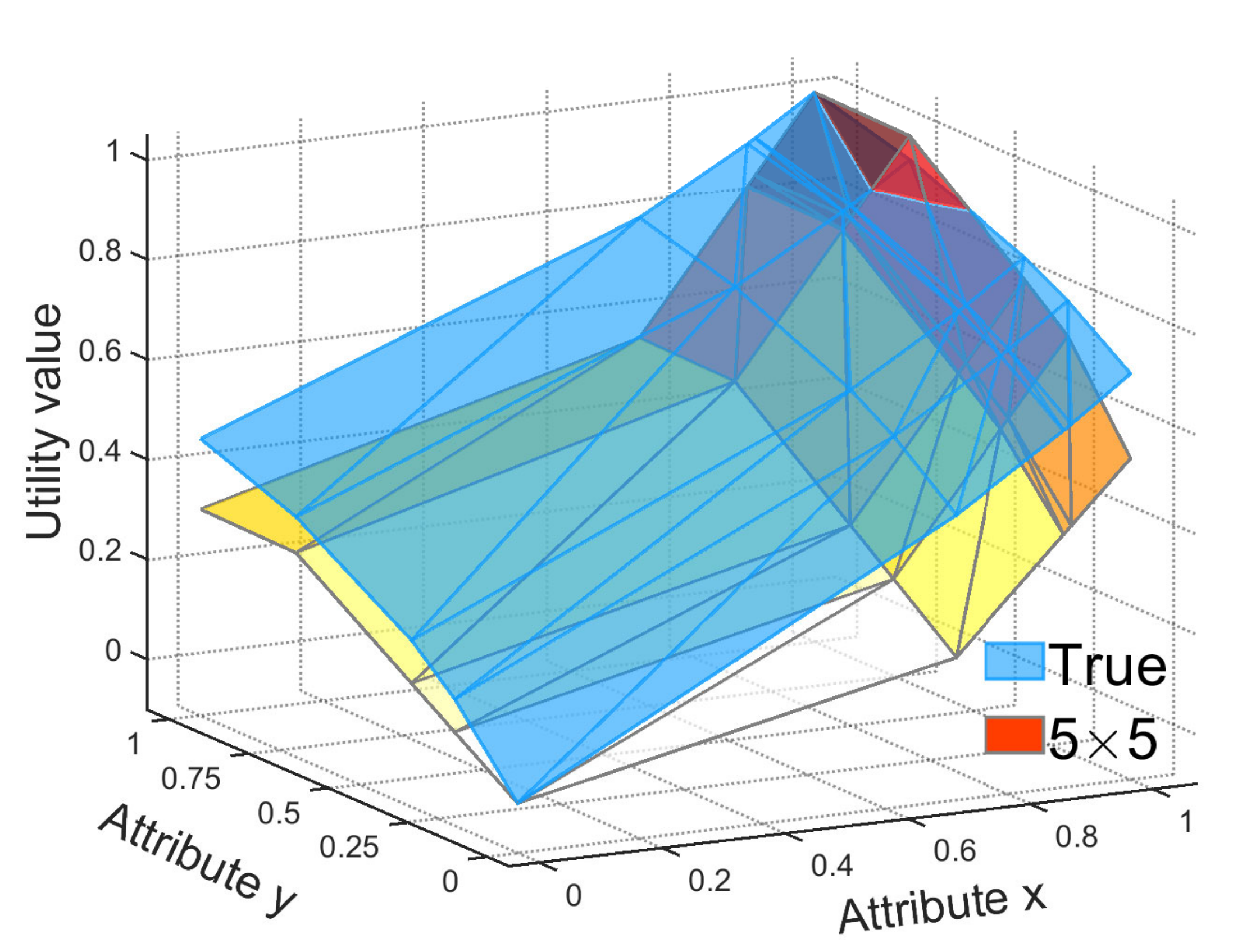}
  }
  \hspace{-0.5em}
  \subfigure{
    \label{subfig-incon-main-ut2}
    \includegraphics[width=0.3\linewidth]{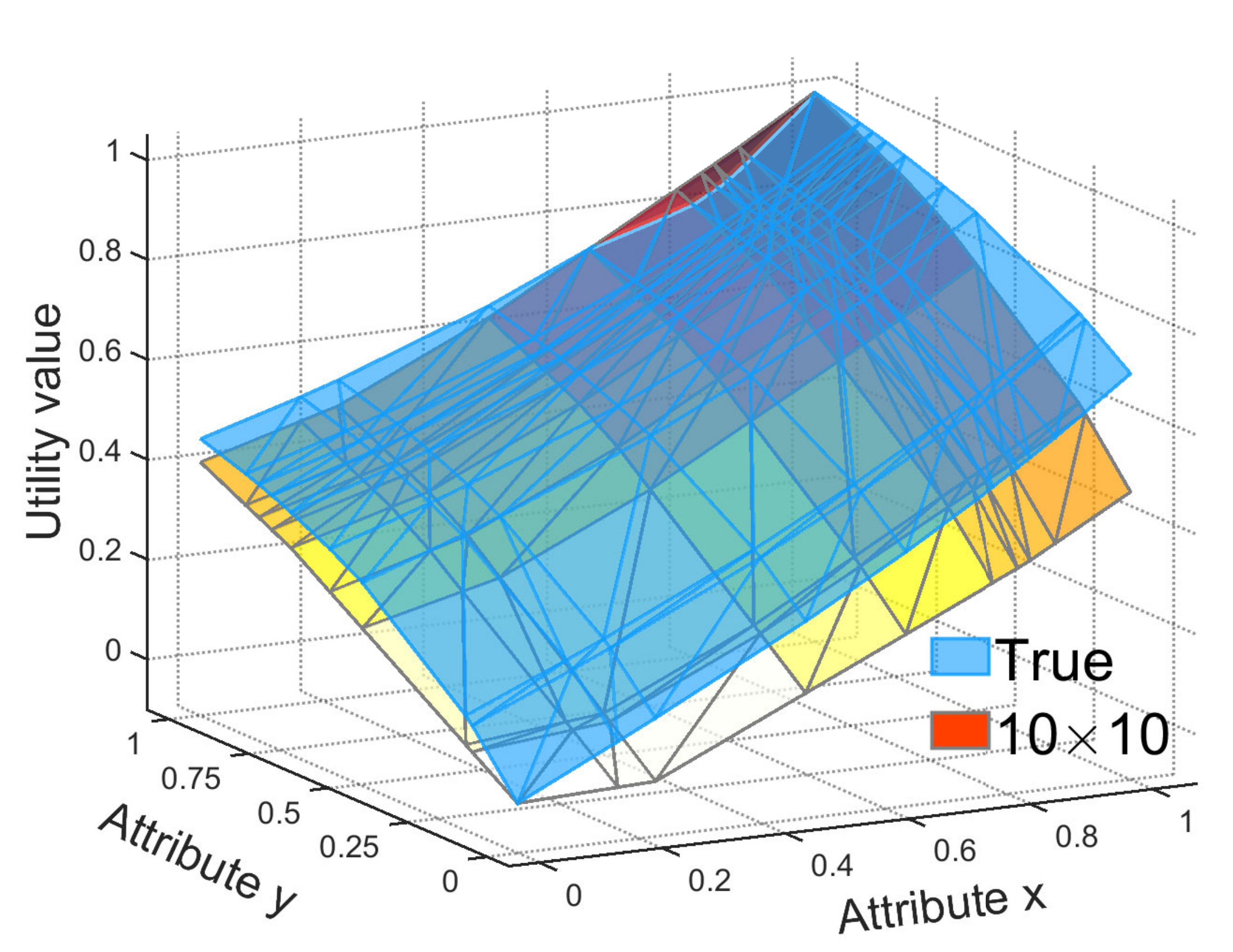}
  }
  \hspace{-0.5em}
  \subfigure{
    \label{subfig-incon-main-ut3}
    \includegraphics[width=0.3\linewidth]{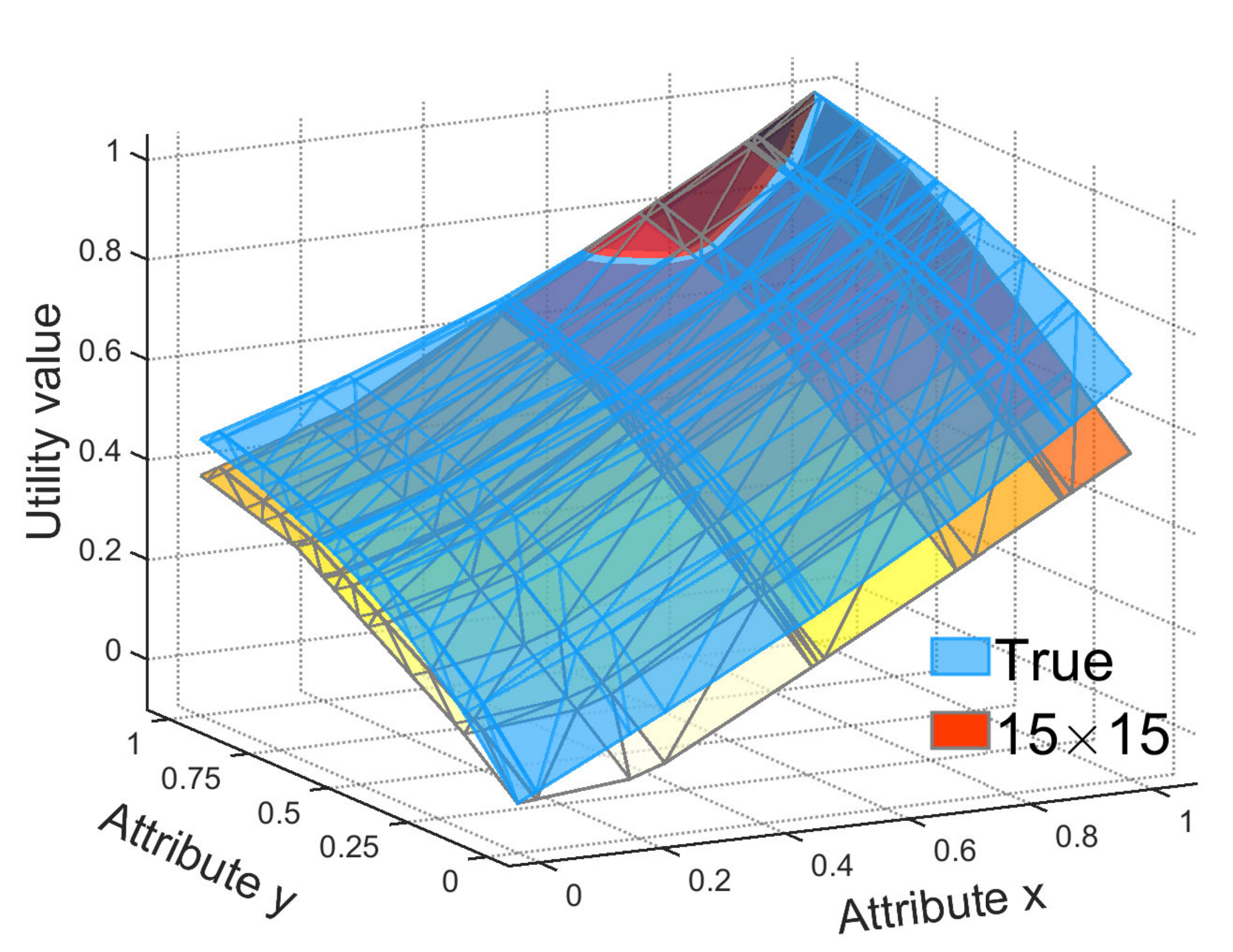}
  }
   \vspace{-1em}
  \captionsetup{font=footnotesize}
  \caption{{\bf Type-1 EPLA}: worst-case utility with $\Gamma=0.5$}
   \vspace{-0.2cm}
  \label{fig-incon-ut-main} 
\end{figure}

\begin{figure}[!ht]
  \centering
  \vspace{-0.5em}
  \subfigure{
    \label{subfig-incon-counter-ut1}
    \includegraphics[width=0.3\linewidth]{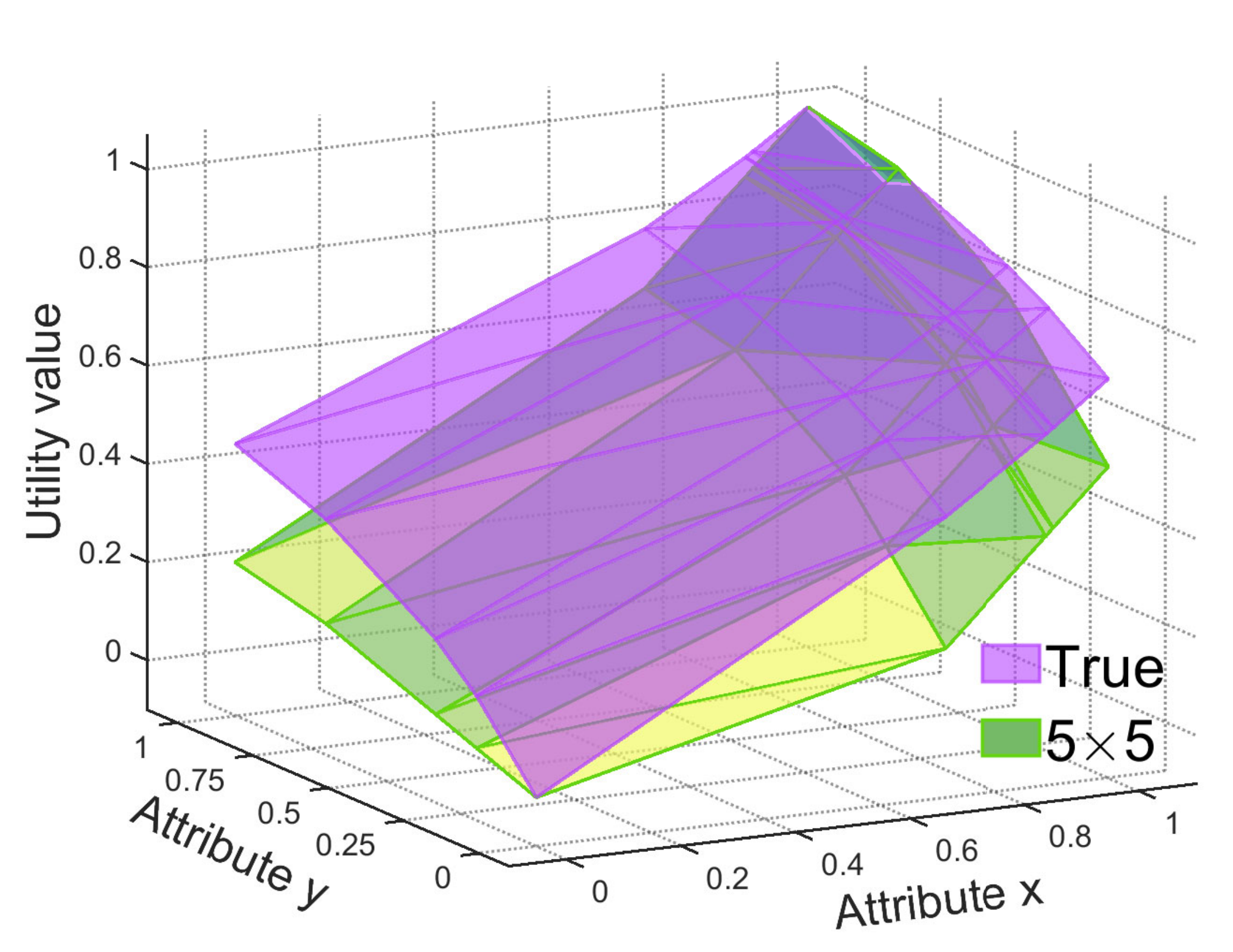}
  }
  \hspace{-0.5em}
  \subfigure{
    \label{subfig-incon-counter-ut2}
    \includegraphics[width=0.3\linewidth]{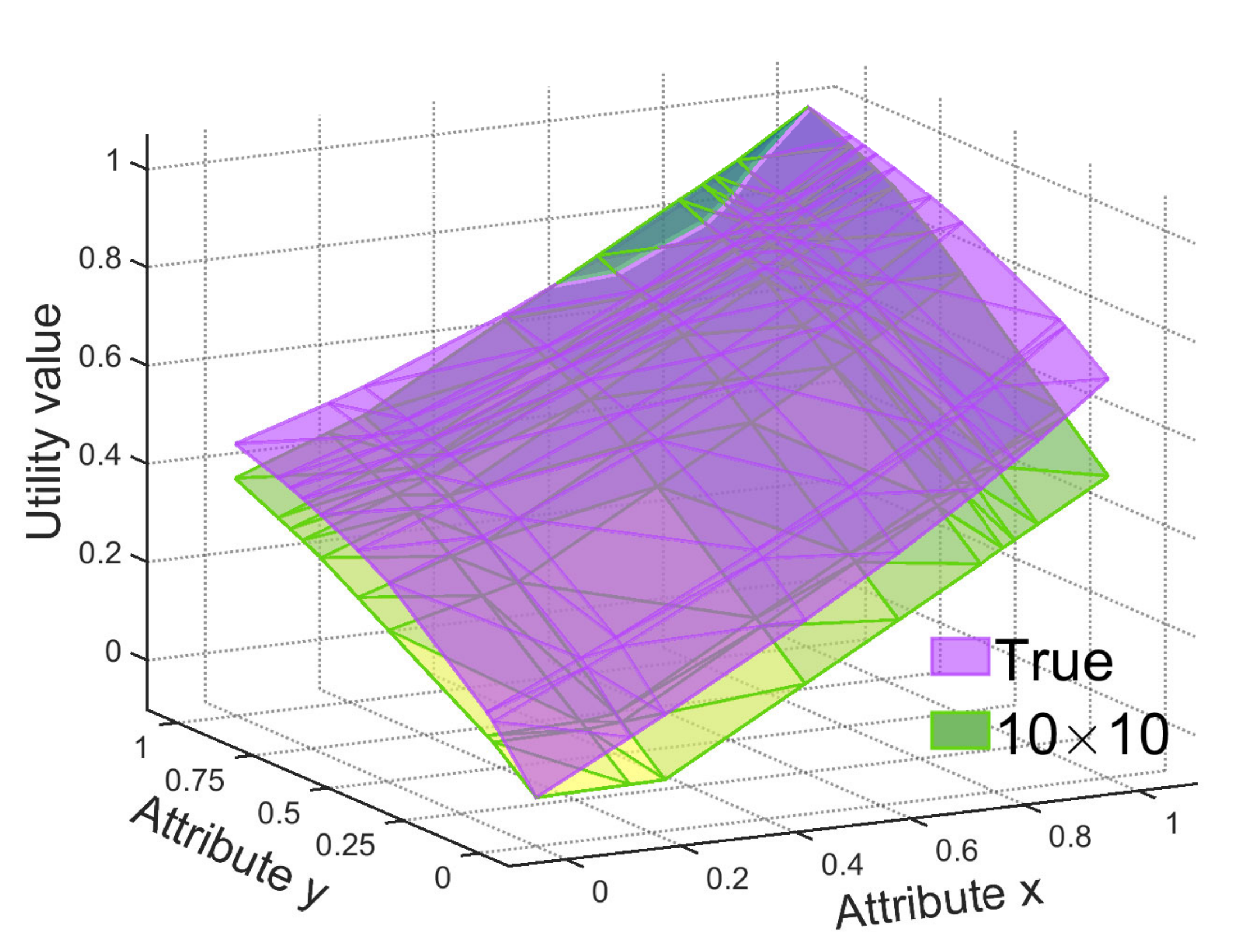}
  }
  \hspace{-0.5em}
  \subfigure{
    \label{subfig-incon-counter-ut3}
    \includegraphics[width=0.3\linewidth]{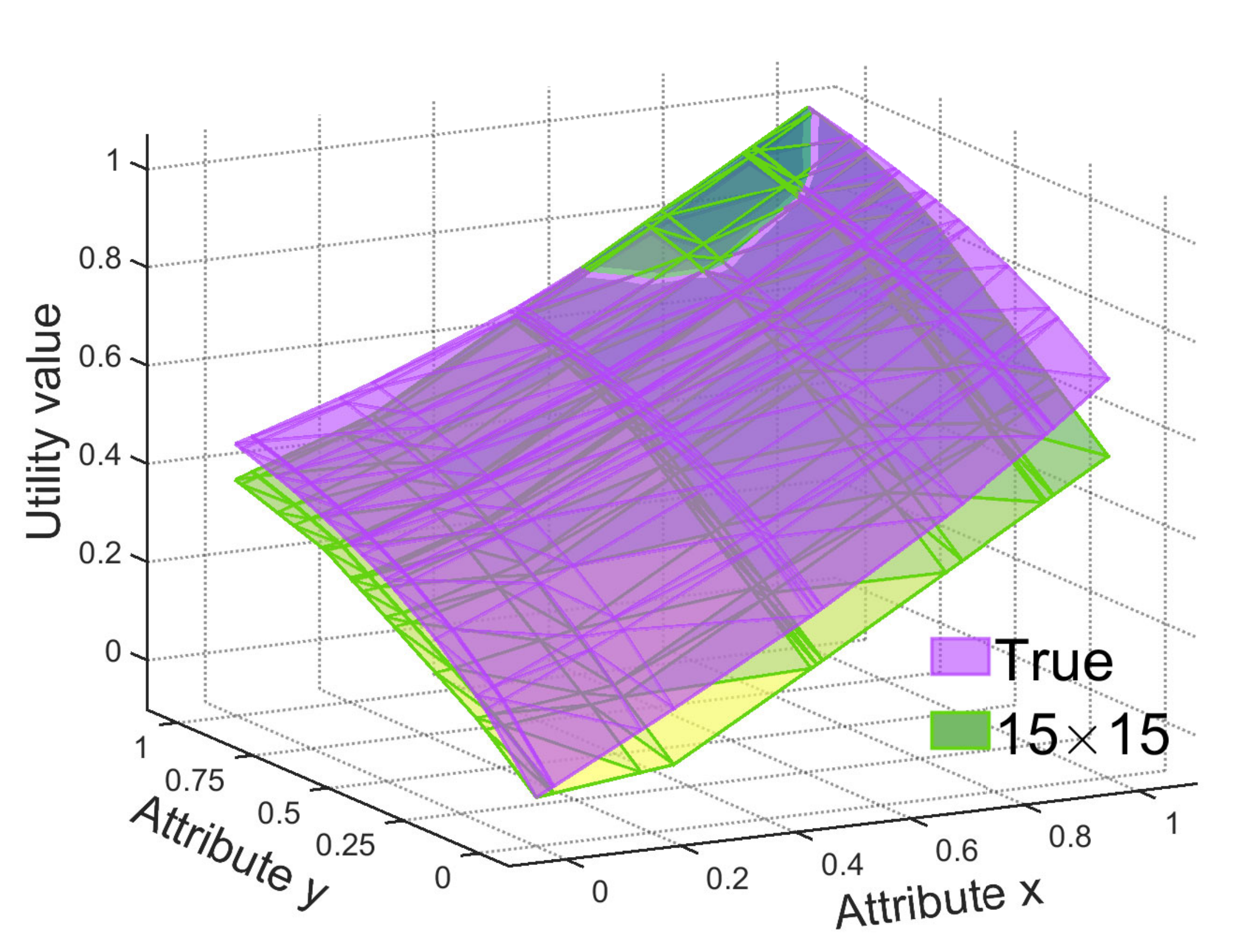}
  }
   \vspace{-1em}
  \captionsetup{font=footnotesize}
  \caption{{\bf Type-2 EPLA}: worst-case utility with $\Gamma=0.5$}
   \vspace{-0.2cm}
  \label{fig-incon-ut-counter} 
\end{figure}

\begin{figure}[!ht]
  \centering
  \vspace{-0.5em}
  \subfigure[]{
    \label{subfig-incon-main}
    \includegraphics[width=0.23\linewidth]{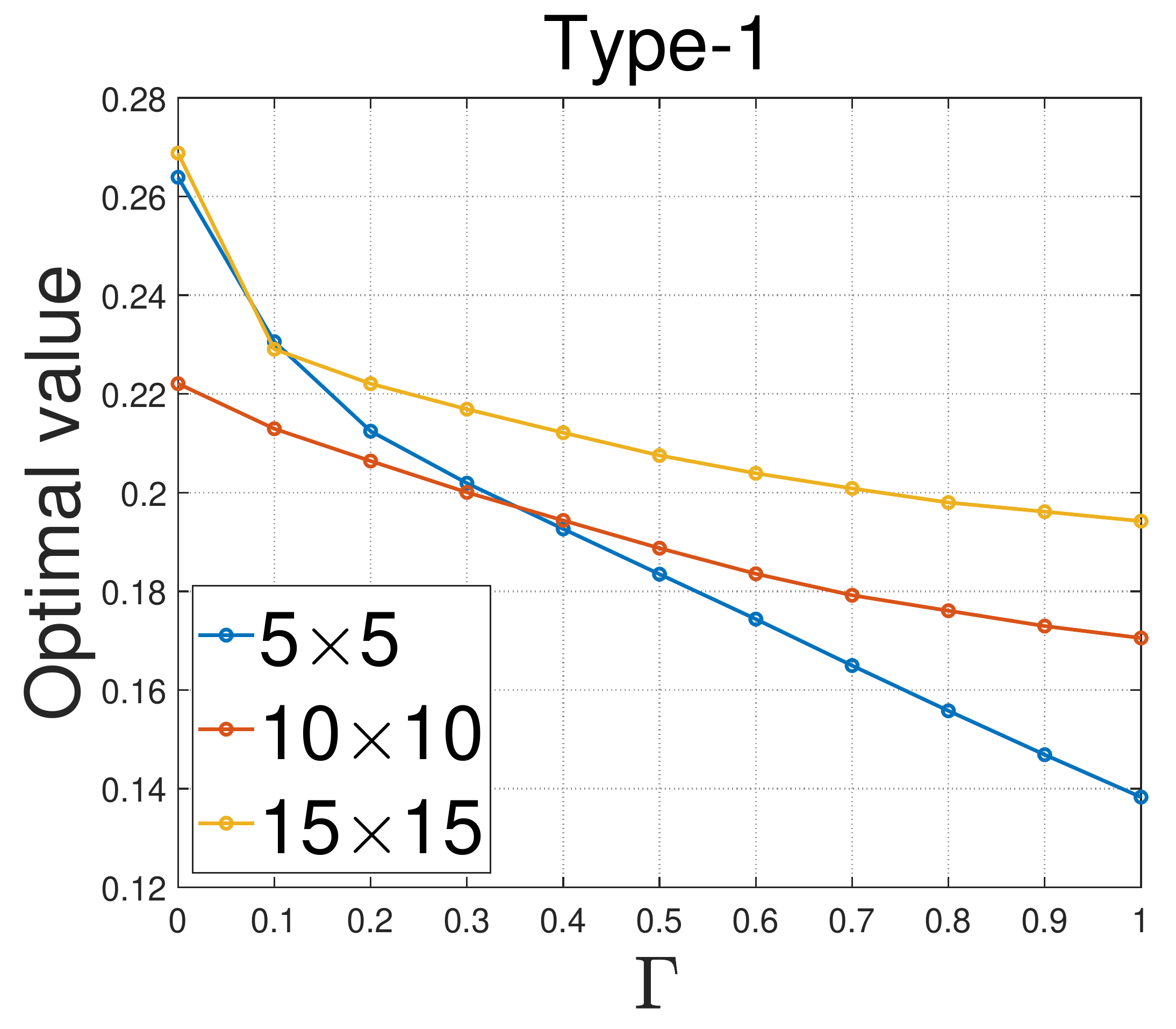}
  }
  \hspace{-1.5em}
  \subfigure[]{
    \label{subfig-incon-counter}
    \includegraphics[width=0.23\linewidth]{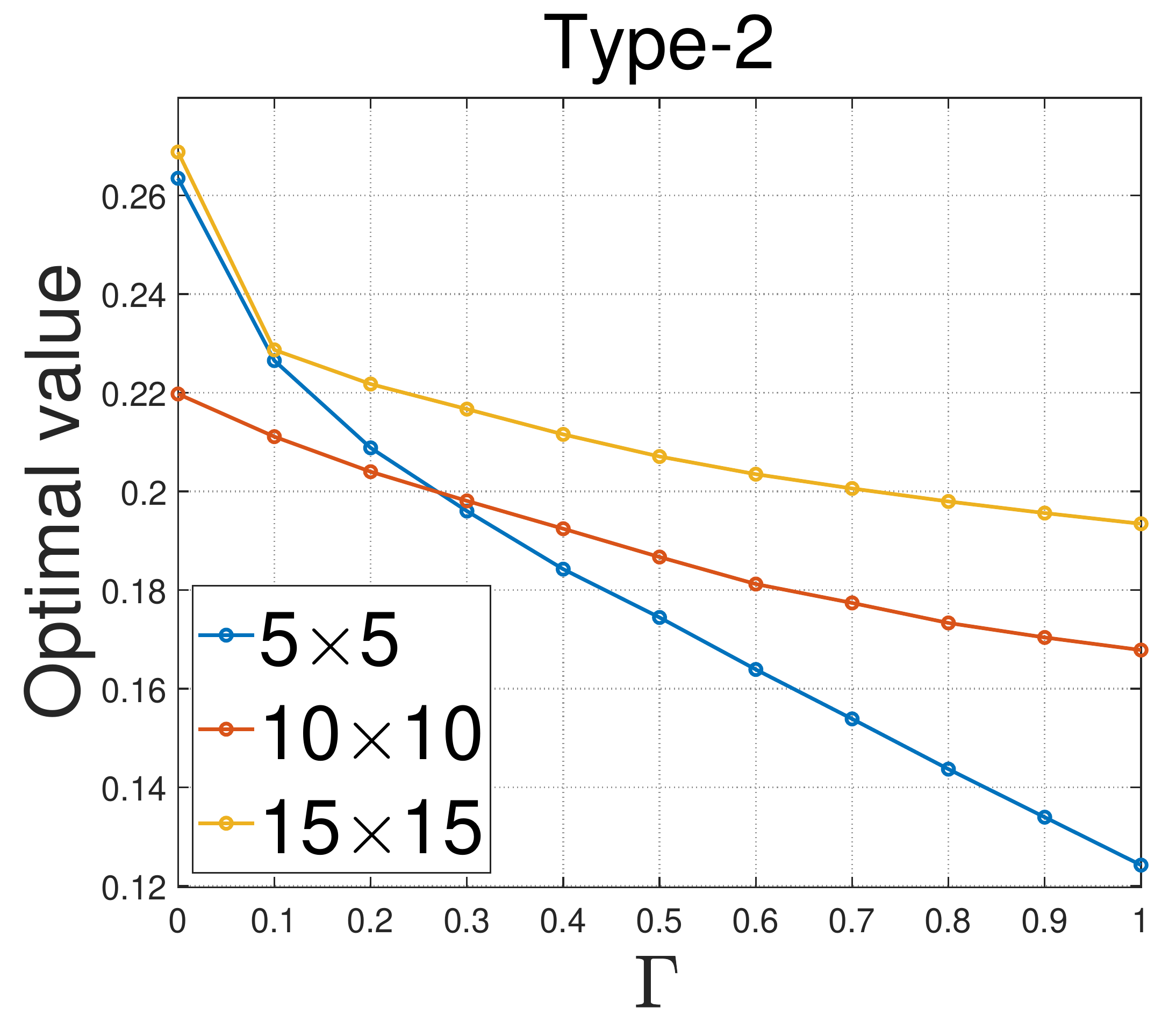}
  }
  \hspace{-1.5em}
  \subfigure[]{
    \label{subfig-responserr-main}
    \includegraphics[width=0.23\linewidth]{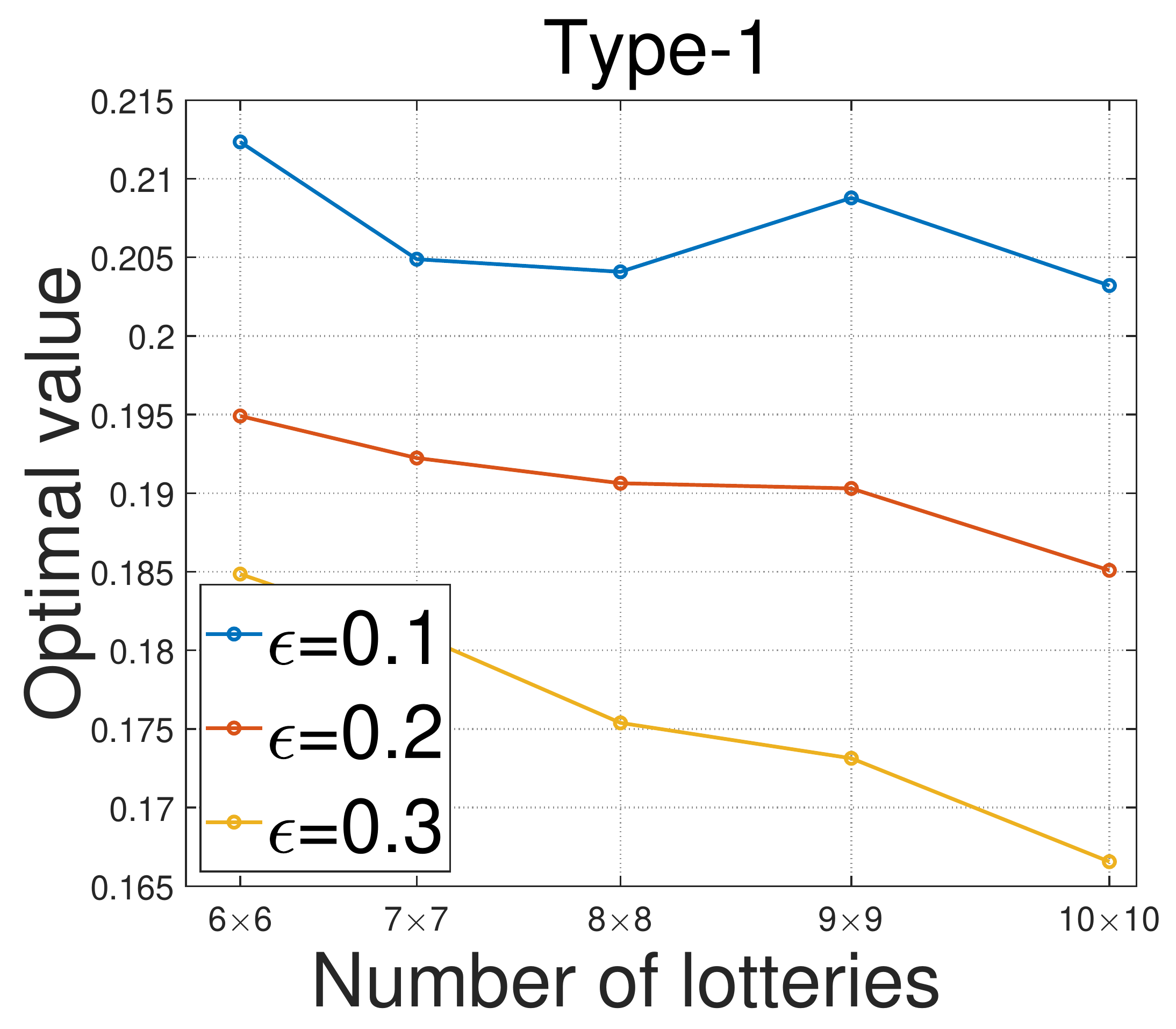}
  }
  \hspace{-1.5em}
  \subfigure[]{
    \label{subfig-responserr-counter}
    \includegraphics[width=0.23\linewidth]{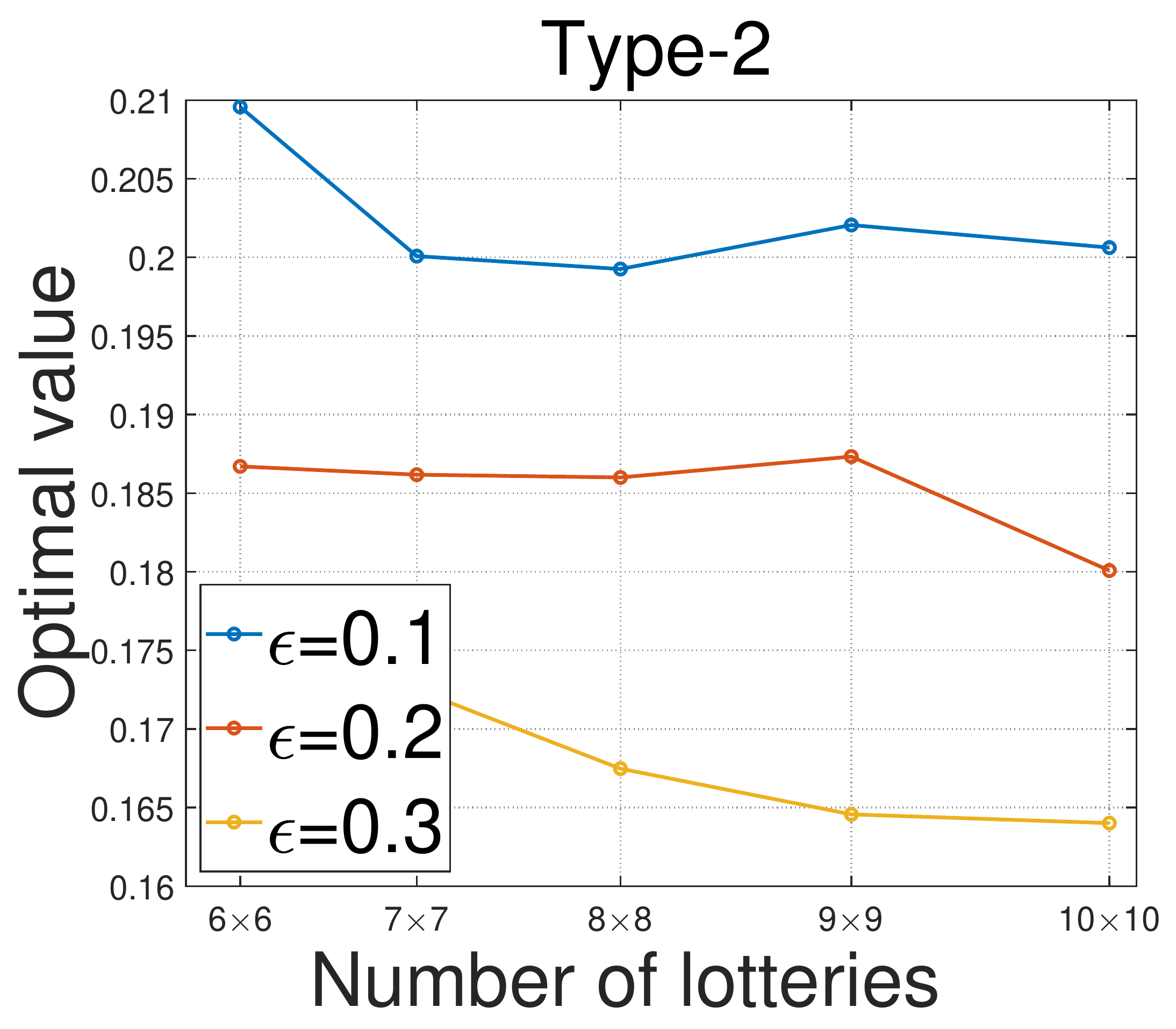}
  }
  \vspace{-1em}
  \captionsetup{font=footnotesize}
  \caption{{\bf EPLA}: the optimal values with total errors and erroneous responses}
   \vspace{-0.5cm}
  \label{fig-incon-ov} 
\end{figure}

\textbf{(i) Limitation on the total quantity of errors. }
We consider 
the rhs of the inequality constraints in the definition of
${\cal U}_N$
to be perturbed by positive constants $\gamma_l$, that is, $\la u_N,\psi_l \ra\leq c_l+\gamma_l, l=1,\ldots,M$.
The perturbation is required for the feasibility of problem (\ref{eq:PRO-N-reformulate}) to hold
when noise  corrupts
the expected utility evaluation when
a comparison is made. 
In other words, the perturbed inequalities accommodate potentially inconsistent responses.
We restrict the total inconsistency by setting $\sum_{l=1}^M \gamma_l\leq \Gamma$,
where $\Gamma$ is the total error to be tolerated.
Figures~\ref{fig-incon-ut-main}-\ref{fig-incon-ut-counter} depict the worst-case utility functions 
and the true utility function.
Figures~\ref{subfig-incon-main}-\ref{subfig-incon-counter} depict the optimal values with $\Gamma$ varies from $0$ to $1$.
As $\Gamma$ increases, the optimal values decrease.
From the figures, we 
find that our PLA approach works very well for this type of inconsistency.

\begin{figure}[!ht]
 \vspace{-0.2cm}
  \centering
  \subfigure{
    \label{subfig-responserr-main-ut1}
    \includegraphics[width=0.3\linewidth]{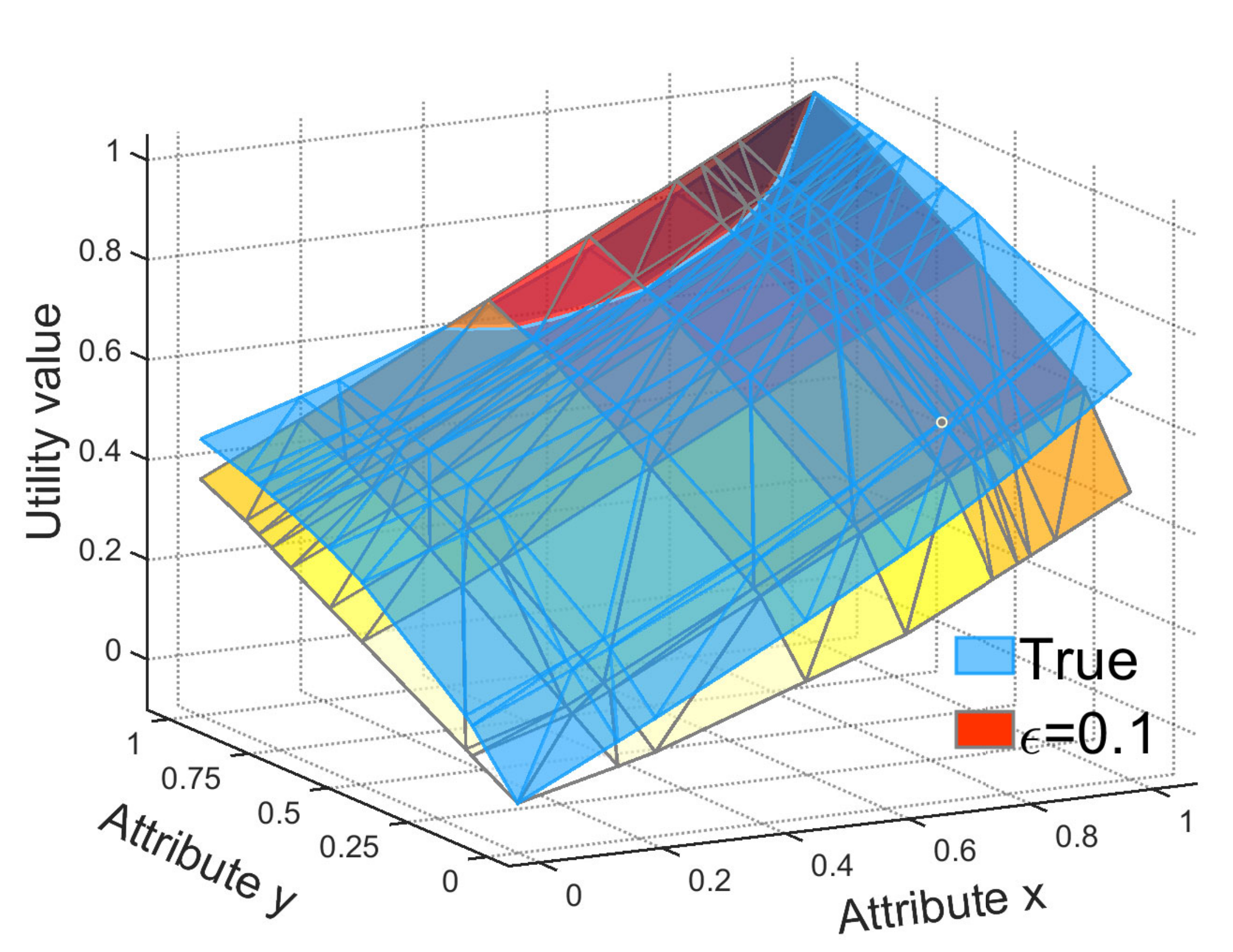}
  }
  \hspace{-0.5em}
  \subfigure{
    \label{subfig-responserr-main-ut2}
    \includegraphics[width=0.3\linewidth]{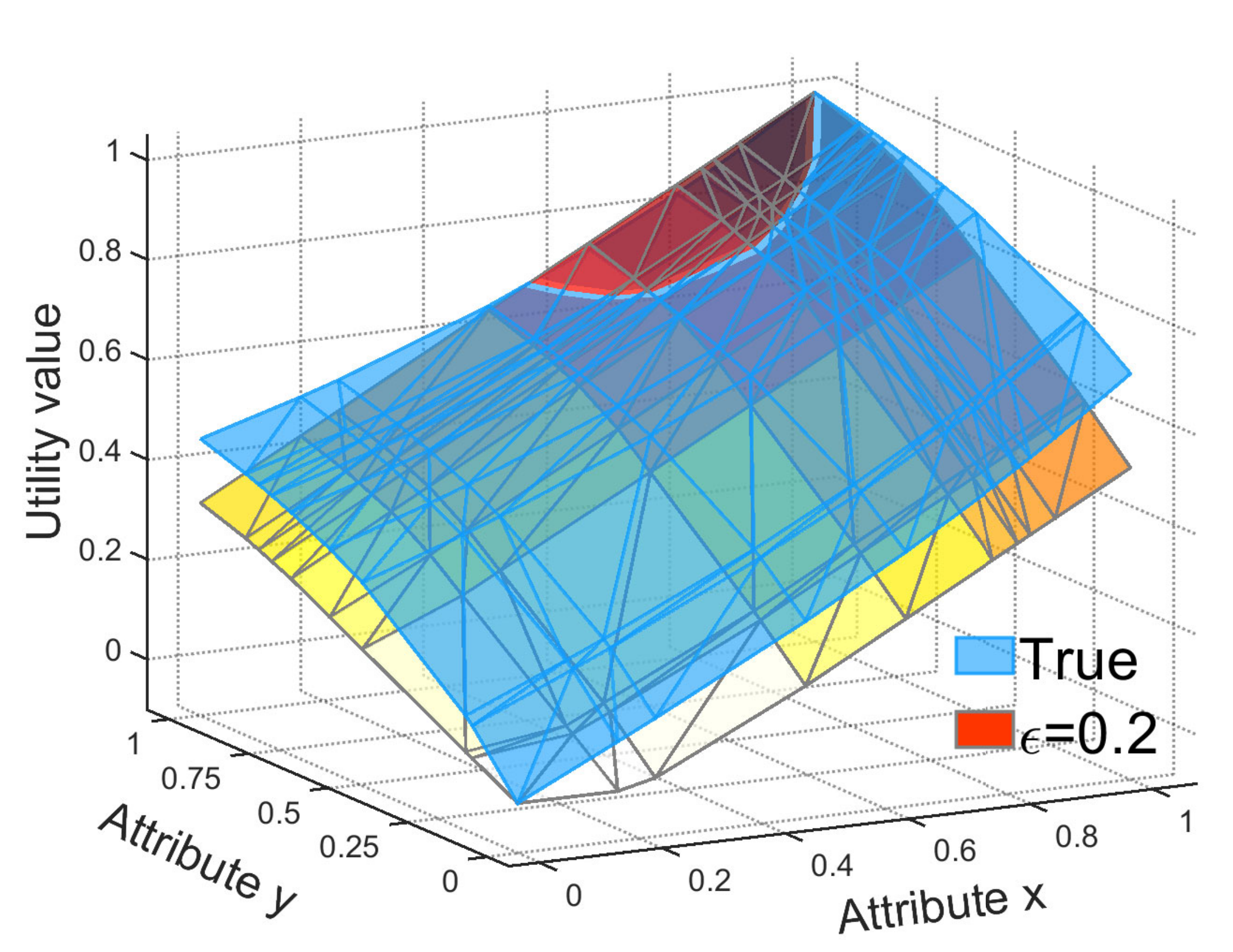}
  }
  \hspace{-0.5em}
  \subfigure{
    \label{subfig-responserr-main-ut3}
    \includegraphics[width=0.3\linewidth]{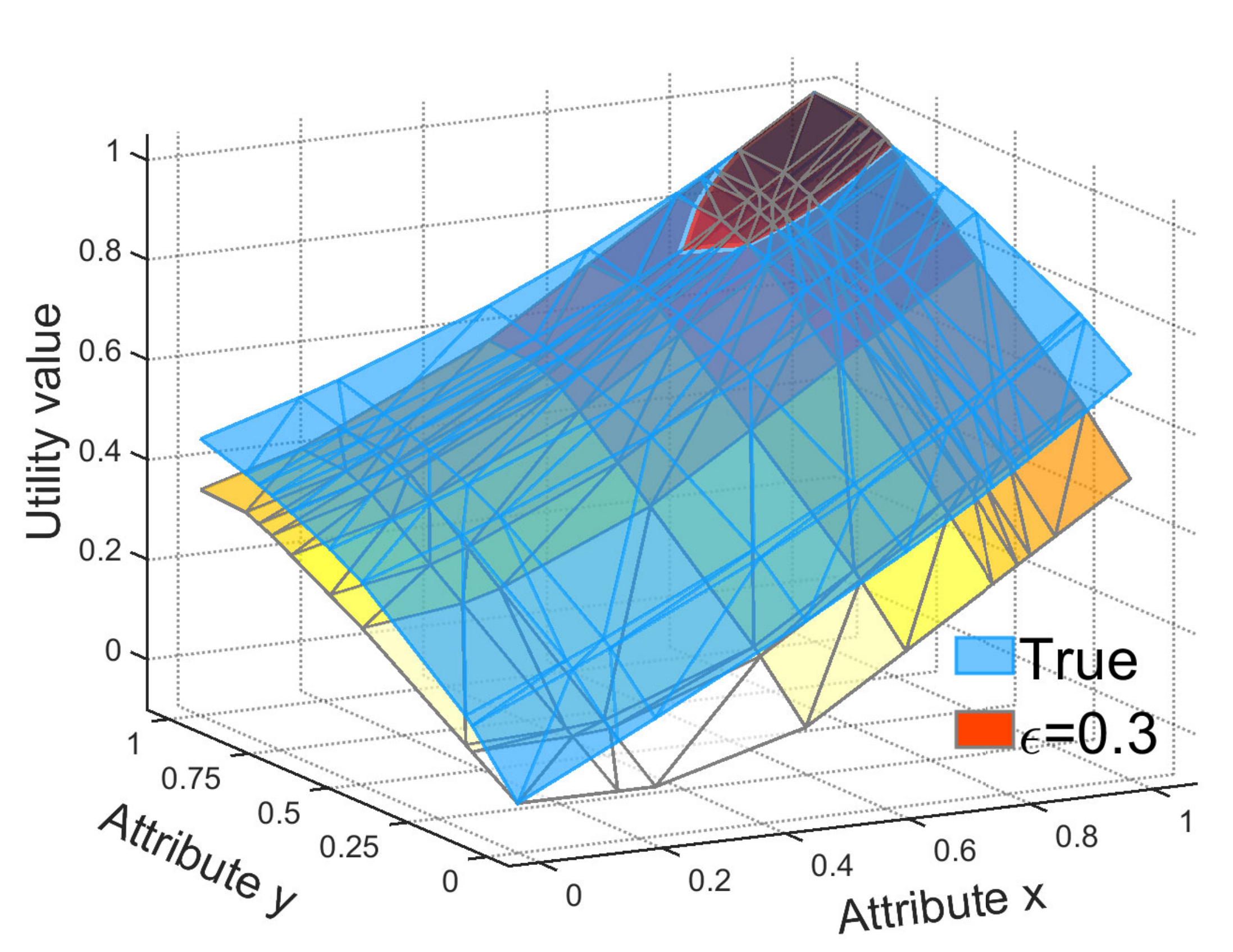}
  }
  \vspace{-1em}
  \captionsetup{font=footnotesize}
  \caption{{\bf Type-1 EPLA}: worst-case utility with $10\times10$ lotteries}
  \vspace{-0.2cm}
  \label{fig-responserr-ut-main}
\end{figure}

\begin{figure}[!ht]
 \vspace{-0.2cm}
  \centering
  \subfigure{
    \label{subfig-responserr-counter-ut1}
    \includegraphics[width=0.3\linewidth]{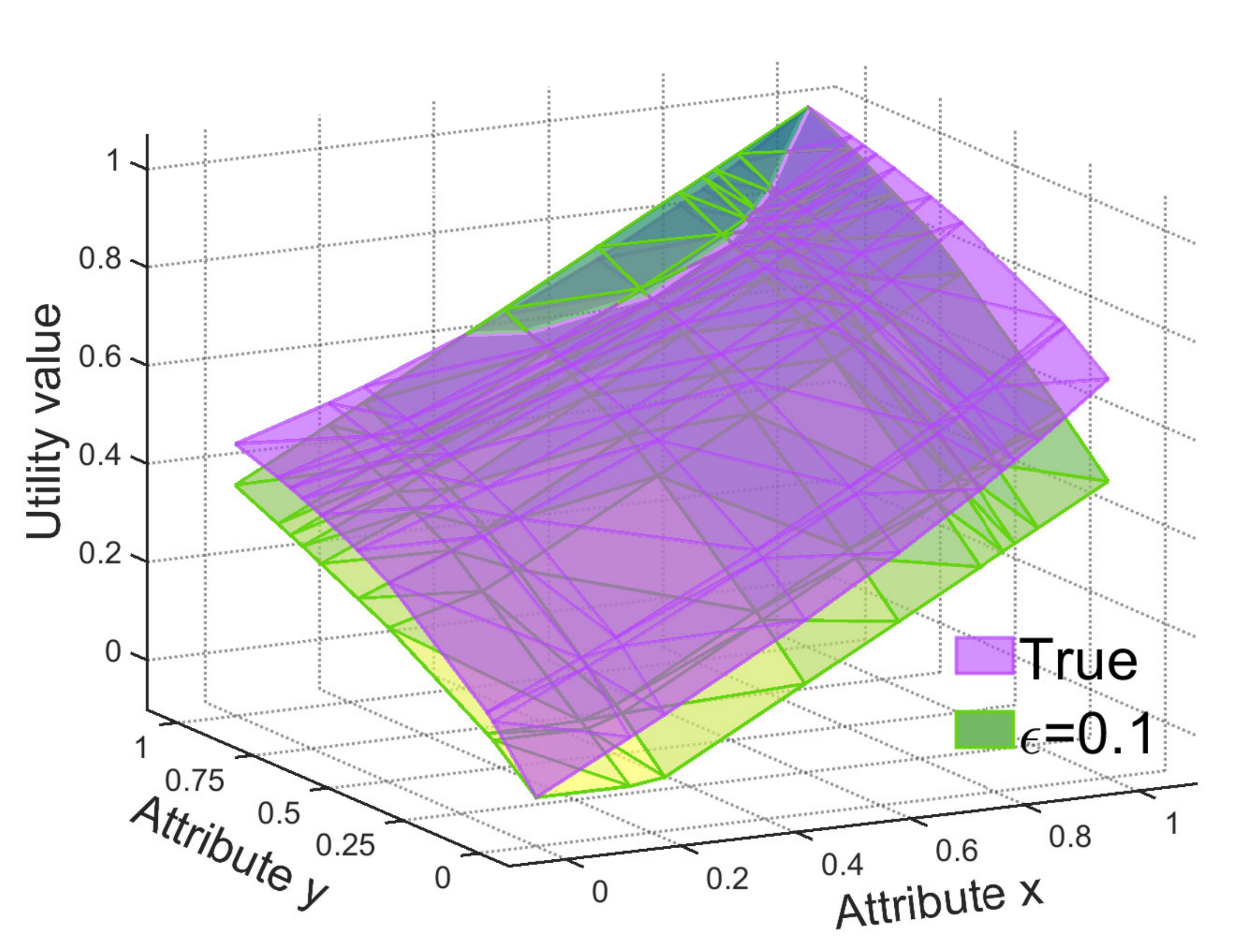}
  }
  \hspace{-0.5em}
  \subfigure{
    \label{subfig-responserr-counter-ut2}
    \includegraphics[width=0.3\linewidth]{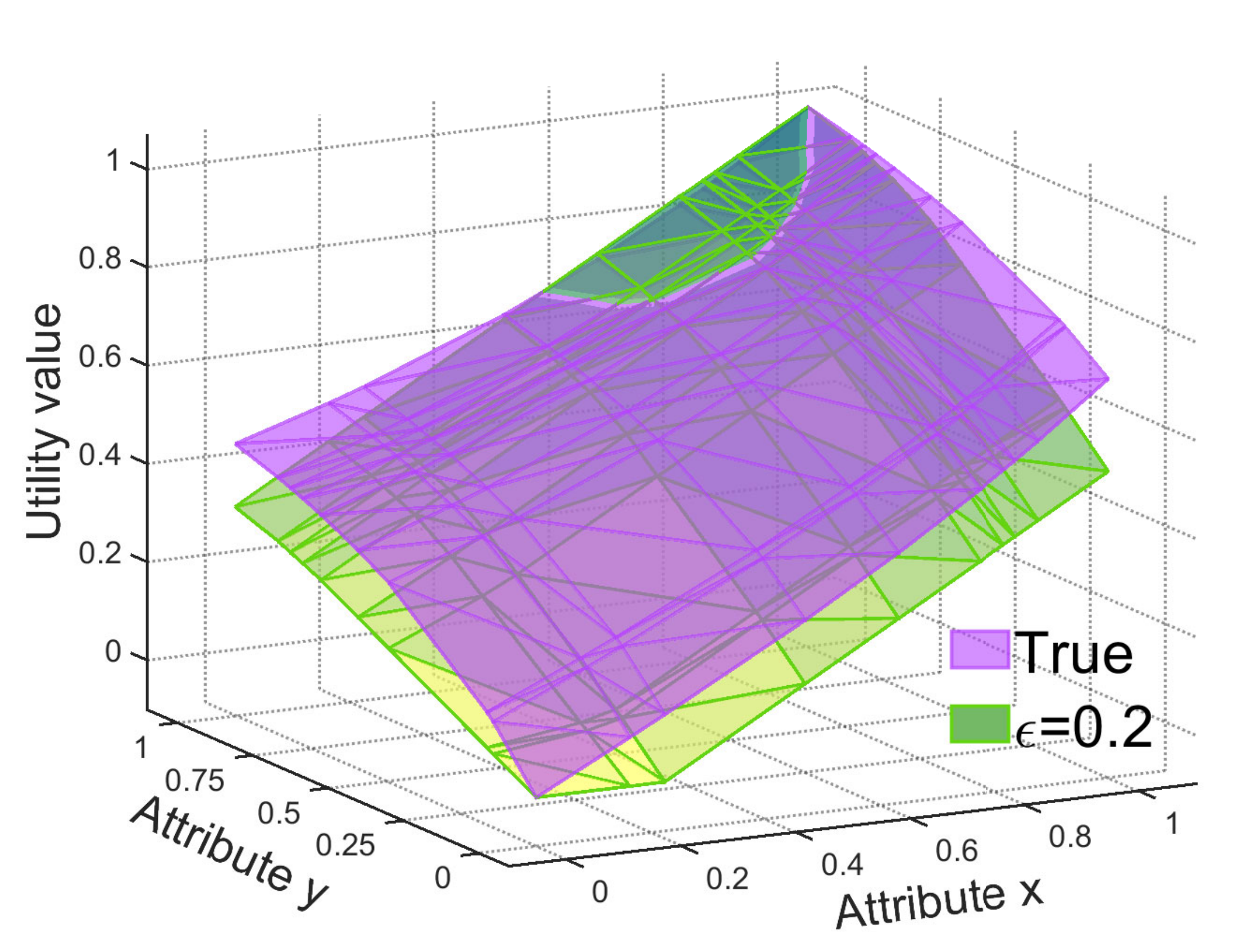}
  }
  \hspace{-0.5em}
  \subfigure{
    \label{subfig-responserr-counter-ut3}
    \includegraphics[width=0.3\linewidth]{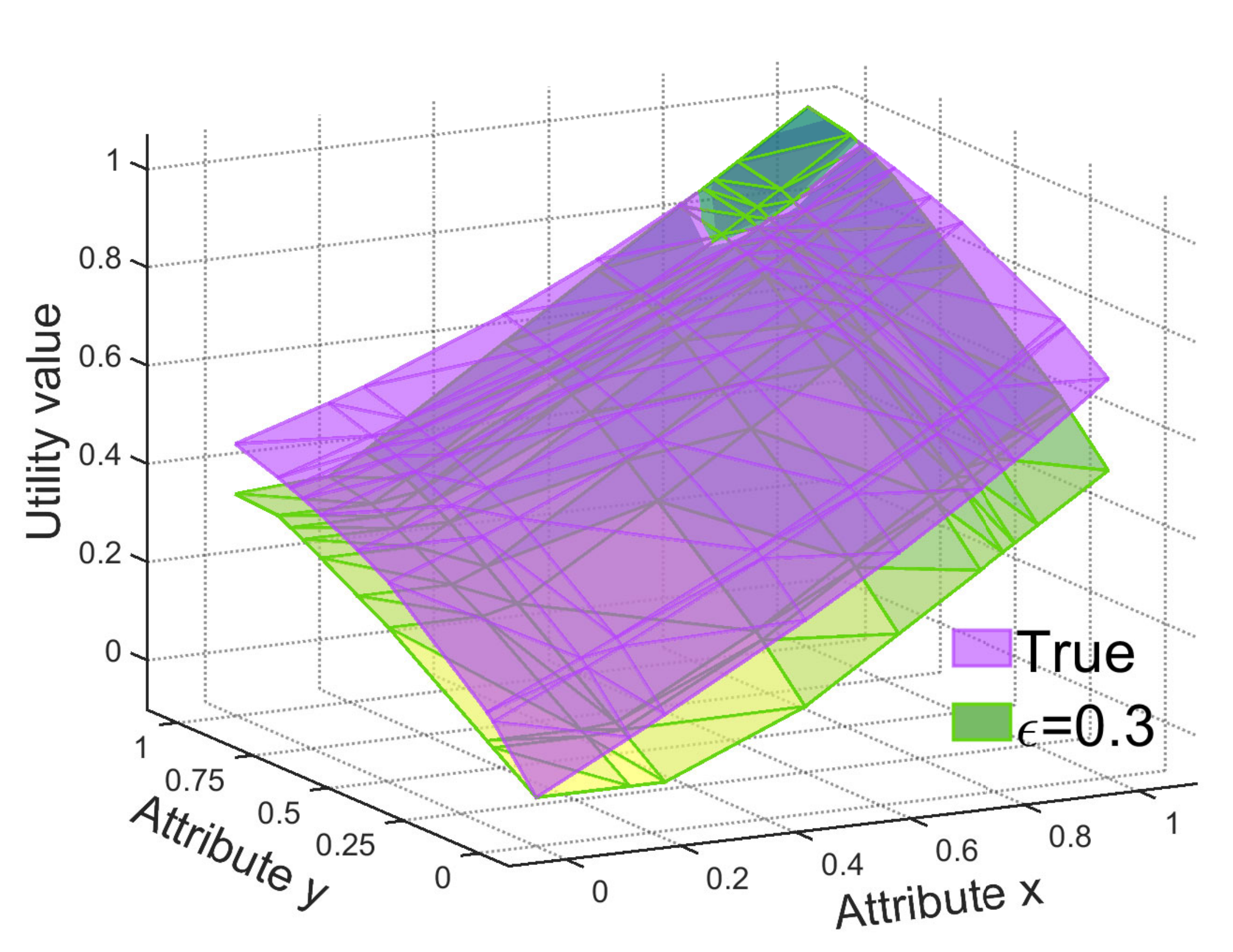}
  }
  \vspace{-1em}
  \captionsetup{font=footnotesize}
  \caption{{\bf Type-2 EPLA:} worst-case utility with $10\times10$ lotteries}
   \vspace{-0.2cm}
  \label{fig-responserr-ut-counter} 
\end{figure}

\textbf{(ii) Limitation on the number of erroneous responses. }
We consider the case that the DM makes mistakes occasionally, that is, 
the DM is mistaken at most $\epsilon M$ of 
lottery comparisons.
We introduce binary variable
$\delta_l$,
which takes
value $1$ if the DM is mistaken about lottery $l$ and $0$ otherwise, and we add the constraint $\sum_{l=1}^M \delta_l\leq \epsilon M$ to limit the total number of mistakes.
If the original comparison is $\bbe_{\mathbb{P}}[u(\bdcz_1^l(\omega))]\geq \bbe_{\mathbb{P}}[u(\bdcz_2^l(\omega))]$, then this condition is replaced by:
\begin{equation*}
    \delta_l \hat{M}+\bbe_{\mathbb{P}}[u(\bdcz_1^l(\omega))]\geq \bbe_{\mathbb{P}}[u(\bdcz_2^l(\omega))] 
    \quad \inmat{and} \quad 
    (1-\delta_l) \hat{M}+\bbe_{\mathbb{P}}[u(\bdcz_2^l(\omega))]\geq \bbe_{\mathbb{P}}[u(\bdcz_1^l(\omega))],
\end{equation*}
where $\hat{M}$ is a large constant (``Big $\hat{M}$'').
These constraints make the inner minimization problem become an MILP.
Figures~\ref{fig-responserr-ut-main}-\ref{fig-responserr-ut-counter} depict the worst-case utility functions, and the gap between 
them and 
the true utility function for Type-1 PLA and Type-2 PLA.
Figures~\ref{subfig-responserr-main}-\ref{subfig-responserr-counter} depict the optimal values with $\epsilon=\{0.1,0.2,0.3\}$.

\section{Concluding remarks}
\label{sec:Concluding remarks}

In this paper, we propose EPLA and IPLA approaches to approximate the true unknown utility function
in the multi-attribute UPRO models and 
demonstrate how the resulting approximate UPRO model 
can 
be solved.
The EPLA approach works only for two-attribute case as it stands because it is complex to derive
an explicit piecewise linear utility function when the utility function has three or more variables. The IPLA is not subject to the limitation of the dimension of the utility function but our numerical test results show that the IPLA-based approach takes considerably longer 
CPU time
to solve
as the numbers of preference elicitation questions and scenarios of exogenous random vector increase. 
This indicates that the formulation is potentially computationally unscalable. It remains an open question as to how to improve the computational efficiency of the IPLA approach. 
For instance, 
in the 
case when $m\geq 4$,
in order to derive IPLA of the utility function, 
we need to 
develop proper triangulation of the hypercube $\bigtimes_{i=1}^{m}[a_{i},b_{i}]$ 
into simplices in $m$-dimensional space.
It will be interesting
to explore such triangulation and
to identify the simplex where the reward function locates efficiently, 
see Hughes and Anderson \cite{hughes1996simplexity} and \cite{COTTLE198225,BROADIE198439} for further study.
Design of questionnaires
to elicit the DM's preference
is another point for potential improvement since our strategy is 
fundamentally based 
on
random  utility split 
scheme in single-attribute PRO models \cite{AmD15}. It 
is worthwhile to explore 
some optimal design strategies
such as in \cite{Vayanos2020}
because in practice, elicitation 
may be time consuming or costly.
Finally,
it will be interesting to explore whether the proposed approaches work more efficiently when the true utility function has some copula structure \cite{abbas2009multiattribute,abbas2013utility}. 
We leave all these for future research. 

\bibliographystyle{siamplain}      
\bibliography{2803-arxiv-main}  

\begin{thebibliography}{10}

\bibitem{abbas2009multiattribute}
{\sc A.~E. Abbas}, {\em Multiattribute utility copulas}, Operations Research,
  57 (2009), pp.~1367--1383.

\bibitem{abbas2013utility}
{\sc A.~E. Abbas}, {\em Utility copula functions matching all boundary
  assessments}, Operations Research, 61 (2013), pp.~359--371.

\bibitem{abbas2005attribute}
{\sc A.~E. Abbas and R.~A. Howard}, {\em Attribute dominance utility}, Decision
  Analysis, 2 (2005), pp.~185--206.

\bibitem{airoldi2011healthcare}
{\sc M.~Airoldi, A.~Morton, J.~Smith, and G.~Bevan}, {\em Healthcare
  prioritisation at the local level: {a} socio-technical approach}, \rm
  {W}orking paper, University of Oxford,  (2011).

\bibitem{andre2007non}
{\sc F.~J. Andr{\'e} and L.~Riesgo}, {\em A non-interactive elicitation method
  for non-linear multiattribute utility functions: Theory and application to
  agricultural economics}, European Journal of Operational Research, 181
  (2007), pp.~793--807.

\bibitem{Ans22}
{\sc J.~Ansari}, {\em On a version of a multivariate integration by parts
  formula for lebesgue integrals}, arXiv preprint arXiv:2203.06772,  (2022).

\bibitem{AmD15}
{\sc B.~Armbruster and E.~Delage}, {\em Decision making under uncertainty when
  preference information is incomplete}, Management Science, 61 (2015),
  pp.~111--128.

\bibitem{ash2000probability}
{\sc R.~B. Ash, B.~Robert, C.~A. Doleans-Dade, and A.~Catherine}, {\em
  Probability and measure theory}, Academic press, 2000.

\bibitem{azaron2008multi}
{\sc A.~Azaron, K.~Brown, S.~Tarim, and M.~Modarres}, {\em A multi-objective
  stochastic programming approach for supply chain design considering risk},
  International Journal of Production Economics, 116 (2008), pp.~129--138.

\bibitem{BeO13}
{\sc D.~Bertsimas and A.~O'Hair}, {\em Learning preferences under noise and
  loss aversion: An optimization approach}, Operations Research, 61 (2013),
  pp.~1190--1199.

\bibitem{BROADIE198439}
{\sc M.~N. Broadie and R.~W. Cottle}, {\em A note on triangulating the 5-cube},
  Discrete Mathematics, 52 (1984), pp.~39--49.

\bibitem{chen2010stochastic}
{\sc A.~Chen, J.~Kim, S.~Lee, and Y.~Kim}, {\em Stochastic multi-objective
  models for network design problem}, Expert Systems with Applications, 37
  (2010), pp.~1608--1619.

\bibitem{chien1977solving}
{\sc M.-J. Chien and E.~Kuh}, {\em Solving nonlinear resistive networks using
  piecewise-linear analysis and simplicial subdivision}, IEEE Transactions on
  Circuits and Systems, 24 (1977), pp.~305--317.

\bibitem{clarkson1933definitions}
{\sc J.~A. Clarkson and C.~R. Adams}, {\em On definitions of bounded variation
  for functions of two variables}, Transactions of the American Mathematical
  Society, 35 (1933), pp.~824--854.

\bibitem{clemen2013making}
{\sc R.~T. Clemen and T.~Reilly}, {\em {M}aking {H}ard {D}ecisions with
  {D}ecision {T}ools}, Pacific Grove, Duxbury, 2013.

\bibitem{COTTLE198225}
{\sc R.~W. Cottle}, {\em Minimal triangulation of the 4-cube}, Discrete
  Mathematics, 40 (1982), pp.~25--29.

\bibitem{DGX22}
{\sc E.~Delage, S.~Guo, and H.~Xu}, {\em Shortfall risk models when information
  on loss function is incomplete}, Operations Research,  (2022),
  \url{https://doi.org/10.1287/opre.2021.2212}.

\bibitem{Dino2017}
{\sc I.~G. Dino and G.~{\"U}{\c{c}}oluk}, {\em Multiobjective design
  optimization of building space layout, energy, and daylighting performance},
  Journal of Computing in Civil Engineering, 31 (2017), p.~04017025.

\bibitem{duncan1977matrix}
{\sc G.~T. Duncan}, {\em A matrix measure of multivariate local risk aversion},
  Econometrica, 45 (1977), pp.~895--903.

\bibitem{DLM10}
{\sc C.~D’Ambrosio, A.~Lodi, and S.~Martello}, {\em Piecewise linear
  approximation of functions of two variables in milp models}, Operations
  Research Letters, 38 (2010), pp.~39--46.

\bibitem{feeny2002multiattribute}
{\sc D.~Feeny, W.~Furlong, G.~W. Torrance, C.~H. Goldsmith, Z.~Zhu, S.~DePauw,
  M.~Denton, and M.~Boyle}, {\em Multiattribute and single-attribute utility
  functions for the health utilities index mark 3 system}, Medical {C}are, 40
  (2002), pp.~113--128.

\bibitem{fishburn1992multiattribute}
{\sc P.~C. Fishburn and I.~H. LaValle}, {\em Multiattribute expected utility
  without the archimedean axiom}, Journal of Mathematical Psychology, 36
  (1992), pp.~573--591.

\bibitem{fliege2014robust}
{\sc J.~Fliege and R.~Werner}, {\em Robust multiobjective optimization \&
  applications in portfolio optimization}, European Journal of Operational
  Research, 234 (2014), pp.~422--433.

\bibitem{GiS02}
{\sc A.~L. Gibbs and F.~E. Su}, {\em On choosing and bounding probability
  metrics}, International Statistical Review, 70 (2002), pp.~419--435.

\bibitem{gonzalez2018utility}
{\sc J.~Gonz{\'a}lez-Ortega, V.~Radovic, and D.~R{\'\i}os~Insua}, {\em Utility
  elicitation}, in Elicitation, Springer, 2018, pp.~241--264.

\bibitem{greco2016multiple}
{\sc S.~Greco, J.~Figueira, and M.~Ehrgott}, {\em Multiple Criteria Decision
  Analysis}, vol.~37, Springer, 2016.

\bibitem{guo2022robust}
{\sc S.~Guo and H.~Xu}, {\em Robust spectral risk optimization when the
  subjective risk aversion is ambiguous: a moment-type approach}, Mathematical
  Programming, 194 (2022), pp.~305--340.

\bibitem{GXZ21}
{\sc S.~Guo, H.~Xu, and S.~Zhang}, {\em Utility preference robust optimization
  with moment-type information structure}, To appear in Operations Research,
  (2023).

\bibitem{haskell2016ambiguity}
{\sc W.~B. Haskell, L.~Fu, and M.~Dessouky}, {\em Ambiguity in risk preferences
  in robust stochastic optimization}, European Journal of Operational Research,
  254 (2016), pp.~214--225.

\bibitem{hildebrandt1963introduction}
{\sc T.~H. Hildebrandt}, {\em Introduction to the Theory of Integration}, Pure
  and Applied Mathematics, Vol. XIII, Academic Press, New York-London, 1963.

\bibitem{hu2018robust}
{\sc J.~Hu, M.~Bansal, and S.~Mehrotra}, {\em Robust decision making using a
  general utility set}, European Journal of Operational Research, 269 (2018),
  pp.~699--714.

\bibitem{hu2015robust}
{\sc J.~Hu and S.~Mehrotra}, {\em Robust decision making over a set of random
  targets or risk-averse utilities with an application to portfolio
  optimization}, IIE Transactions, 47 (2015), pp.~358--372.

\bibitem{hu2022distributionally}
{\sc J.~Hu, D.~Zhang, H.~Xu, and S.~Zhang}, {\em Distributionally preference
  robust optimization in multi-attribute decision making}, arXiv preprint
  arXiv:2206.04491,  (2022).

\bibitem{hughes1996simplexity}
{\sc R.~B. Hughes and M.~R. Anderson}, {\em Simplexity of the cube}, Discrete
  Mathematics, 158 (1996), pp.~99--150.

\bibitem{karni1979multivariate}
{\sc E.~Karni}, {\em On multivariate risk aversion}, Econometrica, 47 (1979),
  pp.~1391--1401.

\bibitem{KDN04}
{\sc A.~B. Keha, I.~R. de~Farias~Jr, and G.~L. Nemhauser}, {\em Models for
  representing piecewise linear cost functions}, Operations Research Letters,
  32 (2004), pp.~44--48.

\bibitem{LeW01}
{\sc J.~Lee and D.~Wilson}, {\em Polyhedral methods for piecewise-linear
  functions {I}: the lambda method}, Discrete Applied Mathematics, 108 (2001),
  pp.~269--285.

\bibitem{levy1991arrow}
{\sc H.~Levy and A.~Levy}, {\em Arrow-pratt measures of risk aversion: The
  multivariate case}, International Economic Review, 32 (1991), pp.~891--898.

\bibitem{liesio2021nonadditive}
{\sc J.~Liesi{\"o} and E.~Vilkkumaa}, {\em Nonadditive multiattribute utility
  functions for portfolio decision analysis}, Operations Research, 69 (2021),
  pp.~1886--1908.

\bibitem{liu2021multistage}
{\sc J.~Liu, Z.~Chen, and H.~Xu}, {\em Multistage utility preference robust
  optimization}, arXiv preprint arXiv:2109.04789,  (2021).

\bibitem{lofberg2004yalmip}
{\sc J.~Lofberg}, {\em Yalmip: A toolbox for modeling and optimization in
  matlab}, in 2004 IEEE international conference on robotics and automation
  (IEEE Cat. No. 04CH37508), IEEE, 2004, pp.~284--289.

\bibitem{maccheroni2002maxmin}
{\sc F.~Maccheroni}, {\em Maxmin under risk}, Economic Theory, 19 (2002),
  pp.~823--831.

\bibitem{Mcs47}
{\sc E.~J. McShane}, {\em Integration}, Princeton University Press, Princeton,
  1947.

\bibitem{meyer2005convex}
{\sc C.~A. Meyer and C.~A. Floudas}, {\em Convex envelopes for edge-concave
  functions}, Mathematical programming, 103 (2005), pp.~207--224.

\bibitem{misener2010piecewise}
{\sc R.~Misener and C.~Floudas}, {\em Piecewise-linear approximations of
  multidimensional functions}, Journal of Optimization Theory and Applications,
  145 (2010), pp.~120--147.

\bibitem{miyamoto1996multiattribute}
{\sc J.~M. Miyamoto and P.~Wakker}, {\em Multiattribute utility theory without
  expected utility foundations}, Operations Research, 44 (1996), pp.~313--326.

\bibitem{richard1975multivariate}
{\sc S.~F. Richard}, {\em Multivariate risk aversion, utility independence and
  separable utility functions}, Management Science, 22 (1975), pp.~12--21.

\bibitem{Rob75}
{\sc S.~M. Robinson}, {\em An application of error bounds for convex
  programming in a linear space}, SIAM Journal on Control, 13 (1975),
  pp.~271--273.

\bibitem{Rom03}
{\sc W.~R{\"o}misch}, {\em Stability of stochastic programming problems},
  Handbooks in operations research and management science, 10 (2003),
  pp.~483--554.

\bibitem{torrance1982application}
{\sc G.~W. Torrance, M.~H. Boyle, and S.~P. Horwood}, {\em Application of
  multi-attribute utility theory to measure social preferences for health
  states}, Operations Research, 30 (1982), pp.~1043--1069.

\bibitem{tsanakas2003risk}
{\sc A.~Tsanakas and E.~Desli}, {\em Risk measures and theories of choice},
  British Actuarial Journal, 9 (2003), pp.~959--991.

\bibitem{tseng1990minimax}
{\sc C.~Tseng and T.~Lu}, {\em Minimax multiobjective optimization in
  structural design}, International Journal for Numerical Methods in
  Engineering, 30 (1990), pp.~1213--1228.

\bibitem{tsetlin2006equivalent}
{\sc I.~Tsetlin and R.~L. Winkler}, {\em On equivalent target-oriented
  formulations for multiattribute utility}, Decision Analysis, 3 (2006),
  pp.~94--99.

\bibitem{tsetlin2007decision}
{\sc I.~Tsetlin and R.~L. Winkler}, {\em Decision making with multiattribute
  performance targets: The impact of changes in performance and target
  distributions}, Operations Research, 55 (2007), pp.~226--233.

\bibitem{tsetlin2009multiattribute}
{\sc I.~Tsetlin and R.~L. Winkler}, {\em Multiattribute utility satisfying a
  preference for combining good with bad}, Management Science, 55 (2009),
  pp.~1942--1952.

\bibitem{Vayanos2020}
{\sc P.~Vayanos, Y.~Ye, D.~McElfresh, J.~Dickerson, and E.~Rice}, {\em Robust
  active preference elicitation}, arXiv preprint arXiv:2003.01899,  (2020).

\bibitem{vielma2015mixed}
{\sc J.~P. Vielma}, {\em Mixed integer linear programming formulation
  techniques}, SIAM Review, 57 (2015), pp.~3--57.

\bibitem{VAN10}
{\sc J.~P. Vielma, S.~Ahmed, and G.~Nemhauser}, {\em Mixed-integer models for
  nonseparable piecewise-linear optimization: Unifying framework and
  extensions}, Operations Research, 58 (2010), pp.~303--315.

\bibitem{von1988decomposition}
{\sc B.~Von~Stengel}, {\em Decomposition of multiattribute expected-utility
  functions}, Annals of Operations Research, 16 (1988), pp.~161--183.

\bibitem{WaX23}
{\sc W.~Wang and H.~Xu}, {\em Preference robust distortion risk measure and its
  application}, Mathematical Finance, 33 (2023), pp.~389--434.

\bibitem{wu2020preference}
{\sc J.~Wu, W.~B. Haskell, W.~Huang, and H.~Xu}, {\em Preference robust
  optimization with quasi-concave choice functions for multi-attribute
  prospects}, arXiv preprint arXiv:2008.13309,  (2020).

\bibitem{WuX22}
{\sc Q.~Wu and H.~Xu}, {\em Preference robust modified optimized certainty
  equivalent}, SIAM Journal on Optimization, 32 (2022), pp.~2662--2689.

\bibitem{young1917multiple}
{\sc W.~Young}, {\em On multiple integration by parts and the second theorem of
  the mean}, Proceedings of the London Mathematical Society, 2 (1917),
  pp.~273--293.

\bibitem{Zakariazadeh2014}
{\sc A.~Zakariazadeh, S.~Jadid, and P.~Siano}, {\em Multi-objective scheduling
  of electric vehicles in smart distribution system}, Energy Conversion and
  Management, 79 (2014), pp.~43--53.

\bibitem{zhang2020preference}
{\sc Y.~Zhang, H.~Xu, and W.~Wang}, {\em Preference robust models in
  multivariate utility-based shortfall risk minimization}, Optimization Methods
  and Software, 37 (2022), pp.~712--752.

\end{thebibliography}



\begin{appendices}
\renewcommand\thefigure{\thesection.\arabic{figure}}
\renewcommand\thetable{\thesection.\arabic{table}}

\setcounter{figure}{0}
\setcounter{table}{0}

\section{Proofs}


\subsection{Proof of Proposition~\texorpdfstring{\ref{prop-uti-N}}{3.1}}
\label{app:proof-uN}

Since $\psi_l$, $l=1,\cdots,M$, take constant values over $T_{i,j}$ for 
$i=1,\cdots, N_1-1$ and
$j=1,\cdots, N_2-1$,
there exist constants $c_{i,j}^l$ such that 
$$
\psi_l(x,y):=\sum_{i=1}^{N_1-1}\sum_{j=1}^{N_2-1} c^l_{i,j} \1_{T_{i,j}}(x,y).
$$
Next, we verify that $u_N(x,y)$
 satisfies the following inequalities: 
 \begin{equation*}
    \int_{T} u_N(x,y)d\psi_l(x,y)\leq c_l, l=1,\ldots,M.
\end{equation*}
By integration in parts (see, e.g., \cite{young1917multiple} and \cite{Ans22}),
\begin{equation*}
\label{eq:parts_integral}
\begin{split}
&\int_T u_N(x,y)d\psi_l(x,y) \\ 
& = u_N(\underline{x},\underline{y})[\psi_l]_{\underline{x},\underline{y}}^{\bar{x},\bar{y}}
+\int_T [\psi_l]_{x,y}^{\bar{x},\bar{y}} d u_N(x,y)
+ \int_X [\psi_l]_{x,\underline{y}}^{\bar{x},\bar{y}}d u_N(x,\underline{y})
+\int_{Y} [\psi_l]_{\underline{x},y}^{\bar{x},\bar{y}} d u_N(\underline{x},y).
\end{split}
\end{equation*}
Since $u_N(\underline{x},\underline{y})=0$, it suffices to calculate the rest three terms at the right hand side of the equation.
Let $\hat{\psi}_l(x,y) :=[\psi_l]_{x,y}^{\bar{x},\bar{y}}$. 
Then by definition (see \cite{young1917multiple} and \cite{Ans22}) 
\begin{equation*}
\begin{split}
     \hat{\psi}_l(x,y)
  &= \psi_l(\bar{x},\bar{y})-\psi_l(x,\bar{y})-\psi_l(\bar{x},y)+\psi_l(x,y) \\
  &= \sum_{i=1}^{N_1-2}\sum_{j=1}^{N_2-2} c^l_{i,j} \1_{T_{i,j}}(x,y) -
    c^l_{N_1-1,N_2-1}\1_{T_{N_1-1,N_2-1}}(x,y).
\end{split}
\end{equation*}
Likewise, we have
\begin{equation*}
\begin{split}
     \psi_{1,l}(x)
  &: = [\psi_l]_{x,\underline{y}}^{\bar{x},\bar{y}} =
  \psi_l(\bar{x},\bar{y})-\psi_l(x,\bar{y})-\psi_l(\bar{x},\underline{y})+\psi_l(x,\underline{y}) \\
  &= \sum_{i=1}^{N_1-1} (c^l_{i,1}-c^l_{i,N_2-1}) \1_{X_i}(x) - c^l_{N_1-1,1}
\end{split}
\end{equation*}
and
\begin{equation*}
\begin{split}
     \psi_{2,l}(y)
  &:= [\psi_l]_{\underline{x},y}^{\bar{x},\bar{y}} =\psi_l(\bar{x},\bar{y})-\psi_l(\underline{x},\bar{y})-\psi_l(\bar{x},y)+\psi_l(\underline{x},y) \\
  &= \sum_{j=1}^{N_2-1} (c^l_{1,j}-c^l_{N_1-1,j}) \1_{Y_i}(y) - c^l_{1,N_2-1}.
\end{split}
\end{equation*}
Consequently, we have
\begin{eqnarray}
\label{eq:Int-by-part-1}
     && \int_T [\psi_l]_{x,y}^{\bar{x},\bar{y}}  d u_{N}(x,y) = \int_T \hat{\psi}_l(x,y) d u_{N}(x,y) \nonumber \\
    &&=  \sum_{i=1}^{N_1-1}\sum_{j=1}^{N_2-1} \int_{T_{i,j}} \lt( \sum_{i=1}^{N_1-2}\sum_{j=1}^{N_2-2} c^l_{i,j} \1_{T_{i,j}}(x,y) - c^l_{N_1-1,N_2-1}\1_{T_{N_1-1,N_2-1}}(x,y) \rt) d u_{N}(x,y) \nonumber \\
    && = \sum_{i=1}^{N_1-2}\sum_{j=1}^{N_2-2} \int_{T_{i,j}}  c^l_{i,j} d u_{N}(x,y) - \int_{T_{N_1-1,N_2-1}} c^l_{N_1-1,N_2-1} d u_{N}(x,y) \nonumber \\
    && = \sum_{i=1}^{N_1-2}\sum_{j=1}^{N_2-2} c^l_{i,j} 
    \lt( u_{N}(x_{i+1},y_{j+1})-u_{N}(x_i,y_{j+1})-u_{N}(x_{i+1},y_j)+u_{N}(x_i,y_j) \rt) \nonumber \\
    && \quad+ c^l_{N_1-1,N_2-1}\lt( u_{N}(x_{N_1},y_{N_2})-u_{N}(x_{N_1-1},y_{N_2-1})-u_{N}(x_{N_1},y_{N_2-1})+u_{N}(x_{N_1-1},y_{N_2}) \rt) \nonumber \\
    && = \sum_{i=1}^{N_1-1}\sum_{j=1}^{N_2-1} 
    c^l_{i,j} \lt( u(x_{i+1},y_{j+1})-u(x_i,y_{j+1})-u(x_{i+1},y_j)+u(x_i,y_j) \rt) \nonumber \\
    && \quad+ c^l_{N_1-1,N_2-1}\lt( u(x_{N_1},y_{N_2})-u(x_{N_1-1},y_{N_2-1})-u(x_{N_1},y_{N_2-1})+u(x_{N_1-1},y_{N_2}) \rt) \nonumber \\
    && = \int_{\underline{x},\underline{y}}^{\bar{x},\bar{y}} \psi_l(x,y) d u(x,y),
\end{eqnarray}
where the third equality 
follows from the definition of 
Lebesgue-Stieltjes
integration given 
that $u_N$ is non-decreasing and 
bounded (see \cite{Mcs47,ash2000probability}).
Likewise, we can show that
\begin{eqnarray}
\label{eq:Int-by-part-2}
    \int_{X}\psi_{1,l}(x) d u_N(x,\underline{y}) = \int_{X}\psi_{1,l}(x) d u(x,\underline{y}) 
\end{eqnarray}
and 
\begin{eqnarray}
\label{eq:Int-by-part-3}
    \int_{Y} \psi_{2,l}(y) d u_N(\underline{x},y) =  \int_{Y} \psi_{2,l}(y) d u(\underline{x},y).
\end{eqnarray}
Combing (\ref{eq:Int-by-part-1})-(\ref{eq:Int-by-part-3}),
we obtain
\begin{equation}
\label{eq-int-u-u-N}
    \int_{T} u_N(x,y)d\psi_l(x,y) = \int_{T} u(x,y)d\psi_l(x,y) \leq c_l, l=1,\ldots,M.
\end{equation}
The proof is complete.
\hfill\Box

\subsection{Proof of Proposition~\texorpdfstring{\ref{prop-int-pl}}{3.2}.}
\label{app:proof-LS}
Since $F$ is a continuous piecewise linear function with two pieces, then
there are only two possibilities that $F$ satisfies the conservative condition (\ref{eq:consevative-condition}) or not.
Without loss of generality, we assume the conservative condition fails.
According to the discussions in 
\cite{Mcs47,ash2000probability},
$F$ generates a 
LS
(outer) measure $\mu_F^*$ defined as
\begin{equation*}
    \mu_F^*((\underline{a}, \bar{a}]\times(\underline{b}, \bar{b}])= F(\bar{a},\bar{b})-F(\underline{a},\bar{b})-F(\bar{a},\underline{b})+F(\underline{a},\underline{b}).
\end{equation*}
By the definition of the
LS
integration, \begin{equation*}
    \int_{\underline{a}, \underline{b}}^{\bar{a}, \bar{b}}\psi(x,y)d F(x,y) = \int_{(\underline{a}, \bar{a}]\times(\underline{b}, \bar{b}]}\psi(x,y)d \mu_F^*.
\end{equation*}
Let $I$ and $II$ denote the triangle regions
in $[\underline{a}, \overline{a}]\times [\underline{b}, \overline{b}]$ above (including) and below (including) 
$AB$ respectively. Let $R=(a,a']\times(b,b']$ be a subset of $I$ or $II$.
Then
\[
F(a,b)+ F(a',b')=2F((a+a')/2,(b+b')/2)=F(a,b')+F(a',b),
\]
because of the linearity of $F$ over the $R$. This implies $\mu_F^*(R)=0$. 
Now we turn to discuss the measure 
over the boundary of $I\cup II$ (
denoted by $\partial (I\cup II)= ((\underline{a},\bar{a}]\times\bar{b})\cup(\bar{a}\times(\underline{b},\bar{b}])$).
For any small constant $\epsilon>0$,
$$
\mu_F^*((\underline{a},\bar{a}]\times\bar{b})\leq \mu_F^*((\underline{a},\bar{a}]\times(\bar{b}-\epsilon,\bar{b}])=F(\underline{a},\bar{b}-\epsilon)-F(\underline{a},\bar{b})-F(\bar{a},\bar{b}-\epsilon)+F(\bar{a},\bar{b}).
$$
By driving $\epsilon$ to zero, 
we obtain
$$
\mu_F^*((\underline{a},\bar{a}]\times\bar{b})\leq \lim_{\epsilon\to0} (F(\underline{a},\bar{b}-\epsilon)-F(\underline{a},\bar{b})-F(\bar{a},\bar{b}-\epsilon)+F(\bar{a},\bar{b}))=0,
$$
which implies $\mu_F^*((\underline{a},\bar{a}]\times\bar{b})=0$.
Likewise, we can also obtain $\mu_F^*(\bar{a}\times(\underline{b},\bar{b}])=0$ and hence $\mu_F^*(\partial (I\cup II))=0$.
Let $\{a_i\}$ and $\{b_i\}$ be two sequences of monotonically increasing numbers such that 
$R_i:=(a_i,a_{i+1}]\times(b_i,b_{i+1}]\subset \inmat{int\,} I$
and $\bigcup_i R_i = I$.
By the property of outer measure, 
$$
\mu_F^*(\inmat{int\,} I)\leq 
\sum_i\mu_F^*(R_i)=0.
$$
This shows $\mu_F^*(\inmat{int\,} I)=0$.
Likewise, $\mu_F^*(\inmat{int\,} II)=0$.
Consequently,
we have $\mu_F^*(I\cup II)=\mu_F^*(I\cap II)$.
Next, let $t\in (\underline{a},\overline{a}]$ and consider the segment $L=(\underline{a},t]\times(\underline{b},y(t)]\cap(I\cap II)$,
we have
\begin{equation*}
    \mu_F^*(L) = \frac{t-\underline{a}}{\bar{a}-\underline{a}} \, \mu_F^*(I\cap II),
\end{equation*}
where $y(t)$ is the linear function representing 
$I\cap II$ (AB).
\begin{equation*}
\begin{split}
    \int_{[\underline{a},\overline{a}]\times [\underline{b}, \overline{b}]}\psi(x,y) & d F(x,y) = \int_{I\cap II} \psi(x,y(x)) d\mu_F^* 
    \\
    & = \frac{\mu_F^*(I\cap II)}{\bar{a}-\underline{a}} \lim_{t\to\bar{a}} \int_{\underline{a}}^t \psi(x,y(x)) dx
    = \frac{\mu_F^*(I\cap II)}{\overline{a}-\underline{a}}\int_{\underline{a}}^{\overline{a}} \psi(x,y(x))d x,
\end{split}
\end{equation*}
where the third equality holds since $\psi(x,y(x))$ is Riemann integrable.
\hfill \Box

\subsection{Proof of Proposition~\texorpdfstring{\ref{prop:single-MILP}}{4.2}}
By introducing dual variables,
we can write down the Lagrange function of the inner minimization problem (\ref{eq:PRO_MILP_eqi})
w.r.t.~${\bm u}$
\bgeq
&& L({\bm u},{\bm \lambda}^1,{\bm \lambda}^2,{\bm \eta}^1, {\bm \eta}^2, {\bm \tau}, \sigma,{\bm \zeta})\\
&& = \sum_{k=1}^K p_k \sum_{i=1}^{N_1} \sum_{j=1}^{N_2} \alpha_{i,j}^k u_{i,j} 
+\sum_{i=1}^{N_1-1}\sum_{j=1}^{N_2} \lambda_{i,j}^1 (u_{i,j}-u_{i+1,j})
+\sum_{i=1}^{N_1}\sum_{j=1}^{N_2-1} \lambda_{i,j}^2 (u_{i,j}-u_{i,j+1})\\
&& \quad +\sum_{i=1}^{N_1-1}\sum_{j=1}^{N_2} \eta_{i,j}^1 (u_{i+1,j}-u_{i,j}-L(x_{i+1}-x_i))
+\sum_{i=1}^{N_1}\sum_{j=1}^{N_2-1} \eta_{i,j}^2
(u_{i,j+1}-u_{i,j}-L(y_{j+1}-y_j))\\
&& \quad
+\sum_{i=1}^{N_1-1}\sum_{j=1}^{N_2-1} \tau_{i,j} 
(u_{i,j}+u_{i+1,j+1}-u_{i,j+1}-u_{i+1,j}) 
 +\sigma(1-u_{N_1,N_2})+\sum_{l=1}^M \zeta_l\sum_{i=1}^{N_1}\sum_{j=1}^{N_2} Q_{i,j}^l,
\edeq
where ${\bm \lambda}^1\in \R^{(N_1-1)\times N_2}_+$,
${\bm \lambda}^2\in \R^{N_1\times (N_2-1)}_+$,
${\bm \eta}^1\in \R^{(N_1-1)\times N_2}_+$,
${\bm \eta}^2\in \R^{N_1\times (N_2-1)}_+$,
$\tau\in \R^{(N_1-1)\times (N_2-1)}_+$,
$\sigma\in \R$ and $\zeta\in \R^M_+$.
We can then derive 
the Lagrange dual formulation and merge it into the outer 
maximization problem to obtain  (\ref{eq:PRO_MILP_single}).
\hfill
$\Box$

\subsection{Proof of Proposition~\texorpdfstring{\ref{prop-d}}{5.1}.}
Case (i). 
$\scrg=\scrg_K$. 
We have
\begin{align}
    & \dd_{\scrg_K} (u,u_{N}) \nonumber \\
    & = \sup_{g\in\scrg_K} \lt|\int_T g(x,y) d u(x,y) - \int_T g(x,y) d u_{N}(x,y)\rt| \nonumber \\
    & \leq \sum_{i=1}^{N_1-1} \sum_{j=1}^{N_2-1} \sup_{g\in\scrg_K} 
    \lt|
    \int_{T_{i,j}} g(x,y) d u(x,y) -  \int_{T_{i,j}} g(x,y) d u_{N}(x,y)
    \rt| \nonumber \\
    & \leq \sum_{i=1}^{N_1-1} \sum_{j=1}^{N_2-1} \sup_{g\in\scrg_K} 
    \lt|
    \int_{T_{i,j}} g(x,y) d u(x,y)
    - \int_{T_{i,j}} g(x_{i},y_{j}) d u(x,y) \rt. \nonumber \\
    &  \hspace{10em} \lt.
    + \int_{T_{i,j}} g(x_{i},y_{j}) d u(x,y) 
    - \int_{T_{i,j}} g(x,y) d u_{N}(x,y)
    \rt| \nonumber \\
    & \leq \sum_{i=1}^{N_1-1} \sum_{j=1}^{N_2-1} \sup_{g\in\scrg_K} 
    \lt( \lt|
    \int_{T_{i,j}} |g(x,y)-g(x_{i},y_{j})| 
    d u(x,y) \rt| \rt. \nonumber \\
    & \hspace{10em} \lt. + \lt | \int_{T_{i,j}} |g(x_{i},y_{j})-g(x,y)| d u_{N}(x,y) \rt|
    \rt) \nonumber \\
    & \leq (\beta_{N_1}^2+\beta_{N_2}^2)^{1/2} \sum_{i=1}^{N_1-1} \sum_{j=1}^{N_2-1} 
    \lt(\lt|
    \int_{T_{i,j}} d u(x,y)\rt|
    +\lt|\int_{T_{i,j}} d u_{N}(x,y) \rt|
    \rt) \nonumber \\
    & = 2(\beta_{N_1}^2+\beta_{N_2}^2)^{1/2} |1-u_N(\underline{x},\bar{y})-u_N(\underline{y},\bar{x})| \leq 2(\beta_{N_1}^2+\beta_{N_2}^2)^{1/2}, \label{eq-u-u-N}
\end{align}
where the last equality holds due to that $u$ and $u_N$ satisfy 
the conservative conditions, which implies that 
$$
\int_{T_{i,j}} d u(x,y) = \int_{T_{i,j}} d u_N(x,y) 
 = u_{i+1,j+1}-u_{i+1,j}-u_{i,j+1}+u_{i,j}\leq 0.
$$
Then 
$$
\sum_{i=1}^{N_1-1} \sum_{j=1}^{N_2-1} \lt| \int_{T_{i,j}} d u(x,y) \rt| = \sum_{i=1}^{N_1-1} \sum_{j=1}^{N_2-1} \lt| \int_{T_{i,j}} d u_N(x,y) \rt| 
=|1-u(\underline{x},\bar{y})-u(\underline{y},\bar{x})|.
$$

Case (ii). $\scrg=\scrg_I$.
We only consider Type 1 PLA. 
Similar arguments can be established 
for Type 2 PLA.
Let $(x,y)\in T_{i,j}$. Consider the 
case that $(x,y)$ lies below the main diagonal, i.e., 
$0\leq\frac{y-y_j}{x-x_i}\leq\frac{y_{j+1}-y_j}{x_{i+1}-x_i}$. Thus
\begin{eqnarray*}
    && |u_{N}(x,y)-u(x,y)| \\
    && = \lt| \lt( 1-\frac{x-x_i}{x_{i+1}-x_i} \rt) u_{i,j}
    +\lt( \frac{x-x_i}{x_{i+1}-x_i} -\frac{y-y_j}{y_{j+1}-y_j} \rt) u_{i+1,j}  +\frac{y-y_j}{y_{j+1}-y_j} u_{i+1,j+1}-u(x,y) \rt| \\
    &&  \leq \lt| \lt(1-\frac{x-x_i}{x_{i+1}-x_i}\rt) (u_{i,j}-u(x,y)) \rt|
    +\lt| \lt(\frac{x-x_i}{x_{i+1}-x_i}-\frac{y-y_j}{y_{j+1}-y_j}\rt) (u_{i+1,j}-u(x,y)) \rt| \\
    && 
    \quad+\lt|\frac{y-y_j}{y_{j+1}-y_j} (u_{i+1,j+1}-u(x,y)) \rt|.
\end{eqnarray*}
Since $u_{i,j}=u(x_i,y_j)$, by the Lipschitz continuity 
of $u$, we have
\begin{equation*}
    |u_{i,j}-u(x,y)| = |u(x_i,y_j)-u(x,y)|\leq L(\beta_{N_1}+\beta_{N_2}).
\end{equation*}
Likewise, we can obtain
$|u_{i+1,j}-u(x,y)|\leq L(\beta_{N_1}+\beta_{N_2})$
and 
$|u_{i+1,j+1}-u(x,y)|\leq L(\beta_{N_1}+\beta_{N_2})$,
which give rise to
\begin{equation}
\label{eq-u-uN}
    |u_{N}(x,y)-u(x,y)|\leq L(\beta_{N_1}+\beta_{N_2}).
\end{equation}
We can obtain the same inequality when $(x,y)\in[x_{i-1},x_i]\times[y_{j-1},y_j]$ and
$\frac{y-y_j}{x-x_i}\geq\frac{y_{j+1}-y_j}{x_{i+1}-x_i}$. 
Summarizing the discussions above, we have 
\begin{eqnarray*}
    \dd_{\scrg_I} (u,u_{N}) &=& \sup_{g\in\scrg} \lt| \int_{\underline{x},\underline{y}}^{x,y} d u(x,y) - \int_{\underline{x},\underline{y}}^{x,y} d u_{N}(x,y) \rt| \\
    &=& \sup_{(x,y)\in T} 
    |u(x,y)-u(x,\underline{y})-u(\underline{x},y)-u_{N}(x,y)+u_{N}(x,\underline{y})+u_{N}(\underline{x},y)| \\
    &\leq& 2L\lt(\beta_{N_1}+\beta_{N_2} \rt), 
\end{eqnarray*}
which implies (\ref{eq-d}).
The proof is complete.
\hfill \Box


\subsection{Proof of Theorem~\texorpdfstring{\ref{thm-erramb}}{5.1}.}
Let $\hat{\alpha}<\alpha$ be a positive number.
Under 
Slater's condition (\ref{eq-sla}), there exists a function $u^0_{N} \in\calu_{N}$ and a positive number $N^0=N_1^0\times N_2^0$ such that 
\begin{equation}
\label{eq-sla-0}
    \la u^0_{N},\bdps \ra -\bdc+\hat{\alpha} \mathbb{B}^M \subset \R_-^M
\end{equation}
for $N\geq N^0$.
The existence follows from Proposition~\ref{prop-d} in that there exists $u^0$ satisfying (\ref{eq-sla}), 
and by (\ref{eq-u-uN})
we can construct a piecewise linear utility function $u^0_{N}$ of $u^0$ such that $u^0_{N}\to u^0$  
under $\|\cdot\|_\infty$ uniformly
as $\beta_{N_i}\to0$, $i=1,2$. 
By applying Lemma~\ref{lem-hof} to $\calu$ under 
Slater's condition (\ref{eq-sla-0}), for any $\tldu\in\scru_N$,
\begin{equation}
\label{eq-bbd}
    \mathbb{D}_{\scrg}(\tldu,\calu_N)\leq \frac{\dd_{\scrg}(\tldu,u^0_{N})}{\hat{\alpha}}
    \|(\la \tldu,\bdps \ra-\bdc)_+\|
\end{equation}
for all $N\geq N^0$. 
Let $u\in\calu$ and $u_{N}$ be defined as in Proposition~\ref{prop-uti-N}.
Then 
\begin{align}
    & \|
    \la u_{N},\bdps \ra 
    - \la u,\bdps \ra 
    \|^2 \nonumber \\
    & = \sum_{l=1}^M \lt| 
    \int_T u_{N}(x,y) d \psi_l(x,y)-\int_T u(x,y) d \psi_l(x,y)
    \rt|^2 \nonumber \\   
    & \leq \sum_{l=1}^M \lt| 
    \int_T |u_{N}(x,y)-u(x,y)|d\psi_l(x,y)
    \rt|^2 \nonumber\\
    & \leq
    L^2 (\beta_{N_1}+\beta_{N_2})^2
    \sum_{l=1}^M \lt| \int_T d \psi_l(t)\rt|^2.
    \label{eq-inpro} 
\end{align}
By the triangle inequality for the pseudo-metric and (\ref{eq-bbd}), we have 
\begin{align}
    \dd_{\scrg}(u,\calu_{N}) &\leq \dd_{\scrg}(u,u_{N})+\dd_{\scrg}(u_{N},\calu_{N}) \nonumber \\
    &\leq \dd_{\scrg}(u,u_{N})+\frac{\dd_{\scrg}(u_{N},u_{N}^0)}{\hat{\alpha}}
    \|(\la u_{N},\bdps \ra-\bdc)_+\| \nonumber \\
    &= \dd_{\scrg}(u,u_{N})+\frac{\dd_{\scrg}(u_{N},u_{N}^0)}{\hat{\alpha}}
    [\|(\la u_{N},\bdps \ra-\bdc)_+\|-\|(\la u,\bdps \ra-\bdc)_+\|] \nonumber \\
    &\leq \dd_{\scrg}(u,u_{N})+\frac{\dd_{\scrg}(u_{N},u_{N}^0)}{\hat{\alpha}}
    \|\la u_{N},\bdps \ra-\la u,\bdps \ra)\| \nonumber \\
    &\leq \dd_{\scrg}(u,u_{N})+\frac{\dd_{\scrg}(u_{N},u_{N}^0)}{\hat{\alpha}} 
    L(\beta_{N_1}+\beta_{N_2}) \lt(\sum_{l=1}^M \lt| \int_T d \psi_l(t)\rt|^2\rt)^{1/2},
    \label{eq-uU} 
\end{align}
where the equality holds due to $u\in\calu$, i.e. $(\la u,\bdps \ra-\bdc)_+=0$ and the last inequality comes from (\ref{eq-inpro}).
In what follows, we turn to estimate $\dd_{\scrg}(u,u_{N})$ and $\dd_{\scrg}(u_{N},u_{N}^0)$ when $\scrg$ have specific form.

Case (i). If $\scrg=\scrg_K$, then $\dd_{\scrg_K}(u_N,u^0_N)\leq \lt((\bar{x}-\underline{x})^2+(\bar{y}-\underline{y})^2\rt)^{1/2}$ and $\dd_{\scrg_K}(u,u_N)\leq 2(\beta_{N_1}^2+\beta_{N_2}^2)^{1/2}$ 
by Proposition~\ref{prop-d}~(i).
Taking supremum w.r.t. $u$ over $\calu$ on both sides of (\ref{eq-uU}), we obtain 
\begin{equation*}
    \mathbb{D}_{\scrg_K}(\calu,\calu_N)\leq 2(\beta_{N_1}^2+\beta_{N_2}^2)^{1/2} + 
    L(\beta_{N_1}+\beta_{N_2})
    \frac{\lt((\bar{x}-\underline{x})^2+(\bar{y}-\underline{y})^2\rt)^2}{\hat{\alpha}} 
    \lt(\sum_{l=1}^M \lt| \int_T d \psi_l(t)\rt|\rt)^{1/2}
\end{equation*}
and hence (\ref{eq-erram-L}) holds since $\mathbb{D}_{\scrg_K}(\calu_{N},\calu)=0$.

Case (ii). 
If $\scrg=\scrg_I$, then
$\dd_{\scrg_I}(u_N,u^0_N)\leq1$ and $\dd_{\scrg_I}(u,u_N)\leq 2L\lt( \beta_{N_1}+\beta_{N_2} \rt)$ by Proposition~\ref{prop-d}~(ii).
Following a similar analysis to Case (i), we obtain (\ref{eq-erram-I}).
\hfill $\Box$


\subsection{Proof of Corollary~\texorpdfstring{
\ref{cor-err-optval-discrete}}{5.1}.}
Since 
$\calu_{N}\subset\calu$ by definition, then
$
\mathbb{D}_{\scrg}(\calu_{N},\calu)
=0
$
for any $\mathscr{G}$.
Thus, it suffices to estimate $\mathbb{D}_{\scrg}(\calu,\calu_{N})$.
For any $u\in \calu$, it follows by Proposition~\ref{prop-uti-N} that
we can construct $u_N$ of Type-1 PLA or Type-2 PLA such that $u_N\in \calu_N$. 
Consequently, 
in the case that $\scrg=\scrg_K$, we have
$$
\dd_{\scrg_K}(u,\calu_{N}) \leq \dd_{\scrg_K}(u,u_{N})+\dd_{\scrg_K}(u_{N},\calu_{N})
= \dd_{\scrg_K}(u,u_{N})\leq 2(\beta_{N_1}^2+\beta_{N_2}^2)^{1/2},
$$
where the last inequality follows from Proposition~\ref{prop-d}~(i).
Hence
$$
\mathbb{H}_{\scrg_K} (\calu,\calu_{N})= \max\{0,
\mathbb{D}_{\scrg_K} (\calu,\calu_{N}) \}= \sup_{u\in\calu} \dd_{\scrg_K}(u,\calu_{N})\leq 2(\beta_{N_1}^2+\beta_{N_2}^2)^{1/2}. 
$$
In the case that $\scrg=\scrg_I$, we have 
$$
\dd_{\scrg_I}(u,\calu_{N}) \leq \dd_{\scrg_I}(u,u_{N}) \leq 2L\lt( \beta_{N_1}+\beta_{N_2} \rt),
$$
where the last inequality follows from 
Proposition~\ref{prop-d}~(ii)
and hence 
$$
\mathbb{H}_{\scrg_I} (\calu,\calu_{N})=
\max\{0,\mathbb{D}_{\scrg_I} (\calu,\calu_{N})\}\leq 2L\lt( \beta_{N_1}+\beta_{N_2} \rt).$$
The proof is complete.
\hfill\Box

\subsection{Proof of Theorem~\texorpdfstring{\ref{thm-optval}}{5.2}.}
Part (i). 
It is well-known that 
\begin{equation*}
    |\vt_{N}-\vt|\leq \max_{\bdz\in Z} \Big{|} \min_{u\in \calu_{N}} \bbe_P [u(\bdf(\bdz,\bdxi))] 
    - \min_{u\in \calu} \bbe_P [u(\bdf(\bdz,\bdxi))] \Big{|}. 
\end{equation*}
Let $\delta$ be a small positive number. 
For any $\bdz\in Z$, we can find 
a $\delta$-optimal solution
$u^{\bdz} \in \calu$ 
and $u^{\bdz}_N \in \calu_N$ such that
$$
\bbe_P [u^{\bdz} (\bdf(\bdz,\bdxi))] \leq \min_{u\in \calu} \bbe_P [u(\bdf(\bdz,\bdxi))]+\delta, \quad
\bbe_P [u^{\bdz}_{N} (\bdf(\bdz,\bdxi))] \geq \min_{u\in \calu_{N}} \bbe_P [u(\bdf(\bdz,\bdxi))].
$$
Combing the above inequalities, we have 
\begin{eqnarray}
    \min_{u\in \calu_{N}} \bbe_P [u(\bdf(\bdz,\bdxi))]-\min_{u\in \calu} \bbe_P [u(\bdf(\bdz,\bdxi))] &\leq& \bbe_P [u^{\bdz}_{N} (\bdf(\bdz,\bdxi))-u^{\bdz} (\bdf(\bdz,\bdxi))] +\delta \qquad \nonumber\\
    &\leq& \sup_{(x,y)\in T} |u^{\bdz}_{N}(x,y)-u^{\bdz}(x,y)| +\delta.
    \label{eq:thm2-proof}
\end{eqnarray}
On the other hand,  for any 
$u,v\in {\cal U}$
\begin{eqnarray*}
    && \sup_{(x,y)\in T} |u(x,y)-v(x,y)| \\
    && \leq \sup_{(x,y)\in T} ( |u(x,y)-v(x,y)-u(x,\underline{y})-u(\underline{x},y)+v(x,\underline{y})+v(\underline{x},y)| \\ 
    && \hspace{5em} +|v(x,\underline{y})+v(\underline{x},y)-u(x,\underline{y})-u(\underline{x},y)| ) \\
    && = \sup_{(x,y)\in T}  \lt( \lt| \int_{\underline{x},\underline{y}}^{x,y} d u(x,y)-\int_{\underline{x},\underline{y}}^{x,y} d v(x,y) \rt| + |v(x,\underline{y})+v(\underline{x},y)-u(x,\underline{y})-u(\underline{x},y)| \rt) \\
    && \leq \sup_{g\in\scrg_I} \lt| \int_T g(x,y) d u(x,y)-\int_T g(x,y) d v(t) \rt| + \sup_{(x,y)\in T} |v(x,\underline{y})+v(\underline{x},y)-u(x,\underline{y})-u(\underline{x},y)|.
\end{eqnarray*}
Since ${\cal U}_N\subset {\cal U}$, by setting $u=u^{\bdz}_{N}$ and $v=u^{\bdz}$, we have 
\begin{eqnarray}
\sup_{(x,y)\in T} |u^{\bdz}_{N}(x,y)-u^{\bdz}(x,y)| 
    &\leq& \dd_{\scrg_I} (u_{N}^{\bdz},u^{\bdz})
   + \sup_{(x,y)\in T} |u_{N}^{\bdz}(x,\underline{y})-u^{\bdz}(x,\underline{y})+u_{N}^{\bdz}(\underline{x},y)-u^{\bdz}(\underline{x},y)| \notag \\
     &\leq& 
    \mathbb{H}_{\scrg_I} (\calu_{N},\calu)+ L(\beta_{N_1}+\beta_{N_2}).
      \label{eq:thm2-proof-1}
\end{eqnarray}
Combining (\ref{eq:thm2-proof})-(\ref{eq:thm2-proof-1}), we obtain
\begin{eqnarray*}
    \min_{u\in \calu_{N}} \bbe_P [u(\bdf(\bdz,\bdxi))]-\min_{u\in \calu} \bbe_P [u(\bdf(\bdz,\bdxi))] 
    \leq \mathbb{H}_{\scrg_I} (\calu_{N},\calu)+ L(\beta_{N_1}+\beta_{N_2}) +\delta.
\end{eqnarray*}
By exchanging the position of $\calu$ and $\calu_{N}$, we can use the same argument to derive
$$
\min_{u\in \calu} \bbe_P [u(\bdf(\bdz,\bdxi))]-\min_{u\in \calu_{N}} \bbe_P [u(\bdf(\bdz,\bdxi))] \leq
\mathbb{H}_{\scrg_I} (\calu,\calu_{N})+ L(\beta_{N_1}+\beta_{N_2}) +\delta.
$$
Since $\delta>0$ can be arbitrarily small, we obtain 
\begin{equation}
\label{eq-dis-vt}
    |\vt_{N}-\vt| \leq \max_{\bdz\in Z} \lt| \min_{u\in \calu} \bbe_P [u(\bdf(\bdz,\bdxi))]-\min_{u\in \calu_{N}} \bbe_P [u(\bdf(\bdz,\bdxi))] \rt| \leq \mathbb{H}_{\scrg_I} (\calu,\calu_{N})+ L(\beta_{N_1}+\beta_{N_2})
\end{equation}
and hence (\ref{eq-err-vt}) follows from (\ref{eq-erram-I}).

Part(ii).
Observe that $\Lambda(\cdot)$ is a non-decreasing function,
thus its generalized inverse is well-defined.
For any $\bdz_{N}^*\in Z_{N}^*$ and $\bdz^*\in Z^*$, 
\begin{eqnarray*}
    \Lambda(d(\bdz_{N}^*,\bdz)) &\leq& v({\bdz_N^*})-\vt^* =v(\bdz_{N}^*)-v(\bdz^*)\leq |v(\bdz_{N}^*)-v_N(\bdz_{N}^*)| + |v(\bdz_{N}^*)-v(\bdz^*)| \\
    &\leq & 2\max_{\bdz\in Z} |v(\bdz)-v_N(\bdz)|.
\end{eqnarray*}
Combining the inequality above with (\ref{eq-err-vt}), we obtain
$$
d(\bdz_{N}^*,Z^*) \leq \Lambda^{-1} \lt( 2\max_{\bdz\in Z} |v(\bdz)-v_N(\bdz)| \rt) \leq \Lambda^{-1}(2 \mathbb{H}_{\scrg_I}(\calu_{N},\calu)+2L(\beta_{N_1}+\beta_{N_2}))
$$
and hence (\ref{eq-err-so}) follows.
\hfill $\Box$
\end{appendices}


\end{document}